\documentclass[11 pt, reqno]{amsart}
\usepackage{hhline}

\usepackage{tikz-cd}
\usepackage{verbatim}
\usepackage{amsmath}
\usepackage{amsxtra}
\usepackage{amscd}
\usepackage{amsthm}
\usepackage{amsfonts}
\usepackage{amssymb}
\usetikzlibrary{arrows}
\usepackage{array}

\makeatletter %%\def\l@subsection{\@tocline{2}{0pt}{1pc}{5pc}{}} 
\def\l@subsection{\@tocline{2}{0pt}{2pc}{6pc}{}} \makeatother

\usepackage[all, cmtip]{xy}
\usepackage{stmaryrd}
\usepackage{xcolor}
\definecolor{dgr}{RGB}{0,47,166}
	
\usepackage{amsthm,amsfonts,amssymb,amsmath,amsxtra}
\usepackage[all]{xy}
\SelectTips{cm}{}
\usepackage{xr-hyper}

\usepackage{setspace}
	\newcommand{\Sh}{\mathrm{Sh}}

\usepackage[pagebackref=true, colorlinks=true, linkcolor=black, citecolor=dgr, urlcolor=black]{hyperref}

\usepackage{verbatim}

\usepackage[margin=2.7cm]{geometry}

\usepackage{mathrsfs}

\RequirePackage{xspace}
\RequirePackage{etoolbox}
\RequirePackage{varwidth}
\RequirePackage{enumitem}
\RequirePackage{tensor}
\RequirePackage{mathtools}
\RequirePackage{longtable}
\RequirePackage{multirow}

\newcommand{\nrm}{e}
\newcommand{\sqD}{\delta} %simpler than \sqrt\Delta, no conflicts in sec.9

\def\<{\langle}	
\def\>{\rangle}
\def\e{\epsilon}

\def\sg{\sigma}
\def\vth{\vartheta}

\newcommand{\sA}{\ensuremath{\mathscr{A}}\xspace}

\newcommand{\sC}{\ensuremath{\mathscr{C}}\xspace}

\newcommand{\sH}{\ensuremath{\mathscr{H}}\xspace}
\newcommand{\sI}{\ensuremath{\mathscr{I}}\xspace}
\newcommand{\sJ}{\ensuremath{\mathscr{J}}\xspace}

\newcommand{\sL}{\ensuremath{\mathscr{L}}\xspace}
\newcommand{\sM}{\ensuremath{\mathscr{M}}\xspace}

\newcommand{\sO}{\ensuremath{\mathscr{O}}\xspace}

\newcommand{\sR}{\ensuremath{\mathscr{R}}\xspace}

\newcommand{\sV}{\ensuremath{\mathscr{V}}\xspace}
\newcommand{\sW}{\ensuremath{\mathscr{W}}\xspace}
\newcommand{\sX}{\ensuremath{\mathscr{X}}\xspace}
\newcommand{\sY}{\ensuremath{\mathscr{Y}}\xspace}
\newcommand{\sZ}{\ensuremath{\mathscr{Z}}\xspace}

\newcommand{\fkm}{\ensuremath{\mathfrak{m}}\xspace}

\newcommand{\fkp}{\ensuremath{\mathfrak{p}}\xspace}

\newcommand{\fkS}{\ensuremath{\mathfrak{S}}\xspace}

\newcommand{\Iw}{{\rm Iw}}
\newcommand{\hT}{T}

\newcommand{\llb}{\llbracket}
\newcommand{\rrb}{\rrbracket}

\newcommand{\nat}{{\natural}}

\newcommand{\G}{\mathrm{G}}
\newcommand{\C}{\mathbf{C}}

\newcommand{\OO}{\mathscr{O}}
\newcommand{\A}{\mathbf{A}}
\newcommand{\Q}{\mathbf{Q}}
\newcommand{\R}{\mathbf{R}}
\newcommand{\Z}{\mathbf{Z}}

\newcommand{\bA}{\mathbf A}

\newcommand{\bC}{\mathbf C}

\newcommand{\bG}{\mathbf G}

\newcommand{\bN}{\mathbf N}
\newcommand{\bQ}{\mathbf Q}

\newcommand{\bZ}{\mathbf Z}

\newcommand{\ff}{ f}
\newcommand{\df}{\dot{f}}
\newcommand{\hI}{ I}
\newcommand{\hJ}{ J}

\newcommand{\BA}{\ensuremath{\mathbb {A}}\xspace}

\newcommand{\BC}{\ensuremath{\mathbb {C}}\xspace}
\newcommand{\BD}{\ensuremath{\mathbb {D}}\xspace}

\newcommand{{\BG}}{\ensuremath{\mathbb {G}}\xspace}

\newcommand{{\BK}}{\ensuremath{\mathbb {K}}\xspace}

\newcommand{\BP}{\ensuremath{\mathbb {P}}\xspace}
\newcommand{\BQ}{\ensuremath{\mathbb {Q}}\xspace}

\newcommand{\BS}{\ensuremath{\mathbb {S}}\xspace}
\newcommand{\BT}{\ensuremath{\mathbb {T}}\xspace}

\newcommand{\BZ}{\ensuremath{\mathbb {Z}}\xspace}

\newcommand{\CG}{\ensuremath{\mathcal {G}}\xspace}

\newcommand{\CJ}{\ensuremath{\mathcal {J}}\xspace}

\newcommand{\CM}{\ensuremath{\mathcal {M}}\xspace}
\newcommand{\CN}{\ensuremath{\mathcal {N}}\xspace}

\newcommand{\CS}{\ensuremath{\mathcal {S}}\xspace}

\newcommand{\CV}{\ensuremath{\mathcal {V}}\xspace}

\newcommand{\RB}{\ensuremath{\mathrm {B}}\xspace}

\newcommand{\RD}{\ensuremath{\mathrm {D}}\xspace}

\newcommand{\RF}{\ensuremath{\mathrm {F}}\xspace}
\newcommand{\RG}{\ensuremath{\mathrm {G}}\xspace}
\newcommand{\RH}{\ensuremath{\mathrm {H}}\xspace}

\newcommand{\RM}{\ensuremath{\mathrm {M}}\xspace}
\newcommand{\RN}{\ensuremath{\mathrm {N}}\xspace}

\newcommand{\RP}{\ensuremath{\mathrm {P}}\xspace}

\newcommand{\RS}{\ensuremath{\mathrm {S}}\xspace}
\newcommand{\RT}{\ensuremath{\mathrm {T}}\xspace}
\newcommand{\RU}{\ensuremath{\mathrm {U}}\xspace}
\newcommand{\RV}{\ensuremath{\mathrm {V}}\xspace}

\newcommand{\RZ}{\ensuremath{\mathrm {Z}}\xspace}

\newcommand{\As}{{\mathrm{As}}}
\newcommand{\ad}{{\mathrm{ad}}}

\DeclareMathOperator{\Aut}{Aut}

\newcommand{\Ch}{{\mathrm{Ch}}}
\DeclareMathOperator{\charac}{char}

\DeclareMathOperator{\diag}{diag}
\newcommand{\disc}{{\mathrm{disc}}}

\DeclareMathOperator{\End}{End}

\DeclareMathOperator{\Gal}{Gal}
\newcommand{\GL}{\mathrm{GL}}

\DeclareMathOperator{\Hom}{Hom}

\renewcommand{\i}{^{-1}}
\newcommand{\id}{\ensuremath{\mathrm{id}}\xspace}
\let\Im\relax
\DeclareMathOperator{\Im}{Im}
\newcommand{\Ind}{{\mathrm{Ind}}}

\DeclareMathOperator{\Ker}{Ker}

\DeclareMathOperator{\Lie}{Lie}
\newcommand{\LNSch}{\ensuremath{(\mathrm{LNSch})}\xspace}

\newcommand{\M}{\mathrm{M}}

\DeclareMathOperator{\Nm}{Nm}

\DeclareMathOperator{\ord}{ord}
\DeclareMathOperator{\tord}{t-ord}

\renewcommand{\Re}{{\mathrm{Re}}}

\newcommand{\reg}{{\mathrm{reg}}}
\DeclareMathOperator{\Res}{Res}
\DeclareMathOperator{\Ros}{Ros}
\newcommand{\rs}{\ensuremath{\mathrm{rs}}\xspace}

\newcommand{\EC}{\ensuremath{{\textup{\'{E}tCorr}}}}

\newcommand{\sform}{\ensuremath{(\text{~,~})}\xspace}

\newcommand{\SL}{{\mathrm{SL}}}
\DeclareMathOperator{\Spec}{Spec}
\DeclareMathOperator{\Spf}{Spf}

\newcommand{\ssm}{\smallsetminus}

\DeclareMathOperator{\Sym}{Sym}

\DeclareMathOperator{\tr}{tr}
\DeclareMathOperator{\Tr}{Tr}

\newcommand{\U}{\mathrm{U}}
\DeclareMathOperator{\uHom}{\underline{Hom}}

\DeclareMathOperator{\vol}{vol}

\newcommand{\one}{\mathbf{1}}

\newcommand{\wt}{\widetilde}
\newcommand{\wh}{\widehat}

\newcommand{\pair}[1]{\langle {#1} \rangle}

\newcommand{\ov}{\overline}
\newcommand{\incl}{\hookrightarrow}
\newcommand{\into}{\hookrightarrow}
\newcommand{\lra}{\longrightarrow}
\newcommand{\imp}{\Longrightarrow}

\newcommand{\bs}{\backslash}
\newcommand{\la}{\langle}
\newcommand{\ra}{\rangle}

\newcommand{\ep}{\varepsilon}
\newcommand{\eps}{\varepsilon}
\newcommand{\lm}{\lambda}
\newcommand{\Lm}{\Lambda}
\newcommand{\vpi}{\varpi}
\newcommand{\vphi}{\varphi}

\providecommand{\twomat}[4]{\left(\begin{array}{cc}#1&#2\\#3&#4\end{array}\right)}
\providecommand{\smalltwomat}[4]{\left(\begin{smallmatrix}#1&#2\\#3&#4\end{smallmatrix}\right)}

\newcommand{\ol}[1]{\overline{#1}{}}
\newcommand{\ul}[1]{\underline{#1}{}}

\newtheorem{theorem}{Theorem}
\newtheorem{proposition}[theorem]{Proposition}
\newtheorem{hypothesis}[theorem]{Hypothesis}
\newtheorem{lemma}[theorem]{Lemma}
\newtheorem {conjecture}[theorem]{Conjecture}
\newtheorem{corollary}[theorem]{Corollary}

\newtheorem{theoA}{Theorem}

\theoremstyle{definition}
\newtheorem{definition}[theorem]{Definition}

\theoremstyle{remark}
\newtheorem{remark}[theorem]{Remark}

\newenvironment{altenumerate}
   {\begin{list}
      {\textup{(\theenumi)} }
      {\usecounter{enumi}
       \setlength{\labelwidth}{0pt}
       \setlength{\labelsep}{0pt}
       \setlength{\leftmargin}{0pt}
       \setlength{\itemsep}{\the\smallskipamount}
       \renewcommand{\theenumi}{\roman{enumi}}
      }}
   {\end{list}}
\newenvironment{altitemize}
   {\begin{list}
      {$\bullet$}
      {\setlength{\labelwidth}{0pt}
	   \setlength{\itemindent}{5pt}
       \setlength{\labelsep}{5pt}
       \setlength{\leftmargin}{0pt}
       \setlength{\itemsep}{\the\smallskipamount}
      }}
   {\end{list}}

\newcommand{\beq}{\begin{equation}\begin{aligned}}
\newcommand{\eeq}{\end{aligned}\end{equation}}

\newcommand{\beqq}{\begin{equation*}\begin{aligned}}
\newcommand{\eeqq}{\end{aligned}\end{equation*}}

\newcommand{\lb}[1]{\label{#1}}

\numberwithin{equation}{subsection}
\numberwithin{theorem}{subsection}
\numberwithin{table}{subsection}

\newcommand{\ot}{\otimes}

\newcommand{\ts}{\times}

\setitemize[0]{leftmargin=*,itemsep=\the\smallskipamount}
\setenumerate[0]{leftmargin=*,itemsep=\the\smallskipamount}

\newcommand{\sto}{\rightarrow}
\renewcommand{\to}{%
   \ifbool{@display}{\longrightarrow}{\rightarrow}%
   }
\let\shortmapsto\mapsto
\renewcommand{\mapsto}{%
   \ifbool{@display}{\longmapsto}{\shortmapsto}%
   }
\newcommand{\hooklongrightarrow}{\mathrel{\mkern 0.5mu\lhook\mkern -3.5mu\relbar\mkern -3mu \rightarrow }}
\newcommand{\inj}{%
   \ifbool{@display}{\hooklongrightarrow}{\hookrightarrow}
   }
\newcommand{\isoarrow}{%
   \ifbool{@display}{\overset{\sim}{\longrightarrow}}{\xrightarrow\sim}%
   }
\newlength{\olen}
\newlength{\ulen}
\newlength{\xlen}
\newcommand{\xra}[2][]{%
   \ifbool{@display}%
      {\settowidth{\olen}{$\overset{#2}{\longrightarrow}$}%
       \settowidth{\ulen}{$\underset{#1}{\longrightarrow}$}%
       \settowidth{\xlen}{$\xrightarrow[#1]{#2}$}%
       \ifdimgreater{\olen}{\xlen}%
          {\underset{#1}{\overset{#2}{\longrightarrow}}}%
          {\ifdimgreater{\ulen}{\xlen}%
             {\underset{#1}{\overset{#2}{\longrightarrow}}}
             {\xrightarrow[#1]{#2}}}}%
      {\xrightarrow[#1]{#2}}
   }
\makeatother
\newcommand{\xyra}[2][]{%
   \settowidth{\xlen}{$\xrightarrow[#1]{#2}$}%
   \ifbool{@display}%
      {\settowidth{\olen}{$\overset{#2}{\longrightarrow}$}%
       \settowidth{\ulen}{$\underset{#1}{\longrightarrow}$}%
       \ifdimgreater{\olen}{\xlen}%
          {\mathrel{\xymatrix@M=.12ex@C=3.2ex{\ar[r]^-{#2}_-{#1} &}}}%
          {\ifdimgreater{\ulen}{\xlen}%
             {\mathrel{\xymatrix@M=.12ex@C=3.2ex{\ar[r]^-{#2}_-{#1} &}}}
             {\mathrel{\xymatrix@M=.12ex@C=\the\xlen{\ar[r]^-{#2}_-{#1} &}}}}}%
      {\mathrel{\xymatrix@M=.12ex@C=\the\xlen{\ar[r]^-{#2}_-{#1} &}}}%
   }
\makeatletter
\newcommand{\xla}[2][]{%
   \ifbool{@display}%
      {\settowidth{\olen}{$\overset{#2}{\longleftarrow}$}%
       \settowidth{\ulen}{$\underset{#1}{\longleftarrow}$}%
       \settowidth{\xlen}{$\xleftarrow[#1]{#2}$}%
       \ifdimgreater{\olen}{\xlen}%
          {\underset{#1}{\overset{#2}{\longleftarrow}}}%
          {\ifdimgreater{\ulen}{\xlen}%
             {\underset{#1}{\overset{#2}{\longleftarrow}}}
             {\xleftarrow[#1]{#2}}}}%
      {\xleftarrow[#1]{#2}}
   }
\renewcommand{\lra}{%
   \ifbool{@display}{\longleftrightarrow}{\leftrightarrow}%
   }
\newcommand{\undertilde}{\raisebox{0.4ex}{\smash[t]{$\scriptstyle\sim$}}}

   \def\XXint#1#2#3{{\setbox0=\hbox{$#1{#2#3}{\int}$}
        \vcenter{\hbox{$#2#3$}}\kern-.5\wd0}}

\usepackage[msc-links, alphabetic]{amsrefs}

\newcommand{\nek}{Nekov{\'a}{\v{r}}}

\usepackage{cases}

\usepackage{blindtext}
\makeatletter
\renewcommand\part{%
   \if@noskipsec \leavevmode \fi
   \par
   \addvspace{4ex}%
   \@afterindentfalse
   \secdef\@part\@spart}

     \def\part{\@startsection{part}{1}%
  \z@{3\linespacing\@plus\linespacing}{2\linespacing}%
  {\centering\normalfont\large\bfseries}}
\makeatother

\title[Gan--Gross--Prasad cycles and derivatives of $p$-adic $L$-functions]{Gan--Gross--Prasad cycles\\and derivatives of $p$-adic $L$-functions}
\author{Daniel Disegni}
\address{Department of Mathematics, Ben-Gurion University of the Negev, Be'er Sheva 84105, Israel}
\address{Aix-Marseille University, CNRS, I2M - Institut de Math\'ematiques de Marseille, campus de Luminy, 13288 Marseille, France}
\address{Institut universitaire de France (IUF)}
\email{daniel.disegni@univ-amu.fr}
\author{Wei Zhang}
\address{Massachusetts Institute of Technology, Department of Mathematics, 77 Massachusetts Avenue, Cambridge, MA 02139, USA}
\email{weizhang@mit.edu}
\thanks{This work was supported by BSF  grants  2018250 and 2022240. D.D. has been supported by ISF grant 1963/20 and ANR grant TIWA-C-SYMPLE. W.Z. has been supported by NSF DMS $\#$1901642, $\#$2401548, and the Simons foundation. }

\begin{document}

\maketitle

\begin{center}
   \small \itshape To the memory of Dick Gross
\end{center}

\begin{abstract} 
We study the $p$-adic analogue of the arithmetic Gan--Gross--Prasad (GGP) conjectures for unitary groups. Let $\Pi$ be a conjugate-selfdual cuspidal automorphic representation of  $\GL_{n}\times\GL_{n+1}$ over a CM field, which is algebraic of  minimal regular weight at infinity. We first show the rationality of  twists of the  ratio of $L$-values of $\Pi$ appearing in the GGP conjectures.  Then, when $\Pi$ is  $p$-ordinary at a prime $p$, we construct a cyclotomic $p$-adic $L$-function $\sL_{p}(\RM_{\Pi})$ interpolating those twists. Finally, under some local assumptions, we prove a precise formula  relating the first derivative  of $\sL_{p}(\RM_{\Pi})$  to the $p$-adic heights of Selmer classes arising from arithmetic diagonal cycles on unitary Shimura varieties.  We deduce applications to the $p$-adic Beilinson--Bloch--Kato conjecture for the motive attached to $\Pi$. All proofs are based on some relative-trace formulas in $p$-adic coefficients. 
\end{abstract}

\tableofcontents

\section{Introduction}\label{s:intro}

The pioneering formulas of Gross--Zagier and  Perrin-Riou,  \cite{GZ, PR}, revealed a remarkable relation between Heegner points and  derivatives of complex and $p$-adic $L$-functions. They  had immediate applications to the (classical and $p$-adic) Birch and Swinnerton-Dyer conjectures,  soon strengthened by the Selmer-group bounds proved by Kolyvagin \cite{koly}.

A ``furtive caress''\footnote{Words borrowed from \cite{Weil}.} between those  formulas 
and  one by Waldspurger  on central $L$-values,  \cite{wald-f}, 
did not escape Gross;
 and in \cite{gross-msri}, he blessed it into a
  representation-theoretic  marriage,  which would blossom in \cite{yzz} (and later $p$-adically in  \cite{D17}).

   The seeds for a new generation were sown in a  paper 
by  Gan, Gross,  and Prasad \cite{GGP}.
Their  influential  work  conjectured a pair of non-vanishing criteria in the context  of  embeddings of unitary groups:
one for automorphic periods, in  terms of Rankin--Selberg $L$-values (generalizing \cite{wald-f}); and one for algebraic cycles in Shimura varieties, in terms of  (complex) $L$-derivatives (the \emph{arithmetic} GGP conjecture, generalizing  \cite{GZ}).

 The conjecture on automorphic periods was refined to an exact formula by Ichino--Ikeda and N. Harris \cite{II, nharris},   and recently proved in this form in  \cite{isolation, BPCZ}. On the other hand, despite considerable progress (see \cite{Z22} for a review), the arithmetic GGP conjecture remains open outside of  cases where it can be reduced to Heegner points \cite{yzz, xue}.
 %\footnote{ The analogous conjecture for orthogonal  groups is treated for $1$-cycles in threefolds \cite{yzz tp}.}

\medskip

The purpose of this work is to formulate and, under  some local assumptions, prove a $p$-adic variant of the arithmetic GGP conjecture.
The result in fact takes the form of a
  precise formula,
  in the spirit of  \cite{PR, D17, II, nharris}.
It has immediate applications to the $p$-adic Beilinson--Bloch--Kato conjecture for the relevant motives, which can be further strengthened by the  Selmer bounds recently established in \cite{LTXZZ, Lai-Ski}. 

 (Indeed, one advantage of working in $p$-adic rather than archimedean coefficients is  that we obtain a non-vanishing criterion in Selmer groups, rather than Chow groups: while the {$p$-adic Abel--Jacobi} map from the latter to the former should be injective,  this is not known beyond {cycles of codimension one}.)

\medskip

In the rest of this introduction, we state our main results, discuss their history and context,  and give some ideas on the proofs.

 In \S~\ref{intro1}, we describe our $p$-adic $L$-function (Theorem \ref{thm B}), preceded by a rationality result for twisted Rankin--Selberg $L$-values (Theorem \ref{thm A}) that should
 be of independent interest. 
 
In \S~\ref{intro2} we state our applications to the  $p$-adic Beilinson--Bloch--Kato conjecture (Theorem \ref{th BBK}; the order of presentation is dictated by ease of exposition rather than logic). In  \S~\ref{intro2.5}  we  define the Gan--Gross--Prasad cycles and state our formula for their $p$-adic heights (Theorem \ref{main thm}).

     In \S~\ref{intro3}, we give a sketch of our methods: inspired by the strategy proposed by  Jacquet--Rallis for the Ichino--Ikeda conjecture \cite{JR},  and by one of us \cite{wei-afl} for the arithmetic GGP conjecture (in archimedean coefficients), we construct a $p$-adic relative-trace formula from which we extract the $p$-adic $L$-function; then, we compare it to another relative-trace formula encoding the $p$-adic heights of GGP cycles.

\subsection{The $p$-adic $L$-function}
 \lb{intro1}
Let $F_{0}$ be a number field, and  denote by $\A$ the ad\`eles of $F_{0}$, by $D_{F_{0}}=\prod_{v\nmid \infty}D_{F_{0, v}}$  the discriminant of $F_{0}$ (here $D_{F_{0,v}}$ is the norm of the different ideal of $F_{0,v}$). Let $F$ be a quadratic extension of $F_{0}$, let  ${\rm c}\in \Gal(F/F_{0})$  be the nontrivial element, and let $\eta\colon F_{0}^{\ts}\bs\bA^{\ts}\to \{\pm 1\}$ be the associated quadratic character. Define a reductive group over $F_{0}$ by 
  $$\RG'\coloneqq  ( \Res_{F/F_{0}} \GL_{n}\ts \Res_{F/F_{0}} \GL_{n+1}) / (\GL_{1}\ts\GL_{1}),$$
where  $\GL_{1}\ts \GL_{1}$ is  the subgroup of  pairs of scalar matrices over $F_{0}$.
Let $\Pi=\Pi_{n}\boxtimes \Pi_{n+1}$ be an (irreducible) automorphic representation of $\RG'(\bA)$. We attach to $\Pi$ a 
Rankin--Selberg $L$-function
 twisted by a character $\chi$ of $F_0^{\ts}\bs \A^{\ts}$ and a product of  Asai $L$-functions\footnote{Throughout the introduction (but  \emph{differently} from the rest of the paper) $L$-functions do not include archimedean factors. See  \cite[{\S~7}]{GGP} for the definition of $ L(s, \Pi_{\nu}, \As^{\pm})$.} 
by
\beqq
L(s, \Pi\ot\chi)& \coloneqq  L(s, \Pi_{n}\ts(\Pi_{n+1}\ot\chi \circ{\rm Nm}_{F/F_{0}})), \\
  L(s, \Pi, \As^{\star}) &\coloneqq  L(s, \Pi_{n}, \As^{(-1)^{n}}) L(s, \Pi_{n+1}, \As^{(-1)^{n+1}});\eeqq
the analogous definition applies to the local factors at places $v$ of $F_{0}$.

Following \cite{isolation}, we say that a cuspidal automorphic representation $\Pi_\nu$ of $\GL_{\nu}(\A)$ (for a positive integer~$\nu$) is \emph{hermitian} if $\Pi_{\nu}\circ {\rm c}\cong \Pi_{\nu}^{\vee}$ and $L(s, \Pi_\nu, \As^{(-1)^\nu})$ is regular at $s=1$. We say that $\Pi=\Pi_{n}\boxtimes \Pi_{n+1}$ is hermitian if $\Pi_{n}$, $\Pi_{n+1}$ are. For such a representation $\Pi$, we define
\beq \lb{sL def}
\sL(s, {\Pi_{v},\chi_{v}})\coloneqq  D_{F_{0,v}}^{n+1}\,  \prod_{i=1}^{n+1} L(i, \eta_{v}^{i}) \cdot  {L(s , \Pi_{v}\ot\chi_{v})\over {L(1, \Pi_{v}, {\rm As}^{\star})}}, \eeq
and
$$\sL(s, {\Pi,\chi})\coloneqq  \prod_{v\nmid \infty}\sL(s, {\Pi_{v},\chi_{v}}).$$
Here,  the abelian factor may be interpreted in terms of  $L$-values of motives of unitary groups (\S~\ref{sec: tam}). 

\subsubsection{Rationality of $\sL$} Assume from now on that $F_{0}$ is totally real and $F$ is CM. Let ${\rm arg}(z)\coloneqq {z/ |z|}$ (a character of $\bC^{\ts}$),  let $\Pi_{\nu, \R}^{\circ}$ be the representation of $\GL_{\nu}(\C)/\GL_{1}(\R)$ induced by the character ${\rm arg}^{\nu -1}\ot {\rm arg}^{\nu -3}\ot\ldots\ot{\rm arg}^{1-\nu}$ of the torus $(\C^{\ts})^{\nu}$,  and define the representation  
$$\Pi_{\infty}^{\circ}=\bigotimes_{v\vert \infty}\Pi_{\R}^{\circ}\coloneqq  \bigotimes_{v\vert\infty} \Pi^{\circ}_{n,\R} \ot \Pi^{\circ}_{n+1, \R}$$
of $\G'(F_{0, \infty})$.  Let us also denote by $\one_{\infty}$ the trivial representation of $\G'(F_{0, \infty})$ over $\Q$. 

Let $\Pi=\Pi^{\infty}\ot \one_{\infty}$ be a representation of $\G'(\A)$ on a characteristic-zero field $L$ (admitting embeddings into $\C$). 
We say that $\Pi$ is a \emph{trivial-weight (algebraic)  cuspidal automorphic representation} if for every $\iota\colon L\into \C$, the representation $\Pi^{\iota}\coloneqq \iota\Pi^{\infty}\ot\Pi_{\infty}^{\circ}$ is an (irreducible) cuspidal automorphic representation of $\G'(\A)$. (It is known that every cuspidal automorphic representation of $\G'(\A)$ over $\C$ such that $\Pi_{\infty}\cong \Pi_{\infty}^{\circ}$ arises in this manner for some number field $L$.)  We say that $\Pi$ is hermitian if $\Pi^{\iota}$ is for some (equivalently, every)~$\iota$.%[[is this independent of $\iota$??? write in text]]

We first prove the following strong rationality property for  the values of $\sL$. For an ideal $m\subset \sO_{F_{0}}$, let $Y(m)_{/\Q}$ be the finite \'etale scheme of characters of $F_{0}^{\ts}\bs \A^{\ts}/F_{0,\infty}^{\ts} (\widehat{\sO}^{\ts}_{F_{0}} \cap 1+m\widehat{\sO}_{F_{0}})$. Let $Y\coloneqq \varinjlim_{m} Y({m})$, the ind-finite scheme over $\Q$ of locally constant characters of  $F_{0}^{\ts}\bs \A^{ \ts}/F_{0,\infty}^{\ts} $.
\begin{theoA} \lb{thm A}Let $\Pi$ be a trivial-weight hermitian cuspidal automorphic representation of $\G'(\A)$ defined over a characteristic-zero field $L$. 
 Then there is an element
\beq \lb{L rat}
\sL(\RM_{\Pi}, \cdot)\in\sO(Y_{L}) .\eeq
such that 
$$ \sL(\RM_{\Pi}, \chi) = {\sL^{}(1/2, \Pi^{\iota}, \chi)\over \ep({1\over 2}, \chi^{2})^{{n+1\choose 2}}}$$
for all $\chi\in  Y_{L}(\C)$ with underlying embedding $\iota\colon L\into \C$. 
\end{theoA}
For the notation `$\RM_{\Pi}$', see Remark \ref{not MPi}.

\begin{remark} For $n=1$, Theorem \ref{thm A} is a variant of a classical result of Shimura \cite{shimura}. 
A conditional proof of the rationality of $\sL(1/2, \Pi, \one)$ for a more general  class of $\Pi$ was   obtained  by Grobner and Lin \cite[Theorem C]{Gro-Lin}. (In fact, their rationality result is also a consequence of the Ichino--Ikeda conjecture, but the method of \cite{Gro-Lin} is different.)  See also \cite{Rag} for a related result, and \cite{GHL} for relations to Deligne's conjecture. Finally, a very general and complete rationality result for Rankin--Selberg $L$-values was recently proved by Li--Liu--Sun, \cite{LLS}.
\end{remark}

\subsubsection{$p$-adic interpolation} \lb{p interp}
Fix from now on a rational prime $p$. For $v\vert p$  a place of $F_{0}$, let  $N^{\circ}_{v}\subset G'_{v}\coloneqq \G'(F_{0, v})$ be the subgroup of  integral unipotent upper-triangular matrices, and let $T^{+}_{v}\subset G_{v}'$ be the monoid of diagonal matrices such that $tN_{v}^{\circ} t^{-1}\subset N_{v}^{\circ}$. 
Let $\Pi$ be a trivial-weight cuspidal  automorphic representation of $\G'(\A)$ over a finite extension $L$ of $\Q_{p}$. We say that $\Pi$ is \emph{$v$-ordinary} if $\Pi_{v}^{N_{v}^{\circ}}$ contains a nonzero vector (necessarily unique up to scalar multiple)
  on which all the operators $U_{t,v} \coloneqq  \sum_{x\in N^{\circ}_{v}/ t N^{\circ}_{v} t^{-1}} [x t]$, for $t\in T_{v}^{+}$, act by units in $\sO_{L}$. We say that $\Pi$ is \emph{ordinary} if it is $v$-ordinary for all $v\vert p$.

For any number field $E$, denote  $\Gamma_{E} \coloneqq  E^{\ts}\bs \A^{\infty \ts}_{E}/ \prod_{w\nmid p} \sO_{E, w}^{\ts}$,
 and let  $$\sY\coloneqq  \Spec \bZ_{p}\llb\Gamma_{F_{0}}\rrb\ot_{\bZ_{p}}\Q_{p}.
$$ 
 We have a natural map $Y(p^{\infty})_{\Q_{p}}\coloneqq \varinjlim_{r}Y(p^{r})_{\Q_{p}} \into \sY$.

If $L'/L$ is a field extension and $S/L$ is an (ind-) scheme, we denote $S_{L'}\coloneqq S\ts_{\Spec L}\Spec L'$. 

\begin{theoA} \lb{thm B} Let $\Pi$ be an ordinary, hermitian, trivial-weight cuspidal automorphic representation of $\G'(\A)$ over  a finite extension $L$ of $\Q_{p}$.  Assume that for each place $v\vert p$ of $F_{0}$ that does not split in $F$, the representation $\Pi_{v}$ is unramified.

There exists a unique function 
$$\sL_{p}(\RM_{\Pi})\in \sO(\sY_{L}) $$
whose  restriction to $Y(p^{\infty})_{L}$ satisfies 
\beqq\lb{eq interp Lp} 
\sL_{p}(\RM_{\Pi})(\chi) =  e_{p}(\RM_{ \Pi\ot \chi})\,  \sL(\RM_{\Pi}, \chi)
\eeqq
where $\sL(\RM_{\Pi})$ is  as in \eqref{L rat}, and  $e_{p}(\RM_{ \Pi\ot \chi}) = \prod_{v\vert p} e(\Pi_{v}, \chi_{v})$ is the product of the  explicit local terms \eqref{e v ord}.
\end{theoA}
\begin{remark} We conjecture that the theorem remains true without the non-ramification condition at nonsplit $p$-adic places. 
\end{remark}
\begin{remark} \lb{Sun}
We say that $\Pi$ is \emph{non-exceptional} if $e_{p}(\RM_{\Pi})\neq0$. By a recent result of Liu and Sun (Proposition \ref{hyp ep}), the factor  $e_{p}(\RM_{\Pi\ot \chi})$ is as conjectured by Coates and Perrin-Riou \cite{coates}; this implies that  
if $\Pi_{v}$ is an irreducible principal series for all $v\vert p$, then $\Pi$ is non-exceptional (see Remark \ref{unr nexc}).
\end{remark}

\begin{remark} Januszewski \cite{janu16} has proven a variant of Theorem \ref{thm B} in a more general context, by the method of modular symbols, recently   substantially improved by Liu--Sun \cite{LS}. Our method is  similar locally at $p$ but very different globally (and at archimedean places), see \S~\ref{sec: rtf an intro} below. \end{remark}

\begin{remark}
Other authors have studied the variation of the above $L$-values (and in fact their `square roots')  in anticyclotomic or more general selfdual $p$-adic families, see \cite{HY, Liu-Ber, liu bessel,  Xenia}.  It is of course  expected that these values can be interpolated  into a function over the entire ordinary deformation space, that specializes to the functions of these works in selfdual subspaces and to  our $\sL_{p}(\RM_{\Pi})$ in the cyclotomic direction. (The case of `two' abelian variables explicitly conjectured in \cite[Hypothesis 7.14]{Liu-Ber} could be achieved by the method of this paper, but we chose not to address it in order to bound the technical aspects.)
\end{remark}

\subsection{On the $p$-adic Beilinson--Bloch--Kato conjectures for Rankin--Selberg motives}   \lb{intro2}
Before moving to discuss our main result, we give its main arithmetic application, which can be stated without much further background.

Let 
$$\Pi = \Pi_{n}\boxtimes \Pi_{n+1}$$
 be a hermitian trivial-weight cuspidal automorphic representation of $\G'(\A)$ over  a finite extension $L$ of $\Q_{p}$.  Denote  by $G_{F}$ the absolute  Galois group of $F$, by $\ol{\Q}_{p}$ an algebraic closure of $L$ and let  $ \rho_{\Pi_{\nu}, \ol{\Q}_{p}} \colon G_{F}\to \GL_\nu({\ol{\Q}_{p}})$ be the semisimple representation attached to $\Pi_{\nu}$  by the global Langlands correspondence (as described in \cite[Theorem 1.1]{Ca}).  Assuming that $\ep(\Pi)\coloneqq \ep(1/2, \Pi_{n}^{\iota}\ts \Pi_{n+1}^{\iota})=-1$ for any (equivalently, all) $\iota\colon L\into \C$, 
 we construct a continuous representation
 \beq \lb{rho Pi intro}
\rho_{\Pi}\colon G_{F}\to  \GL_{n(n+1)}(L) \eeq
whose base-change $\rho_{\Pi}\ot_{L}{\ol{\Q}_{p}}$ is isomorphic, up to semisimplication, to 
 $ \rho_{\Pi_{n},{\ol{\Q}_{p}}}\ot   \rho_{\Pi_{n+1}, \ol{\Q}_{p}}(n)$ (Remark \ref{construct rho}).
It satisfies $\rho_{\Pi}^{\rm c}\cong\rho_{\Pi}^{*}(1)$, where  $\rho^{\rm c}(g)\coloneqq \rho (c^{-1}gc) $ for any lift $c\in G_{F}$ of ${\rm c}$.

The Beilinson--Bloch--Kato (BBK) conjecture relates the dimension of the Bloch--Kato Selmer group $$H^{1}_{f}(F, \rho_{\Pi})$$ to the order of vanishing of $\sL(s,\Pi^{\iota})$ at $s=1/2$, for any $\iota\colon L\into\C$. Assuming that  $\Pi$ is ordinary, we can consider a variant in terms of $${\rm ord}_{\chi=\one}\sL_{p}(\RM_{\Pi})\coloneqq \sup\, \{r \ |\  \sL_{p}(\RM_{\Pi})\in \mathfrak{m}_{\one}^{r}\subset \sO(\sY_{L})\},$$ where ${\frak m}_{\one}$ is the ideal of functions vanishing at $\chi=\one$. We prove the following.

\begin{theoA} \lb{th BBK}
Let $\Pi$ be an ordinary, hermitian, trivial-weight cuspidal automorphic representation of $\G'(\A)$ over  a finite extension $L$ of $\Q_{p}$. Assume that $\ep(\Pi)=-1$, and that the following further conditions are satisfied:
 \begin{itemize}
 \item $F/F_{0}$ is unramified; in particular, $F_{0}\neq \Q$;
\item  all places $v|2$ are split in $F/F_0$;
\item $p>2n$ if $n>1$;
 \item for every place $v\vert p$ of $F_{0}$, we have that $v$ splits in $F$ and $\Pi_{v}$ is unramified;
\item  for every place $v$ of $F_{0}$ that  splits in $F$, at least one of $\Pi_{n,v}$ and $\Pi_{n+1,v}$ is unramified;
\item for every place $v$ of $F_{0}$ that is inert in $F$, each of $\Pi_{n,v}$, $\Pi_{n+1,v}$ has conductor at most~$1$ and trivial central character,\footnote{Note that  if  $\Pi_{\nu, v}$  has  conductor~$0$ (i.e., it is unramified) and it is the local component at the inert place $v$ of a hermitian representation, then its central character is necessarily trivial. The notion of conductor is reviewed in \S~\ref{sss:cond01}.} and moreover if $\Pi_{n, v}$ has conductor~$1$ then $\Pi_{n+1,v}$ also has conductor~$1$.
 \end{itemize}
 Then 
\beq \lb{BBK1} {\rm ord}_{\chi=\one}\sL_{p}(\RM_{\Pi})=1\ \Longrightarrow\ \dim_{L} H^{1}_{f}(F, \rho_{\Pi})\geq1.\eeq
If moreover $p$ is an admissible prime for $\Pi$ in the sense of \cite[Definition 8.1.1]{LTXZZ}, then 
\beq
 \lb{BBK2}
 {\rm ord}_{\chi=\one}\sL_{p}(\RM_{\Pi})=1\ \Longrightarrow\  \dim_{L} H^{1}_{f}(F, \rho_{\Pi})=1.\eeq
\end{theoA}

For a comment on the reason for the condition on the conductor, see Remark \ref{chao rmk} below.

This result is  a consequence of a non-vanishing criterion for certain explicit elements of $H^{1}_{f}(F, \rho_{\Pi})$ arising as classes of algebraic cycles, which we describe in the rest of this section. The stronger \eqref{BBK2} follows from combining that criterion with the Selmer bounds of  \cite{LTXZZ}.
In particular, in this case we have that $ H^{1}_{f}(F, \rho_{\Pi})$ is generated by the class  of an algebraic cycle -- a result analogous to the finiteness of the $p^{\infty}$-torsion of the Tate--Shafarevich group of an elliptic curve.  

By \cite[Theorem 1.2.4 (3)]{LTX}, if we start with a hermitian, trivial-weight cuspidal automorphic representation $\Pi_{0}$ of $\G'(\A)$  over a number field  $L_{0}$ that satisfies all of the above assumptions except the ones involving $p$, then all but finitely many rational primes are admissible\footnote{The admissibility condition used in \cite{LTX} is slightly different, but the result applies to the admissibility notion of \cite{LTXZZ} thanks to \cite[footnote to Definition 5.12]{LTX} and either \cite[Corollary D.14]{LTXZZ} or \cite[Lemma 5.1.4]{LTX}.}
 for $\Pi$. Alternatively, one may appeal to the Selmer bound of \cite{Lai-Ski} instead of \cite{LTXZZ},  under different  conditions on~$p$.

\medskip

\begin{remark}\lb{history}
The history of theorems of type \eqref{BBK1} 
consists of several works for similar $2$-dimensional Galois representations over CM fields (starting with \cite{PR} and continuing with \cite{N95, kobayashi, shnidman, D17, nonsplit, univ}),
  together with  a very recent result by Y. Liu and one of us for a family of higher-dimensional representations \cite[Theorem 1.7]{DL}. The result \eqref{BBK2} appears to be the first one of its kind for higher-dimensional representations, \emph{ex aequo} with the main result of \cite{dd-euler} building on \cite{DL}; previously, only  $2$-dimensional cases were known, based on  generalizations of \cite{koly}.
  \end{remark}

\begin{remark} \lb{not MPi} Our notation (and the definitions going back to \eqref{sL def}) suggest that one may think of $\sL(\RM_{\Pi})$, $\sL_{p}(\RM_{\Pi})$ as attached to the virtual motive $\RM_{\Pi}$ over $F_{0}$ whose $p$-adic realization is (up to abelian factors) $$\RM_{\Pi, p}\coloneqq  ({\rm Ind}_{G_{F}}^{G_{F_{0}}} \rho_{\Pi})\ominus {\As}^{\star}(\rho_{\Pi}).$$ 
Here, $\As^{\star}(\rho_{\Pi})=\As^{\star}(\rho_{\Pi_{n}}) \oplus \As^{\star}(\rho_{\Pi_{n+1}})$ with the factors  defined by 
\beqq
\As^{\pm}(\rho_{\Pi_{\nu}})\colon G_{F_{0}} &\to \GL(L^{\nu}\ot_{L} L^{\nu}) \\
 G_{F}\ni g &\mapsto \rho_{\Pi_{\nu}}(g) \ot \rho_{\Pi_{\nu}}^{\rm c}(g),\\
c &\mapsto (x\ot y\mapsto \pm y\ot x)
\eeqq
and the sign $\star=(-1)^{\nu}$ on the $\nu$-factor.

Then the $p$-adic BBK conjecture would rather relate   ${\rm ord}_{\chi=\one}\sL_{p}(\RM_{\Pi})$ with 
$$\dim_{L} H^{1}_{f}(F_{0}, {\rm Ind}_{G_{F}}^{G_{F_{0}}}\rho_{\Pi}) - \dim_{L} H^{1}_{f}(F_{0}, \As^{\star}(\Pi)).$$ 
The first term equals $\dim_{L}  H^{1}_{f}(F, \rho_{\Pi})$. Under our assumption that $\Pi$ is hermitian,  $\As^{\star}(\Pi_{\nu})$ coincides with the adjoint representation defined in the opening paragraphs of \cite{NT} (cf. \cite[Proposition 7.4]{GGP}). By the results obtained there and in \cite{thorne-adjoint}, under some irreducibility assumptions on $\rho_{\Pi_{\nu}}$, we have $H^{1}_{f}(F_{0},  \As^{\star}(\rho_{\Pi}))=0$.
\end{remark}

\begin{remark} 
Theorem \ref{th BBK} and Theorem \ref{main thm} below
rely on a decomposition of the tempered part of the cohomology of unitary Shimura varieties (Hypothesis \ref{hyp coh}), which is expected to be proven in a sequel to \cite{KSZ}. (At a more basic level, we also freely use the results of \cite{mok, KMSW} on automorphic representations of unitary groups.) \end{remark}

\begin{remark} 
 For $\Pi$ satisfying the conditions of  Theorem \ref{th BBK}, let  $\Sigma(\Pi)$ be the set of inert places $v$ of $F_0$ such that $\Pi_{\textup{even},v}$ has conductor~$1$, where $\Pi_{\textup{even}}\coloneqq \Pi_{2\lfloor{(n+1)/2}\rfloor}$. Then by Corollary \ref{coro root n all},
we have 
 $$\ep(\Pi)=(-1)^{\vert\Sigma(\Pi)\vert}.$$
 (In particular, if $n$ is even and $\eps(\Pi)=-1$, then $\Sigma(\Pi)$ is non-empty.)
\end{remark}

In the next subsection we describe, after some preliminaries, the construction of the Selmer classes of interest and our formula relating those to the derivative of $\sL_{p}(\RM_{\Pi})$ (Theorem \ref{main thm}).

\subsection{The $p$-adic arithmetic Gan--Gross--Prasad conjecture}   \lb{intro2.5}

The cycles of interest arise from Shimura varieties attached to certain unitary group. We start by describing the representation-theoretic background.

\subsubsection{Incoherent unitary groups and their representations}\lb{sec incoh} 
Fix a totally positive $\nrm\in F_0^\ts$. (For the purposes of the Introduction, the reader may simply assume $\nrm=1$ without loss of generality, see Remark  \ref{indep nrm} below.)

For a place $v$ of $F_{0}$,  denote by $\sV_{v}$ the set of isomorphism classes of pairs $V_{v} = (V_{n,v}, V_{n+1,v})$ of (non-degenerate) $F_{v}/F_{0,v}$-hermitian spaces over $F_{v}$, where $V_{n,v}$  has rank $n$ and $V_{n+1,v} = V_{n,v} \oplus F_{v} u$ where the special vector $u$ has hermitian norm~$\nrm$.
Let $\sV$ be the set of isomorphism classes of pairs $(V_n, V_{n+1})$ of $F/F_0$-hermitian spaces where $V_n$ has rank $n$ and $V_{n+1}=V_n\oplus F u$ with $u$ of hermitian norm $\nrm$; we have obvious localization maps from $\sV$ to $\sV_v$.

 Let $\sV^{\circ}$ be the set of collections $V=(V_{v})_{v}$  with $V_{v}\in \sV_{v}$ such that  $V_{n, v}$ is positive-definite for all archimedean places, and   for all but finitely many places $v$, the Hasse(--Witt) invariant for the pair $V_{v}$ of hermitian spaces defined by  
\beq\lb{eq: HW}\epsilon(V_{v})\coloneqq \epsilon(V_{{\rm even},v})
\eeq
equals $+1$. Here $V_{{\rm even},v}$ denotes $V_{2\lfloor{(n+1)/2}\rfloor,v}$ and the Hasse invariant for a hermitian space $V_{\nu,v}$ is defined by
\beq\lb{eq: Hasse}
\epsilon(V_{\nu,v})\coloneqq\eta_{v}( (-1)^{\nu\choose 2} \det V_{\nu,v}),\quad \nu\geq 1.
\eeq

We say that $V\in \sV^{\circ}$ is \emph{coherent} if 
there exists a (unique up to isomorphism) element still denoted $V\in \sV$, whose $v$-localization is $V_{v}$.
This holds  if and only if $\epsilon (V)\coloneqq \prod_{v} \epsilon(V_{v})$ equals $+1$. (Thus we identify coherent collections in $\sV^\circ$ with a subset of $\sV$.) When $\epsilon (V)$ equals $-1$, we refer to $V$ as an \emph{incoherent} pair of $F/F_{0}$-hermitian spaces. 
For $V\in \sV^{\circ}\cup \sV$, we denote by 
\beq \lb{unit}
\RH_{v}^{V_{v}}\coloneqq  \RU(V_{n,v}) \subset \G_{v}^{V_{v}} \coloneqq  \RU(V_{n,v})\ts  \RU(V_{n+1,v}), \eeq
(where the embedding is diagonal), by $H_{v}^{V_{v}}\subset G_{v}^{V_{v}} $ their $F_{0, v}$-points. When $V$ is coherent, these are localizations of  unitary groups  $\RH^{V}\coloneqq  \RU(V_{n})\into \G^{V}\coloneqq \RU(V_{n})\ts \RU(V_{n+1})$ over $F_{0}$. When $V$ is incoherent, we fix  models $V_{\nu, \A}$ of $\prod_v V_{\nu, v}$ over $\A$ such that $V_{n+1, \A} =V_{n,\A} \oplus \A_F u$. {That is, $V_{\nu,\A}$ is an $\A_F/\A$-hermitian space such that $V_{\nu, \A}\ot_\A \prod_v F_{0, v}\cong \prod_v V_{\nu,v}$; this is equivalent to choosing selfdual hermitian lattices $\Lm_{\nu, v}\subset V_{\nu, v} $ at all but finitely many places, up to equivalences at all but finitely many places. The condition $V_{n+1, \A} =V_{n,\A} \oplus \A_F u$ amounts to requiring that $\Lm_{n+1, v} $ is equivalent to $\Lm_{n, v}\oplus \sO_F u$ for all but finitely many  places $v$.}
We still use the notation
$$\RH^{V}\subset \G^{V}$$ 
for  the collections \eqref{unit}, which we refer to as \emph{incoherent} unitary groups over $F_{0}$, and we denote $\G^{V}(\A^{S})=\prod'_{v\notin S}G_{v}^{V_{v}}$, where the restricted product is with respect to the stablizers of $\Lm_{n, v} \ts \Lm_{n+1, v}$. (Different choices of adèlic models will lead to isomorphic incoherent adèlic groups.)

 In \S~\ref{sec: meas0}, for each $V_{v}\in \sV_{v}$, we fix Haar measures  $dh_{v}$ on $H_{v}=H_{v}^{V_{v}}$ such that (i) if  $v$ is finite, $dh_{v}$ is $\Q$-valued; (ii)  if $v$ is archimedean and $V_{v}$ is positive definite, 
 $\vol(H_{v}, dh_{v})\in \Q^{\ts}$; 
 (iii) if $V\in \sV^{\circ}$ is coherent, $\prod_{v}dh_{v}$ is the Tamagawa measure on $\RH^{V}(\A)$. We also have measures $dg_{v}$ on $G_{v}=G_{v}^{V_{v}}$ with the analogous properties.

\medskip

Suppose that $V\in \sV^{\circ}$ is incoherent. If $v$ is a place of $F_{0}$ non-split in $F$, we let $V(v)\in \sV$ be the `$v$-nearby' coherent pair such that: $V(v)_{w}= V_{w}$ if $w\neq v$; and    $V(v)_{v}\in \sV_{v}$ is the unique element different from $V_{v}$ if $v$ is non-archimedean, and the element such that $V(v)_{n, v}$ has signature $(n-1, 1)$ if $v$ is archimedean.  
For a characteristic-zero field $L$ and an incoherent $\G=\G^{V}$, a \emph{discrete automorphic representation of $\G(\A)$ over $L$ trivial at infinity} is a representation $\pi=\pi^{\infty}\ot \one_{\infty}$ of $\G(\A)$ over $L$,  such that for every $\iota\colon L\into \C$, every $v\vert \infty$, {and some irreducible admissible representation $\pi_{v}^{\circ}$ of $G^{V(v)}_v$ with the same infinitesimal character as the trivial representation},
the complex representation  of $\G^{V(v)}(\A)$
$$\pi^{\iota, (v)}\coloneqq  \iota \pi^{v}\ot \pi_{v}^{\circ}$$
is irreducible, discrete and automorphic. We say that $\pi$ is  cuspidal, tempered (at all finite places), or stable if each $\pi^{\iota, (v)}$ is (where we call $\pi^{\iota, (v)}$ {stable} if its  base-change  to $\G'(\A)$ is cuspidal).
If $\pi$ is cuspidal, tempered, and stable, then the  base-change of $\pi^{\iota, (v)}$ is necessarily of the form $\Pi^{\iota}$ for a (unique) trivial-weight hermitian cuspidal automorphic representation $\Pi$ over $L$, which we call the base-change of $\pi$ and denote by ${\rm BC}(\pi)$.

\subsubsection{Arithmetic diagonal cycles} 
When $V$ is incoherent, we attach to  $\G=\G^{V}$ a tower of Shimura varieties $(X_{K})_{K\subset \G(\A^{\infty})}$ over $F$ of dimension $2n-1$,\footnote{ Note that Shimura {\em varieties} in this paper are generally {\em Deligne--Mumford stacks} unless the level subgroup is small enough. We use this name for historical reasons; stacky phenomena will not  play any significant role.}, and to  $\RH=\RH^{V}$ a tower of Shimura varieties $(Y_{K'})_{K'\subset \RH(\A^{\infty})}$ over $F$ of dimension $n-1$. {(These Shimura varieties are obtained by fixing an auxiliary archimedean place of $F$, and switching the hermitian spaces to the nearby ones; see \S~\ref{ss:AD} and \S~\ref{incoh Sh}.)} They are proper provided that $F_{0}\neq \Q$, a condition that we  assume for the rest of the introduction.

The embedding $\jmath\colon \RH(\A)\to \G(\A)$ induces a morphism of Shimura varieties still denoted by $\jmath$.
Consider the (well-defined) normalized fundamental class $[Y]^{\circ}\coloneqq \lim_{K'}\vol(K')[Y_{K'}]\in\varprojlim_{K'} \Ch^{0}(Y_{K'})_{\Q}$ and  the \emph{arithmetic diagonal cycle} $\jmath_{*}([Y]^{\circ})\in\varprojlim_{K} {\rm Ch}^{n}(X_{K})_{\Q}$ (where $\Ch^{i}(Z)_{\Q}$ denotes the Chow group of codimension-$i$ cycles on $Z$ with rational coefficients). The $p$-adic  absolute cycle class of $\jmath_{*}([Y]^{\circ})$ can be projected to an element
 $$Z\in H^{1}_{f}(F, M^{\rm temp})$$
where 
$$M^{\rm temp}=\varprojlim_{K}  M_{K}^{\rm temp}, \qquad M_{K}^{\rm temp}\coloneqq
H^{2n-1}_{\text{\'et}}(X_{K, \ol{F}}, \Q_{p}(n))^{\rm temp},$$
 with the superscript `temp' referring to the tempered part of cohomology (see \S~\ref{sec: hecke gal coh}).

\subsubsection{Gan--Gross--Prasad cycles} \lb{ggp cycle def}
Let $\pi$ be a stable,\footnote{If $\pi$ is only assumed to be tempered but not stable, we can still define $Z_{\pi}$ with values in the Selmer group of a certain direct summand of~$\rho_\Pi$ (see \S~\ref{ssec: ggp cycle}).}
cuspidal automorphic representation of $\G(\A)$ trivial at infinity, over some finite extension $L$ of $\Q_{p}$;
 let $\Pi={\rm BC}(\pi)$. According to Hypothesis \ref{hyp coh},  there is an injective map 
$$\pi \to \Hom_{\Q_{p}[G_{F}]} ( M^{\rm temp}, \rho_{\Pi}),$$
well-defined uniquely up to scalar multiples.
We identify $\pi$ with the image of this map, and define the \emph{Gan--Gross--Prasad functional}
\beq \lb{Zpi intro}
Z_{\pi}\colon \pi &\to H^{1}_{f}(F, \rho_{\Pi})\\
\phi&\mapsto Z_{\pi}(\phi)\coloneqq  \phi_{*}Z.
\eeq
We call elements in its image \emph{Gan--Gross--Prasad cycles}.

 \subsubsection{The $p$-adic arithmetic Gan--Gross--Prasad  conjecture}
 By construction, we have
$$Z_{\pi}\in 
 \Hom_{\RH^{V}(\A)}(\pi, L)\ot_{L} H^{1}_{f}(F, \rho_{\Pi}).$$ 
 
 The space $\Hom_{\RH^{V}(\A)}(\pi, L)$ is known to be of dimension~$0$ or~$1$; in the latter case,  $\pi$ is  said to be \emph{distinguished}. By the local Gan--Gross--Prasad conjecture proved in \cite{BP-ggp1, BP-ggp2}, for a given representation $\Pi$ over $L$ as in Theorem \ref{thm A}, there exists a unique (up to isomorphism) pair $(V, \pi)$ where $V\in \sV^{\circ}$ and $\pi$ is a tempered representation of $\G^{V}(\A)$ as above that is distinguished and satisfies  ${\rm BC}(\pi)=\Pi$. Moreover, $\pi$ can be defined over $L$, and $V$ is incoherent if and only if $\ep(\Pi)=-1$ (see \S~\ref{sec: ggp rat}). 
 
 The following is a $p$-adic analogue of  the arithmetic Gan--Gross--Prasad conjecture \cite[Conjecture 27.1]{GGP} for unitary groups.
 \begin{conjecture} \lb{conj AGGP} Let $\Pi$ be a representation as in Theorem \ref{thm B}.
Assume that $\ep(\Pi)=-1$ and that $\Pi$ is non-exceptional. Let  $(V,\pi)$ be the unique pair  with $\pi$ tempered and distinguished and ${\rm BC}(\pi)=\Pi$ as in the previous paragraph. The following conditions are equivalent:
 \begin{enumerate}
 \item ${\rm ord}_{\chi=\one}\sL_{p}(\RM_{\Pi})=1$;
  \item $Z_{\pi}\neq 0.$
 \end{enumerate}
 \end{conjecture}

 \begin{remark}\lb{C and D} According to the $p$-adic BBK conjecture, both conditions are also equivalent to
 \begin{enumerate}
\item[{\em (3)}] $\dim_{L}H^{1}_{f}(F, \rho_{\Pi})=1$. 
\end{enumerate}
 The implication {\em (2)} $ \Longrightarrow$ {\em (3)} is \cite[Theorem 1.1.9]{LTXZZ}  or \cite[Theorem 1.4]{Lai-Ski} 
 (each under suitable conditions on $p$; see remark (1) following \cite[Theorem 1.4]{Lai-Ski} for a comparison of the two sets of conditions).
 \end{remark}

 \begin{remark} \lb{indep nrm}
 Conjecture \ref{conj AGGP} is independent of the choice of the norm $\nrm$ of the special vector $u$ that was fixed at the beginning of \S~\ref{sec incoh}. 
Indeed, by rescaling the hermitian form we may always assume $\nrm$ to be $1$ without changing the associated unitary groups (up to isomorphisms); yet it is more convenient to allow a general $\nrm$ when we consider the relevant Shimura varieties.
 \end{remark}

 As a refinement of Conjecture \ref{conj AGGP},  we prove (under some conditions) a formula that `measures' the product $Z_{\pi} \ot Z_{\pi^{\vee}}$ in terms of the derivative of $\sL_{p}(\RM_{\Pi})$; in order to state it, we need to define some pairings.

\subsubsection{Dualities} 
Continue with the setup of \S~\ref{ggp cycle def}. 
Fix a non-degenerate pairing 
$$ \la\, , \, \ra_{\Pi}\colon \rho_{\Pi}\ot_{L}\rho_{\Pi^{\vee}}\to L(1),$$
and for a compact open subgroup $K\subset \G(\A^{\infty})$, let $\la \, , \, \ra_{K} \colon M_{K}^{\rm temp}\ot M_{K}^{\rm temp}\to L(1)$ be the pairing induced by Poincar\'e duality. Then we (well-)define a pairing 
\beq \lb{pair pi}
(\, , \, )_{\pi} \colon \pi\ot \pi^{\vee}\to L\eeq
 by $(\phi, \phi')_{\pi}\coloneqq \vol(K)^{-1}\,   \phi\circ u_{K}( \phi'^{*}(1))$ for  any $K\subset \G(\A^{\infty})$ fixing $\phi$, $\phi'$. Here, $\phi'^{*}(1)\colon \rho_{\Pi^{\vee}}^{*}(1)\to M_{K}^{{\rm temp},*}(1)$ is the transpose, the volume uses the measure $\prod_{v} dg_{v}$, and $u_{K}\colon M_{K}^{{\rm temp},*}(1)\to M^{\rm temp}_{K}$ is the isomorphism induced by $ \la\, , \, \ra_{K}$.

\subsubsection{Invariant functionals}
If $\pi$ is distinguished, there is a canonical generator 
$$\alpha\in \Hom_{\RH^{V}(\A)}(\pi, L)\ot_{L} \Hom_{\RH^{V}(\A)}(\pi^{\vee}, L)$$ defined as follows. 
Pick a factorization $(\, ,  \, )_{\pi}= \prod_{v\nmid \infty}(\, ,  \, )_{\pi_{v}}$, where each factor is a pairing on $\pi_{v}\ot\pi_{v}^{\vee}$. 
Then  $\alpha$ is defined on factorizable elements  $\phi=\ot_{v\nmid \infty}\phi_{v}$, $\phi'=\ot_{v\nmid \infty}\phi'_{v}$
by the product of absolutely convergent integrals
\beq \lb{alpha intro}
\iota\alpha(\phi, \phi^{'})\coloneqq  \vol(H_{\infty}^{V}, dh_{\infty})\cdot \prod_{v\nmid \infty} \sL(1/2, \iota\Pi_{v})^{-1} \int_{H_{v}} \iota(\pi(h)\phi, \phi')_{\pi}\, dh_{v},\eeq
where  $\iota\colon L\into \C$ is any embedding, $\vol(H_{\infty}^{V}, dh_{\infty})=\prod_{v\vert \infty}\vol(H_{v}^{V_{v}}, dh_{v})\in \Q^{\ts}$, 
and almost all factors are equal to~$1$.

\subsubsection{$p$-adic heights and main result} \lb{ht intro} Assume that $\Pi$ is ordinary. Then $\rho_{\Pi}$ is Panchishkin-ordinary in the sense of \cite{N93} (recalled in \S~\ref{pair def}). By \nek's theory (see \cite{N93} or \S~\ref{sec ht}), the pairing  $\la\, , \, \ra_{\Pi}$ and the natural projection $\lm\colon \Gamma_{F}\to \Gamma_{F_{0}}$ induce a height pairing 
$$h_{\pi}\colon  H^{1}_{f}(F, \rho_{\Pi}) \ot_{L}  H^{1}_{f}(F, \rho_{\Pi^{\vee}})\to \Gamma_{F_{0}}\hat{\ot}L.$$

For $\sL\in \sO(\sY)_{L}$, set $$\partial \sL\coloneqq  [\sL-\sL(1)]\in T^{*}_{\one}\sY_{L}=\frak{m}_{\one}/\frak{m}_{\one}^{2}\ot_{\Q_{p}} L =\Gamma_{F_{0}}\hat{\ot} L.$$
The following is a $p$-adic analogue of  the refined arithmetic Gan--Gross--Prasad conjecture (cf. \cite[Conjecture 5.1]{xue}), in the spirit of the Ichino--Ikeda refinement of the usual Gan--Gross--Prasad conjecture. The case $n=1$ is essentially equivalent to the $p$-adic Gross--Zagier formula as in \cite{D17}.

 \begin{conjecture} \lb{conj AGGP II} 
  Let $V\in \sV^{\circ}$ be an incoherent pair, and let $\pi$ be a distinguished, stable,  ordinary, cuspidal automorphic representation of $\G^{V}(\A)$,  trivial at infinity, over a finite extension $L$ of $\Q_{p}$. Let $\Pi\coloneqq {\rm BC}(\pi)$ and assume that it is ordinary and non-exceptional. 
   Then  for all $\phi\in \pi$, $\phi'\in \pi^{\vee},$ we have
 $$ 
h_{\pi} \left(Z_{\pi}(\phi), Z_{\pi^{\vee}}(\phi')\right)
 =  e_{p}(\RM_{\Pi})^{-1} \cdot {1\over 4} \partial \sL_{p}  (\RM_{\Pi}) 
\cdot \alpha(\phi, \phi')$$
in $\Gamma_{F_{0}}\hat{\ot}L$. 
 \end{conjecture}
\begin{remark}
This conjecture implies the direction $(1)\Longrightarrow (2)$ in Conjecture \ref{conj AGGP}; the converse implication is reduced to the conjectural non-degeneracy of $h_{\pi}$. \end{remark}
 
We have the following theorem, confirming the above refined conjecture in certain cases.
\begin{theoA}\lb{main thm}
Conjecture \ref{conj AGGP II} holds if we  further assume that:
 \begin{itemize}
 \item $F/F_{0}$ is unramified; in particular, $F_{0}\neq \Q$;

\item  all places $v|2$ are split in $F/F_0$;
\item $p>2n$ if $n>1$;
 \item for every place $v\vert p$ of $F_{0}$, we have that $v$ splits in $F$ and $\pi_{v}$ is unramified;

\item  for every finite place $v$ of $F_0$ that  splits in $F/F_0$, at least one of $\pi_{n,v}$ and $\pi_{n+1,v}$ is unramified;
\item  for every finite place $v$ of $F_0$ that is inert in $F/F_0$, one of the following conditions holds:
\begin{itemize}
{\item $\e(V_v)=1$  and $\pi_{n,v}$, $\pi_{n+1,v}$ are both unramified;
    \item $\e(V_v)=(-1)^n$ and $\pi_{n,v}$ is unramified, $\pi_{n+1,v}$ is  almost unramified and  not unramified;}
    \item $\e(V_v)=-1$ and the representations $\pi_{n,v}$,  $\pi_{n+1,v}$ are both almost unramified but not unramified.
\end{itemize}

 \end{itemize}
\end{theoA}

 Here, for an inert finite place $v$ of $F_0$ and a $\nu$-dimensional $F_v/F_{0,v}$-hermitian space $V_{\nu,v}$, a representation of  $\U(V_{\nu,v})$  is called {\em almost unramified} if it has a non-zero vector fixed by the stabilizer of a vertex lattice of type $1$ or $\nu-1$; see \cite{Liu} for the case where $\nu$ is even and \S~\ref{sss:almostun} in general.

\begin{remark}
Besides the $p$-adic Gross--Zagier results mentioned in Remark \ref{history}, the only other $p$-adic height formula in the literature is the recent \cite[Theorem 1.8]{DL}.
 While our setup and global approach to the proof are different,  a  theorem on $p$-local heights in \cite{DL} is essential for us.
\end{remark}

\subsection{$p$-adic relative trace formulas and the  proofs} \lb{intro3}
Our approach to Theorem \ref{main thm}  is based on the comparison of a pair of  relative-trace formulas with $p$-adic coefficients, analogously to  the approach proposed by one of us \cite{wei-afl} over archimedean coefficients. 
In fact,  Theorem \ref{thm A} and Theorem \ref{thm B} are also proved by constructing rational and $p$-adic relative-trace formulas. 
We give a brief overview; unexplained terminology will be defined in the main body of the paper.

\subsubsection{Rationality} \lb{intro rat} Let us first explain the proof of Theorem \ref{thm A}.  For each $\chi\in Y(\C)$, we have a Jacquet--Rallis relative-trace distribution  
$$I(-, \chi)\colon \sH(\G'(\A) , \C)\to \C$$
on the Hecke algebra for $\G'$.
For a `regular' $f'\in \sH(\G'(\A), \C)$, it admits a spectral and a geometric expansion
\beq\lb{rat rtf intro}
\sum_{\Pi}{1\over 4} \sL(1/2, \Pi, \chi)\prod_{v}I_{\Pi_{v}} (f'_{v}, \chi_{v}) =I(f',\chi)= \sum_{\gamma\in \RB'_{}(F_{0})} I_{\gamma}(f',\chi),\eeq
where: $\Pi$ ranges over isomorphism classes of cuspidal representations of $\G'(\A)$; the $I_{\Pi_{v}}$ are local relative  characters;  the variety  $\RB'_{/F_{0}}  = \RH'_{1} \backslash \G' /  \RH'_{2}$ for certain reductive subgroups $\RH_{1}', \RH_{2}'\subset \G'$; and the $I_{\gamma}$ are products of local orbital integrals.  

The only possible sources of irrationality in the right-hand side of \eqref{rat rtf intro} are essentially the archimedean orbital integrals. However, there is a particularly well-behaved class of $f'_{\infty}\in \sH(\G'(F_{0, \infty}))$ (and corresponding $f'\in \sH(\G'(\A))$), the so-called (rational)  Gaussians, whose orbital integrals are controlled. Building on \cite{isolation}, we are able to show that for $\Pi$ as in Theorem \ref{thm A}, there exist  $L$-rational Gaussians $f'$ annihilating every automorphic representation of $\G'(\A)$  but   $\Pi$. Moreover, we need to show that one can pick $f'$ to be `regular' (that is, supported on suitably regular elements for the group action  of $\RH'_{1} \times \RH'_{2}$): this could be quickly done by invoking the results of \cite[Appendix A]{wei-fourier}, but we do it in a more explicit way as described in \S~\ref{intro test}.
 Then  the rationality of $\sL(1/2, \Pi, \chi)$ can be deduced from \eqref{rat rtf intro}. 

\subsubsection{$p$-adic analytic distribution} \lb{sec: rtf an intro}
We have a $p$-adic variant of $I(-, \chi)$, that we describe at first in a slightly idealized form.
For any  `convenient' subgroup $K_{p}'\subset \G'(F_{0,p})$, we construct a distribution
$$\sI =\sI_{K_{p}'}\colon \sH(\G'(\A^{p}))^{\circ}_{K'_{p}, \rm  rs, qc} \to \sO(\sY) $$
on  a certain space of  regularly supported, $\Q_{p}$-rational Gaussian elements of the Hecke algebra away from~$p$.  It admits a spectral and a geometric expansion
\beq\lb{rtf an intro}
\sum_{\Pi}
{1\over 4}
\sL_{p}(\RM_{\Pi} , \chi) \prod_{v\nmid p} \sI_{\Pi_{v}} (f_{v}', \chi_{v})=  \sI(f'^{p}, \chi)=
\int_{\RB'_{\rm rs}(F_{0})} \sI_{\gamma}(f'^{p},\chi) \, dI^{\ord}_{\gamma, K_{p}',p}, \eeq
where
$\Pi$ ranges over 
 representations as in Theorem~\ref{thm B} with nontrivial $K_{p}'$-invariants; 
the  $\sI_{\Pi_{v}}$,  $\sI_{\gamma}$ are $\sO(\sY)$-valued relative characters and  orbital integrals, respectively; and finally,    $dI_{-, K_{p}',p}^{\ord}$ is a certain generalized Radon measure on the rational points of  $ \RB'_{\rm rs}\subset \RB'$, the open subvariety of regular semisimple orbits. 
In fact, we construct $\sI$ from its geometric expansion, and prove Theorem \ref{thm B} by extracting  $\sL_{p}(\RM_{\Pi})$ from $\sI$.
\begin{remark}
 This appears to be a new method for constructing $p$-adic $L$-functions. Let us linger on the archimedean input:  while previous works relied on the non-vanishing of zeta integrals for explicit cohomological test vectors (as proved  by Sun in  \cite{Sun}), we use instead the  `spectral matching' property (proved by Beuzart-Plessis \cite{BP-Planch}), which  relates the value of  $I_{\Pi_{v}}$ on Gaussians with  relative characters  of  constant Hecke measures on a definite unitary group, whose computation is trivial.
\end{remark}
  Under some conditions on $K'_{p}$, we can relax the conditions of regularity on $f'^{p}$ by using the recent work of Lu \cite{Lu}; then the orbital integrals corresponding to non-semisimple orbits need an interesting regularization featuring Deligne--Ribet $p$-adic $L$-functions.

 We note that Urban \cite[\S~6]{urban-e} has constructed a $p$-adic  Arthur--Selberg trace formula; it would be interesting to compare or combine our two approaches.

\subsubsection{The derivative} \lb{sec: rtf'  intro}
For suitable inputs $f'^{p}$, we then have a similar expansion for the derivative of $\sI$.  We will be especially interested in those  $f'^{p}$ that  `purely match' an element $f^{p}\in \sH(\G^{V}(\A^{p}))$ for some incoherent $V$, in the following sense. We have a `matching of orbits' map for all places $v$
\beq \lb{match intro}
\ul{\delta} \colon \RB_{\rs}(F_{0,v}) \to  \bigsqcup_{V'_{v} \in \sV_{v}} \RH^{V'_{v}}(F_{0,v})\bs \G^{V'_{v}}(F_{0,v})/ \RH^{V'_{v}}(F_{0,v}) \eeq
with image the set of regular semisimple orbits on the right hand side. 
The matching condition on $f^{p}=\ot_{v\nmid p} f_{v}$, $f'^{p}=\ot_{v\nmid p}f'_{v}$ is that, defining unitary-group   orbital integrals  by
$$J_{\delta}(f_{v})=\int_{\RH^{V_{v}}(F_{0,v})^{2}} f_{v}(h^{-1} \gamma h')\, dh dh',$$
we should have, for all $v\nmid p$, that
 $I_{\gamma}(f'_{v},\one)=J_{\ul{\delta}(\gamma)}(f_{v})$ if ${\ul{\delta}(\gamma)}$ belongs to the component indexed by $V_{v}$ in  \eqref{match intro}, 
 and $I_{\gamma}(f'_{v},\one)=0$ otherwise.

For such $f'^{p}$, we have $\sI(f'^{p}, \one)=0$ and the  $\Gamma_{F_{0}}\hat{\ot} \Q_{p}$-valued expansions
\beq\lb{rtf' intro}
\sum_{\Pi}
{1\over 4}
 \partial \sL_{p}(\RM_{\Pi} )\prod_{v\nmid p} \sI_{\Pi_{v}} (f_{v}', \chi_{v}) 
= \partial \sI(f'^{p})
= \int_{\RB'_{\rm rs}(F_{0})} \partial \sI_{\gamma}({f}'^p) \, dI_{\gamma, K_{p}', p}^{\ord}\eeq
for the derivative. Moreover
$$\partial \sI_{\gamma} (f'^{p}) =\sum_{v\nmid p\infty \text{ nonsplit in $F$}}  I_{\gamma}(f'^{v p})\,  \partial \sI_{\gamma}(f'_{v})$$
with  $I_{\gamma}(f'^{v p})=\sI_{\gamma}(f'^{v p}, \one)$. The $v$-component of the sum can be nonzero only if $\gamma$ matches an orbit $\delta$ of $ \RH^{V(v)}(\A^{p})\bs \G^{V(v)}(\A^{p})/  \RH^{V(v)}(\A^{p} )$ for the nearby coherent  pair $V(v)\in \sV^{\circ}$.

In practice, unless $K_{p}'$ is suitably symmetric, we are only able to prove the geometric expansion in \eqref{rtf an intro} after specialization at a $\chi\in Y(p^{\infty})$, and with a generalized Radon measure $I_{-, K_{p}', p}^{\ord}(\chi_{p})$ depending on $\chi_{p}$; nevertheless we can show that \eqref{rtf' intro} still holds with   $I^{\ord}_{-, K_{p}', p}\coloneqq I^{\ord}_{-, K_{p}', p}(\one)$.

\subsubsection{Arithmetic distribution} \lb{intro ar d}  Let $V\in \sV^{\circ}$ be incoherent, $\G=\G^{V}$. For a convenient subgroup $K_{p}\subset \G(F_{0,p})$, we define another  $\Gamma_{F_{0}}\hat{\ot} \Q_{p}$-valued 
distribution on a suitable subset of $\sH(\G(\A^{p}))$ by 
$$\sJ_{K_{p}}(f^{p}) = h(Z^{\ord}_{K_{p}}\,  T(f^{p}), Z_{K_{p}}^{\ord}),$$
where $Z^{\ord}_{K_{p}}$ is an ordinary modification of the arithmetic diagonal cycle in level $K_{p}$, and $h$ is a limit of height pairings on the Selmer group of the  tempered, ordinary part of $ H^{2n-1}(X_{K^{p}K_{p},\ol{F}_{0}}, \Q_{p}(n))$.

 When the cycles have disjoint support on the generic fiber, the $p$-adic height pairing admits an expansion $h=\sum_{v\nmid \infty} h_{v}$ into local height pairings. The disjointness is guaranteed if $f$ has regular support at some place.
 
By results in \cite{DL, LL1}, the local height pairing at a place $v$ away from $p$    
is related  to the arithmetic intersection pairing on a regular $v$-integral model, at least after applying  suitable Hecke correspondences to the cycles, and under some vanishing condition for the absolute cohomology of the model (upon localizing at a non-Eisenstein ideal).  After a base change, for suitable levels  we may use the models constructed in the previous work of Rapoport, Smithling and the second author \cite{RSZ3,RSZ4}; here,  a technical difficulty is to verify the vanishing of cohomology in the case of non-trivial level structure, 
 as required  in order to treat both the places of ramification of $\Pi$ in Theorem \ref{main thm}, and the places fostering the regular support condition.  Once this is settled:
 \begin{itemize}
 \item for \emph{split} places away from $p$, we can show that  the local arithmetic intersection numbers  vanish, by refining an argument of \cite{wei-afl, RSZ3};
  \item for \emph{inert} places $v$ (thus away from $p$), by results in \cite{wei-afl} and \cite{RSZ3}, the local arithmetic intersection numbers admit geometric expansions over the orbits $\delta$ for $V(v)$, whose terms  are products of local orbital integrals $J_{\delta}(f_{v'})$ ($v'\neq v$) and arithmetic intersection numbers $\sJ_{\delta}(f_{v})$ in a certain $v$-adic  Rapoport--Zink space.
 \end{itemize}
On the other hand, the contribution of  \emph{$p$-adic} places
 vanishes: this is proved by     a variant of an argument of Perrin-Riou,
which  in our higher-dimensional case relies  on a recent   foundational result of Y. Liu and the first author in \cite{DL}.

{\begin{figure}[h] \lb{fig k}
$$\includegraphics[width=.9\textwidth]{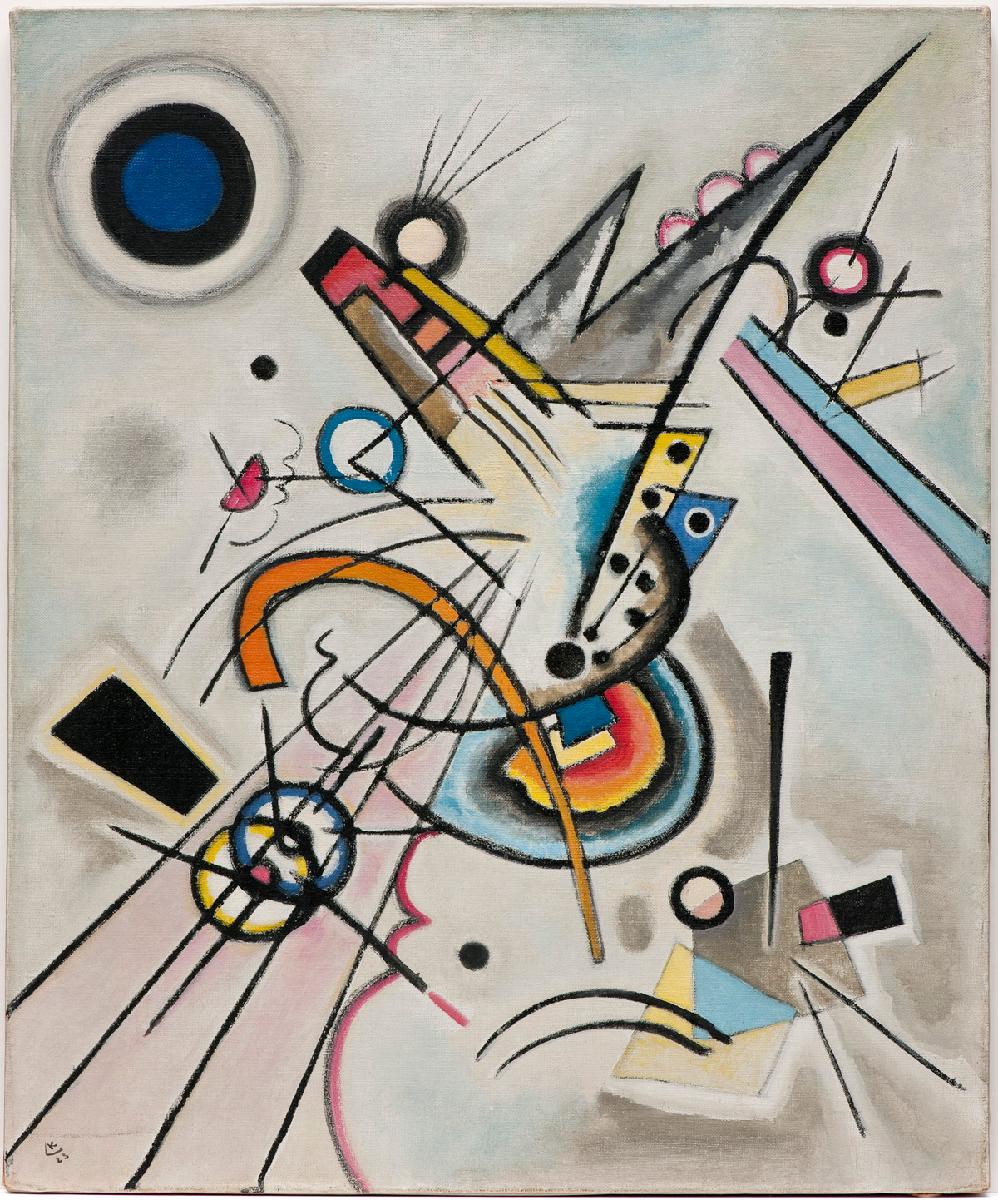} 
$$
\caption{Wassily Kandinsky, \emph{Diagonal}, 1923.
 }
\end{figure} }

 The distribution arising from arithmetic diagonals is thus shown to admit  a spectral and a geometric expansion  (depicted in Figure 1)
\beqq \sum_{\pi}
 \sJ_{\pi}(f^{p})
&= \sJ_{K_{p}} (f^{p})\\
&= \int_{\RB'_{\rs}(F_{0})} \sum_{v\nmid p\infty \ \rm nonsplit}
\one_{V(v)}(\gamma) J^{vp}_{\ul{\delta}(\gamma)}(f^{vp}) \sJ_{\ul{\delta}(\gamma),v}(f_{v}) \, d I^{\ord}_{\gamma,p, K'_{p}},
\eeqq
where:
$\pi$ ranges over
equivalence classes of automorphic representations as in Theorem~\ref{main thm}; the geometric expansion is pulled back to   $\RB_{\rs}'$ via the `matching of orbits' map $\ul{\delta}$, and $\one_{V(v)}$ is the indicator function of those orbits matching one on $\G^{V(v)}$; and finally,  $ d I^{\ord}_{\gamma,p, K'_{p}}$ is as in \eqref{rtf' intro}.

\subsubsection{Comparison}  Theorem \ref{main thm} is eventually deduced from the spectral sides of an equality
\beq\lb{comp intro}
 \sJ_{K_{p}} (f^{p}) = \partial \sI_{K_{p}'}(f'^{p}) \eeq
for suitable  matching $f^{p}$, $f'^{p}$. 

We  prove  \eqref{comp intro} by comparing the geometric expansions. By the definitions of local  matching of Hecke elements (which can be globally assembled thanks to the Fundamental Lemma \cite{Yun, BP-FL}), orbital integrals on either side are the same, thus we are reduced to  identities
\beq\lb{afl intro}
\sJ_{\ul{\delta}(\gamma),v}(f_{v})= \partial \sI_{\gamma}(f'_{v}) 
\eeq
for inert places $v$.
For the spherical $f_{v}'$, $f_{v}$, the identity \eqref{afl intro} is the Arithmetic Fundamental Lemma proposed by one of us \cite{wei-afl} and then proved in \cite{afl-pf1, afl-pf2}; for certain $f_{v}'$, $f_{v}$ of maximal parahoric level,  \eqref{afl intro}  is the arithmetic transfer conjecture  recently proved by Z. Zhang \cite{zhiyu}.

\begin{remark}\lb{chao rmk} We point out the main obstacle to removing the condition of our representations being at worst almost unramified at inert places from our main theorems. The condition  comes from working with Shimura varieties at ``almost selfdual" levels (namely, for vertex-parahoric subgroups of type $1$ or $\nu-1$).
    Although Z. Zhang's result on the  arithmetic transfer conjecture  \cite{zhiyu}
holds in greater generality 
    (at maximal parahoric levels), we can only show the vanishing result for absolute cohomology of integral models alluded to in \S~\ref{intro ar d} in almost selfdual levels,
    see Proposition \ref{thm:van}. The proof of that proposition
      relies on a refined understanding of the cohomology of the irreducible components of the special fiber, which currently seems only available at almost selfdual levels. The generalization to other maximal parahoric levels seems a very interesting yet challenging question.  
      
     Even when restricting to vertex parahorics, our results miss the case where $\pi_{n,v}$ is almost unramified (but not unramified) and $\pi_{n+1,v}$ is unramified. For a discussion of the difficulties, see Remark \ref{rmk dongwen}.
\end{remark}

\subsubsection{Construction of test Gaussians} \lb{intro test}  In order to deduce Theorem \ref{main thm} from the comparison, we need to pick suitable matching $f^{p}$, $f'^{p}$ that annihilate all terms in the spectral expansions but those corresponding to $\pi$, $\Pi$; then  we may use the  comparison in  \cite{BP-Planch, BP-comparison} of the functionals $I_{\Pi_{v}}$ with corresponding ones, $J_{\pi_{v}}$, that are related to the local components of $\alpha=\eqref{alpha intro}$.

The most challenging requirement for the Gaussian $f'^{p}$ is that  the relative character $\ot_{v\nmid p\infty} I_{\Pi_{v}} (f^{'p})$ should not vanish, its non-regular-semisimple orbital integrals should vanish,  while at the same time its level should be controlled in order to allow working with nice integral models on the arithmetic side. This turns out to be a rather hard semi-local problem, which is solved by an explicit construction of a pair of elements $f_{v, \pm}'$ of Iwahori level (to be used at a pair of split non-$p$-adic places), and two explicit local computations: one on the spectral side, which is 
  Proposition \ref{hyp ep} (deduced from a result of Liu--Sun); and one on the geometric side, which is part of Proposition \ref{chi1 along}, whose proof occupies the entire \S~\ref{app orb p}.   It is curious to note that $f_{v, +}'$ also occurs in the construction of the $p$-adic relative-trace formula (and in fact, this is how we discovered it).  
 
\subsubsection{Organization of the paper}
After some preliminaries in \S~\ref{sec: prelim}, this paper is divided into two parts and an epilogue. In Part \ref{part1}, we  construct the  analytic distribution $\sI$ and prove the associated RTF, as well as  Theorems \ref{thm A} and \ref{thm B}. In Part \ref{part 2}, we construct the distribution  $\sJ$ and prove the associated RTF. In the epilogue, we  compare the two RTFs  to prove Theorems \ref{main thm} and \ref{th BBK}. More details on the contents of the two parts are provided at the beginning of each.

\subsection*{Acknowledgements}  We would like to thank Yifeng Liu for many helpful discussions, and especially for providing us with the material of \S~\ref{sec: prove Gauss}.  We are also grateful to Dongwen Liu, Weixiao Lu, Binyong Sun, and Lei Zhang for sharing and discussing with us their  works \cite{LS, Lu, LZ}; to Ryan Chen, Michael Harris, Chao Li, Yu Luo,  Eric Urban, and Zheng Wang for their comments on a draft of the paper; to the referee for their comments and suggestions; and to SLMath for its hospitality to both of us during the Spring 2023 semester on “Algebraic cycles, $L$-values, and Euler systems”, when part of this work was done.

\addtocontents{toc}{\medskip}

\section{Notation and preliminaries}\lb{sec: prelim} 
\subsection{Basic notation}  We set up some  notation to be used throughout the paper unless otherwise noted.
\subsubsection{Fields} \lb{sss tau}
We  denote by $F\supset F_{0}$ a quadratic extension of number fields, as in the introduction, and by ${\rm c}\in \Gal(F/F_{0})$ the conjugation. We  denote by $\A$ the ad\`eles of $F_{0}$. 
 From \S~\ref{sec rat} on, we will assume that $F_{0}$ is totally real and $F$ is CM.

 We denote 
 by 
 $$\eta\colon F_{0}^{\ts}\bs \A^{\ts}\to \{\pm 1\} $$
  the quadratic character associated with $F/F_{0}$. We fix an auxiliary element $\tau\in F^\ts$ such that ${\tau}^{\rm c}=-\tau$, and an extension $\eta'\colon F^{\ts}\bs \A_{F}^{\ts}\to \C^{\ts}$ of $\eta$.

If $F'$ is a number field and $S$ is a finite set of places of $F'$, we denote by $F'_{S}=\prod_{v\in S} F'_{v}$, and by $\A_{F'}^{S}=\prod'_{v\notin S} F_{v}'$. If $F''\subset F'$ is a subfield and $\ell$ is a place of $F''$, for notational purposes we identify $\ell$ with the set of places of $F'$ above $\ell$.

{\subsubsection{Identity matrix}
We denote by $1_\nu$ the identity matrix of rank $\nu$.}

\subsubsection{Coefficient fields}\lb{coeff} Throughout this paper,  a  \emph{coefficient field} is a field admitting embeddings into $\C$.

\subsubsection{Points of schemes}\lb{points of sch} If $Y$ is an (ind-)scheme of (ind-)finite type over a field $L$, by a slight abuse of notation we write `$y\in Y$' to indicate  that $y$ is a point of the underlying topological space of $Y$. The residue field at such a point $y$ is denoted by $L(y)\supset L$, and $y$ is identified with an element of $Y(L(y))$.

\subsubsection{Galois groups} If $E$ is a perfect field, we denote by $\ol{E}$ an algebraic closure (if choosing a particular closure is important, we will specify it in the text), and by $G_{E}\coloneqq \Gal(\ol{E}/E)$ the associated absolute Galois group.

\subsubsection{$L$-functions}
In the rest of the paper (\emph{unlike} in  the introduction), all global  $\zeta$- and  $L$-functions valued in the complex numbers are complete including the archimedean factors (this also includes the ratio of $L$-functions $\sL(1/2,\Pi,\chi)$). 
 If $L^{S}(s)$ is a global $L$-function, we denote by $$L^{ S, *}(s_{0})$$ its leading term at $s=s_{0}$.

\subsubsection{Groups}\label{sss:gp}
We now recall the groups under consideration in this paper.
  We denote by $\bG_{m}=\Spec\Q[T^{\pm1}]$ the multiplicative group over $\Q$. If $\G$ is a (usually, group-) scheme over a global field $F_{0}$ and $v$ is a place of $F_{0}$, we denote $G_{v}\coloneqq \G(F_{0, v})$ with its $v$-adic topology. We also denote 
$$[G]=\G(F_{0})\bs \G(\A).$$

For $*=\emptyset, 0$ (where in this type of context, `$\emptyset$' will always mean `no subscript') and $\nu\in \bN$,   let 
$\RG'_{\nu, *}\coloneqq  \Res_{F_{*}/F_{0}}\GL_{\nu}$. We consider 
\beq\lb{G' def}
\RG'\coloneqq  \RG'_{n}/{\RG}'_{1,0} \ts  \RG'_{n+1}/\G_{1, 0}',
\eeq
where $\G_{1,0}'$ is the $F_{0}$-split center of $\G'_{\nu}$, and its subgroups
$$j_{1}\colon \RH_{1}'\coloneqq  \RG'_{n}\incl \RG',$$ where $j_{1}(h)\coloneqq  [(\diag(h, 1), h)]$, and
$$j_{2}\colon \RH'_{2}\coloneqq   \RG'_{n, 0}/{\RG}'_{1, 0} \ts  \RG'_{n+1, 0}/\RG'_{1,0}\incl \RG',$$
where $j_{2}$ is induced by $F_{0}\incl F$.

\lb{coh incoh G} 
We fix an element $\nrm \in F_0^\ts$; when $F_{0}$ is totally real and $F$ is CM,  we assume that $\nrm$ is totally positive. 
We denote by $\sV$ the set of isomorphism classes of pairs $V=(V_{n}, V_{n+1}=V_{n}\oplus Fu)$ of $F/F_{0}$-hermitian spaces with $(u,u)=\nrm$. For a place $v$ of $F$ (or $v=\infty$), we denote by $\sV_v=\sV_{e,v}$ the set of isomorphism classes of pairs $V=(V_{n,v}, V_{n+1,v}=V_{n,v}\oplus F_vu)$ of $F_v/F_{0,v}$-hermitian spaces with $(u,u)=\nrm$. 
 When $F_{0}$ is totally real and $F$ is CM,  
 we denote by $V_{\infty}^{\circ}=(V_{v}^{\circ})_{v\vert\infty}$ the pair over $F_\infty$ such that each $V_{n,v}$ is positive-definite, and by
$\sV^{\circ} =\sV^\circ_e\subset \prod_v \sV_v$ the set of (coherent or incoherent) pairs $(V_{v})_v$ such that $V_{v}=V_{v}^{\circ}$ for all $v\vert\infty$ defined as in \S~\ref{sec incoh}. We partition
$$\sV^{\circ}=\sV^{\circ , +}\sqcup \sV^{\circ, -}, $$
where $V\in \sV^{\circ, \epsilon}$ if and only if $\epsilon(V)=\epsilon$. (Then $\sV^{\circ, +}$ may be  identified with  the intersection of $\sV$ and $\sV^\circ$ inside $\prod_v \sV_v$.)
 
For $V\in \sV$ (or $V\in \sV^\circ$ if $F_{0}$ is totally real and $F$ is CM), we use the notation $\RH^{V}$, $\G^{V}$ introduced in \S~\ref{sec incoh} for the (products of) unitary groups associated to $V$.

\subsubsection{Relative selfdual hyperspecial subgroups} \lb{sss rel sd}
  Let $v$ be a finite place of $F_0$, let $V\in \sV_v$, and let $G_v=G_{v}^{V_{v}}$.
  A \emph{relative selfdual hyperspecial subgroup}  of $ G_{v}$ is one of the form  $K_{v} =K_{n,v}\ts K_{n+1,v}$  where each $K_{\nu, v}$ is the  stabilizer of a selfdual lattice $\Lm_{\nu,v}$ such that $\Lm_{n+1,v}\cap V_{n,v}= \Lm_{n,v}$.  If $v$ splits in $F$,   relative selfdual hyperspecial subgroups always exist. If $v$ is inert, they exist if and only if $\e(V_v)=1$ and $v(\nrm)$ is even.

\subsubsection{Vertex parahoric subgroups}\label{sss:vert} 
Let $v$ be a finite place of $F_0$ unramified in $F$, and let $W$ be a $F_v/F_{0,v}$-hermitian space of dimension $\nu\geq 1$. 
An $\sO_{F,v}$-lattice $\Lambda_v$ in $W$ is called a vertex lattice if we have 
 \begin{equation*}
   \Lambda_v\subset \Lambda_v^\ast \subset\varpi_v\i\Lambda_v ,
\end{equation*}
where $\Lambda_v^\ast$ is the dual lattice with respect to the hermitian form on $W_v$. Here $\varpi_v$ denotes a uniformizer of $F_{0,v}$.
The type $t(\Lambda_v)$ of a vertex lattice $\Lambda_v$ is by definition the rank of the free  $\sO_{F_v}/(\varpi_v)$-module $ \Lambda_v^\ast /\Lambda_v$. In particular, a selfdual lattice is of type $t=0$. The stabilizer of a vertex lattice (of type $t$) under the action of $\U(W)(F_{0,v})$ is called a vertex parahoric subgroup (of type $t$). (If $v$ is ramified, we can make the same definitions after replacing $\vpi_v$ by a uniformizer of $F_v$.)

Let now $v$ be inert. For $W$ to admit a vertex lattice of type $t$, it is necessary and sufficient for its Hasse invariant to be $(-1)^t$. Moreover, rescaling the hermitian form on $W$ by $\vpi_v$ induces the identity on the unitary group and exchanges vertex parahorics of type $t$ and $\nu -t$.
We will be interested in the following situation. Let $\e=\e_v\in \{0, 1\}$ have the same parity as $v(\nrm)$, where $\nrm\in F_0^\ts$ is the hermitian norm of the special vector $u$ as fixed above.   Let $V_v\in \sV_v$. We fix a vertex lattice  $\Lambda_{n,v}\subset V_{n,v}$ of type $0\leq t \leq n$ and let $ \Lambda_{n+1,v}=\Lambda_{n,v}\oplus \sO_F{\varpi_v^{-\lfloor v(\nrm)/2\rfloor}u}\subset V_{n+1,v}$.
Then $\Lambda_{n+1,v}$ is a vertex lattice of type $t+\epsilon$.  Let $K_{\nu,v}$ be the stabilizer of  $\Lambda_{\nu,v}$ for $\nu=n, n+1$. We then call $K_v=K_{n,v}\times K_{n+1,v}$ a vertex parahoric subgroup of type $(t,t+\epsilon)$. The case of type $(0,0)$ gives  a  relative selfdual hyperspecial subgroup.

\subsection{Measures}  \lb{sec: meas0}
Let $F_{0}$ be a number field,  and let $D=D_{F_{0}}$ be the absolute value of its discriminant. 
 Fixing a nontrivial character $\psi\colon F_{0}\bs \A\to \C^{\ts}$, we denote by  $dx=\prod_{v}d x_{v}$ the selfdual measure  on $\A$ with respect to $\psi$; it satisfies $ \vol(F_{0}\bs \A , dx)=1$. For a finite place $v$, let $d_{v}$ be a generator of the different ideal of $F_{0,v}$ and let $D_{v}\coloneqq  |d_{v}|^{-1}$;
 if $\Ker (\psi_{v}) = d_{v}^{-1}\sO_{F_{0,v}}$,
  then  $\vol(\sO_{F_{0, v}}, dx_{v}) = D_{v}^{-1/2}$. 
 We have $D=\prod_{v\nmid \infty}D_{v}$, and for a finite set of places $S$ of $F_{0}$  we define $D^{S}\coloneqq  \prod_{v\nmid S\infty} D_{v}$.

\subsubsection{Tamagawa measures}  \lb{sec: tam}
If $\G $ is a reductive group over a local or global field $E$, we denote by  $M_{\G}$  the Artin--Tate motive attached to (the quasi-split inner form of) $\G$ by Gross  \cite{gross-mot}. If $E$ is a local field, let
$$
 \Delta_{\G}^{} \coloneqq   D_{v}^{\dim \G/2} 
 L(M_{\G}^{\vee}(1)).$$
Then the abelian term in \eqref{sL def} (including the factor $D_{F_{0,v}}^{n+1}$)  equals $\Delta_{G^{V_{v}}}^{}/\Delta_{H^{V_{v}}}^{}$ where $H^{V_{v}}$, $G^{V_{v}}$ are as in \eqref{unit} (for any $V_{v}\in \sV_{v}$.)

Assume from now on that $E$ is the global field $F_{0}$. For a finite set   $S$  of places  of $F_{0}$,  let 
$$
\Delta_{\G}^{S}\coloneqq  (D^{S})^{\dim \G/2} \, L^{ S, *}(M_{\G}^{\vee}(1), 0).$$
If $\G$ is a unitary group in $\nu$ variables, this only depends on $S$ and  $\nu$, explicitly
$$
\Delta_{\G}^{S} =  (D^{S})^{\nu^2/2} \, \prod_{j=1}^\nu
L^{ S}(1, \eta^j).$$
In particular, the definition can be applied to incoherent unitary groups as well.

Let $\omega$ be any non-zero top-degree invariant differential form on $\G$. 
We denote by $$d_{\omega}g_{v}\coloneqq  |\omega|_{v}$$ 
 its  modulus  with respect to $dx_{v}$ (see \cite[\S I.4]{Oes}, \cite[\S 10.1.7]{Bou}),  a Haar measure on $\G(F_{0,v})$. We define
 $$d^{\nat}g_{v}\coloneqq 
 \Delta_{\G, v} \, d_{\omega}g_{v} 
 $$
Then for all finite places $v$ and any open compact subgroup $K_{v}\subset G_{v}$, we have $\vol(K_{v}, d^{\natural}g_{v})\in \Q^{\ts}$.
Moreover, if $G_{v}$ is unramified and $K_{v}$ is hyperspecial, we have $\vol(K_{v}, d^{\natural}g_{v})=1$.  
The  Tamagawa measure on $\G(\A)$ is
\beq \lb{eq: tam m}
dg&\coloneqq \Delta_{\G}^{ -1} \prod_{v}d^{\natural} g_{v}.
\eeq
The definition also makes sense when $\G(\A)$ is an 
incoherent unitary group (or a product thereof). 

\subsubsection{Variants}
We define a variant 
\beq\lb{dgv}
dg_{v}=\begin{dcases} d^{\natural}g_{v}  &\text{if $v\nmid \infty$}\\
\Delta_{\G}^{ -1} d^{\natural}g_{v}=\Delta_{\G}^{\infty,  -1} d_{\omega}g_{v} & \text{if $v=\infty$}
\end{dcases}
\eeq 
so that $dg=\prod_{v}dg_{v}$.
The `rationale' for this choice is the following.
\begin{lemma}\lb{gross fact} If $G_{\infty}$ is compact, then $\vol(G_{\infty}, dg_{\infty})$ is rational. 
\end{lemma}
\begin{proof} We say that two measures $\mu$, $\mu'$ are commensurable if $\mu=c\mu'$ for some $c\in \Q^{\ts}$. 
 Let $\mu\coloneqq \prod_{v} \mu_{{v}}$ be the measure on $\G(\A)$ considered in \S~9 of \cite{gross-mot}, to which all citations in this proof will refer. The measure $\mu$   is nonzero by Propositions 9.4, 9.5. For almost all finite $v$, $\mu_{{v}}=dg_{v}$; for all finite $v$, $\mu_{{v}}$ gives rational volume to compact open subgroups (equation (5.2)), hence it is commensurable with $dg_{v}$; and   $ \mu$ is commensurable with $dg$ (Theorem 9.9). It follows that $dg_{\infty}$ is commensurable with $\mu_{\infty}$, which  (again by equation (5.2)), gives rational volume to $G_{\infty}$. 
\end{proof}

We also consider a different measure, for comparison with some of the literature (notably \cite[\S~2]{wei-fourier}). Let $\RZ$ be the center of $\G$,  let $\G^{\rm ad}\coloneqq  \G/\RZ$;
put
$$\zeta_{\G, v}(1)\coloneqq  D_{v}^{-\dim \RZ/2}\Delta_{\RZ, v}, \qquad 
\zeta^{S, *}_{\RG}(1)\coloneqq  D^{-\dim\RZ/2}\Delta_{\RZ}^{S},$$ sot that $D_{v}^{\dim \RZ/ 2}\zeta_{\G, v}(1)\Delta_{\G^{\rm ad}, v}=\Delta_{\G, v}$.
 Then we set 
$$d^{*}g_{v}\coloneqq \zeta_{\G,v}(1)^{} \, d_{\omega}g_{v}, \qquad d^{*}g=\prod_{v}d^{*}g_{v}$$
so that  
$$
dg=
\zeta_{\G, v}^{*}(1)^{-1} \prod_{v}d^{*}g_{v}$$
and for finite $v$, $dg_{v}= D_{v}^{\dim \RZ/2}\Delta_{\G^{\ad},v} d^{*}g_{v}$.

\subsubsection{Local and incoherent measures}\lb{sec: meas}
The global measures do not depend on $\omega$, but the local ones do (up to multiplying $\omega$ by a scalar of absolute value~1). We fix the following explicit choices:
\begin{itemize}
\item if $\G=\GL_{\nu}$, we take\footnote{We do not specify the order of the variables in the wedge product, which leaves a sign ambiguity of no consequence for the associated Haar measures.} 
 $$\omega\coloneqq \det(g)^{-\nu}\, {\wedge}_{i, j}dg_{ij}.$$
\item if $\G$ is a (product of) unitary groups over a local or a global field, we fix $\omega$ as in \cite[\S~2]{wei-aut}. 
\end{itemize}
If $\G$ is a (product of) incoherent unitary groups, we then get a measure on $\G(\A)$ by \eqref{eq: tam m}, with a factorization $dg=\prod_{v}dg_{v}$ as in \eqref{dgv}.

\subsection{Hecke algebras}
\lb{sec hecke}
Let $\G$ be a reductive group over a number field $F_{0}$, let $v$ be a finite  place of $F_{0}$ or $v=\infty$,  and let $L$ be a characteristic-zero field, with $L=\C$ if $v=\infty$ is archimedean and $G_{\infty}$ is not compact. We denote by $\CS(G_{v},L)$ the space of Schwartz functions on $G_{v}$ valued in $L$:  when $v$ is finite, this is the same as the space of smooth compactly supported $L$-valued functions on $G_v$, whereas when $v=\infty$
is archimedean this is the space of smooth functions  $\dot f\colon G_\infty\to \C$ such that $D\dot f$ is bounded for every algebraic differential operator $D$ on $G_\infty$.
We denote by
 $\sH(G_{v}, L)$   the space of \emph{Schwartz measures} on $G_{v}$: those are measures of the form $\df\, dg_{v}$ where $\df\in \CS(G_{v}, L)$ and $dg_{v}$ is the Haar measure fixed above. The field $L$ will be omitted when it is unimportant or understood from context.   For an open compact $K_{v}\subset G_{v}$, we denote 
 $$e_{K_{v}}\coloneqq  {1\over \vol(K_{v}^{} , dg_{v})}
\one_{K_{v}^{}} \, dg_{v}$$
for any Haar measure $dg_{v}$; it is an idempotent in $\sH(G_{v})$. We denote by $\sH(G_{v})_{K_{v}}$ the subspace of bi-$K_{v}$-invariant measures.

For $S$ a finite set of places of  $F_{0}$, and $\G$ denoting either $\G'$ or $\G^{V}$ for some $V\in \sV^{\circ}\cup \sV$, we consider the Hecke algebra
  $$\sH(\G(\A^{S}))\coloneqq  \bigotimes_{v\notin S}' \sH(G_{v}),$$
 where the restricted tensor product is with respect to   
\beq\lb{unit hecke}
\ff_{v}^{\circ}\coloneqq   e_{K_{v}^{\circ}}.
\eeq
Here,  when $\G=\G'=  \eqref{G' def}$, we take 
$K_{v}^{\circ}\subset G_{v}'$ to be the image of $\GL_{n}(\sO_{F_{ v}}) \ts \GL_{n+1}(\sO_{F, v})$.
When $\G=\G^V$, for all but finitely many $v$ we take $K_v^\circ$ to be the relative selfdual hyperspecial $K_{v}^{\circ}\subset G_{v}$ stabilizing the product of the selfdual lattices $\Lm_{\nu,v}$ fixed when defining $V_{\nu,\A}$ (\S~\ref{sec incoh}). 

If $K^{S}=\prod_{v}K_{v}\subset \G(\A^{S})$ is an open compact subgroup, with $K_{v}=K_{v}^{\circ}$ for all but finitely many $v$,  we denote $e_{K^{S}}\coloneqq  \prod_{v}e_{K_{v}}$, and
 we say that an element $f\in \sH(\G(\A^{\infty}))$ is \emph{supported in the set $S$} if we can write $f=f_{S}\ot \ot_{v\notin S\infty}f_{v}^{\circ}$ for some $f_{S}\in 
\sH(G_{S})$.
 
For $f\in \sH(G_{v})$, we   denote $f^{\vee}(x)\coloneqq f(x^{-1})$. 
 We denote by  $\star$ the convolution operation $f_{1}\star f_{2}(x)\coloneqq \int_{G_{v}} f_{1}(xg)f_{2}(g^{-1})$, and sometimes omit this symbol.

 \subsubsection{Convention.} We stipulate that groups and Hecke algebras act on locally symmetric spaces, Shimura varieties, and their homology and algebraic cycles {on the right}; on automorphic forms {on the left}.

\subsection{Local base-change and distinction}
Let $v$ be a place of $F_{0}$. If $v$ is nonarchimedean,  $G^{?}_{v}$ is the group of $F_{0,v}$-points of a reductive group over $F_{0,v}$, and $L$ is a coefficient field, 
we say that an absolutely irreducible (that is, $\pi\ot_{L}\ol{L}$ is irreducible) smooth admissible representation $\pi$ of $G^{?}_{v}$ over $L$ is \emph{tempered} if $\pi\ot_{L, \iota}\C$ is tempered in the usual sense for every $\iota\colon L\into \C$. When we are working over $L=\C$ and there is a risk of ambiguity, we will say `absolutely tempered' for the notion just introduced.

Let $\sV_{v}$ be the set of isomorphism classes of pairs $V_{v}=(V_{n,v}, V_{n+1,v}=V_{n,v}\oplus F_{0,v}u)$ of hermitian spaces over $F_{v}/F_{0,v}$ {where $\pair{u,u}=e$ (hence $\sV_{v}$ depends on $e$).} 
Let ${\rm Temp}(G_{v}^{?})(L)$ be the set of $\Gal(\ol{L}/L)$-orbits of isomorphism classes of irreducible tempered representations of $G_{v}^{?}$ over $\ol L$.

\subsubsection{Local base-change} \lb{loc BC}
Let $v$ be a place of $F_{0}$.
Thanks to  \cite{mok, KMSW}, we have a local base-change map ${\rm BC}$ from complex irreducible admissible representations of $G_{v}^{V_{v}}$  to complex irreducible admissible representations of $G_{v}'$, whose definition is recalled in \cite[\S~2.10]{BP-Planch}. It has the following properties:
\begin{enumerate}
\item it restricts to a map
\beq\lb{local BC} 
{\rm BC}\colon {\rm Temp}(G_{v}^{V_{v}})(\C)\to {\rm Temp} (G_{v}')(\C);
\eeq
\item by Lemma \ref{equiv bc} below, it is rational in the sense that it yields a map
$${\rm BC}\colon {\rm Temp}(G_{v}^{V_{v}})(L)\to {\rm Temp} (G_{v}')(L)$$
for any characteristic-zero field $L$;
\item  when $v$ splits in $F$, we simply have 
\beq \lb{split bc}
{\rm BC}(\pi)= \pi\boxtimes \pi^{\vee}\eeq if we identify $G_{v}^{V_{v}}\cong G _{n, 0,v}'\ts G_{n+1, 0,v}'$ for the unique $V_{v}\in \sV_{v}$;
\item\lb{bc p3} when $G_{v}^{V_{v}}=U(n)\ts U(n+1)$ (the compact unitary group) over $\R$, the preimage of $\Pi^{\circ}_{\R}$ under \eqref{local BC} consists of the trivial representation only;
\item\lb{bc p4} when $G_{v}^{V_{v}}=U(n-1, 1)\ts U(n,1)$ over $\R$, the preimage of $\Pi^{\circ}_{\R}$ under \eqref{local BC}
consists of the $n(n+1)$ discrete series representations having the Harish-Chandra parameter $\{{1-\nu\over 2}, {3-\nu\over 2}, \ldots, {\nu-1\over 2}\} $ on the $U(\nu-1, 1)$-component. (See \cite[Proposition C.3.1]{LTXZZ}.)
\end{enumerate}
If $v$ is non-archimedean and $\pi_{v}$, respectively $\Pi_{v}$, is a representation of $G_{v}^{V_{v}}$, respectively $G_{v}'$, over a coefficient field $L$, we will write ${\rm BC}(\pi_{v})=\Pi_{v}$ if ${\rm BC}(\iota\pi_{v})=\iota\Pi_{v}$ 
for every embedding $\iota\colon L\into \C$. 

{For a given $\Pi_v\in {\rm Temp}(G_v')(L)$,  the set of pairs $(V_v, \pi_v)$ where $V_v\in \sV_v$ and $\pi_v\in   {\rm Temp}(G_v^{V_v})(L)$ satisfies ${\rm BC}(\pi_v)=\Pi_v$ is called the \emph{Vogan $L$-packet} of $\Pi_v$.}

\subsubsection{Hermitian representations and local root numbers}
\lb{sss tau psi}
We will say  that a tempered representation $\Pi_{v}$ of  $G_{v}'$ is \emph{hermitian} if  the space $\Hom_{H_{2,v}'}(\Pi_{v}, \eta^{n}_{v}\boxtimes \eta_{v}^{n-1})$ is nonzero. 
By the local Flicker--Rallis conjecture proved  by Matringe, Mok, and others (see \cite[\S 3.1]{anan} and references therein), a representation  $\Pi_{v}$ over $\C$ is hermitian if and only if it is in the image of base-change for some $V_{v}\in \sV_{v}$.

 Let $\Pi_{v}$ be a hermitian tempered representation of $G_{v}'$ over a coefficient field $L$; we restrict to $L=\C$ if $v$ is archimedean. 
 {Recall the fixed the totally imaginary  $\tau\in F$ from \S~\ref{sss tau}, and put $\tau \psi_{F}(z)\coloneqq \psi_{F}(\tau z)$. The local root number (see  \cite[\S5]{GGP} for the definition when $L=\C$)
\beq\lb{def local rn}
\eps(\Pi_v)=\eps(\Pi_{v}^\iota, \tau\psi_{F_v}) \in {\pm 1} \eeq 
is independent of the choice of $\iota\colon L\into \C$ (if $v$ is finite), cf. \cite{ces}. Moreover,  
\cite[Theorem 6.1.3]{GGP} shows that it is also independent of $\tau$.}
 Later in Lemma \ref{lem: loc ep}  we will compute the local root number for a certain class of representations when $v$ is non-archimedean. Here we give the local root number in the archimedean case $v\mid\infty$ for our special representation $\Pi_v^\circ$ of  interest. The $L$-parameter of $\Pi_{\nu,v}$ is of the form 
$$\bigoplus_{i=1}^{\nu} \mu_v^{\nu+1-2i}, \quad\mu_v(z)=\frac{z}{|z|}.
$$
Using \cite[Proposition 2.1]{GGP2} we obtain 
\beq\lb{eq:eps v=infty}
\eps(\Pi_v)=\prod_{1\leq i\leq n, 1\leq j\leq n+1}\eps( \mu_v^{n+1-2i+(n+2-2j)})=(-1)^{n+1 \choose 2}.
\eeq

 \subsubsection{Distinction and the local Gan--Gross--Prasad conjecture} \lb{sss: dist}
Let $v$ be a place of $F_{0}$ and let $L$ be a field of characteristic zero; we restrict to $L=\C$ if $v$ is archimedean and $V_{v}$ is not definite. 
We say that a tempered representation $\pi$ of $G_{v}^{V_{v}}$ over $L$ is \emph{distinguished} if the space $\Hom_{H_{v}^{V_{v}}}(\pi_{v}, L)$ is nonzero;  by the {multiplicity-one} result of 
\cite{AGRS}, this space is one-dimensional if nonzero.  It is clear that  distinction is a $\Gal(\ol{L}/L)$-invariant property. We denote by
$${\rm Temp}({H_{v}^{V_{v}}}\bs G_{v}^{V})(L)\subset {\rm Temp}(G_{v}^{V})(L) $$
the subset of orbits of distinguished tempered representations.

 The following fundamental result is 
the local Gan--Gross--Prasad conjecture for unitary groups.

\begin{proposition}\lb{local ggp} Let $\Pi_{v}$ be a hermitian tempered representation of $G_{v}'$ over a coefficient field $L$; we restrict to $L=\C$ if $v$ is archimedean. 

There exists a unique pair $(V_{v}, \pi_{v})$ with $V_{v}\in \sV_{v}$ and $\pi_{v}\in {\rm Temp}(H^{V_{v}}_{v}\bs G_{v}^{V_{v}})(L)$ such that $\Pi_{v}={\rm BC}(\pi_{v})$.
{Moreover, the Hasse invariant of  $V_{v}$ (defined by \eqref{eq: HW}) is determined by the local root number:
$$
\epsilon(V_v)=\eps(\Pi_v).
$$
} 
\end{proposition}
{In the context of the lemma, we will refer to  $(V_v, \pi_v)$ as the \emph{distinguished element} in the Vogan $L$-packet of $\Pi_v$.}
\begin{proof} 
 If $L=\C$, this is proved  by Beuzart-Plessis in \cite{BP-ggp1, BP-ggp2}. In general (if $v$ is finite), we may pick an embedding $\iota\colon L\into \C$ and apply the result to $\iota\Pi_{v}$ to obtain a pair $(V_{v}, \pi_{v}^{\C})$. By uniqueness, $\pi_{v}^{\C}$ is isomorphic to its ${\rm Aut}(\C/\iota L)$-conjugates.
\end{proof}

In the special archimedean case considered in \eqref{eq:eps v=infty}, it is straightforward to check that 
${n+1 \choose 2}\equiv { 2\lfloor{(n+1)/2}\rfloor \choose 2}\, (\!\!\!\mod 2)$. On the other hand, if $V_v$ is the pair of positive definite hermitian spaces, its  Hasse invariant of the pair is by definition (cf. \eqref{eq: HW}) equal to $(-1)^{ 2\lfloor{(n+1)/2}\rfloor \choose 2}$.   Thus  the identity $\epsilon(V_v)=\eps(\Pi_v)$ is verified directly.

\subsection{Automorphic base-change}

\subsubsection{Rational spaces of automorphic representations} \lb{rat sp aut}
The following discussion is based on \cite[Th\'eor\`eme 3.1.3]{clozel}.
Let $L$ be a  coefficient field,  and let $\Pi=\Pi^{\infty}\ot \one_{\infty}$ be an absolutely irreducible  representation of $\G'(\A)$  over an $L$-vector space. We say that $\Pi$ is cuspidal automorphic of trivial weight if for every (equivalently, some) embedding $\iota\colon L\into \C$, the representation $\Pi^{\iota}\coloneqq  \iota \Pi^{\infty}\ot \Pi_{\infty}^{\circ}$ is cuspidal and automorphic. Every cuspidal automorphic representation $\Pi_{\C}$ of $\G'(\A)$ such that $\Pi_{\C, \infty}\cong \Pi_{\infty}^{\circ}$ arises as $\Pi^{\iota}$ for some $\Pi$ defined over a number field, there is a smallest such number field $\Q(\Pi)\subset \C$, and a representation $\Pi $ over $\Q(\Pi)$, unique  up to $\Q(\Pi)$-isomorphism, such that $\Pi^{\iota}\cong \Pi_{\C}$.

Denote by $\wt\sC(\G')(L)$ the set of isomorphism classes of trivial-weight cuspidal automorphic representations defined over $L$, and by $\sC(\G')(L)\coloneqq  \wt\sC(\G')(\ol{L})/ G_{L}$, where we recall that $G_{L}\coloneqq \Gal(\ol{L}/L)$.  
 
\begin{lemma}\lb{rat aut rep} For every  coefficient  field  $L$, the   natural map 
$$\wt\sC(\G')(L) \to \sC(\G')(L)$$
is a bijection.
\end{lemma}
\begin{proof} This follows from the above discussion and \cite[Proposition 3.1, Th\'eor\`eme 3.1.3]{clozel}.
\end{proof}

\subsubsection{On the Ramanujan conjecture} 
Consider the class 
of automorphic representations $\Pi$ of $\G'(\A)$ that are isobaric sums
\beq \lb{bc isob}   
\Pi = \Pi_1 \boxplus \cdots \boxplus \Pi_r 
\eeq
where the $\Pi_i$ are pairwise non-isomorphic discrete automorphic representations of $\GL_{r_i}(\A_F)\ts \GL_{r_i'}(\A_F)$ satisfying $\Pi_i^{\rm c}\cong \Pi_i^\vee$. {Such a  $\Pi$ is called \emph{isobaric} if each $\Pi_i$ is cuspidal.}

\begin{proposition} \lb{ramanujan}
Let $\Pi$ be a complex automorphic representation of $\G'(\A)$ of the form \eqref{bc isob}. 
\begin{enumerate}
    \item \lb{ram 1} If $\Pi$ is isobaric, $F_0$ is totally real, $F$ is CM, and $\Pi_\infty$ has the same infinitesimal character as an algebraic representation, then $\Pi$ is absolutely tempered at every finite place.
    \item If $\Pi$ is not isobaric, then $\Pi$ is not tempered at any finite place.
\end{enumerate}
\end{proposition}
\begin{proof} Part (1) is proved by Caraiani \cite{Ca}, completing previous works of Shin and other authors. (Caraiani assumes that $\Pi$ itself is cuspidal, but temperedness is preserved under isobaric sums.)  Part (2) follows from the  description of the discrete spectrum of $\GL_N$ in \cite{MW}.
\end{proof}

\subsubsection{Ramakrishnan's automorphic Tchebotarev theorem} 
We will  use the  following  special case of \cite[Theorem A]{Ram}.

\begin{proposition} \lb{Ram-prop}
Let $\Pi$, $\Pi'$ be two isobaric automorphic representations of $\GL_{n}(\A_{F})$. Assume that $\Pi_{w}\cong \Pi_{w}'$ for all but finitely many primes $w$ of $F$ split over $F_{0}$. Then $\Pi\cong \Pi'$. 
\end{proposition}

\subsubsection{Base change} \lb{sec: G BC}
Let $V\in \sV$ or, if $F_{0}$ is totally real and $F$ is CM, let $V\in \sV\cup \sV^{\circ}$.  Let $\G=\G^{V}$, $\RH=\RH^{V}$. 
For a coefficient field $L$,
denote by 
$$\wt\sC(\G)(L) \quad\supset\quad \quad\wt\sC({\RH\bs\G})(L)$$
 the set of isomorphism classes  of tempered\footnote{That is, tempered at all finite places.} cuspidal automorphic representations of $\G(\A)$ that are trivial at infinity, and its subset of representations that are $\RH(\A)$-distinguished. (If $V\in \sV^{\circ, -}$, the relevant notions of automorphic representation of $\G(\A)$ were defined in \S~\ref{sec incoh}.)  We also put
$$\sC(\G)(L)\coloneqq  \wt\sC(\G)(\ol L) /G_{L}  \quad\supset\quad \quad\sC({\RH\bs\G})(L)\coloneqq \quad\wt\sC({\RH\bs\G})(\ol L)/G_{L}.$$

Since representations  in $\wt\sC(\G')$, $\wt\sC(\G)$ always  have a model over a number field (by \cite[Th\'eor\`eme 3.13]{clozel} for $\G'$ and \cite[Corollary 2.13]{Har13} for $\G$), we may view $\sC(\G')$,  $\sC(\G)$ and $\sC(\RH\bs \G)$ as ind-finite schemes over $\Q$.

\begin{definition}\lb{BC def}
Let $V\in \sV$, and let $\G=\G^{V}$. Let $\pi$ be a (complex) {discrete} automorphic representation of $\G(\A)$,
and let $\Pi$ 
be an automorphic representation of $\G'(\A)$ of the form \eqref{bc isob}.
We say that $\Pi$ is an \emph{automorphic base-change of $\pi$}, and write $\Pi\cong {\rm BC}(\pi)$, if for all but finitely many places $v$ of $F_{0}$ such that $\pi_v$ is unramified, we have  $\Pi_{v}\cong {\rm BC}(\pi_{v})$ for the local base-change of \S~\ref{loc BC}. (By strong multiplicity one for $\GL_\nu$, this is unique if it exists.)
\end{definition}

{
\begin{remark}  \lb{rmk bc new} We have the following results on automorphic base-change.

\begin{enumerate}
\item By Shin's result in \cite[Theorem A.1]{WG} together with \cite[footnote to Proposition 3.2.8]{LTXZZ}, if $F_0$ is totally real, $F$ is CM, $V\in \sV$ and $\pi$ is a discrete automorphic representation of $\G^V(\A)$ such that $\pi_\infty$ has the same infinitesimal character as an algebraic representation, then $\pi$ admits an automorphic base-change. (This is the only case of interest for this paper; the full general theory is developed in \cite{mok, KMSW}.)
\item\lb{rmk bc 2} If $\pi$ is tempered at every place, then  by \cite[Theorem 1.7.1]{KMSW} we have ${\rm BC}(\pi)_v\cong {\rm BC}(\pi_v)$ for every place $v$ of $F_0$.  
\end{enumerate}
\end{remark}
}

\subsubsection{Rationality of base-change}
The associated local and automorphic complex base-change maps enjoy some nice equivariance properties.

\begin{lemma}\lb{equiv bc} 
Let $V\in \sV$ and let $\G=\G^V$.  
    \begin{enumerate}
        \item \lb{equiv loc} Let $v$ be a finite place of $F_0$, let $G_v=\G(F_{0, v})$. For every irreducible absolutely tempered   representation $\pi_v$ of $G_v$ over $\C$ and every $\sg\in \Aut(\C/\Q)$, we have ${\rm BC}(\pi_v^\sg) \cong {\rm BC}(\pi_v)^\sg$.
        \item \lb{equiv glob} Assume that $F_0$ is totally real and $F$ is CM. For every complex absolutely tempered automorphic representation $\pi$ of $\G(\A)$ and every $\sg\in \Aut(\C/\Q)$, we have ${\rm BC}(\pi^\sg) \cong {\rm BC}(\pi)^\sg$. 
    \end{enumerate}
\end{lemma}
\begin{proof}
First we note that if $v$  splits, say $v\sO_F =w\ol{w}$, then part \eqref{equiv loc} is clear  from \eqref{split bc}.
   The analogue of relation \eqref{equiv loc}  also holds for the base-change of characters of $U(1)_{F_v/F_{0,v}}$, as ${\rm BC}(\chi_v) (t)=\chi_v(t/t^{\rm c})$. 

For part \eqref{equiv glob}, by the previous paragraph and the compatibility of local and global base-change at split places, we have that ${\rm BC}(\pi)_v^\sg$ is the automorphic base-change of $\pi_v^\sg$ for every split place~$v$. By Proposition \ref{Ram-prop}, this implies  ${\rm BC}(\pi)^\sg \cong  {\rm BC}(\pi^\sg)$. 

Finally, consider the general case of  part \eqref{equiv loc}. At least if $G_v$ is quasi-split this is proved as \cite[(5.3)]{ST} by a local argument. We follow Gan's alternative argument sketched in Remark 5.4 \emph{ibid.} Since the operation of taking Langlands subquotients of a parabolic induction is compatible with base-change and equivariant under  automorphisms of $\C$ (for the latter, we renormalize the operation as in \cite[pp. 92, 132-133]{clozel}), it suffices to consider the case where $\pi_v$ is a discrete series. The operation of twisting by a character also enjoys the same compatibilities, so that by the result for characters, we may harmlessly replace $\pi_v$ by a twist. 

Consider a CM quadratic extension $\dot F/\dot F_0$ with a place $\dot v$ of $\dot F_0$ such that $(\dot F_{\dot v}/\dot F_{0, \dot v}) \cong (F_v/F_{0,v})$, and a product $\dot \G$ of (say, totally definite)  $\dot F/\dot F_0$-unitary groups  in $n$, $n+1$ variables such that $\dot G_{\dot v}\cong G_v$. By \cite[Theorem 5.7]{Shin-pl}, we can find a trivial-weight  automorphic representation $\dot \pi$ of $\dot \G$ such that  $\dot \pi_{\dot v} $  is isomorphic to a twist of $\pi_v$ and $\dot \pi_{w}$ is supercuspidal for some other place $w$; this implies that   $\dot \pi$ is cuspidal and stable and then, by Proposition \ref{ramanujan} \eqref{ram 1} applied to ${\rm BC}(\pi)$,  absolutely tempered. By part \eqref{equiv glob},  we have ${\rm BC}(\dot\pi^\sg) \cong{\rm BC}(\dot\pi)^\sg$ for every $\sg\in \Aut(\C/\Q)$. By the compatibility of local and global base-change from Remark \ref{rmk bc new}~\eqref{rmk bc 2}, this implies the statement of part \eqref{equiv loc} for a twist of $\pi_v$. Since base-change is compatible with character twisting, by part \eqref{equiv loc} for $U(1)_{F_v/F_{0,v}}$ we conclude that the statement holds for $\pi_v$ as well. 
\end{proof}

Suppose now that  $F_{0}$ is totally real and $F$ is CM. Let $V\in \sV^{\circ}$, let $\G=\G^{V}$, and let $L$ be a characteristic-zero field. Let $\pi$ be a cuspidal automorphic representation of $\G(\A)$ over $L$ which is trivial at infinity. (If $V\in \sV^{\circ, +}$, this simply means that $\pi=\pi^\infty\ot\one_\infty$ and $\pi\ot_{L, \iota} \C$ is irreducible, cuspidal, and automorphic for every $\iota\colon L\into \C$. If $V\in \sV^{\circ, -}$, this was defined in \S~\ref{sec incoh}.)
Let $\Pi$ be a trivial-weight cuspidal automorphic representation of $\G'(\A)$ over $L$. We say that $\Pi$ is the cuspidal \emph{automorphic base-change} of $\pi$, and write $$\Pi\cong {\rm BC}(\pi),$$ if for every $\iota\colon L\into \C$ and every finite place $v$, we have $\iota\Pi_{v}\cong{\rm BC}(\iota\pi_{v})$.
 We say that  $\pi$ is \emph{stable} if it admits a cuspidal automorphic base-change over $L$;
we denote by 
$$\wt\sC(\G)(L)^{\rm st} \subset \wt\sC({\G})(L),  
\qquad
\sC(\G)(L)^{\rm st} \subset \sC({\G})(L)
$$
the subsets consisting of (orbits of) isomorphism classes of  representations that are stable.
 
By the definitions and Lemma \ref{equiv bc},
 the stability condition is Galois-invariant, so that the above definition makes sense {and  is compatible with Definition \ref{BC def} if $L=\C$ and $V\in \sV^{\circ, +}$.}

\subsubsection{Hermitian automorphic representations as the image of base-change} \lb{sec: ggp rat}

\begin{lemma}\lb{her bc} 
Let $\Pi$  be a complex cuspidal automorphic representation of  $\G'(\A)$.
The following are equivalent:
\begin{enumerate}
\item $\Pi$  is hermitian;
\item for every $V \in \sV$ such that $\G^{V}$ is quasi-split at all  places,
  there exists a cuspidal automorphic representation $\pi'$ of $\G^{V}(\A)$
   such that $\Pi\cong {\rm BC}(\pi')$.
\end{enumerate}
\end{lemma}
\begin{proof}
That (1) implies (2) is the automorphic descent of \cite{GRS}. The converse implication is proved in \cite{mok}.
\end{proof}

\begin{proposition} \lb{cor C her}
Assume that $F_0$ is totally real and $F$ is CM.
\begin{enumerate}
    \item For $\Pi\in \sC(\G')$ and $\e\in\{\pm 1\}$, the condition that\footnote{Note that the root number below is necessarily a sign.} 
    $$\text{$\Pi^\iota$ is hermitian and the central root number  $\eps(\Pi)\coloneqq \varepsilon(\Pi^\iota_{n}\ts \Pi^\iota_{n+1}, 1/2)=\epsilon$}$$
    is independent of the choice of an embedding $\iota\colon \Q(\Pi)\into \C$.
    \item Denote by $\sC(\G')^{\rm her, \e}$ the sub-ind-scheme consisting of those $\Pi$ satisfying the above condition, and let 
    $$\sC(\G')^{\rm her} \coloneqq \sC(\G')^{\rm her, +}\sqcup \sC(\G')^{\rm her, -}\quad \subset \quad \sC(\G').$$
 The base-change map gives an isomorphism of $\Q$-ind-schemes
\beq\lb{BC rat st}
{\rm BC}\colon \bigsqcup_{V\in \sV^{\circ}} \sC(\RH^{V}\bs \G^{V})^{\rm st} \to \sC(\G')^{\rm her},
\eeq
such that 
$$\sC(\G')^{\rm her, \epsilon}={\rm BC}\left(  \bigsqcup_{V\in \sV^{\circ, \epsilon}} \sC(\RH^{V}\bs \G^{V})^{\rm st}  \right).
$$
\end{enumerate}
\end{proposition}
{Similarly to  the local case, given $\Pi\in \sC(\G')^{\rm her}$, we will refer to its preimage $(V, \pi)$ under \eqref{BC rat st} as the \emph{distinguished element} in the global Vogan $L$-packet of $\Pi$.}
\begin{proof} 
Let $\Pi\in \sC(\G')$,  pick any embedding $\iota\colon \ov{\Q(\Pi)}\into \C$, assume that $\Pi^\iota$ is hermitian,   and let $\e\coloneqq  \varepsilon(\Pi^\iota_{n}\ts \Pi^\iota_{n+1}, 1/2)$. Note that $\Pi^\iota$ is absolutely tempered at every finite place by Proposition \ref{ramanujan}. 

Let $(V', \pi'_\iota)$ be as provided by Lemma \ref{her bc} (2). By that lemma applied to any conjugate of $\pi'_\iota$ by an automorphism of $\bC$ and Proposition \ref{equiv bc}, we see that $\Pi^{\iota'}$ is also hermitian for any $\iota'\colon \Q(\Pi)\into \C$. (This is probably an overly roundabout way of proving this assertion.)

By \cite[\S~26, discussion of Question (1)]{GGP} and Proposition \ref{local ggp},  we can find an element  $V\in \sV^{\circ, \epsilon}$ and a representation $\pi\in \sC(\RH^V\backslash\G^V)(\ov{\Q(\Pi)})$ such that $\iota {\rm BC}(\pi_{v})\cong \iota\Pi_v$ for every finite place $v$; moreover, $(V_v,\iota\pi_v)$ and $(V'_v, \pi_{\iota,v}')$ agree at all but finitely many non-archimedean places~$v$. By the special cases of base-change for real groups stated in \S~\ref{loc BC} \eqref{bc p3}-\eqref{bc p4},  the equivariance of base-change of Lemma \ref{equiv bc}, and the absence of global obstructions from \cite[Theorem 1.7.1]{KMSW} (for generic packets), we have ${\rm BC}(\pi) =\Pi$ and $\pi$ lies above a $\Q(\Pi)$-point of $\sC(\RH^V\backslash\G^V)$. This completes the proof of the proposition.
\end{proof}

\subsection{Relative traces}

\subsubsection{Definitions} 

\begin{definition} \lb{rtf defi}
 Let $L$ be a normed field. A \emph{relative-trace datum}
 $D=(\Pi_{1},  \widehat\Pi_1, \Pi_{2},  \widehat\Pi_2; \vth,\beta,  T)$ consists of:
\begin{itemize}
\item topological $L$-vector spaces $\Pi_{1}\subset \widehat\Pi_1$, $\Pi_{2}\subset \widehat\Pi_2$, where  $\Pi_i$ is dense in $\widehat\Pi_i$ for $i=1,2$;
\item a  bilinear form $\vth\colon \widehat \Pi_{1} \ot  \widehat \Pi_{2}\to L$;
\item a bilinear form $\beta\colon  \widehat \Pi_{1}\ot  \widehat \Pi_{2}\to \Gamma$, where $\Gamma$ is a finite-dimensional $L$-vector space;
\item a map $T\colon \Pi_{1} \to  \widehat\Pi_{1}$,
\end{itemize} 
satisfying:
\begin{itemize}
\item for $i=1,2$ we can write $\Pi_{i}=\varinjlim_{\lm\in \Lambda} \Pi_{i, \lm}$ as a filtered direct limit of finite-dimensional $L$-vector spaces and injective maps, in such a way that:
\begin{enumerate}
    \item[(a)] for every $\lm \in \Lambda$, $\vth_{|\Pi_{1, \lm} \ot \Pi_{2, \lm}}$ is a perfect pairing;
    \item[(b)] there exist idempotents $e_{i, \lm}\colon \Pi_i\to \Pi_{i, \lm}$ with $e_{i, \lm\vert \Pi_{i, \lm}}=\id$ such that $\vth(e_{1,\lm}v_1, v_2)= \vth(v_1, e_{2, \lm}v_2)$ for all $v_1\in \Pi_1$, $v_2\in \Pi_2$.
\end{enumerate}
\end{itemize}
Let us say that a basis $\{\phi\}$ of $\Pi_{1}$ is \emph{filtered} if there is a  presentation $\Pi_{1}=\varinjlim_{\lm\in \Lambda} \Pi_{1, \lm}$ with the above properties, such that $\{\phi\}\cap \Pi_{1, \lm}$ is a basis of $ \Pi_{1, \lm}$ for all $\lm\in \Lambda$;
if this is the case we denote by $\{\phi^{\vee}\}$ the basis of $\Pi_{2}$ whose restriction to $\Pi_{2, \lm}$ is the $\vth$-dual basis of $\{\phi\}\cap \Pi_{1, \lm}$, which exists by conditions (a) and (b) above.

We define  the \emph{trace of $T$ relative to $\beta$,  $\vth$} to be
\beq \lb{eq rtf def}
\Tr_{\vth}^{\beta}(T)\coloneqq \sum_{\phi}\beta(T\phi, \phi^{\vee}),\eeq
provided the sum is absolutely convergent and independent of the choice of a filtered  basis $\{\phi\}$ of $\Pi_{1}$.  

We say that the datum $D$ is \emph{discrete} when 
the spaces $\Pi_i =\widehat\Pi_i$ for $i=1, 2$ and the operator $T$ has finite rank, so that the sum \eqref{eq rtf def} has only finitely many nonzero terms.  In this case, we will omit listing $\widehat \Pi_{1}$ , $\widehat \Pi_{2}$ in $D$.
\end{definition}
\begin{remark} If $\Gamma=L$ and $\beta=\vth$, we recover the usual notion of trace. In the examples of interest to us:
\begin{itemize}
\item when $L$ is $\C$, $\widehat{\Pi}_1\cong \widehat{\Pi}_2 $ will be the space underlying a separable Hilbert space,  $\vth$ will be the bilinear form associated to the hermitian form, and $\Pi_i$ will be the span of a Hilbert basis;
\item 
when $L$ is not $\C$, the data will be discrete;
\item
 we will have $\beta=h\circ (P_{1}\boxtimes P_{2})$ for some  linear functionals $P_{i}\colon \widehat \Pi_{i}\to S_{i}$ valued in an $L$-vector space $S_{i}$, and some bilinear form $h\colon S_{1}\ot S_{2}\to \Gamma$. (In fact, in the first part of the paper we will only consider $S_{1}=S_{2}=\Gamma=L$, and $h$ equal to the multiplication map.) 
\end{itemize}
\end{remark}

\subsubsection{Relations between different relative traces} 
Let $D=(\Pi_{1},  \widehat\Pi_1, \Pi_{2},  \widehat\Pi_2; \vth,\beta,  T)$  be a relative-trace datum, and   let $\alpha_{2}\in\End_{L}(\Pi_{2})$. 
Let $\mu\colon \Lambda \to \Lambda $ be a strictly increasing function with cofinal image such that $\alpha_{2}(\Pi_{2, \lm})\subset \Pi_{2, \mu(\lm)}$. We define the \emph{$\vth$-transpose} of $\alpha_{2}$ to be the unique $\alpha_{2}^{\vth}\in \End_{L}(\Pi_{1})$ whose restriction to $\Pi_{1, \mu(\lm)} $ is the transpose of $\alpha_{2|\Pi_{2, \lm}}$ for the restriction of $\vth$.  

\begin{lemma} \lb{different RT} Let $D=(\Pi_{1}, \Pi_{2};  \vth, \beta, T)$ and $D'=(\Pi_{1}' , \Pi_{2}' ; \vth', \beta', T')$ be discrete relative-trace data  as in Definition \ref{rtf defi}.  
\begin{enumerate} 
\item \lb{drt0} Suppose that $\Pi_{1}=\Pi_{1}'$, $\Pi_{2}=\Pi_{2}'$, $\vth=\vth'$, $T=T'$, and 
$$\beta' =\beta\circ (1\boxtimes \alpha_{2})$$
 for some $\alpha_1\in \End_L(\Pi_1)$, $\alpha_{2}\in\mathrm{Isom}_{L}(\Pi_{2})$. Then 
$$ \Tr_{\vth'}^{\beta'}(T')
=
\Tr_{\vth}^{\beta}(\alpha_1 T\alpha_{2}^{\vth, -1})
$$
where $\alpha_{2}^{\vth}\in \End_{L}( \Pi_{1})$ is the $\vth$-transpose of $\alpha_{2}$.
\item\lb{drt1} Suppose that $\Pi_{2}=\Pi_{2}'$ and there is an $L$-isomorphism $\alpha_{1}\colon \Pi'_{1}\to \Pi_{1}$ for which
$$\vth'=\vth \circ (\alpha_{1}\boxtimes {\rm id}), \quad \beta'=\beta \circ (\alpha_{1}\boxtimes {\rm id}), \quad  T' = {\alpha_{1}^{-1}}T\alpha_{1}.$$
Then 
$$\Tr_{\vth'}^{\beta'}(T')=\Tr_{\vth}^{\beta}(T).$$
\item\lb{drt2} Suppose that for $i=1$, $2$ the spaces $\Pi_{i}'\subset \Pi_{i}$ are part of  direct sum decompositions $\Pi_i=\Pi_i'\oplus \Pi_i''$ such that $\Pi_i'$ is $\vth$-orthogonal to $\Pi_{3-i}''$  and that $\vth'\coloneqq \vth_{|\Pi_{1}'\ot \Pi_{2}'} $ is a perfect pairing (in the sense that it satisfies the condition of Definition \ref{rtf defi}). Suppose that $T(\Pi_{1})\subset \Pi_{1}'$.  If $\beta'= \beta_{|\Pi_{1}'\ot \Pi_{2}'}$ and $T'= T_{|\Pi_{1}'}$, then 
$$\Tr_{\vth}^{\beta}(T)= \Tr_{\vth'}^{\beta'}(T').$$
\end{enumerate}
\end{lemma}
The proof is elementary linear algebra.

\part{$p$-adic $L$-functions and the analytic relative-trace formula} \lb{part1}
 
We study Rankin--Selberg $L$-functions and the related Jacquet--Rallis relative-trace formulas in a sequence of  contexts. In \S~\ref{sec: JR}, we review the theory in complex coefficients. In \S~\ref{sec rat}, we construct a Jacquet--Rallis RTF in rational coefficients  and at the same time prove Theorem~\ref{thm A} on the rationality of twisted Rankin--Selberg $L$-values. The construction  relies in particular on the existence of suitable Gaussians, obtained from a refinement of the results of \cite{isolation}. 
 In \S~\ref{sec: prtf}, we construct an RTF in $p$-adic coefficients and at the same time prove Theorem \ref{thm B} on the existence of $p$-adic $L$-functions. The construction  relies on some  local theory and in particular on a suitable family of  explicit test Hecke measures at $p$-adic places: the theory is developed on the spectral side in
 \S~\ref{sec: loc rtf} (whose centerpiece is an explicit calculation of Liu--Sun) and on the geometric side in \S~\ref{app orb p} (whose centerpiece is a new explicit calculation of {orbital integrals for the aforementioned  test Hecke measures}).

\section{Jacquet--Rallis relative-trace formulas} \lb{sec: JR}
 We consider the traces of Hecke operators relative to two period functionals and the Petersson inner product on automorphic forms for $\G'$, and compare (the resulting local terms) with a parallel relative-trace distribution for $\G$.  The substance of this section is not new, rather it recalls some related work done by previous authors, particularly \cite{JR,wei-aut,isolation}. We omit detailed discussions of convergence issues, for which we refer to \cite{wei-aut} or \cite[Appendix A]{BP-comparison}.
 
 In \S~\ref{sec:31}, we define the periods on $\G'$ and the associated relative-trace distribution~$I$.  In  \S~\ref{sec:32} and  \S~\ref{sec:33}, we give the spectral and the geometric expansions of~$I$. In  \S~\ref{sec:34}, we define the  relative-trace distribution $J$ on $\G$, and in  \S~\ref{sec:35} we state the fundamental results on  the local comparison  of the two distributions. {In \S~\ref{ss:loc char} we consider ``almost unramified representations", a class of representations with the mildest ramification, as well as their local root numbers and local relative characters.}
 
\subsection{Period functionals and the distribution}  \lb{sec:31}
Let $\mathfrak z$ be the center of the universal enveloping algebra of ${\rm Lie}\,  G'_{\infty}$, and fix a maximal compact subgroup $K_\infty\subset G'_\infty$.
Let $\sA(\G')$ be the space of smooth, $\frak z$-finite, $K_\infty$-finite, moderate-growth automorphic forms on $\G'(\A)$;  let $\sA_{\rm cusp}(\G')$ be its cuspidal subspace, and let $L^2_{\rm cusp}(\G'(F_0)\backslash\G'(\A))$ be the closure of  $\sA_{\rm cusp}(\G')$
in $L^2(\G'(F_0)\backslash\G'(\A))$. We endow $L^2_{\rm cusp}(\G'(F_0)\backslash\G'(\A))$ 
  with the bilinear Petersson product
  $$\vartheta^{}(\phi_{1}, {\phi}_{2})\coloneqq 
 \int_{[\RG']}  \phi_{1}(g) {\phi}_{2}(g) \, dg.
$$

\subsubsection{Period functionals} We define two functionals on $L^2_{\rm cusp}(\G'(F_0)\backslash\G'(\A))$.

For $\chi\in Y(\C)$,  the ($\chi$-twisted) \emph{Rankin--Selberg period} is the functional
$$P^{}_{1, \chi}(\phi)\coloneqq  \int_{[\RH'_{1}]} \phi(h_{1}) \chi(h_{1})\, dh_{1}.$$
where  $\chi(h_{1})\coloneqq  \chi(N_{F/F_{0}}\det h_{1})$,

The \emph{Flicker--Rallis period} is the functional 
$$P^{}_{2}(\phi)\coloneqq  \int_{[\RH_{2}']}\phi(h_{2})\eta(h_{2})\, dh_{2}, $$
where 
$\eta(h_{2}) \coloneqq  \eta( \det (h_{n})^{n+1}\det( h_{n+1})^{n})$  if $h_{2}= ([h_{n}],  [h_{n+1}])$.

\subsubsection{Relative-trace distribution} \lb{sec I def}
 We  say that  $f'\in \sH(\G'(\A))$  is \emph{quasicuspidal} if $R(f')$ sends $\sA(\G')$ to $L^2_{\rm cusp}(\G'(F_0)\backslash\G'(\A))$ (cf. \cite[Definition 3.2]{isolation}), and we denote by $\sH(\G'(\A))_{\rm qc}$ the space of quasicuspidal Hecke measures.

 \begin{definition} We define a  relative-trace distribution on $\sH(\G(\A))_{\rm qc}\ts Y(\C)$   by 
$${\hI} ({\ff'}, \chi)\coloneqq  C\cdot  \Tr_{\vth}^{P_{1, \chi}\ot P_{2}}(R({\ff'})),$$
where  the constant
\beq\lb{def C}
C\coloneqq  {\Delta_{\G}^{}\over \Delta_{\RH}^{2}} \, {\Delta_{\RH_{1}'}^{} \Delta_{\RH_{2}'}\over \Delta_{\G'}^{}}
\eeq
is motivated by rationality considerations.

\end{definition}

We note that the above definition does fit within the setup of  Definition \ref{rtf defi}: we may write 
$$\sA_{\rm cusp}(\G')=\varinjlim_{(K,\tau,  \frak a)} \sA_{\rm cusp}(\G')^{K,\tau,  \frak a=0} \qquad \subset L^2_{\rm cusp}(\G'(F_0))$$
  as  $K$ varies among  compact open subgroups of $\G'(\A^{\infty})$, $\tau$ varies among finite-dimensional representations of $K_\infty$,
and   $\frak a$ varies among finite-codimensional  ideals in $\frak z$; here the subscript $\tau$ refers to forms whose $K_\infty$-Fourier transform has zero $\tau'$-component for every irreducible $\tau'$ that does not occur in  $\tau$. The existence of idempotents $e_{(K, \tau, \frak a)}$ with image  $\sA_{\rm cusp}(\G')^{K,\tau,  \frak a=0}$ is clear, and the relative trace is well-defined by (the proof of) \cite[Theorem 2.3]{wei-fourier}. See also \cite[Proposition 2.8.4.1]{BPCZ} for a more general result in a framework similar to ours.

\medskip

In the next two subsections we discuss the two expansions of $I$: a spectral expansion, in terms of automorphic representations, and a geometric expansion, in terms of orbits  (double-cosets).

\subsection{Spectral expansion}   \lb{sec:32}
Let $\Pi$ be a cuspidal  automorphic representation of $\G'(\A)$, which by multiplicity one  we may and do identify with a subspace of $\sA_{\rm cusp}(\G')$.
  We define a distribution on $\sH(\G'(\A))$ by
$$\hI_{\Pi}(\ff', \chi)\coloneqq C\cdot  \Tr_{\vartheta_{\Pi}}^{P_{1,\Pi, \chi}\otimes P_{2, \Pi^{\vee}}} (\Pi(\ff')),$$
where we use subscripts to indicate the restriction of period functionals and Petersson product to (the $L^2$ closures of) $\Pi$, $\Pi^{\vee}$, $\Pi\ot \Pi^{\vee}$.

We define some local periods, in order to factorize $I_{\Pi}$. 

 \subsubsection{Whittaker models and rational structures}  Let $\psi\colon F_{0}\bs \A\to \C^{\ts}$ be a nontrivial character, and let $$\psi_{F}\coloneqq  \psi({1\over 2} \Tr_{F/F_{0}}(\cdot))\colon F\bs \A_{F}\to \C^{\ts}.$$ 
  We denote by ${\rm N}_{\nu}\subset \GL_\nu$ the  upper-triangular unipotent subgroup, and  inflate $\psi_{F}$ to a character of ${\rm N}_{\nu}(\bA_F)$ by $\psi_{F, \nu}(u)=\psi_{F}(\sum_{i=1}^{\nu-1} u_{i, i+1} )$.

  Let $\Pi_{\nu}$ be an automorphic representation of $\RG_{\nu}(\A)$.  
 Its $\psi$-Whittaker model $\sW_{\psi}(\Pi_{\nu})$ is the image of the map 
\beq\sW\colon \Pi &\to C^{\infty}({\rm N}_{\nu}(\A) \bs\GL_{\nu}(\A), \psi_{F, \nu})\\
 \phi &\mapsto W_{\phi} (g)\coloneqq \int_{{\rm N}_{\nu}(\A)} \phi(ug) \overline{\psi}_{F, \nu}(u) \, du.
 \eeq
 The $\psi$-Whittaker model of $\Pi=\Pi_{n}\boxtimes \Pi_{n+1}$ is $\sW_{\psi}(\Pi)=\sW_{ \overline{\psi}}(\Pi_{n})\boxtimes \sW_{\psi}(\Pi_{n+1})$; it has a $\G'(\A)$-factorization $ \sW_{ \psi}(\Pi)= \bigotimes_{v}\sW_{ \psi_{v}}(\Pi_{v})$.

We now consider rational structures, along the lines of  \cite[\S~ 3.2]{RS}.  Let $v$ be a finite place of $F_{0}$ with underlying rational prime $\ell$, and suppose that $\Pi_{v}$ is a smooth irreducible admissible representation of $G'_{v}$  over a  subfield $L\subset \C$.  For $\sg \in \Aut(\C/L)$, let $a_{\sg}\in \Z_{\ell}^{\ts}$ be its image under the composition
$$\Aut(\C/L)\to \Gal(L(\mu_{\ell^{\infty}})/L)\to \Z_{\ell}^{\ts}$$
of the restriction and the cyclotomic character. Let $t_{\sg, \nu}\coloneqq \diag(a_{\sg}^{\nu-1}, \ldots, 1)$ and let $t_{\sg}\coloneqq (t_{\sg, n}, t_{\sg, n+1})\in G'_{v}$.
 Then we may define an action of $\Aut(\C/L)$ on $\sW_{\psi_{v}}(\Pi_{v}\ot_{L}\C)$ by 
 \beq\lb{act W}
 W^{\sg}(g)\coloneqq  \sg( W(t_{\sg}^{-1}g)); \eeq
 we will denote by $\sW_{\psi_{v}}(\Pi_{v})$ the space of $\Aut(\C/L)$-invariants; it is an $L[G_{v}']$-module satisfying  $\sW_{\psi_{v}}(\Pi_{v})\ot_{L}\C\cong \sW_{\psi_{v}}(\Pi_{v}\ot_{L}\C)$ (see \cite[Lemma 3.2]{RS}).
 
 \subsubsection{Factorizations of the periods and Petersson product} \lb{sec fact}
 For the following factorization results, see \cite[\S~3]{wei-aut} and references therein.
Let $\epsilon'_{\nu}(\tau)\coloneqq \diag(\tau^{\nu+\epsilon -1},\tau^{\nu+\epsilon -2}, \ldots, \tau^{\epsilon -1})\in \GL_{\nu}(F_{v})$, where $\epsilon \in \{0, 1\}$ has the same parity as $\nu$.

For $W=W_{n}\ot W_{n+1}\in \sW_{\psi,v}(\Pi_{v})$, define\footnote{Our definition of $P_{2, v}$ differs form the one of \cite{wei-aut} by the factor $\ep({1\over 2}, \eta, \psi)^{n+1\choose 2}$, cf.  \S~\ref{sp ma} below.}
\beq\lb{local periods}
P_{1, \Pi_{v}, \chi_{v}}(W) &\coloneqq {\ep({1\over 2} , \chi_{v}^{2}, \psi_{v})^{{n+1\choose 2}}\over  L(1/2, \Pi_{v} \ot\chi_{v})} \, 
\int_{\RN_{n}(F_{v}) \bs\GL_{n}(F_{v})} W(j_{1}(h_{1})) \chi_{v}(h_{1})\, d^{\natural}h_{1} ,\\
P_{2,\Pi_{v}}(W) &\coloneqq {\ep({1\over 2}, \eta_{v}, \psi_{v})^{n+1\choose 2}
\over  L(1, \Pi_{v} , \As^{-\star})}\,  P^{\sharp}_{2, \Pi_{n, v} }(W_{n})  P^{\sharp}_{2, \Pi_{n+1, v} }(W_{n+1}), \qquad \text{where}\\
P^{\sharp}_{2, \Pi_{\nu, v}}(W_{\nu}) &\coloneqq \displaystyle{\int\limits_{\RN_{\nu-1}(F_{0,v})\bs \GL_{\nu-1}(F_{0,v})}} W_{\nu}\left(\twomat {\epsilon'_{\nu-1}(\tau)  h_{2, \nu-1}} {}{}1\right)\, \eta_{v} (\det h_{2,\nu-1})^{\nu-1}\, d^{\natural}h_{2, \nu-1}
\eeq
and  $L(1, \Pi_{v}, \As^{-\star})= \prod_{\nu=n}^{n+1}L(1, \Pi_{\nu, v}, \As^{(-1)^{\nu-1}})$. 

For $W\in \sW_{\psi_{v}}(\Pi_{v})$, $W^{\vee}\in \sW_{\overline{\psi}_{v}}(\Pi_{v}^{\vee})$, define
\beqq
\vth_{\Pi_{v}}(W, W^{\vee}) &\coloneqq 
\prod_{\nu=n}^{n+1}  L(1, \Pi_{\nu, v}\ts \Pi_{\nu, v}^{\vee})^{-1} \\
& \int\limits_{\RN_{\nu-1}(F) \bs \GL_{\nu-1}(F)} W_{\nu}\left(\twomat { g_{ \nu-1}} {}{}1\right)W_{\nu}^{\vee}\left(\twomat {g_{\nu-1}}{}{}1\right)\, d^{\natural}g_{\nu-1}.
\eeqq

\begin{remark}\lb{aa1}
With our normalizations, when all the data are unramified and $W(1)=W^{\vee}(1)=1$, we have  
\beqq P_{1, \Pi_{v}, \chi_{v}}(W) =P_{2,\Pi_{v}}(W)= \vth_{\Pi_{v}}(W, W^{\vee})=1.\eeqq
Moreover,  the three functionals are rational in the  following sense.
\end{remark}

\begin{lemma} \lb{aa1 L}
 If $\Pi_{v}$ is defined over  a subfield $L\subset \C$, for every $\sg\in \Aut(\C/L)$ we have
\beq \lb{aa1 eq}
P_{1,\Pi_{v}, \chi_{v}^{\sg}}(W^{\sg})&=\sg P_{1, \Pi_{v}, \chi_{v}}(W), \\
\quad P_{2,\Pi_{v}}(W^{\sg})&=\sg P_{2,\Pi_{v}}(W),\\
\quad \vth_{\Pi_{v}}(W^{\sg}, W^{\vee, \sg})& =\sg \vth_{\Pi_{v}}(W, W^{\vee}).\eeq
\end{lemma}
\begin{proof}
First note that  $ L(1, \Pi_{\nu, v}\ts \Pi_{\nu, v}^{\vee})$ and $\As^{(-1)^{\nu-1}})$ belong to $L$, and  $L(1/2, \Pi_{v} \ot\chi_{v}^{\sg})=\sg L(1/2, \Pi_{v} \ot\chi_{v}) $ for every $\sg\in \Aut(\C/L)$. This is immediate from writing those $L$-values in terms of $L$-values \emph{at integers} obtained from the  Weil--Deligne representations $r_{\nu, v}$ over $L$  associated with $\Pi_{\nu, v}$ by the local Langlands correspondence: using the normalization $r_{\nu, v}= \sL_{\nu, F_{v}}(\Pi_{\nu, v})$ of  \cite[\S~2.3]{Shin} (to which we also refer for background), we have for instance  
$$L(1/2, \Pi_{v} \ot\chi_{v})=L(1/2 -(1-(n+1)n)/2, r_{n,v}\ot r_{n+1,v}\ot \chi_{v}).$$

It then remains to verify the analogue of \eqref{aa1 eq} for the products of the integrals with the epsilon factors.
This follows by applying \eqref{act W} and a change of variable.
\end{proof}

We may now state the factorization (see \cite[\S 3]{wei-aut} and \S~\ref{sec comp wei} below): for any  $\phi\in \Pi$ with factorizable $\psi$-Whittaker function $W=\ot_{v}W_{v}\in \sW_{\psi}(\Pi)$, and $\phi^{\vee}\in \Pi^{\vee}$ with factorizable $\overline{\psi}$-Whittaker function $W^{\vee}=\ot_{v}W^{\vee}_{v}\in \sW_{\overline{\psi}}(\Pi^{\vee})$, we have
\beq\lb{fct}
P_{1,\Pi,  \chi}(\phi) &={ L(1/2, \Pi\ot \chi)\over  \Delta_{\RH_{1}'}^{} \cdot  \ep({1\over 2},  \chi^{2}, \psi)^{{n+1\choose 2}}} \prod_{v}P_{1, \Pi_{v}, \chi_{v}} (W_{v}) \\
P_{2, \Pi}(\phi^{\vee}) &= { n(n+1)L^{*}(1, \Pi, \As^{-\star}) \over  \Delta_{\RH_{2}'}}  \prod_{v}P_{2,\Pi_{v}} (W^{\vee}_{v}) \\
\vth_{\Pi}(\phi, \phi^{\vee}) &= {4 n(n+1)L^{*}(1, \Pi \ts \Pi^{\vee}) \over  \Delta_{\G'}^{}} \prod_{v}\vth_{\Pi_{v}}(W_{v}, W_{v}^{\vee}),
\eeq
where $4$ is the Tamagawa number of $\G'$. In the factorization of $P_{2}$, we have used that $\ep({1\over 2}, \eta, \psi)=1$.

\subsubsection{Local relative character} We define a  character
 $$I^{}_{\Pi_{v}}   (f_{v}', \chi_{v}) = I_{\Pi_{v}}  (f_{v}', \chi_{v}, \psi_{v}) =  \Tr_{\vth_{\Pi_{v}}}^{P_{1,\Pi_{v},\chi_{v}}\otimes P_{2, \Pi^{\vee}_{v}}}
( R(f_{v}'))$$
  on $\sH(\G'_{v})$.

\subsubsection{Comparison with the normalization of \cite{wei-aut}} \lb{sec comp wei} 
Let 
$$\wt{\G}'\coloneqq  \G_{n}'\times \G_{n+1}',$$ and let us identify representations of $G_{v}'$ with representations of $\wt{G}_{v}'$ whose central character is trivial on $(F_{0,v}^{\ts})^{2}$.  Let  ${\rm p}\colon \wt{G}'_{v}\to G_{v}'$ be the projection, and let 
\beq
{\rm p}_{*} \colon  \CS(\wt{G}'_{v}) & \to \CS(G_{ v}') \\
\df' &\mapsto \left(g= [\wt{g}] \mapsto \int_{F_{0, v}^{\ts, 2}} \df'(z\wt{g}) \, d^{*}_{}z\right).
\eeq
In \cite{wei-aut}, one considers a distribution $\wt{I}^{}_{\Pi_{v}}$  on $\CS(\wt G'_{v})$ (denoted there by $I^{\natural}_{\Pi_{v}}$), and a global distribution $\wt{I}_{\Pi}$ on $\CS(\wt\G'(\A))$ (denoted there by $I_{\Pi}$).\footnote{Strictly speaking only $\chi_{v}=\one_{v} $ is considered in \cite{wei-aut}, but the definition remains valid in our more general case too. When this is again the case in the rest of the paper, we will simply cite \cite{wei-aut} without repeating this remark.} Suppose that $f'=\ot_{v}f'_{v}\in \sH(\G'(\A))$ and  $\df' =\ot_{v}\df'_{v}\in \CS(\wt{\G}'(\A))$ are related by
\beq \lb{f and df} f'_{v}= {\rm p}_{*}(\df'_{v}) \, d^{*}_{}g_{v}.\eeq
 Then 
\beq\lb{compare I wei}
I_{\Pi_{v}}^{}(f_{v}, \chi_{v}) & =   {\Delta_{\RH_{1}'^{\ad}, v}^{} \Delta_{\RH_{2}'^{\ad}, v}^{} \over   \Delta_{\G^{\ad}, v}^{}}  
\cdot \left(\ep(\textstyle{1\over 2}, \eta_{v}, \psi_{v}) \ep(\textstyle{1\over 2}, \chi_{v}^{2}, \psi_{v})\right)^{{n+1\choose 2} }\, 
\cdot  \wt{I}^{}_{\Pi_{v}}(\df_{v}, \chi_{v}),\\ 
\hI_{\Pi}(f, \chi) &= {1 \over \vol([\RZ_{\G'}], d^{*}z)}   {\zeta_{\G'}^{*}(1)\over \zeta_{\RH'_{1}}^{*}(1) \zeta_{\RH'_{2}}^{*}(1)}\,  \wt{I}_{\Pi}(\df, \chi) 
= {1\over 4}  {1\over \zeta_{\RH'_{1}}^{*}(1) \zeta_{\RH'_{2}}^{*}(1)}\,  \wt{I}_{\Pi}(\df, \chi), 
\eeq
where the factor $\vol([\RZ_{\G'}], d^{*}z)=4L(1, \eta)^{2}= 4\zeta_{\G'}^{*}(1)$ accounts for the fact that the Petersson 
 product in \cite{wei-aut} is defined via integration on $[\G'^{\ad}]=[\wt{\G}'^{\ad}]$ and not $[\G']$.

\subsubsection{Factorization of the relative character} Define 
$$\sL(1/2, \Pi, \chi)\coloneqq  {\Delta_{\RG}^{}\over \Delta_{\RH}} {  L(1/2 , \Pi\ot\chi)\over {L(1, \Pi, {\As}^{\star})}},$$
which agrees with the definition made in the introduction as noted in \S~\ref{sec: tam}.
 We use the analogous notation relative to the constituents $\Pi_{v}$, $\chi_{v}$ for $v$ a finite place of $F_{0}$ or $v=\infty$.

\begin{proposition}\lb{factor I}
For all $\ff'=\ot_{v}\ff'_{v}\in \sH(\G'(\A))$,  there is a factorization
\beqq
\hI_{\Pi}(\ff', \chi)
&=  {1\over 4 }
{ \sL(1/2, {\Pi,\chi}) 
\over \Delta_{\RH}\cdot \ep(\textstyle{1\over 2} , \chi^{2})^{{n+1\choose 2}} }
 \prod_{v} \hI_{\Pi, v}^{}(f'_{v}, \chi_{v}).
\eeqq
\end{proposition}
\begin{proof} 
Using \eqref{compare I wei},  the factorization in \cite[Proposition 3.6]{wei-aut} is equivalent to 
$$ \hI_{\Pi}(\ff', \chi)= C{1\over 4} {\Delta_{\G'}\over \Delta_{\RH'_{1}} \Delta_{\RH'_{2}}} 
\cdot
 \ep({\textstyle{1\over 2}} , \chi^{2})^{-{n+1\choose 2}} 
\cdot 
{  L(1/2 , \Pi\ot\chi)\over {L(1, \Pi, {\As}^{\star})}}
\cdot  \prod_{v} \hI^{}_{\Pi_{v}}(\ff'_{v}, \chi_{v}).$$
By the definition of $C$ in \eqref{def C} and of $\sL$, this is equivalent to the asserted formula. (Equivalently, the factorization follows from \eqref{fct}.)
\end{proof}

\subsubsection{Spectral expansion} 
 Let $\ff' \in \sH(\G'(\A))$ be   quasicuspidal, and let $\chi\in Y(\C)$.
By   \cite[Prop. 4.1]{isolation} (where it is assumed that $\chi=1$, but the proof extends to the general case),  we have
\beq\lb{eq sp exp}
\hI(\ff', \chi)  = \sum_{\Pi} \hI_{\Pi} (\ff', \chi) 
\eeq
where the sum,   running over the cuspidal hermitian automorphic representations of $\G'(\A)$, is absolutely convergent.

\subsection{Geometric expansion}   \lb{sec:33}
The distribution $I$ also admits an expansion as a sum of orbital integrals, which we review.

\subsubsection{Orbit varieties} \lb{sec: orb'}
 Let 
$${\RB}' \coloneqq  \RH'_{1}\bs \G'/\RH'_{2}$$ 
be the categorical quotient, which is an affine variety over $F_{0}$, cf. \cite{wei-fourier}. 
Let 
\beqq  {\RS}\coloneqq  \{\gamma \in \G'_{n+1}\ \ |\  \gamma {\gamma}^\mathrm{c}=1_{n+1}\} .\eeqq
The maps 
\beq \lb{g *}  g=(g_{n}, g_{n+1})\mapsto g_{\star} \coloneqq  g_{n}^{-1}g_{n+1}, \qquad  g_{n+1} \mapsto g_{n+1}{g}_{n+1}^{\rm c, -1}\eeq
induce maps and isomorphisms 
$$s\colon \wt\G' \to \G'_{n+1}/\G'_{n+1, 0} \cong \RS, \qquad \RB'\cong \G'_{n,0} \bs \G'_{n+1}/ \G'_{n+1,0}\cong  \G'_{n,0} \bs \RS.$$
The second map in \eqref{g *} also  yields  a bijection on $F'$-points $\G'_{n+1}(F') /\G'_{n+1,0}(F') \to \RS(F')$ for any field $F'\supset F_{0}$.

Representing a point of $\RB'$ by a matrix in $\RS$, written in  block decomposition with the upper-left block of size $n\ts n$, the  invariant map
\beq\lb{eq inv}
 {\rm inv}\colon \RB' &\to \Res_{{F}/{F_{0}}} \BA^{2n+1}\\
\twomat Abcd &\mapsto ((\Tr(\wedge^{i}A))_{i=1}^{n}, (cA^{j-1} b)_{j=1}^{n}, d)
\eeq
gives an embedding into affine space.

\subsubsection{Regular, plus-regular, and semisimple orbits}\lb{sss:reg}
We define three functions  on $\RS$ (hence on $\wt\G'$) by\footnote{In the first formula, and in similar situations later, the external transpose turns a row vector whose entries are row vectors into a column vector whose entries are row vectors.}
\beq \lb{disc rs}
D^+(s)&\coloneqq \det (e_{n+1}^{\rm t} ,e_{n+1}^{\rm t} s,\dots, e_{n+1}^{\rm t} s^{n})^{\rm t} \\
D^-(s)&\coloneqq \det (e_{n+1}, se_{n+1},\dots, s^n e_{n+1}) \\
D(s)&\coloneqq  \det ((e_{n+1}^{\rm t} s^{i+j}e_{n+1})_{0\leq i, j\leq n}) = D^{+}(s)D^{-}(s),
\eeq
where $e_{n+1}= (0, \ldots, 0, 1)^{\rm t}\in F^{n+1}$. We denote by $\wt\G'_{\rs}\subset \wt\G'_{\reg^\pm} \subset \wt\G'$  the open subschemes  defined, respectively,  by  $D\neq 0$ and $D^\pm\neq 0$, 
and by 
$$\G'_{\rs}\subset \G'_{\reg^\pm} \subset \G'$$
the respective images in $\G'$; thus $\G'_{\rs}=\G'_{\reg^{+}}\cap \G'_{\reg^{-}}$.  The function $D$ descends to $\RB'$ and we denote by $\RB'_{\rs}$ its non-vanishing locus, whose preimage in $\G'$ is $\G'_{\rs}$.

\begin{remark} \lb{invol +-}
The involution $g^{\diamond}\coloneqq  g^{\rm c, -1, \rm t} $ on $\wt\G'$ satisfies $D^{\pm}(s(g^{\diamond}))=D^{\mp}(s(g))$, and it descends to $\G'$; in particular, it swaps $\G'_{\reg^{+}}$ and $\G'_{\reg^{-}}$. 
\end{remark}

 Let $F'\supset F_{0}$ be a field. 
An $\RH'_{1}(F')\times\RH'_{2}(F')$-orbit in $\G'(F')$ is said to be \emph{regular} if its stabilizer is trivial; \emph{semisimple} if the orbit is Zariski-closed. 

The regular semisimple orbits in $\G'(F')$ are in bijection with  $\RB'_{\rs}(F')$, and the preimage  in $\G'(F')$ of any $\gamma\in \RB'_{\rs}(F')$ consists of a single orbit. The preimage of a general $\gamma\in \RB'(F')$  contains finitely many regular orbits (but possibly infinitely many orbits), of which exactly one belongs to $\G'_{\reg^{+}}(F')$ and exactly one belongs to $\G'_{\reg^{-}}(F')$
 (these two coincide precisely when $\gamma\in \RB'_{\rs}(F')$).  We will call the elements in  $\G'_{\reg^+}$ (respectively  $\G'_{\reg^-}$) {\em plus-regular}  (respectively {\em minus-regular}). We refer to \cite[\S~2.4]{Lu} for more details.

\subsubsection{Local  orbital integrals} 
 Let $v$ a place of $F_{0}$, and let  $\gamma' \in G_{\reg^{\pm}, v}'$. Then for all $\chi_{v}\colon F_{0,v}^{\ts}\to \C^{\ts}$, we define the orbital integral
\beq \lb{orb int def}
\hI^{\sharp}_{\gamma'}( \ff'_{v}, \chi_{v}) & \coloneqq 
 \int_{H'_{1,v}}\int_{H'_{2,v}} \ff'_{v}(h_{1,v}^{-1}\gamma'  h_{2,v}) \chi_v(h_{1,v}) \eta(h_{2,v}) \,  {d^{\natural}{h_{1,v}}d^{\natural}h_{2,v}\over d^{\natural}g_{v}},\eeq
  where we recall that $f'_{v}/d^{\natural}g_{v}$ is a function. If $\gamma' \in G_{\rs, v}'$ or $f'_{v}$ is supported in the regular locus of $G_{v}'$, the integral is absolutely convergent. In general, Lu proved that the integral is convergent when $\chi$ is the product of a unitary character and $|\cdot |^{s}$ for $\pm\Re(s)<-1$ \cite[Lemma 5.14]{Lu}, and  gave the following regularization.

\begin{proposition}[Lu] \lb{lu reg}
Let $\sR_{v, 0}^{\pm}$ be the set of functions, on the space of  smooth characters of $F_{0,v}^{\ts}$, of the form
\beq \lb{L factor}
\chi_{v}\mapsto\prod_{j=1}^{m}L(\pm1\mp j, (\chi_{v}\eta_{v})^{\mp j}\circ N_{F_{0}'/F_{0}}),
\eeq
for varying integers $m\geq1$ and finite field extensions $F_{0}'/F_{0}$; 
and let $\sR_{v}^{\pm}$ be the set of finite products of functions in $\sR_{v,0}^{\pm}$. 

Let $\gamma' \in G_{\reg^{\pm}, v}'$. Then one can define an element 
$$L_{\gamma'}\in \sR_{v}^{\pm},$$
 such that the following hold. 
\begin{enumerate}
\item $L_{\gamma'}$ only depends on the $H_{1,v}\ts H_{2,v}'$-orbit of $\gamma'$, and it   equals  $1$ if and only if $\gamma'\in G'_{\rs, v}$,
  \item\lb{lu aa1} For unramified data (the precise meaning is given in \cite[Prop. 5.12 (4)]{Lu}), we have 
  \beq \lb{I unram}
  \hI_{\gamma'}^{\sharp}( \ff'_{v}, \chi_{v})=L_{\gamma'}(\chi_v).
  \eeq 
  \item \lb{lu reg 3} Define a  normalized orbital integral by 
   \beq \lb{n orb int def}
  \hI_{\gamma'}^{}( \ff'_{v}, \chi_{v})\coloneqq \frac{\hI_{\gamma'}^{\sharp}( \ff'_{v}, \chi_{v})}{L_{\gamma'}(\chi_v)}.
  \eeq
  Then for every fixed character $\chi_{v}^{\circ}$ of $F_{0, v}^{\ts}$:
  \begin{enumerate}
  \item \lb{lu reg 3a} if $v$ is archimedean, the function $s\mapsto \hI_{\gamma'}^{}( \ff'_{v}, \chi_{v}^{\circ}|\cdot|^{s})$ extends to an entire function on $\C$;
   \item \lb{lu reg 3b} if $v$ is non-archimedean,  the function
    $$\chi\mapsto \hI_{\gamma'}^{}( \ff'_{v}, \chi_{v}^{\circ}\chi)$$
    on the space of  unramified characters of $F_{0,v}^{\ts}$ is a polynomial in $\chi(\varpi_v)$, $\chi(\varpi_v)^{-1}$, whose coefficients belong to the field of rationality of $\ff'$ and $\chi_0$.
  \end{enumerate}
  \end{enumerate}
\end{proposition}

\begin{proof}
All references in this proof point to \cite{Lu}. The definition of $L_\gamma$ is given in \S~5.3, and its properties follow from the results of \S\S~5.1-5.2 and Remark 5.13. Parts \eqref{I unram} and \eqref{lu reg 3a} are part of Proposition 5.12.
The assertion that the normalized orbital integrals in part  \eqref{lu reg 3b} are Laurent polynomials is  Proposition 5.16. (Strictly speaking Lu considers  $\chi_v^\circ$ to be trivial  and $\chi=\vert \cdot \vert^s$, but his proof applies to our context as well.) 
Finally, the rationality statement on the coefficients of part \eqref{lu reg 3b}  can be verified by tracing through the proof of  Proposition 5.16., which uses operations that preserve the field of coefficients to  reduce to the consideration of  Tate integrals.
\end{proof}

\begin{remark} 
In most of \cite{Lu},  all regular orbits are treated, and we expect that all the results of the present paper involving plus-regular orbital integrals could be extended to general  regular orbits as well. Nevertheless, for our purposes in the general construction of the $p$-adic relative-trace formula, we will only need to consider plus-regular orbital integrals.
For this reason,  we will  restrict our attention to plus-regular orbits, which  introduces some simplifications. 
 (In fact, for the applications in the  proofs of our main theorems, we will not  need to consider any  regularized divergent orbital integrals; however we consider the more general $p$-adic relative-trace formula to be of independent interest.)
\end{remark}
  
  Let  $v$ be a place of $F_{0}$, and let $\gamma'\in G'_{\reg^{+},v}$. For $h_{1}\in H_{1,v}'$, $h_{2}'\in H_{2, v}'$, we have $I_{h_{1}\gamma' h_{2}} (-, \chi_{v})= \chi_{v}(h_{1})\eta_{v}(h_{2}) I_{\gamma'}(-, \chi_{v})$.  We then add  a renormalization factor to the orbital integral so that, when $\chi_{v}=\one$, it only depends on the orbit   of $\gamma'$.
   Let $\eta'\colon F^{\ts}\bs \A_{F}^{\ts}\to \C^{\ts}$ be a character such that $\eta'_{| \A^{\ts}}=\eta$. 
   With the notation $\gamma_{\star }'$ and $s=s(\gamma')$ as in \eqref{g *}, 
     we define  
\beq\lb{eq: tr f}
\kappa_{v}(\gamma')&\coloneqq  \eta'\left(\det(\gamma'_{\star})^{\epsilon} \det s^{-(n+\epsilon)/2} D^{+}(s)\right)\\
&=
 \eta'\left(\det(\gamma'_{\star})^{\epsilon} \det s^{-(n+\epsilon)/2}\det (e_{n+1}^{t},e_{n+1}^{t}s, \ldots, e_{n+1}^{t}s^{n})^{\rm t}  \right)
 \eeq
where $\epsilon \coloneqq 0$ if $n$ is even, $\epsilon \coloneqq 1$ if $n$ is odd.  
This equals  the transfer factor  denoted by $\Omega_{v}$ in  \cite[(4.12)-(4.13)]{wei-aut} (cf. \S~\ref{compare geo match} below), and it   satisfies
$$\kappa_{v}(h_{1}\gamma' h_{2})=
\eta_{v}(h_{2})\kappa_{v}(\gamma', \chi_{v}),$$ and 
\beq\lb{prod kappa}
\prod_{v}\kappa_{v}(\gamma')=1\eeq
 for all $\gamma'\in \G'_{\reg^+}(F_{0})$.  We also record the following rationality property.

\begin{lemma} \lb{kappa i} Let { $\gamma'\in G'_{v, \reg^+}$}.  Then $\kappa_{v}(\gamma')$ is a square root of  $\eta_{v}(-1)^{-{n+1\choose 2}}$.
\end{lemma}
\begin{proof} With notation as above (but dropping the apex from $ \gamma_{\star}'$  for lightness), write $n=2m+\epsilon$ and let 
 $$a\coloneqq  \det (e_{n+1}^{t},\ldots, e_{n+1}^{t}s^{n})^{\rm t}  \det s^{-(m+\epsilon)}  \det(\gamma_{\star})^{\epsilon}
 = \det (e_{n+1}^{t} s^{-m-\epsilon}, \ldots, e_{n+1}^{t} s^{m} )^{\rm t}   \det \gamma_{\star}^{\epsilon},
  $$
which satisfies $\kappa_{v}(\gamma')=\eta_{v}'(a)$. 
Using ${s}^{\rm c}=s^{-1}$ and $\gamma_{\star}^{\rm c}=s ^{-1}\gamma_{\star}$, we find
  $${a}^{\rm c}=(-1)^{n+1\choose 2} \det  (e_{n+1}^{t} s^{-m-\epsilon}, \ldots, e_{n+1}^{t} s^{m} )^{\rm t}  \det \gamma_{\star}=(-1)^{n+1\choose 2} a$$
  where $(-1)^{n+1\choose 2}$ is the sign of the longest permutation on ${n+1}$ elements. The assertion of the lemma follows.
\end{proof}
We now define, for  $\gamma \in B'_{v}$,
\beq\lb{eq orb on B}
L_\gamma (\chi_{v})&\coloneqq L_{\gamma'}(\chi_{v})\\
\hI^{\sharp}_{\gamma}(f_{v}', \chi_{v})& \coloneqq  \kappa_{v}(\gamma') \,\hI^{\sharp}_{\gamma'}(f_{v}', \chi_{v})\\
\hI^{}_{\gamma}(f_{v}', \chi_{v})& \coloneqq  \kappa_{v}(\gamma') \,\hI^{}_{\gamma'}(f_{v}', \chi_{v}).
\eeq
for any $\gamma'$ in the unique plus-regular orbit above $\gamma$. When $\chi_{v}=\one$, it is straightforward to check that the right hand side is independent of the choice of $\gamma'$; in general, our notation is somewhat abusive, but the ambiguity is canceled out in the global context as discussed next.

\subsubsection{Global orbital integrals}  
Let $\sR_{0}^{+}$ be the set of functions on Hecke characters of $F_{0}$ of the form 
\beq\lb{pre DR}\chi\mapsto\prod_{j=1}^{m}L(1- j, (\chi\eta)^{- j}\circ N_{F_{0}'/F_{0}}),
\eeq
for varying integers $m\geq1$ and finite field extensions $F_{0}'/F_{0}$, and  let $\sR^{+}$ be the set of finite product of functions in $\sR_{0}^{+}$. 

For any $\gamma\in \RB'(F_{0})$, by \cite[\S6]{Lu} we can define an element $L_{\gamma}= \prod_{v}L_{\gamma_{v}}$ in $\sR$. For 
  any   $\ff'\coloneqq  \ot_{v} \ff'_{v}
   \in\sH(\G'(\A))$, and any character  $\chi\in Y(\C)$, we put
\beq\lb{prod orb int}
\hI_{\gamma}(\ff', \chi) &\coloneqq  
C{\Delta_{\G'} \over \Delta_{\RH_{1}} \Delta_{\RH_{2}'}}  \,  L_{\gamma}(\chi)\prod_{v}\hI^{}_{\gamma}(\ff'_{v}, \chi_{v})  
={\Delta_{\G}\over \Delta_{\RH}^{2}}  L_{\gamma}(\chi)\prod_{v }\hI_{\gamma}^{}(\ff'_{v}, \chi_{v}),
\eeq
where all but finitely many factors equal~$1$ (the finite set of exceptions depends on $\gamma$); we take the orbital integrals in the product to be defined as in \eqref{eq orb on B} 
 by means of a common \emph{rational} plus-regular  lift  $\gamma' \in \G'(F_{0})$ of $\gamma$, which ensures that the product is well-defined. 
  When $\gamma\in \RB'_{\rm rs}(F_{0})$, it is clear that we have
\beqq\lb{loc orb lc}
\hI_{\gamma}(\ff', \chi) &=  C\cdot \int_{\RH'_{1}(\A)}\int_{\RH'_{2}(\A)} \ff'_{}(h_{1}^{-1}\gamma'  h_{2}) \chi(h_{1}) \eta(h_{2}) \,  {d^{}{h_{1}}d^{}h_{2}\over d^{}g_{}}.\eeqq

\subsubsection{Comparison with the normalization of \cite{wei-aut}} In  \cite[\S 4]{wei-aut}, one considers the distribution on Hecke functions (and not measures) on $\wt G'_{\rs, v}$, defined by 
$$\wt{I}_{\gamma'}(\df'_{v} , \chi_{v}) \coloneqq 
 \int_{H_{1,v}} \int_{ H_{2,v}} \df'_{v}(h_{1,v}^{-1}\gamma' h_{2,v}) \chi(h_{1,v}) \eta(h_{2,v}) \,  {d^{*}{h_{1,v}}d^{*}h_{2,v}}$$
 (and denoted there by $O(\gamma, f')$),
 and the global analogue 
 $$\wt{I}_{\gamma'}(\df' , \chi) \coloneqq  \prod_{v} \wt{I}_{\gamma'}(\df'_{v} , \chi_{v}).$$
 If $f'=\ot_{v}f'_{v}\in \sH(G'_{v})$ and $\df'=\ot_{v}\df'_{v}\in \CS(\wt{G}'_{v})$ are related by \eqref{f and df},
then 
\beq\lb{compare I gamma} I^{}_{\gamma}(f'_{v}, \chi_{v}) &= \kappa_{v}(\gamma') \,  {\Delta_{\RH_{1}'^{\ad}, v}^{} \Delta_{\RH_{2}'^{\ad}, v}\over   \Delta_{\G'^{\ad}, v}} \, 
\wt{I}_{\gamma'}(\df'_{v}, \chi_{v}),\\
\hI_{\gamma}(\ff', \chi)  &= C {\zeta_{\G'}^{*}(1) \over \zeta_{\RH_{1}'}^{*} (1)\zeta_{\RH_{2}'}^{*} (1) }\prod_{v} \wt{I}_{\gamma'}(\df'_{v} , \chi).
\eeq

We will also denote by
\beq \lb{pr Hecke mu}
{\rm p}_{*} \colon  \sH(\wt{G}_{v}) & \to \sH(G_{ v}) 
\eeq
the pushforward map of Hecke measures.

\subsubsection{Relative-trace formula for ${\hI}$} We describe the spectral and geometric expansions of $\hI$.
 For $S$ a finite set of places of $F_{0}$ and $? \in \{ \rs, \mathrm{reg}^{\pm}\}$, a Hecke measure $f'^{S} \in \sH(\G'(\A^{S}))$  is said to have \emph{locally ?-support} if it is in the span of those $\ot_{v\notin S} f'_{v}$ such that for some place $v$, $f'_{v}$ is supported on $G'_{?, v}$.  We introduce a weaker notion.
 
\begin{definition}\lb{def: weak support}
We say that the Hecke measure $f'^{S}\in  \sH(\G'(\A^{S}))$ has \emph{globally $?$-support} if it belongs to the subspace 
spanned by those pure tensors $\ot_{v \notin S} f_{v}'$ such that for every $\gamma'\in \G'(F_{0})-\G'_{?}(F_{0})$, there is some $v\notin S$ such that $f_{v}'$ vanishes on the $H_{1,v}\times H_{2,v}$-orbit of  $\gamma'$. 

For the sake of readability we will  also write \emph{globally regular semisimple} (respectively, \emph{plus-regular, minus-regular}) \emph{support} for \emph{globally $\rs$-} (respectively, $\mathrm{reg}^{+}$-, $\mathrm{reg}^{-}$) \emph{support}.
\end{definition}

\begin{proposition}\lb{exps rtf} Let $\ff' \in \sH(\G'(\A))$ be   quasicuspidal with globally plus-regular support.
 Then for every character $\chi\in Y(\C)$, we have
\beqq
\sum_{\Pi} \hI_{\Pi} (\ff', \chi)  = \hI(\ff', \chi) =\ \sum_{\gamma \in {\RB}'(F_{0})} \hI_{\gamma}(\ff', \chi).
\eeqq
where both sums are absolutely convergent,  the first one running over the cuspidal hermitian automorphic representations of $\G'(\A)$. 
\end{proposition}

\begin{proof} 
The spectral expansion was already observed in \eqref{eq sp exp}.
The geometric expansion is \cite[Theorems 3.1, 6.1]{Lu}; the definition of the summands in \emph{loc. cit.} contains extra terms corresponding to regular but non-plus-regular orbits, but those vanish by our assumption on $f'$. 
\end{proof}

 Let us motivate Definition \ref{def: weak support}, and at the same time make more precise part of the outline in \S~\ref{intro test}. 
 For the proof of Theorem \ref{main thm}, we will need to  consider an RTF, similar to that of Proposition \ref{exps rtf}, in which we can only understand the terms corresponding to regular semisimple global orbits.\footnote{In fact, there are some technical difficulties in even defining that RTF when non-regular-semisimple orbits may be involved, see Proposition \ref{prop:de I}.}   We thus need to input an element $f' \in \sH(\G'(\A))$ annihilating the contribution of all non-regular semisimple orbits, and enjoying some other favorable properties.
  Whereas we are not able to find an $f'_v$ with rs-support and controlled level for any finite place $v$ -- hence no global $f’$ with rs-support in the usual sense --, for split places $v$ we can find an $f’_{v, \pm}$ with $\text{reg}_\pm$-support, controlled level, and pleasing spectral properties (Proposition \ref{test fin} \eqref{test fin 3}). This suffices to produce a global $f’$ with \emph{globally} regular semisimple support by taking components $f’_{v_+, +}f’_{v_-, -}$ at two different split places $v_+, v_-$.

 \subsection{Relative traces for unitary groups}  \lb{sec:34}
 We review the Jacquet--Rallis RTF for unitary groups.

\subsubsection{Orbit spaces} Let $v$ be a finite  place of $F_{0}$ or $v=\infty$, and
 recall from \S~\ref{coh incoh G} the set $\sV_{v}$ of  pairs of hermitian spaces. For $V\in \sV_{v}$, let
\beq \lb{BV}
B_{ v, V}  \coloneqq  H_{v}^{V} \bs G_{v}^{V} / H_{v}^{V} ,\eeq
which is isomorphic to the quotient of  $U(V_{n+1})$ by the adjoint action of $U(V_{n})$ via the map $[(g_{n}, g_{n+1})]\mapsto [g_{n} g_{n+1}^{-1}]$. 
Differently from the linear case, $B_{ v, V} $ is a \emph{subset}  (open for the $v$-adic topology if $v$ is non-archimedean) of the set of $F_{0,v}$-points of $\RB_{v}^{V} \coloneqq  \RH^{V}_{v}\bs\G^{V}_{v}/ \RH^{V}_{v}$.  Similar to \S~\ref{sss:reg}, we say that $g\in G_{v}^{V}$ is {\em regular}  (for the $H_{v}^{V}\ts H_{v}^{V}$-action) if its stabilizer is trivial; \emph{semisimple} if  its orbit is closed. We denote by $G^{V}_{\rs, v}\subset G^{V}_{v}$ the (Zariski-open) subset of  regular semisimple elements and by  $B_{\rs, v, V}$ its image in $B_{v, V}$.  When $v=\infty$  and $V_{\infty}=V_{\infty}^{\circ}$ is the positive 
definite pair, every orbit is semisimple, and we denote $B_{\infty}^{\circ}\coloneqq   B_{\infty, V_{\infty}^{\circ}} = B_{\rs, \infty, V_{\infty}^{\circ}}$.

Consider now the global case, and let $V\in \sV$. We similarly define $\RG_{\rs}^{V}\subset \RG^{V}$ to be the sub-group-scheme of those $g$ with closed orbit and trivial stabilizers for the $\RH^{V}\ts \RH^{V}$-action. For uniformity of notation, we denote  by
\beq \lb{Brs not}
\RB_{\rs}(F_{0})_{V}\subset \RB_{\rs}^{V}(F_{0})\eeq
the image of $\RG^{V}_{\rs}(F_{0})$.

\subsubsection{Local distributions} \lb{loc J def}
Let $\delta \in B_{\rs, v, V}$ and let  $\delta'\in G_{\rs, v}^{V}$ be a preimage of $\delta$. We define
a  local orbital-integral distribution  $ J_{\delta,v}= J_{\delta,v}^{V}$ on the Hecke algebra of $G_{v}=G_{v}^{V}$ by
\beq
\lb{J delta v}
J_{\delta,v}^{}(f_{v}) & \coloneqq  \int_{H_{v}} \int_{H_{v}} f_{v}(x^{-1} \delta' y)\, {d^{\nat}x d^{\nat}y \over d^{\nat}g}.\eeq

For $\pi_{v}$ a representation of $G_{v}=G_{v}^{V}$, we define a relative  character ${J}^{}_{\pi_{v}}= {J}_{\pi_{v}}^{V}$ on $\sH(G_{v})$ by
\beq\lb{J pi v}
{J}^{}_{\pi_{v}}(f_{v}) & \coloneqq 
\sL(1/2, {\rm BC}(\pi_{v}))^{-1} 
 \int_{H_{v}}  \Tr_{\pi_{v}}( \pi_{v}(h) \pi_{v}(f)) \, {d^{\nat}h}.
 \eeq
 It vanishes unless $\pi_v$ is distinguished.

By our choices of measures, for all finite $v$, if $f_{v}$ is $L$-valued (for some subfield $L\subset\C$) then so are $J_{\delta,{v}}(f_{v})$, $J_{\pi_{v}}(f_{v})$.

\subsubsection{Comparison with the normalization of \cite{wei-aut}} Let $f_{v}\in \sH(G_{v})$ and $\df_{v}\in \CS(G_{v})$ be related by 
\beq\lb{f and df G}
f_{v}= \df_{v} d^{*}g.
\eeq
\begin{enumerate}

\item Let $\wt{J}^{}_{\pi_{v}}$ be the relative  character on $\CS(G_{v})$  defined in \cite[(1.8)]{wei-aut} (using the measure $d^{*}h_{v}$ on $H_{v}$ as in \S~4 \emph{ibid.}), and denoted there by $J_{\pi_{v}}^{\natural}$. Then 
\beq\lb{compare J pi Wei}
 {J}^{}_{\pi_{v}}(f_{v})= D_{v}^{1/2} \Delta_{\RH_{v}^{\ad}}^{}
  \, \wt{J}^{}_{\pi_{v}}(\df_{v}),\eeq
\item
Let $ \wt{J}_{\delta}$ be the orbital integral distribution on $\CS(G_{v})$  defined in \cite[(4.2)]{wei-aut}, and denoted there by $O(\delta, \cdot )$. Then 
\beq\lb{compare J delta}
J^{}_{\delta}(f_{v})= {(\Delta_{\RH^{\ad}, v})^{2}\over \Delta_{\G^{\ad}, v}} \, \wt{J}_{\delta}(\df_{v}).
\eeq
\end{enumerate}

\subsubsection{Global relative-trace formula} Let now $V\in \sV$ be coherent, and let $\G=\G^{V}$.
 Let $\vth\colon \sA_{\rm cusp}(\G)\otimes \sA_{\rm cusp}(\G)\to \C$ be the Petersson inner product (with respect to the Tamagawa measure on $[\G]$), and consider the $\RH$-period
\beq \label{H per}
P=P^{V}\colon  \sA_{\rm cusp}(\G) &\to \C \\
\phi&\mapsto\int_{[\RH]} \phi(h)\, dh
\eeq

We define the following distributions on (subspaces of) $\sH(\G(\A))$:
\begin{itemize}
\item let  $\sH(\G(\A))_{\rm qc}\subset \sH(\G(\A))$ be the  quasicuspidal subspace (defined as in \S\ref{sec I def}). For $f\in \sH(\G(\A))_{\rm qc}$, we define
$$J(f)\coloneqq  \Tr_{\vth}^{P\ot P}(R(f));$$
\item let  $\pi$  be a cuspidal automorphic representation of $\G(\A)$. For $f\in \sH(\G(\A))$,  we define
$$J_{\pi}(f)\coloneqq  \Tr_{\vth_{\pi}}^{P_{\pi}\ot P_{\pi}}( \pi(f));$$
\item let $\delta\in \RB_{\rm rs}(F_{0})$. For $f=\ot_{v}f_{v}\in \sH(\G(\A))$, we define
\beq \lb{J delta global} J_{\delta} (f) \coloneqq 
{\Delta_{\G}\over (\Delta_{\RH})^{2}}\,
\prod_{v} J_{\delta, v}^{}(f_{v}) 
= \prod_{v\nmid \infty} J_{\delta,v}(f_{v}) \cdot J^{\circ}_{\delta,\infty}(f_{\infty}).\eeq
\end{itemize}

Analogously to Proposition \ref{exps rtf}, we have the 
 spectral and geometric expansions (\cite[Proposition A.2.1]{BP-comparison})
\beq\lb{J rtf} \sum_{\pi } J_{\pi}(f)= J(f)=
\sum_{\delta\in \RB_{\rm rs}(F_{0})} J_{\delta} (f),\eeq
valid whenever $f\in \sH(\G(\A))$ is quasicuspidal with globally regular semisimple support (in the analogous sense as to Definition \ref{def: weak support}), where the first sum runs over cuspidal automorphic representations of $\G(\A)$.

However, unlike the factorization 
$$\hI_{\Pi}(\ff', \chi)
=  {1\over 4 }
\displaystyle{ \sL(1/2, {\Pi,\chi}) 
\over \Delta^{}_{\RH}\cdot \ep(\textstyle{1\over 2} , \chi^{2})^{{n+1\choose 2}} }
 \prod_{v} \hI_{\Pi, v}^{}(f'_{v}, \chi_{v})$$
of Proposition \ref{factor I}, 
 the analogous factorization 
$$J_{\pi}= {1\over 4   \Delta^{}_{\RH}} \sL({{1/ 2}} , \Pi, \one)
\prod_{v} J_{\pi_{v}}$$
for a stable cuspidal tempered representation $\pi$  of $\G(\A)$
 is highly nontrivial, and equivalent to the Ichino--Ikeda conjecture for unitary groups \cite[Conjecture 1.1]{wei-aut}, whose proof  is completed in \cite{isolation}. The proof, which we briefly review in \S~\ref{sec: II}
 below (for expository purposes), goes through a  comparison of local orbital integrals $I_{\gamma,v}$ and $J_{\delta,v}$ and of local relative characters $I_{\Pi_{v}}$ and $I_{\pi}$.
 We first review the main results on the  local  comparison, which are equally  important in the arithmetic setting.

 \subsection{Comparison of the local distributions}
 \lb{sec:35}

\subsubsection{Spectral matching}\lb{sp ma}
 Let $v$ be a place of $F_{0}$.
For  $V\in \sV_{v}$ and $\pi^{V}_{v}\in {\rm Temp}(G_{v}^{V})$, define  a spectral transfer factor
\beq\lb{tr f}\kappa({\pi^{V}_{v}})= \kappa(\pi^{V}_{v}, \psi_{v}, \tau) 
\coloneqq   \eta_{v}'((-1)^{n+1}\tau)^{n+1\choose 2} \cdot \eta_{v}({\rm disc}(V_{n}))^{n}\cdot \omega_{\pi^{V}_{v}}(-1);\eeq
this is the same  as in \cite[Conjecture 4.4]{wei-aut} with the correction of \cite[Remark 5.52]{BP-Planch}, up to a factor $\ep({1\over 2}, \eta_{v}, \psi)^{n+1\choose 2}\, $.\footnote{To compare the last factor in \eqref{tr f}
 with \cite{wei-aut},  recall that $\omega_{\Pi_{v}}(z)= \omega_{\pi_{v}}(z/z^{c})$, so that $\omega_{\pi_{v}}(-1)=\omega_{\Pi_{v}}(\tau)$.  The absence of the factor $\ep({1\over 2} ,\eta_{v}, \psi_{v})^{n+1\choose 2}$, which  cancels out its presence in our local Flicker--Rallis period $P_{2, \Pi_{v}}$, is helpful in Lemma \ref{rat kappa}.}

 Let $S$ be a finite set of places of $F_{0}$.
Denote $\sV_{S}\coloneqq \prod_{v\in S} \sV_{v}$; for $V=(V_{v})_{v\in S}\in \sV_{S}$, denote ${\rm Temp}({(H_{S}^{V})}\bs G_{S}^{V})=\prod_{v\in S}{\rm Temp}({(H_{v}^{V_{v}})} \bs G_{v}^{V_{v}})$; for $\pi_{S}^{V}\in {\rm Temp}(G_{S}^{V})$, set $\kappa(\pi_{S}^{V})\coloneqq \prod_{v\in S}\kappa(\pi_{v})$ and $\hJ_{\pi_{S}}\coloneqq \ot_{v\in S} \hJ_{\pi_{v}}$. For $\Pi_{S} \in {\rm Temp}(G_{S}')\coloneqq \prod_{v\in S}{\rm Temp}(G_{v}')$, let $I_{\Pi_{S}}\coloneqq \ot_{v\in S} I_{\Pi_{v}}$. 

\begin{definition}
We say that Hecke measures $\ff'_{S}\in \sH(G_{S}')$ and $(\ff^{V}_{S})_{V\in \sV_{S}}\in \prod_{V\in \sV_{S}}\sH(G_{S}^{V})$ \emph{match spectrally}  if  for all $V\in \sV_{S}$ and all $\pi^{V}_{S}\in {\rm Temp}({H_{S}^{V}}\bs G^{V}_{S})$, 
we have
\beq \lb{mat sp}
\hI^{}_{{\rm BC}(\pi^{V}_{S})}(\ff'_{S}, \one)  = \kappa(\pi^{V}_{S}) \hJ_{\pi_{S}}(\ff_{S}^{V}).\eeq
\end{definition}

\subsubsection{Geometric matching}\label{sss:orb match}

Let us first recall the matching of orbits for $\G'$ and $\G$; for the details, see \cite[\S~2.1]{wei-afl}.  Let $V\in \sV_{v}$. Orbits $\gamma \in B'_{\rs,v }$ and $\delta \in B_{\rs, v, V}$ are said to \emph{match} if a lift  $\gamma'\in \RS_{v}\subset \GL_{n+1}(F_{v})$  of $\gamma$ and a lift $\delta'\in {U}(V_{n+1})\subset \GL_{n+1}(F_{v})$ of $\delta$ are conjugate for the adjoint action of $\GL_{n}(F_{v})$. (This notion is independent of the choices of the lifts and of the basis of $V_{n+1}$ giving the embedding 
${U}(V_{n+1})\subset \GL_{n+1}(F)$.) The matching relation defines a bijection (an isomorphism of $F_{0,v}$-manifolds if $v$ is non-archimedean)
\beq \lb{dec Brs}
\ul{\delta}
\colon
B'_{\rs, v} \cong \bigsqcup_{V\in \sV_{v}} B_{\rs, v,V}.\eeq
If $S$ is a finite set of places of $F_{0}$, by taking products we obtain a matching bijection
\beqq 
\ul{\delta}
\colon
B'_{\rs, S}\cong  \bigsqcup_{V\in \sV_{S}} B_{\rs, S,V}.\eeqq
where $B'_{\rs, S}\coloneqq  \prod_{v\in S} B'_{\rs, v}$, $B_{\rs, S, V}\coloneqq  \prod_{v\in S} B_{\rs, v, V_{v}}$.

For the number field $F_{0}$ and for $V\in \sV$, with the notation of  \eqref{Brs not} we have an analogous bijection 
\beq\lb{dec Brs gl}
\ul{\delta}
\colon
\RB'_{\rs}(F_{0}) \cong \bigsqcup_{V\in \sV} \RB_{\rs}(F_{0})_{V}
\eeq
 compatible with \eqref{dec Brs}

We  will also denote by $B'_{\rs, v, V} $ (respectively $B'_{\rs, S, V}$, $B'_{\rs}(F_{0})_{V}$), 
  the preimage of $B_{\rs,v,V}$ (respectively $B_{\rs, S, V}$, $B_{\rs}(F_{0})_{V}$) under 
$\ul{\delta}$.

\begin{definition} Let $S$ be a finite set of places of $F_{0}$. 
We  say that Hecke measures $\ff_{S}'\in \sH(G_{v}')$ and $(\ff_{S}^{V})_{V} \in \prod_{V\in \sV_{S}} \sH (G^{V}_{S})$  \emph{match geometrically}  if
\beq\lb{mat ge}
\hI_{\gamma, S}(\ff_{S}', \one_{S})= \hJ_{\delta, S}(\ff^{V}_{S})
\eeq
whenever $\gamma\in B'_{\rs, S}$ and $\delta\in B_{\rs,S,V}$ match.
\end{definition}

\subsubsection{Comparison with the normalization of \cite{wei-aut}}
 \lb{compare geo match} 
Let $v$ be a place of $F_{0}$.  Suppose that  $f_{v}'$ is related to $\df_{v}'$  as in \eqref{f and df} and $f_{v}$ is related to $\df_{v}$ as in  \eqref{f and df G}. Let 
\beq\lb{cv}
c_{v}\coloneqq 
   {\Delta_{\RH_{1}'^{\ad}, v} \Delta_{\RH_{2}'^{\ad}, v}\over   \Delta_{\G'^{\ad}, v}}
\cdot   {\Delta_{\G^{\ad}, v} \over\Delta_{\RH^{\ad}, v}^{2}}.
\eeq
\begin{enumerate}
\item By \eqref{compare I wei}, \eqref{compare J pi Wei}, the Hecke measures
$f'_{v}$ and $(f_{v}^{V})$ match spectrally if and only if 
$ c_{v}  \df_{v}' $ and $(\df_{v}^{V})$
match spectrally in the sense (\cite[Conjecture 4.4 and last equation on p. 566]{wei-aut}) that 
 $$\wt{I}_{\Pi_{v}}^{}(c_{v}\df_{v}') 
 =\kappa(\pi^{V}_{v}) {D_{v}^{-\dim \RG/2}\Delta_{\RG, v}\over D_{v}^{-\dim \RH/2}\zeta_{\RH, v}(1) \Delta_{\RH, v}} \cdot \wt{J}^{}_{\pi_{v}}(\df_{v}^{V}).$$
\item
By  \eqref{compare I gamma} and  \eqref{compare J delta}, the Hecke measures
$f'_{v}$ and $(f_{v}^{V})$ match geometrically if and only if 
$c_{v} \df' $ {and} $ (\df_{v}^{V})$
match  geometrically in the sense of \cite[(4.14)]{wei-aut}, namely 
$$\kappa_{v}(\gamma')^{-1} \, \wt{I}_{\gamma'}(c_{v}\df', \one)=\wt{J}_{\delta}(\df_{v})$$  
for all matching pairs of orbits $(\gamma, \delta)$. 
\end{enumerate}

\subsubsection{Main results on the local comparisons}
Each of the  following propositions is a  deep result.

\begin{proposition}[Equivalence of   spectral and geometric matching] \lb{equiv match} 
Let $S$ be a finite set of places of $F_{0}$. The pairs 
 $f_{S}'\in \sH(G_{S}')$ and  $(f_{S}^{V})\in \prod_{V\in \sV_{S}}\sH(G^{V}_{S})$   match spectrally if and only if they match geometrically.
\end{proposition}
\begin{proof}  The proof of \cite[Lemma 4.9]{isolation},  based on  \cite{BP-Planch, BP-comparison}, applies. (As noted in \cite[Remark 4.10]{isolation}, in general this relies on \cite{mok, KMSW}.)   Note that by  \S~\ref{compare geo match}, the  comparisons of matchings in \emph{loc. cit.}, whose conventions are inherited from \cite{wei-aut}, are compatible with ours. 
\end{proof}

From now on we will simply say that $f'_{S}$ and $(f_{S}^{V})$ \emph{match} when they match spectrally and geometrically. For a fixed $V\in \sV_{S}$, we say that $f_{S}'$ \emph{purely matches} $f_{S}^{V}$ if it matches $(f_{S}^{V}, (0^{V'})_{V'\neq V})$.

 \begin{proposition}[Existence of matching {\cite[Theorem 2.6]{wei-fourier}}] 
 Let $v$ be a finite place of $F_{0}$. For every $f'\in \sH(G_{v}')$, a matching $(f^{V})\in \prod_{V\in \sV_{v}}$ exists; conversely, for every $(f^{V})\in \prod_{V\in \sV_{v}}$, a matching  $f'\in \sH(G_{v}')$ exists.
 \end{proposition}

 We will need the following two explicit matching statements.

\begin{proposition}[Fundamental Lemma {\cite{Yun, BP-FL}}] \lb{JRFL} Let $v$ be a finite place of $F_{0}$ that is unramified in $F$,   let $V\in \sV_{v}$ be the unramified pair of hermitian spaces, and recall the unit Hecke measures from \eqref{unit hecke}. 

The unit Hecke measure $f_{v}'^{\circ}$ on $G_{v}'$ purely matches the unit measure  $f_{v}^{\circ}$ on $G_{v}^{V}$.
\end{proposition}

\begin{proposition}[Vertex parahoric transfer, {\cite[Theorem 1.1]{zhiyu}}] \lb{zhiyu tr}
    Let $v$ be an inert place of $F_0$ with odd residue characteristic, let $\e=\e_v\in \{0,1\}$ have the same parity as $v(\nrm)$, let $V_v\in \sV_v$, and let $K_v=K_{n, v} \ts K_{n+1, v}\subset G_v^{V_v}$ be the stabilizer of a vertex lattice $\Lm=\Lm_n\ts \Lm_{n+1}$ of type $(t, t+\e)$ (\S~\ref{sss:vert}).
    Let $K_v'\subset G_v'$ be the stabilizer of the lattice chain $(\Lm\subset \Lm^\vee)$. Then $e_{K_v}$ purely matches $e_{K_v'}$.
\end{proposition}

A matching result for archimedean places will be proved in \S~\ref{sec E gauss}. Finally, we will also need to note the following (relatively easy) fact.
\begin{lemma}[{\cite[Proposition 2.5]{wei-fourier}}] \lb{match split} Let $v=w \ol{w}$ be a split place of $F_{v}$, and identify $G_{v}\cong G_{n,0}'\ts G_{n+1, 0}'$. 
Then $f_{v}'={\rm p}_{*}(f_{w}'\ot f_{\ol{w}}') \in \sH(G_{v}')$ matches $f_{v}\coloneqq  f_{w}\star f_{\ol w}^{\vee} \in \sH(G_{v})$.
\end{lemma}

\subsection{Almost unramified representations} \lb{ss:loc char}

For this subsection, we fix a coefficient field $L$. 
 All representations of (products of) general linear groups under our consideration will be over~$L$.

\subsubsection{Representations with conductor at most $1$}
\lb{sss:cond01}
We now allow $F_v$ to be any non-archimedean local field of characteristic zero. 

For $\Pi_{\nu,v}$  a generic irreducible representation of $\GL_\nu(F_{v})$ {over $\C$},
we recall that  Jacquet, Piatetski-Shapiro, and Shalika \cite{jpps-c} have defined the conductor $c(\Pi_{\nu,v})$ of $\Pi_{\nu,v}$: it is the smallest integer $c\geq 0$ such that  $\Pi_{\nu,v}^{K'_{\nu,v}(c)}\neq 0$, where $K'_{\nu,v}(c)\subset \GL_\nu(\sO_{F_v})$ is the subgroup  consisting of those $k$ whose reduction modulo $\varpi_v^c$ lies in the mirabolic subgroup (that is, $k\mod \varpi_v^c$ fixes the row vector $(0, \ldots, 0, 1)\in (\sO_{F_v}/\vpi_v^c)^\nu$). {Moreover, they proved that for $c=c(\Pi_{\nu, v})$, the space  $\Pi_{\nu,v}^{K'_{\nu,v}(c)}$ of \emph{new vectors} is $1$-dimensional.

Now let $\Pi_{\nu,v}$ be  an irreducible representation of $\GL_\nu(F_{v})$ over a coefficient field~$L$. We say that $\Pi_{\nu, v}$ is generic if it is associated with a product of unlinked segments in the sense of \cite[Theorem 9.7]{Zel}; by \emph{loc. cit.}, if $L=\C$ this condition is equivalent to the usual notion. Suppose that  $\Pi_{\nu,v}$ is generic.
By the uniqueness of new vectors over $\C$ and \cite[Lemme I.1]{wald-a}, the value $c(\Pi_{\nu, v}\ot_{L, \iota} \C)$ is independent of the choice of an embedding $\iota L\into \C$; we again call this value the conductor of~$\Pi_{\nu,v}$.}

\begin{lemma}\label{lem: c=1}
Let $\Pi_{\nu,v}$ be a generic irreducible representation of $\GL_\nu(F_{v})$ with conductor~$1$. {Assume that the coefficient field $L$ is algebraically closed.} Then the $L$-parameter of  $\Pi_{\nu,v}$ (a Weil--Deligne representation over $L$) must have one of the following two forms:
\begin{itemize}
\item
a direct sum of (one dimensional) characters  
$$\bigoplus_{i=1}^\nu\chi_i,
$$ where exactly one of them is ramified with conductor~$1$, {and the other ones are unramified};
\item a direct sum 
$${\rm St}_2\otimes\chi+ \bigoplus_{i=1}^{\nu-2}\chi_i,
$$ where all $\chi$ and $\chi_i$ are unramified characters, and
 ${\rm St}_2$ is the Weil--Deligne representation corresponding to the Steinberg representation of $\GL_2(F_v)$ (namely, ${\rm St}_2$ is trivial on the Weil group ${\rm W}_{F_{v}}$ of $F_v$ and the standard 2-dimensional representation  ${\rm Std}$ of $\SL_2$).
 \end{itemize}
 
\end{lemma}

\begin{proof}
 Suppose that the Weil--Deligne representation is of the form 
$$
M=\bigoplus_{i\geq 0} M_i \otimes \Sym^i  
$$
where $\Sym^i= \Sym^i({\rm Std})$   is the irreducible $\SL_2$-representation of dimension $i+1$, and where $M_i$ is a semisimple complex representation of the Weil group ${\rm W}_{F_{v}}$. By \cite[formula after equation (10) on p.441]{GR},  the conductor of $M$ equals
$$
\sum_{i\geq 0}[(i+1)\, a(M_i)+i\,\dim M_i^{I_{F_v}}]
$$
where $a(M_i)$ is the Artin conductor of $M_i$, and $I_{F_v}$ denotes the inertia subgroup of ${\rm W}_{F_v}$.
Note that, by the local Langlands correspondence for $\GL_\nu(F_{v})$, the conductor of $\Pi_{\nu,v}$ is equal
to the one defined for the Weil--Deligne representation in \cite[\S2.2]{GR}. We note that $$
a(M_i)\geq \dim M_i/M_i^{I_{F_v}}
$$
by \cite[last line on p.435]{GR}, hence 
$$
a(M_i)+\dim M_i^{I_{F_v}}\geq \dim M_i.
$$
It follows that, for the conductor to be equal to~$1$, we must have $M_i=0$ for all $i\geq 2$.
There are two cases:

Case (1): $M_1=0$. Then we have $1=a(M_0)\geq \dim M_0/M_0^{I_{F_v}}$, and hence $M_0^{I_{F_v}}=M_0$ or 
$M_0=M_0^{I_{F_v}}+ \chi$ for a $1$-dimensional character 
$\chi$. In the former case $M_0^{I_{F_v}}=M_0$, we have $a(M_0)=0$; in the latter case $a(M_0)=a(\chi)$, which gives the first case in the assertion of the lemma.

Case (2): $\dim M_1=1$. Then $a(M_1)=0$, hence $M_1$ is unramified; then $1=\dim M_1^{I_{F_v}}=\dim M_1 $. We must have $a(M_0)=0$, and hence $M_0$ is unramified (and semisimple) as desired.
\end{proof}

\subsubsection{Local root numbers}
We now return to the setup where $F/F_0$ is a CM quadratic extension. Let $v$ be an inert finite place of $F_{0}$.
We now calculate the local root number $\eps(\Pi_{v}) =\eqref{def local rn}$  for representations   $\Pi_v=\Pi_{n,v}\boxtimes\Pi_{n+1,v}$ with the mildest ramification.
We denote 
$$\Pi_{\textup{even},v}\coloneqq   \Pi_{2 \lfloor(n+1)/2\rfloor,v},\quad \Pi_{\textup{odd},v}=\Pi_{2\lfloor n/2\rfloor+1,v}.
$$
\begin{lemma}\label{lem: loc ep}
Let $\Pi_v=\Pi_{n,v}\boxtimes\Pi_{n+1,v}$ be generic hermitian irreducible representations of $G'_v=\GL_{n}(F_v)\times \GL_{n+1}(F_v)$. Assume that  both $\Pi_{n,v}, \Pi_{n+1,v}$ have trivial central characters and  conductor at most~$1$. 
Then
$$\eps(\Pi_{v})=-1$$
if and only if $\Pi_{\textup{even},v}$ has conductor~$1$.
\end{lemma}

\begin{proof}
{After extending scalars through an arbitrary embedding, we may and do assume that the coefficient field is $\C$.}
For the sake of lightness we write $\psi$ for $\tau\psi_{F,v}$  (cf. \S\ref{sss tau psi}), and we follow the conventions in \cite[\S4, 5, 6]{GGP} (except that they use $\sigma$ to denote the conjugation ${\rm c}$). Let $M$ (resp., $N$) denote the $L$-parameter of  $\Pi_{\textup{even},v}$ (resp., $\Pi_{\textup{odd},v}$). Then $M$ and $N$ are conjugate-selfdual representations of the Weil--Deligne group ${\rm WD}_{F_v}={\rm W}_{F_v}\times\SL_2(\C)$  of $F_v$ with sign $b(M)=-1 $ and $ b(N)=1$.

First we assume that $c(\Pi_{\textup{even},v})=0$. Then $M$ is unramified and we may write it as $P+\,^{\rm c} P^\vee$. We have 
$$\eps(M\otimes N,\psi)
=\eps(P\otimes N+\,^{\rm c} P^\vee\otimes N,\psi)
=\eps(P\otimes N+\,^{\rm c} (P\otimes N)^\vee,\psi)
=1,
$$ where the last step follows from  \cite[Proposition 5.1~(2)]{GGP}.

It remains to consider the case  $c(\Pi_{\textup{even},v})=1$. Since it has the trivial central character and the sign $b(M)=-1$, by Lemma \ref{lem: c=1},
 its $L$-parameter has the form $M={\rm St}_2+(P+\,^{\rm c} P^\vee)$. (Note that here we are using the unramifiedness of $F_v/F_{0,v}$ to see that only the trivial character of $F_v^{\times}$ is unramified and conjugate-selfdual with sign $+1$.) A similar argument using \cite[Proposition 5.1~(2)]{GGP} shows
 $$\ep(M\otimes N,\psi)
=\ep({\rm St}_2\otimes N, \psi).
$$ 
If $N$ is unramified, then it is necessarily of the form $\one +(Q+\,^{\rm c} Q^\vee)$ where $\one$ denotes the trivial one-dimensional representation of ${\rm WD}(F_v)$. Then we have 
$$\eps(M\otimes N,\psi)
=\eps({\rm St}_2\otimes \one, \psi)
=-1.
$$ 
If $N$ has trivial central character and  conductor~$1$ (and it is conjugate-selfdual of sign $b(N)=1$), it is necessarily of the form ${\rm St}_2\otimes \mu+\one
+(Q+\,^{\rm c} Q^\vee)$
where $\mu$ is the unique extension of the character $\eta_{F_v/F_{0,v}}$ to an unramified character of $F_v^\times$. Then we have
$$\eps(M\otimes N,\psi)
=\eps({\rm St}_2\otimes {\rm St}_2\otimes \mu, \psi)  \eps({\rm St}_2\otimes \one ,
\psi)
=-\eps({\rm St}_2\otimes {\rm St}_2\otimes \mu, \psi).
$$ 
Now it is a routine exercise to check using the definition of the local root number (e.g., in the beginning of \cite[\S5]{GGP}) that $\eps({\rm St}_2\otimes {\rm St}_2\otimes \mu, \psi)=1$. The proof is complete.
\end{proof}

\begin{remark}When $v$ is ramified,
one can similarly compute the local root number when the conductors are at most~$1$; however, the result will depend on more refined information of the $L$-parameter: the parameter  
$M$ may be of the form $M={\rm St}_2\otimes\chi+(P+\,^{\rm c} P^\vee)$ where $\chi: F_v^\times\to L^\times$ is either of the two unramified characters of order two.
\end{remark}

\begin{corollary}\lb{coro root n all}
    Assume that $F/F_0$ is unramified, and let $\Pi  =\Pi_n\boxtimes \Pi_{n+1}\in \sC(\G')^{\rm her}$ be a representation such that,
   for every inert finite place $v$ of $F_0$, both $\Pi_{n,v}$ and $\Pi_{n+1, v}$ have trivial central characters and conductor at most~$1$.        
 Let $\Sigma(\Pi)$ be the set of inert places $v$ of $F_0$ such that $\Pi_{\textup{even},v}$ has conductor~$1$.
 Then 
  $$\ep(\Pi)=(-1)^{\vert \Sigma(\Pi)\vert}.$$
\end{corollary} 
\begin{proof} Since $F/F_0$ is unramified, we have  $1=\eta(-1)=(-1)^{[F_0:\Q]}$, which shows that $[F_0:\Q]$ is even; therefore, the local root numbers (cf. \eqref{eq:eps v=infty}) at archimedean places multiply  to~$1$. Then the result follows from Lemma \ref{lem: loc ep}.
\end{proof}

 \subsubsection{Almost unramified representations} \lb{sss:almostun}
Let $v$ be an inert place in $F/F_0$.  Assume that the coefficient field $L$ is algebraically closed. An irreducible representation $\pi_{\nu,v}$ of $U(V_{\nu,v})$ is called {\em almost unramified} if it has a non-zero  vector invariant under the stabilizer of a vertex lattice of type $1$ or $\nu-1$. 
An irreducible representation   $\Pi_{\nu,v}$ of $\GL_\nu(F_v)$ is  called {\em almost unramified} if it is the base change of an almost unramified representation of a unitary group. {The notion is invariant under the action by automorphisms of~$L$; therefore it passes to Galois quotients.}

\begin{lemma}\label{lem: a unram}
Assume that the  residue characteristic of $v$ is odd. Let $\Pi_{\nu,v}$ be a tempered hermitian representation of $\GL_\nu(F_v)$. The following conditions are equivalent:
\begin{enumerate}
\item $\Pi_{\nu,v}$ is almost unramified.
\item The central character of $\Pi_{\nu,v}$ is trivial and the conductor of $\Pi_{\nu,v}$ is~$1$ (resp., at most~$1$) when $\nu$ is even (resp., odd).
\end{enumerate}
\end{lemma}
\begin{proof}
We first prove that
$ (1)\imp(2)$. {Let $V_{\nu, v}$ be a hermitian space of dimension $\nu$ with  a vertex lattice $\Lambda_{\nu,v}$ of type $1$ (the case of type $\nu-1$ is similar and omitted here), and let $K^{[1]}_{\nu,v}$ be the stabilizer of $\Lambda_{\nu,v}$. Let $\pi_{\nu,v}$ be a tempered almost unramified representation of  $U(V_{\nu,v})$, with  $\pi_{\nu,v}^{K^{[1]}_{\nu,v}}\neq 0$. Note that the Hasse invariant of $V_{\nu, v}$ is $-1$, hence when $\nu$ is even the unitary group $U(V_{\nu,v})$ is non-quasi-split.} Let $\Pi_{\nu,v}$ be the base change of $\pi_{\nu,v}$.  Since the center of $U(V_{\nu,v})$ is contained in $K^{[1]}_{\nu,v}$,  the central character of  $\pi_{\nu,v}$ must be trivial, hence so is the  central character of $\Pi_{\nu,v}$.  

Next we show that the conductor of $\Pi_{\nu,v}$ is at most~$1$, or equivalently, $\Pi_{\nu,v}^{K_{\nu, v}'}
\neq 0$
for some parahoric subgroup  $K_{\nu,v}'\subset \GL_\nu({F_v})$ of type $(\nu-1, 1)$ (that is,  $K_{\nu,v}'$  is conjugate to $\sO_{F, v}^\ts K'_{\nu, v}(1)$ where $K'_{\nu, v}(1)$ is an in \S~\ref{sss:cond01} and $\sO_{F, v}^\ts$ is the center of $\GL_\nu(\sO_{F_v})$).
{We may and do extend scalars via an arbitrary embedding of the coefficient field into~$\C$}. By \cite[\S5.1]{Dang}  we can choose the following auxiliary data:
\begin{itemize}
\item
 a hermitian subspace $V_{\nu-1, v}\subset V_{\nu,v}$ of dimension $\nu-1$ with a selfdual lattice $\Lambda_{\nu-1,v} \subset \Lambda_{\nu,v}$, with the stabilizer of $\Lambda_{\nu-1,v}$ denoted by $ K_{\nu-1}^{[0]}$, and
\item a tempered unramified  representation $\pi_{\nu-1, v}$ of $U(V_{\nu-1, v})$, whose base change we denote by $\Pi_{\nu-1,v}$,
 \end{itemize}
 such that the local relative character 
 $$J_{\pi_{\nu-1, v}\times \pi_{\nu,v}} (e_{ K_{\nu-1}^{[0]}\times K^{[1]}_{\nu,v} })\neq 0.$$ (Here $\nu$ plays the role of $n+1$ in our usual notation.) By Proposition  \ref{zhiyu tr}, the test measure $e_{ K_{\nu-1}^{[0]}\times K^{[1]}_{\nu,v} }$ purely matches $e_{K'_{\nu-1,v} \ts K'_{\nu,v}}$ 
 where $K'_{\nu-1,v}\subset\GL_{\nu-1}({F_v})$ is a maximal compact subgroup 
 and $K'_{\nu,v}\subset\GL_{\nu}({F,v})$ is a parahoric of tpye $(\nu-1, 1)$.
 Hence by Proposition \ref{equiv match}, 
 the local relative character $$I_{\Pi_{\nu-1,v}\times \Pi_{\nu,v}}( e_{K'_{\nu-1,v}\times K'_{\nu,v}})\neq 0.$$ It follows that $\Pi_{\nu,v}^{K_{\nu,v}'} \neq 0$, as desired.

 When $\nu$ is even, $\pi_{\nu,v}$ is a representation of a non-quasi-split unitary group, hence its base change $\Pi_{\nu,v}$ cannot be unramified. Therefore $\Pi_{\nu,v}$ must have conductor exactly~$1$.

Let us now prove that
 $ (2)\imp(1)$. The case where the conductor is equal to~$1$ follows from \cite[Theorem 1.1]{Dang} (the case where $\nu$ id odd also follows from \cite{AOY}). We note that in \cite{Dang} it is assumed that the residue field has at least $\nu$ elements; this is used to establish explicit smooth transfers when the conductor is $>1$ and is not needed when the conductor is $\leq 1$. Now let $\nu$ be 
odd and  let  $V_{\nu,v}$ be a hermitian space with Hasse invariant $+1$. Then $V_{\nu,v}$ admits lattices of types $0$ and $\nu-1$, with  stabilizers denoted by $K_{\nu,v}^{[0]}$ and $K_{\nu,v}^{[\nu-1]}$. Let $\pi_{\nu,v}$ be a tempered unramified representation of $U(V_{\nu,v})$ whose base-change is $\Pi_{\nu, v}$.  Then the intertwining Hecke operator ${\rm I}_{\nu}^\circ$ in \cite[Proposition B.4.3~(1)]{LTXZZ} acting on the one-dimensional space $\pi_{\nu,v}^{K_{\nu,v}^{[0]}}$ is a non-zero scalar. Since the operator ${\rm I}_{\nu}^\circ$ factors through  a map $ \pi_{\nu,v}^{K_{\nu,v}^{[0]}}\to \pi_{\nu,v}^{K_{\nu,v}^{[\nu-1]}}$, the target space must be non-zero, hence $\pi_{\nu,v}$ is also almost unramified. 
\end{proof}

 \begin{lemma}\lb{lem: L-packet}
 Assume that the  residue characteristic of $v$ is odd. Let $\Pi_{\nu,v}$ be a tempered hermitian representation of $\GL_{\nu}(F_v)$. Assume that $\Pi_{\nu,v}$ is almost unramified and the coefficient field $L$ is algebraically closed.
\begin{itemize}
\item When $\nu$ is even, if $\e(V_{\nu,v})=-1$ there is a unique almost unramified representation $\pi_{\nu,v}$ of $U(V_{\nu,v})$ with base change $\Pi_{\nu,v}$ (whereas there are no almost unramified representations of the quasi-split unitary group due to the non-existence of vertex lattices of type $1$, $\nu-1$ in a hermitian space with Hasse invariant $1$).
\item When $\nu$ is odd, for every hermitian  space $V_{\nu,v}$ of dimension $\nu$, there is a unique almost unramified representation $\pi_{\nu,v}$ of $U(V_{\nu,v})$ with base change $\Pi_{\nu,v}$. Moreover, the representation  $\pi_{\nu,v}$ is generic. 
    \end{itemize}
\end{lemma}

\begin{proof}
{We may assume that $L=\C$ after fixing an arbitrary embedding.}
The even $\nu$ case follows from \cite[Theorem 1.1]{Dang} (the assumption in \emph{loc. cit.} that $\vert \sO_{F_{0},v}/(\varpi_v)\vert \geq \nu$  is not needed in the special case of conductor $1$). The odd $\nu$ case follows from \cite[Theorem 1.1]{AOY}. Note that in the odd $\nu$ case the two non-isomorphic classes of hermitian spaces define isomorphic unitary groups.
\end{proof}

\subsubsection{Distinction and non-vanishing of local relative characters}

{For the rest of this section we assume that $v$ is an inert place of $F_{0}$ with odd residue characteristic.}

Recall from \S~\ref{sec incoh} the fixed element $\nrm\in F_0^\ts$ based on which we have the set  $\sV_v$ of pairs of hermitian spaces (up to isomorphism). Recall also that an element $V_v\in \sV_v$ is uniquely determined by its Hasse invariant $\e(V_v)\in\{\pm 1\}$ defined in \eqref{eq: HW}.
We will say that a representation $\pi_v=\pi_{n,v}\boxtimes\pi_{n+1,v}$ of $G_v^{V_v}$ is almost unramified if both  $\pi_{n,v}, \pi_{n+1,v}$ are almost unramified.
\begin{lemma}\lb{lem: GGP au}
Let $\Pi_v=\Pi_{n,v}\boxtimes\Pi_{n+1,v}$ be a tempered hermitian irreducible representation of $G'_v$ where both $\Pi_{n,v}$ and $\Pi_{n+1,v} $ have conductors at most $1$ and trivial central characters. 
{Let $(V_v, \pi_v)$ be the unique distinguished element in the Vogan $L$-packet of $\Pi_v$ (\S~\ref{sss: dist}).

\begin{itemize}
\item If $\Pi_{{\rm even},v}$ is unramified, then $V_{v}$ has Hasse invariant $\e(V_v)=1$ and $\pi_v$  is the unique element in the Vogan $L$-packet such that $\pi_{{\rm odd},v}$ is almost unramified.
\item If $\Pi_{{\rm even},v}$ has conductor $1$, then $V_{v}$ has Hasse invariant $\e(V_v)=-1$, and $\pi_v$  is the unique almost unramified element in the Vogan $L$-packet. 
\end{itemize}}
\end{lemma}
 {Here, recall that  the Vogan $L$-packet is a Galois orbit of sets of pairs $(V_v, \pi_v)$ where $\pi_v$ is an isomorphism class of irreducible representations of $G_v^{V_v}$ over an algebraic closure of $L$, and the Galois action is by automorphisms on the coefficients of $\pi_v$. The statement still makes sense since the notion of almost unramified is Galois-invariant (and in fact, we obtain that the  Galois orbit associated to the distinguished element is a singleton).
}

\begin{proof}
{We may and do extend scalars to $\C$.}
The proof is a continuation of that of Lemma \ref{lem: loc ep}, and we will use the same notation.  By that lemma and Proposition \ref{local ggp}, we already know the Hasse invariant of $V_v$.  We apply  the explicit form of the local GGP conjecture  proved by Beuzart-Plessis to determine which representation in the Vogan $L$-packet is distinguished. Recall from \cite[\S17]{GGP} that the unique distinguished representation corresponds to the character $\chi: A_M\times A_N\to\{\pm1\}$, $\chi(a,b)=\chi_N(a)\chi_M(b)$ for $\chi_M,\chi_N$ defined in \cite[\S6]{GGP}.

First we consider the case where $\Pi_{{\rm even},v}$  is unramified.  The $L$-parameter of $\Pi_{{{\rm even}},v}$ (resp. $\Pi_{{\rm odd},v}$) is of the form $M=P+\,^{\rm c}P^\vee$ (resp. $N={\bf 1}+(Q+\,^{\rm c} Q^\vee)$ or $N={\rm St}_2\otimes\mu+{\bf 1}+(Q+\,^{\rm c}Q^\vee)$), and its component group $A_M$ is trivial (resp. $A_N\simeq\Z/2\Z$ or  $A_N\simeq \Z/2\Z\times\Z/2\Z$). It follows that the character of the component group satisfies
$$\chi(1,a)=\ep(M\otimes N^a,\psi)=1
$$ for all $a$ in the component group $A_N$ of $N$ (here and below, `$1$' in $\chi(1,a)$ denotes the identity element in the component group). Hence, the distinguished representation $\pi_v=\pi_{n,v}\boxtimes\pi_{n+1,v}$ is the unique generic representation in the Vogan $L$-packet, where $\pi_{n,v}$ is unramified and $\pi_{n+1,v}$ is the almost unramified one by Lemma \ref{lem: L-packet}.

We now consider the case where $\Pi_{{\rm even},v}$  has conductor $1$. Then the $L$-parameter of $\Pi_{{\rm even},v}$ (resp. $\Pi_{{\rm odd},v}$) is of the form $M={\rm St}_2+(P+\,^cP^\vee)$ (resp.,  $N={\bf 1}+(Q+\,^{\rm c}Q^\vee)$ or $N={\rm St}_2\otimes\mu+{\bf 1}+(Q+\,^{\rm c}Q^\vee)$), with its component group $A_M\simeq\Z/2\Z$ (resp.,  $A_N\simeq\Z/2\Z$ or $A_N\simeq \Z/2\Z\times\Z/2\Z$). The case $N={\bf 1}+(Q+\,^{\rm c}Q^\vee)$ is similar to the first case and we omit the detail. When $N={\rm St}_2\otimes\mu+{\bf 1}+(Q+\,^{\rm c}Q^\vee)$,  we have 
  \begin{align*}
  \chi(1,a_1)&=\ep(M\otimes {\rm St}_2\otimes\mu,\psi)=1, \\
  \chi(1,a_2)&=\ep(M\otimes {\bf 1},\psi)=-1,
  \end{align*} for $a_1\in A_N$ ($a_2\in A_N$) corresponding to ${\rm St}_2\otimes\mu$ (resp. $\bf 1$). The above character of $A_N$ is $\chi(1,a)=(-1)^{\dim N^a}$, hence it corresponds to the generic representation, or equivalently, the almost unramified representation by Lemma \ref{lem: L-packet}.  
\end{proof}

{
\begin{remark}\lb{rmk: GGP au}
From Lemma \ref{lem: GGP au}, we deduce the following  statements (in the same context, but not covering all cases)  that will be directly used in the proof of Theorem \ref{th BBK}.
\begin{itemize}
\item If $\Pi_{n,v}$  is unramified and $\Pi_{n+1,v}$ has conductor $1$,   then $V_{v}$ has Hasse invariant $\e(V_v)=(-1)^n$ and $\pi_v$ is the unique representation in the Vogan $L$-packet such that $\pi_{n+1,v}$ is almost unramified and not unramified.  
\item If $\Pi_{n,v}$ has conductor $1$, and $n$ is even or $\Pi_{n+1,v}$ has conductor $1$, then $V_{v}$ has Hasse invariant $\e(V_v)=-1$ and $\pi_v$  is the unique almost unramified representation in the Vogan $L$-packet. 
\end{itemize}
\end{remark}
}

It follows from Lemma \ref{lem: GGP au} that the local relative character\footnote{{See the comment after the next lemma for the exact meaning of this over a general  coefficient field.}} of the indicated $\pi_v$ (as a distribution) is nonzero. We will need  the following pair of  \emph{explicit} non-vanishing statements. Let $\e_v\in\{0, 1\}$ have the same parity as $v(\nrm)$, and recall from \S\ref{sss:vert} the notion of vertex parahoric subgroups of type $(t,t+\e_v)$.  Let $V_v\in \sV_v$ and let $G_v=G_v^{V_v}$.

\begin{lemma}[Dang, D. Liu -- L. Zhang] \lb{lem dang}
  Let $\Pi_{v}=\Pi_{n,v}\boxtimes \Pi_{n+1,v}$ be a tempered hermitian irreducible representation of $G'_v$ such that $\Pi_{n,v}$, $\Pi_{n+1,v}$ have trivial central character and conductors $c_{n}$, $c_{n+1}$ as in one of the entries of the following table.
  
  \bigskip
  
  \begin{center}
 \begin{tabular}{| c | c | c  |} 
 \hline
Conductors $ (c_{n}, c_{n+1})$ & $\e_{v}$ & Type of $K_{v}$   \\ [0.5ex] 
 \hline
$ (0, 0)$ or $(0, 1)$ & $0$ & $(n,n)$  \\ 
 \hline
 $ (0, 0)$ or $(0, 1)$ & $1$ & $(0,1) $ \\ 
 \hline
 $ (1, 1)$ & $0$ & $(1,1)$  \\ 
 \hline
 $ (1, 1)$  & $1$& $(n-1,n)$  \\  
 [0,5ex] 
 \hline
\end{tabular}
\end{center}

  \bigskip

 Let $(V_v, \pi_v)$ be the distinguished element in the Vogan $L$-packet of $\Pi_v$ (cf. Lemma \ref{lem: GGP au}).
 Let $f_{v}=e_{K_{v}}$ where $K_{v}\subset G_{v}$ is a vertex parahoric subgroup of type given by the above table (depending on the conductors and on $\e_{v}$).
 % \WZ{e.g.,} type $(n,n)$ if conductors $ (0, 0)$ or $(0, 1)$ and $\e_v=0$, type  $(0, 1)$ if conductors $ (0, 0)$ or $(0, 1)$ and $\e_v=1$.
 Let $f_v'=e_{K_v'}$, where $K_v'\subset G_v'$ is related to $K_v$ as in Proposition \ref{zhiyu tr}.
Then 
$$
I_{\Pi_v}(f'_v)=J_{\pi_{v}}(f_{v})\neq 0.
$$
\end{lemma}

{(The above, as well as the conclusion of the next conjecture, may be interpreted as holding after taking any embedding of the coefficient field into $\C$. For the general definitions of relative characters over arbitrary coefficient fields, see Proposition \ref{rat rtf} \eqref{rat sph ch} and
Remark \ref{rat matching} below.)}

\begin{proof} The equality follows from Proposition \ref{zhiyu tr}.

In the cases with $c_{n}=0$, the non-vanishing in the case where $\e_{v}=1$   is a special case of \cite[Theorem 1.2]{Dang}; the  case with $\e_{v} =0$ is reduced to the case of $\e_{v}=1$ by simultaneously rescaling the hermitian forms by $\varpi_v$.

The case $(c_{n}, c_{n+1})=(1, 1)$ is \cite[Theorem 5.1]{LZ}.
\end{proof}

\begin{remark}  \lb{rmk dongwen}
The case of Lemma \ref{lem dang} concerning conductors $(1,1)$ was conjectured in a preprint version of this paper.\footnote{See Conjecture 3.6.10 in \url{https://arxiv.org/pdf/2410.08401v2}.} There, inspired by the second case of Remark \ref{rmk: GGP au},  we had also  suggested that the same conclusion would hold for conductors $(1,0)$ when $n$ is even; however, Dongwen Liu has recently communicated to us a possible counterexample in the latter case. 

Thus the {\em explicit} non-vanishing of relative characters in the case  where $\Pi_{n,v}$ has conductor~$1$ and $\Pi_{n+1,v}$  is unramified remains an interesting topic for further research.
 That situation is closely related to the context  of  \cite[\S~13.1]{LRZ}.
  However, neither an {\em explicit} transfer nor an arithmetic transfer statement  is known in that case.
 \end{remark}

\section{Rationality}\lb{sec rat}
This section is dedicated to  establishing the rationality of our $L$-values, namely Theorem \ref{thm A} from the introduction (Theorem \ref{rat L} below), and a rational relative-trace formula (Proposition \ref{rat rtf}).

From now on, we assume that $F_{0}$ is totally real and $F$ is CM.

 In \S\ref{ss: arch} we deal with the archimedean computations using  Gaussian test measures.  In  \S\ref{sec: rat st} we state the rationality theorem and the rational relative-trace formula, and prove some easy parts.  In \S\ref{sec: test} we study the existence of suitable Hecke measures: the non-archimedean components rely on later results from \S\S~\ref{sec: loc rtf}-\ref{app orb p}; the (mixed) archimedean components are given in \S\ref{sec: prove Gauss} using an argument provided by Yifeng Liu, refining the technique of isolating cuspidal representations in \cite{isolation}. In \S\ref{last 422 pf}, we finish the proof of  Proposition \ref{rat rtf} and Theorem \ref{rat L}. In \S\ref{sec: II} we recall an outline of the proof of (a special case of) the  Ichino--Ikeda--Harris conjecture. Logically this is not needed for this paper, but it will make the proof of our main Theorem \ref{main thm}  in \S~\ref{sec:11}  easier to understand.

\subsubsection*{Notation} For a locally compact and totally disconnected group $G$ with a fixed Haar measure $dg$ and a ring $R$, from now on we denote by  $\sH(G, R)$   the 
module of smooth compactly supported $R$-multiples of $dg$.
(Thus the object denoted by $\sH(G)$ in the previous section will henceforth be denoted by $\sH(G,\C)$.)

\subsection{Archimedean theory}
\lb{ss: arch}
We define some rational variant of the archimedean distributions of the previous section. 
Denote $G_{\infty}^{\circ}\coloneqq G_{\infty}^{V_{\infty}^{\circ}}$, $H_{\infty}^{\circ}\coloneqq H_{\infty}^{V_{\infty}^{\circ}}$, $B_{\infty}^{\circ}\coloneqq  B_{\infty, V^{\circ}_{\infty}} = B_{\rs, \infty, V^{\circ}}$, $B_{\rs, \infty}'^{\circ}\coloneqq B'_{\rs, \infty, V^{\circ}_{\infty}}$.

\subsubsection{A product of transfer factors} 
Let
\beqq \kappa(\one_{\infty})\coloneqq \prod_{v\vert \infty}\kappa(\one_{v})\eeqq
 be the product of \eqref{tr f} for  the trivial representation of the positive-definite group $G_{\infty}^{\circ}$. 
 
\begin{lemma} \lb{rat kappa}
For each $\gamma'\in \RG'_{\reg^+}(F_{0, \infty})$, we have
${ \kappa_{\infty}(\gamma') \kappa(\one_{\infty})}\in 
\{\pm 1\}$. 
\end{lemma}
\begin{proof} By Lemma \ref{kappa i}, the first factor is a square root of $(-1)^{-{n+1\choose 2} [F_{0}:\Q]}$; so are $\eta'_{\infty}(\tau)^{n+1\choose 2}$ and, hence, the second factor.
\end{proof}

\subsubsection{Distributions on $\sH(G'_{\infty}, \C)$}
For any tempered representation $\Pi_{\infty}$ of $G_{\infty}'$, any $\gamma\in B'_{ \infty}$, and any $f'\in  \sH(G'_{\infty}, \C)$, we define 
\beq\lb{def I circ}
I_{\Pi_{\infty}}^{\circ} (f'_{\infty}, \chi_{\infty })&\coloneqq  
{1 \over \kappa(\one_{\infty})\Delta_{\RH}^{}} \sL(1/2, {\Pi_{\infty},\chi_{\infty}}) \, {I}^{}_{\pi_{\infty}}(f_{\infty},\chi_{\infty}),
\\
\hI_{\gamma}^{\circ}( \ff'_{\infty}, \chi_{\infty})&\coloneqq {\Delta_{\G}^{}\over \Delta_{\RH}^{2}} L_{\gamma} (\chi_{\infty}) \hI_{\gamma}^{}( \ff'_{\infty}, \chi_{\infty}),\eeq
where the last two terms  were defined in \eqref{eq orb on B}.
Then the factorizations of Proposition \ref{factor I} and of \eqref{prod orb int}   are equivalent to
\beq\lb{factor I rat}
\kappa(\one_{\infty})^{-1} \hI_{\Pi}(\ff', \chi)
&= {1\over 4}{ \sL^{\infty}(1/2, \Pi,\chi) 
\over  \ep(\textstyle{1\over 2} , \chi^{2})^{{n+1\choose 2}} }
\,    \prod_{v\nmid \infty} \hI_{\Pi_{v}}(\ff'_{v}, \chi_{v}) \cdot I_{\Pi_{\infty}}^{\circ} ( f'_{\infty}, \chi_{\infty}),\\
\hI_{\gamma}(\ff', \chi) 
&= L^{\infty}_{\gamma}(\chi) \prod_{v\nmid\infty }\hI_{\gamma}(\ff'_{v}, \chi_{v})\cdot  \hI_{\gamma}^{\circ}( \ff'_{\infty}, \chi_{\infty}).
\eeq

\subsubsection{Distributions and special elements in $\sH(G^{\circ}_{\infty} , \C)$}
For any $V\in \sV_{\infty}=\prod_{v\vert \infty}\sV_{v}$, every representation $\pi^{V}_{\infty}$ of $G_{\infty}^{V}$, and every {$\delta\in B_{\infty, V}$ (note that $G_{\infty}^{V}$ is compact and hence every orbit is semisimple), we define  variants of $J_{\pi^{V}_{\infty}}$ and $J_{\delta, \infty}$ by

\beq \lb{def J circ}
{J}^{\circ}_{\pi_{\infty}}(f_{\infty})  &\coloneqq  
  \int_{H^{V}_{\infty}}  \Tr_{\pi_{\infty}}( \pi_{\infty}(h) \pi_{v}(f_{\infty})) \, {d^{} h} 
  = {1 \over 
\Delta_{\RH}^{}} \sL(1/2, {\rm BC}(\pi_{\infty})) \cdot {J}^{}_{\pi_{\infty}}( f_{\infty}), \\
J^{\circ}_{\delta, \infty}(f_{\infty}) &\coloneqq    \int_{H^{V}_{\infty}} \int_{H^{V}_{\infty}} f_{v}(x^{-1} \delta' y)\, {d^{}x d^{}y \over d^{}g} =
{\Delta_{\G} \over \Delta_{\RH}^{2}} \, J^{}_{\delta,\infty}(f_{\infty}).\eeq

Then the  matching relations \eqref{mat sp} and, respectively, \eqref{mat ge} for $S=\{v\vert \infty\}$ are equivalent to  
\beq\lb{match circ}
\hI^{\circ}_{{\rm BC}(\pi^{V_{}}_{\infty})}(\ff'_{\infty} )&= {\kappa(\pi^{V_{}}_{\infty}) \over \kappa(\one_{\infty})}\hJ^{\circ}_{\pi^{V_{}}_{\infty}}(\ff_{\infty}) \\
\hI^{\circ}_{\gamma, \infty}(\ff_{\infty}', \one_{\infty})&=\hJ^{\circ}_{\delta, \infty}(\ff^{V_{}}_{\infty}).
\eeq

\begin{lemma}\lb{J infty compute} 
Let
\beq \lb{f inf}
f_{\infty}^{\circ}\coloneqq  \vol(G_{\infty}^{\circ}, dg)^{-1} dg \quad \in \sH(G_{\infty}^{\circ}, \Q).
\eeq 
Then:
\begin{enumerate}
\item 
for  all $\pi_{\infty}\in {\rm Temp}(G_{\infty}^{\circ})$, we have
 \beq\lb{J circ}
J^{\circ}_{\pi_{\infty}}(f_{\infty}^{\circ})=\begin{dcases} \vol(H_{\infty}^{\circ})\coloneqq \vol(H_{\infty}^{\circ}, dh_{\infty}) & \text{if $\pi_{\infty}\cong \one$}\\ 0 &\text{otherwise};\end{dcases}\eeq
\item \lb{J infty compute b}
for all  $\delta \in G^\circ_{\infty}$, we have
$$J^{\circ}_{\delta} (f_{\infty}^{\circ}) = \vol(B_{\infty}^{\circ})^{-1}
: = {\vol(H_{\infty}^{\circ}, dh_{\infty})^{2}\over \vol(G_{\infty}^{\circ}, dg_{\infty})}.$$
\end{enumerate}
Moreover, both of the above values are rational.
\end{lemma}
\begin{proof} The calculation is immediate. The rationality follows from  Lemma \ref{gross fact}.
\end{proof}

\subsubsection{Gaussians} 
Let $f_{\infty}^{\circ}= \eqref{f inf}$ be the standard Hecke measure on $G_{\infty}^{\circ}$. For each characteristic-zero field $L$, we put $\sH(G^{\circ}_{\infty}, L)^{\circ}\coloneqq Lf_{\infty}^{\circ}$.

For $L$ a subfield of $\C$, we denote by
$$\sH(G'_{\infty}, L)^{\bullet}\subset \sH(G_{\infty}', \C)$$ the preimage
of $L f_{\infty}\subset \sH(G^{\circ}_{\infty}, L)^{\circ}$ under pure matching;  if $S$ is a finite set of non-archimedean places of $F_{0}$, we also put 
$ \sH(\G'(\A^{S}), L)^{\bullet}\coloneqq  \sH(\G'(\A^{S\infty}), L)\ot_{L}\sH(G'_{\infty}, L)^{\bullet}$.  By Proposition \ref{gauss E} below, the pure matching map
$${\rm tr} \colon \sH(G'_{\infty}, L)^{\bullet}\to  \sH(G^{\circ}_{\infty}, L)^{\circ}$$
is surjective (here $\rm tr$ stands for ``(smooth) transfer"). We put 
$$ \sH(G'_{\infty}, L)^{\circ}\coloneqq  \sH(G'_{\infty}, L)^{\bullet}/\Ker({\rm tr}),$$
we extend the definition to any characteristic-zero field $L$ by $\sH(G'_{\infty}, L)^{\circ}\coloneqq  \sH(G'_{\infty}, \Q)^{\circ}\ot_{\Q}L$, and we extend the notion of matching by linearity. Elements of  $\sH(G'_{\infty},L)^{\circ}$ are called $L$-rational \emph{Gaussians}. If $L$ is a subfield of $\C$,
we also refer to an element $f_{\infty}'\in\sH(G'_{\infty},L)^{\bullet}$ as a Gaussian; we say that $f_{\infty}'$ is nontrivial if its image in $\sH(G'_{\infty},L)^{\circ}$ is nonzero.

  If $S$ is a finite set of non-archimedean places of $F_{0}$, we put 
\beqq
\sH(G'_{S\infty}, L)^{\circ}&\coloneqq  \sH(G'_{S}, L)\ot_{L}\sH(G'_{\infty}, L)^{\circ}, \\
\sH(\G'(\A^{S}), L)^{\circ} &\coloneqq  \sH(\G'(\A^{S\infty}), L)\ot_{L}\sH(G'_{\infty}, L)^{\circ},
\eeqq
and refer to the elements of those spaces as Gaussians too.

\begin{proposition}[Existence of  Gaussians]
\lb{gauss E}  The space $\sH(G'_{\infty}, \Q)^{\circ}$ is nonzero.
\end{proposition}
\begin{proof} This follows from \cite[Proposition 4.11]{isolation}, or from the recent explicit construction by Mihatsch--Sankaran--Yang \cite{MSY}.
\end{proof}

\begin{lemma}\lb{I infty reg} 
Let $\ff'$  be a Gaussian matching  $f_{\infty}^{\circ}= \eqref{f inf}$. Then for any $\gamma\in B'_\infty$, we have $I_{\gamma}^{\circ}(\ff',\one)\in\Q$.
\end{lemma}

\begin{proof}
We recall from
\eqref{def I circ}
$$
\hI_{\gamma}^{\circ}( \ff', \one)\coloneqq {\Delta_{\G}^{}\over \Delta_{\RH}^{2}} L_{\gamma} (\one) \hI_{\gamma}^{}( \ff', \one), 
$$
where from \eqref{eq orb on B},
\beq\notag
\hI^{}_{\gamma}(\ff', \one) &=  \kappa(\gamma') \,\hI^{}_{\gamma'}(\ff', \one)
  = \kappa(\gamma') \frac{\hI_{\gamma'}^{\sharp}( \ff', \chi)}{L_{\gamma'}(\chi)}.
\eeq
Here $\gamma'$ is any plus-regular element above  $\gamma\in B'_\infty$, and $\kappa(\gamma')$ is the local transfer factor.

Recall also 
 the orbital integral  $J^{\circ}_{\delta}$ in the unitary side \eqref{def J circ}.  
If $\gamma$ is regular semisimple, the lemma follows immediately from   the rationality of $J^{\circ}_{\delta}(f_{\infty}^{\circ})$ (Lemma \ref{J infty compute}~\eqref{J infty compute b}) and the matching relation \eqref{match circ}.
Though the matching relation is defined only using regular semisimple orbits, the definition implies non-trivial relations for non-regular-semisimple orbital integrals. We  record the result of Lu \cite[Thm. 7.9, Remark 7.10]{Lu} comparing the local orbital integrals.  Let $f'\in \sH(G'_{\infty}, \C)$  purely match an element $f\in \sH(G_{\infty}^{\circ},\C)$.
 If $\gamma\in B'$, then  
\beq\lb{orb reg=ss}
 L_{\gamma} (\one) ^{-1}\hI_{\gamma}^{\circ}( \ff', \one)=\sum_{\delta} c_{\delta} J^{\circ}_{\delta}(\ff),
\eeq
where the sum runs over all {\em semisimple} orbits in the {compact group} $G_{\infty}^{\circ}$ with image $\gamma\in B'_\infty$, and  
$$
c_{\delta} = \prod_{W\in \sW(\gamma)} c_{W},
$$
 where the set    $\sW(\gamma)$  
 and  the constants $ c_{W}$ will be recalled next.  The set  $\sW(\gamma)$  
is a finite set  of positive definite $\C/\R$-hermitian spaces $W$ defined in {\em loc. cit.}, and it
 can be described as follows: the stabilizer of any semisimple $\delta$ matching $\gamma$ is isomorphic to the product of the compact unitary groups $U(W)$ for $W\in \sW(\gamma)$. For  $W$ of dimension $n'$, by
  \cite[\S7.4 on the Lie algebra, and (7.15) and Remark 7.10 on the group]{Lu} we have
\beq\lb{eq Lu factor}
c_{W} = &\vol^{\nat}(U(n',\R))^{-1} \cdot \eta_{\C/\R}(\disc(V'))^{n'+1} \\
&\cdot \ep(1/2,\eta_{\C/\R}, \psi)^{n'(n'+1)/2} 
\cdot\prod_{j=1}^{n'}\ep(1-j, \eta_{\C/\R}^j , \psi)^{-1} .
\eeq
Here
 $\vol^{\nat}(U(n',\R))$ is the volume  of the compact unitary group $U(n',\R)$ for the normalized measure $d^{\nat}g$ of \S~\ref{sec: tam}, which is the measure in \cite[\S~7.0.1]{Lu}  (for a suitable differential $\omega$). The formula for the constant $ c_{W}$ differs slightly from \cite{Lu} due to a few different conventions between ours and those in {\it loc. cit.}:
\begin{itemize}
 \item the factor ${\Delta_{\G}^{}\over \Delta_{\RH}^{2}}$ appears on both  the $\GL$ and the unitary side, and hence our notion of matching is equivalent to that of \cite{Lu};
 \item Theorem 7.15 in  \cite{Lu} is expressed in terms of the normalized orbital integral, and this results in the factor $L_{\gamma} (\one) ^{-1}$ on the left hand side of \eqref{orb reg=ss};
\item when defining $\hI_{\gamma}(\ff', \one)$, we only consider the plus-regular element  $\gamma'$ with image $\gamma\in B'$ and our notation has already included the transfer factor (our transfer factor is the plus-transfer factor in \cite{Lu});
 \item our orbital integral in the unitary side is taken over the full group $H^{\circ}_{\infty}\times H^{\circ}_{\infty}$, whereas in  {\it loc. cit.} the integral is taken over the quotient of the full group by the stabilizer: this results into the volume factor in \eqref{eq Lu factor};
 \item In \cite[(7.15)]{Lu} the additional factor disappears because our choice of $\gamma'$ above $\gamma$ is plus-regular {and the formula in \cite[Thm. 7.9]{Lu} simplifies to \cite[Remark 7.10, (7.16)]{Lu}.} 
\end{itemize}

The factors in the second line in \eqref{eq Lu factor} are signs,  hence lie in $\Q^\times$. { 
By definition, the $L$-factor $ L_{\gamma} (\one) $ is the product}
$$
 L_{\gamma} (\one)  = \prod_{W\in \sW(\gamma)}\prod_{j=1}^{\dim W}L(1-j, \eta_{\C/\R}^j).
$$
Therefore, to show that $\hI_{\gamma}^{\circ}( \ff', \one)\in\Q$, from \eqref{orb reg=ss} and \eqref{eq Lu factor}  it suffices to show that 
for all $n'\geq 1$, the product 
\begin{align}\lb{eq rational}
 \vol^{\nat}(U(n',\R))^{-1} \cdot   \ep(1/2,\eta_{\C/\R}, \psi)^{n'(n'+1)/2}
 \cdot \prod_{j=1}^{n'}L(1-j,  \eta_{\C/\R}^j) \ep(1-j, \eta_{\C/\R}^j , \psi)^{-1}
\end{align}
lies in $\Q^\times$.

By Tate's thesis (e.g. \cite[\S3.2]{Ta}), the standard choice of $\psi(x)=e^{2\pi i x}$  gives 
$$
\ep(s,\eta^a_{\C/\R}, \psi)=i^a\in \BC
$$ 
for all $s\in\bC$ and $a\in\{0,1\}$.  In particular, 
we have
$$
 \ep(1/2,\eta_{\C/\R}, \psi)^{n'(n'+1)/2}=i^{n'(n'+1)/2}.
$$
and
$$
\prod_{j=1}^{n'}\ep(1-j, \eta_{\C/\R}^j , \psi)^{-1}=i^{-\lfloor \frac{n'+1}{2}\rfloor}.
$$
We note that $n'(n'+1)/2\equiv \lfloor (n'+1)/2\rfloor\mod 2$, and hence 
\beq\label{ep prod n'}
\prod_{j=1}^{n'}\ep(1-j, \eta_{\C/\R}^j , \psi)^{-1}\cdot \ep(1/2,\eta_{\C/\R}, \psi)^{n'(n'+1)/2}=\pm1.
\eeq

Next we note for $a\in\{0,1\}$,
 $$
 L(s,\eta^a_{\C/\R})=L(s+a,\one)=\pi^{-(s+a)/2}\Gamma((s+a)/2)
 $$ and we have its special values 
 $$
L(1-j,\eta^j_{\C/\R})=\begin{cases}\pi^{-(1-i)/2}\Gamma((1-i)/2), & j\geq 0 \text{ even},
\\
\pi^{-(1-j+1)/2}\Gamma((1-j+1)/2), & j\geq 0 \text{ odd}.
\end{cases}
$$
 In both cases we have
\beq\lb{eq L 1-i}
L(1-j,\eta^j_{\C/\R})\in\pi^{\lfloor j/2\rfloor}\cdot \Q^\times,\quad j\geq 1.
\eeq
Similarly,
\beq\lb{eq L i}
L(j,\eta^j_{\C/\R})\in\pi^{-\lfloor (j+1)/2\rfloor}\cdot \Q^\times,\quad j\geq 1.
\eeq

We compute the volume $\vol^{\nat}(U(n',\R))$ of the compact unitary group. Denote by $\vol(U(n',\R))$ the volume under the unnormalized measure $d_{\omega}g$ of \S~\ref{sec: tam}, then 
\begin{align*}
\vol^{\nat}(U(n',\R))&=\prod_{j=1}^{n'}L(j,\eta^j_{\C/\R})\cdot \vol(U(n',\R))\\
&=\prod_{j=1}^{n'}L(j,\eta^j_{\C/\R})\prod_{j=1}^{n'}\vol(S^{2j-1})\\
\text{(by \eqref{eq L i})}\quad &
\in \prod_{j=1}^{n'}\pi^{-\lfloor (j+1)/2\rfloor}   \pi^{j} \cdot\Q^\times
\\&= \prod_{j=1}^{n'}\pi^{\lfloor j/2\rfloor}  \cdot\Q^\times
\end{align*}
where $\vol(S^{2j-1}) \in \pi^{j}\Q^{\times}$ is the usual volume of the unit sphere of dimension~$2j-1$.
Combining this with \eqref{eq L 1-i}, we have
\beq\lb{eq vol times L}
 \vol^{\nat}(U(n',\R))^{-1} \prod_{j=1}^{n'}L(1-j, \eta_{\C/\R}^j)\in\Q^\times.
 \eeq
Therefore the rationality of \eqref{eq rational} follows from \eqref{ep prod n'} and \eqref{eq vol times L}, and
the lemma follows from this and the rationality of $J^{\circ}_{\delta}(f_{\infty}^{\circ})$ (Lemma \ref{J infty compute}~\eqref{J infty compute b}).
\end{proof}

\subsection{Rationality statements} \lb{sec: rat st}
We state the main results of this section, whose proofs will be completed in \S~\ref{4pf}.

\subsubsection{Rationality of $L$-values}
Recall from Proposition \ref{cor C her} the ind-finite scheme  $\sC(\G')^{{\rm her}}$ 
(over $\Q$) of trivial-weight  cuspidal hermitian automorphic representations of $\G'(\A)$.
The following is  Theorem \ref{thm A} from the introduction.

\begin{theorem} \lb{rat L} 
Let $\Pi \in \sC(\G')^{{\rm her}}$, 
 and let $L=\Q(\Pi)$ be its field of definition. There is a function $$\sL(\RM_{\Pi}, \cdot)\in \sO(Y_{L})$$ 
such that for each $\chi \in Y_{L}(\C)$ with underlying embedding $\iota\colon L\into \C$,
$$\sL(\RM_{\Pi}, \chi)={ \sL^{\infty}(1/2,{\Pi^{\iota}, \chi})\over \ep({1\over 2}, \chi^{2})^{{n+1\choose 2}} }.$$
\end{theorem}

\subsubsection{Special Hecke algebras} \lb{rha}
Let $L$ be a coefficient field (in the sense, defined in \S~\ref{coeff}, that $L$ admits embeddings into $\C$),
let $S$ be a finite set of finite places of $F_{0}$, and let $?\in \{\rs, \mathrm{reg}^{+}, \emptyset\}$.
 We denote by
 $$\sH(\G'(\A^{S}) , L)^{\circ}_{K_{S}, ?, \rm qc}\subset \sH(\G'(\A^{S}) , L)^{\circ}$$
the space of Gaussian measures $f'^{S}$ with globally $?$-support  (Definition \ref{def: weak support}; there is no condition if $?=\emptyset$) such that for every $\iota\colon L\into \C$, some preimage  $f'^{S, \iota}e_{K_{S}}\in \sH(\G'(\A) , \iota L)^{\bullet}$ of $\iota f'^{S}e_{K_{S}}$ is quasicuspidal.

If $\Pi\in \sC(\G')(L)$ and $\chi\in Y(L)$, we say that a Hecke measure 
 $f'^{S}\in \sH(\G'(\A^{S}) , L)^{\circ}$
is \emph{adapted to $(\Pi^{}, \chi^{}, K_{S})$} if  \begin{itemize}
\item
 $(\ot_{v\notin S}I_{\Pi_{v}})(f'^{S}, \chi^{S})\neq 0$, and \item  for every $\iota\colon L\into \C$, some preimage $f'^{S, \iota}e_{K_{S}} \in \sH(\G'(\A) , \iota L)^{\bullet} $ of $\iota f'^{S}e_{K_{S}}$ 
 sends $\sA(\G')$ into (the image in $\sA(\G')$ of) $\Pi$. 
 \end{itemize}
 We  denote by 
 $$
 \sH(\G'(\A^{S}) , L)^{\circ}_{K_{S}, ?, \Pi, \chi}
 $$ 
the space of  Gaussians   with globally   $?$-support
that  are adapted to $(\Pi, \chi, K_{S})$. When $\chi=\one$ we omit it from the notation.

\subsubsection{Rational relative-trace formula} \lb{sec 423} We introduce a variant of the distribution $I$ and its expansions. From now on,  we change the notation for the distributions $I_{?}$  of the previous section by  appending a superscript `$\C$', thus writing   $I_{?}^{\C}$ in place of $I_{?}$; we also write $L_{?}^{\C}$ for the abelian complex $L$-functions attached to orbits. 

We introduce some further notation. 
For a  finite place $v$ of $F_{0}$ and an ideal $m\subset \sO_{F_{0,v}}$, let $Y_{v}(m)=\Spec \Q[F_{0,v}^{\ts}/(\sO_{F_{0,v}}^{\ts}\cap 1+m\sO_{F_{0,v}})]$, viewed as the space of characters of the group within square brackets. Let $Y_{v}\coloneqq  \varinjlim_{r} Y_{v}(v^{r})$. For the sake of uniformity, we will denote $\sH(G_{v}', L)^{\circ}\coloneqq  \sH(G_{v}', L)$ if $v\nmid \infty$, and $Y_{\infty}\coloneqq \Spec\Q$. 

In the rest of the paper, unless otherwise noted all products `$\prod_{v}$' run over the union of the set of finite places $v$ of $F_{0}$  and   $\{v=\infty\}$.  If $\sH$ is a Hecke algebra over  a field $L$ and $Y$ is an ind-scheme over  $L$, an $L$-linear functional $D\colon \sH\to \sO(Y)$  will be  called a \emph{distribution}.

\begin{proposition}
\lb{rat rtf}
 Let $L$ be a coefficient field.
There exist: 
\begin{enumerate}
\item   \lb{rat sph ch} 
for each finite place $v$ of $F_{0}$ and for $v=\infty$, and for  each tempered irreducible admissible representation $\Pi_{v}$ of $G_{v}'$ over $L$, a  distribution 
$$I_{\Pi_{v}} \colon \sH(G_{v}', L)^{\circ}\to \sO( Y_{v, L})$$
uniquely  characterized by  
$$\iota I_{\Pi_{v}}(f_{v}', \chi_{v} )=
\begin{cases}
 I^{\C}_{\iota\Pi_{v}}(\iota f_{v}', \chi_{v}) & \text{if $v\nmid \infty$}\\
 I_{\Pi_{\infty}}^{\circ, \C}(\iota f_{\infty}', \chi_{v})& \text{if $v= \infty$}
    \end{cases}
  $$
for each $\chi_{v}\in Y_{v, L}(\C)$ with underlying embedding $\iota \colon  L \into \C$;
\item  \lb{rat sp exp} for each representation $\Pi\in \sC(\G')^{{\rm her}}(L)$, 
 a   distribution 
$$I_{\Pi}  \colon \sH(\G'(\A),L)^{\circ}\to \sO(Y_{L})$$
defined on factorizable elements $f'=\ot_{v\nmid \infty}f_{v}'\ot f_{\infty}'$ by  
\beq \lb{fact I Pi rat}
I_{\Pi}(f', \chi) &= {1\over 4} \sL(\RM_{\Pi}, \chi)\cdot  \prod_{v} I_{\Pi_{v}}(f_{v}', \chi_{v});
\eeq
\item \lb{rat orb int reg}
for each finite place $v$ of $F_{0}$ and for $v=\infty$,  and each  $\gamma\in B'_{v}$, a distribution $$\hI_{\gamma, v}\colon \sH(G'_{v}, L)^{\circ} \to \sO(Y_{v,L[\sqrt{-1}]}) $$
uniquely characterized by 
$$ \iota \hI^{}_{\gamma,v}(f'_{v}, \chi_{v})=\begin{cases}
  \hI_{\gamma,v}^{\C}(\iota f'_{v}, \chi_{v})& \text{if $v\nmid \infty$}\\
  \hI_{\gamma,v}^{\circ, \C}(\iota \iota f'_{v}, \chi_{v})& \text{if $v= \infty$}
    \end{cases}
  $$
for each
$\chi_{v}\in Y_{v,L[\sqrt{-1}]}(\C)$ with underlying embedding $\iota\colon L[\sqrt{-1}]\into \C$;
\item  \lb{rat orb int 2 reg} 
for each  $\gamma\in \RB'(F_0)$,
\begin{enumerate}
\item  \lb{rat orb int 2 reg a} 
an element  $L_{\gamma}\in \sO(Y)$, uniquely  characterized by $L_{\gamma}(\chi)=L_{\gamma}^{\infty, \C}(\chi)$ (the $L$-function without archimedean local $L$-factors) for every $\chi\in Y(\bC)$,
\item  \lb{rat orb int 2 reg b} a distribution 
 $$\hI_{\gamma}^{}= \kappa(\one_{\infty}^{})^{-1}\cdot  L_{\gamma}\cdot 
 \, \prod_{v} \hI_{\gamma,v}^{}\colon \sH(\G'(\A), L)^{\circ}\to\sO(Y_{L})$$
where the product is locally finite; 
\end{enumerate}

\item \lb{rat rtf part} a distribution 
$$\hI^{}  \colon  \sH(\G'(\A), L)^{\circ}_{\reg^+, \rm qc}\to \sO(Y_{L})$$
 admitting the spectral and geometric expansions
$$
 \sum_{\Pi\in \sC(\G')^{\rm her}} \hI_{\Pi}=
 \hI^{} = \sum_{\gamma \in \RB'(F_{0})}
\hI^{}_{\gamma} $$
where both sums are locally finite (that is, for every input $f'$, both expansions of $I(f')$ are finite sums). 
\end{enumerate}

\end{proposition}

\begin{remark} It should be possible to interpret the rational distribution $\hI^{}$ as the inner product of  analogues of $P_{1, \chi}$, $ P_{2}$ in the rational Betti homology (in complementary degrees) of the symmetric space for $\G'$. 
\end{remark}

We prove   Proposition \ref{rat rtf}   \eqref{rat sph ch}-\eqref{rat orb int 2 reg}; the proof of part   \eqref{rat rtf part} is deferred to \S~\ref{last 422 pf}.

\begin{proof}[Proof of Proposition \ref{rat rtf}  \eqref{rat sph ch}-\eqref{rat orb int 2 reg}]
We need to show the existence of various distributions.

\smallskip

\paragraph{\em Archimedean distributions}
Suppose that $f_{\infty}'$ is an $L$-rational Gaussian  matching  $f_{\infty}=c f_{\infty}^{\circ}\in \sH(G^{\circ}_{\infty}, L)$. Then by Lemma \ref{J infty compute}, Lemma \ref{I infty reg}, and \eqref{match circ}, 
we may define
\beqq
I_{\Pi_{\infty}}(f'_{\infty}, \one)&\coloneqq  \begin{cases} c \vol(H_{\infty}^{\circ}) &\text{if $\Pi_{\infty}\cong \Pi_{\infty}^{\circ}$}\\ 0 &\text{otherwise},
\end{cases} \\
I_{\gamma}(f'_{\infty})&\coloneqq  \begin{cases}
c I_{\gamma}^{{\circ}}(\ff'^{\circ},\one) 
&\text{if $\gamma\in  B_{ \infty}'^{\circ}$}, \\ 0 &\text{otherwise},
\end{cases} 
\eeqq
for any Gaussian $\ff'^{\circ} $ matching $\ff^{\circ}$, where the orbital integral is rational by Lemma \ref{I infty reg}.

\smallskip

\paragraph{\em Orbital integrals}
Suppose $v$ is non-archimedean. Then part  
 \eqref{rat orb int reg}
 of Proposition \ref{rat rtf}  follows from Proposition \ref{lu reg}  \eqref{lu reg 3b} together with Lemma \ref{kappa i}.

 Part \eqref{rat orb int 2 reg a}  is a well-known rationality theorem of Klingen and Siegel \cite{siegel}. Part  \eqref{rat orb int 2 reg b} then follows  from part \eqref{rat orb int reg} and Proposition \ref{lu reg} \eqref{lu aa1}, together with  \eqref{prod kappa} and Lemma \ref{rat kappa} for the elimination of $\sqrt{-1}$ from the field of rationality.

 \smallskip

\paragraph{\em Local relative character} It suffices to show that $P_{1, \Pi_{v}, \chi_{v}}(W_{v})$, $P_{2,\Pi_{v}}(W_{v}')$ and $\vartheta_{\Pi_{v}}(W_{v}, W_{v}')$ are polynomial functions in the values of $\chi_{v}$, with coefficients in the field $L[W_{v}, W_{v}']$ generated over $L$ by the values of  $W_{v}$, $W_{v}'$.
By Lemma \ref{aa1 L}, for each fixed $\chi_{v}$ those expressions belong to the field generated over $L$ by the values of $W_{v}$, $W_{v}'$, $\chi_{v}$. Then we only need to prove that the function $\chi_{v}\mapsto P_{1, \Pi_{v}, \chi_{v}}(W_{v})$ belongs to $\sO(Y_{v, L[W_{v}]})$.  Let $Y_{v}'$ be the ind-finite scheme over $L$ of  smooth characters of $\sO_{F_{0,v}}^{\ts}$; then $\chi\mapsto \chi_{v|\sO_{F,0,v}^{\ts}}$ gives  an exact sequence of ind-group-schemes $1\to Y_{v}^{\circ}\to Y_{v}\to Y_{v}'\to 1$ where $Y_{v}^{\circ}\cong {\bf G}_{m, L}$ parametrizes unramified characters of $F_{0,v}^{\ts}$. Thus locally we may reduce to proving the desired result when $\chi_{v}$ is restricted to $Y_{v}^{\circ}$ at the cost of replacing  $\Pi_{v}$ by (one of locally finitely many) ramified twists. In this case, that $P_{1, \Pi_{v}}$ is a  polynomial in $Y_{v}^{\circ}$ and the values of $W_{v}$  is  one of the main results of \cite{JPSS}, whose proof considers unramified characters of the form $|\cdot|_{v}^{s}$ but goes through in our context.
\end{proof}

\subsection{Test Hecke measures} \lb{sec: test}

We now give some key results asserting the existence of suitable Hecke measures.

\subsubsection{Test measures at finite places} \lb{def reg}
Let $v$ be a finite place of $F_{0}$ and let $L$ be a coefficient field.
  A character $\xi' =\xi'_{1}\boxtimes \cdots \boxtimes \xi'_{\nu}\colon (F_{w}^{\times})^{\nu}\to \C^{\ts}$ is called \emph{regular} if the characters $\xi_{i}'$ are pairwise distinct. A \emph{regular principal series} representation of $G_{v}'$ is a representation $\Pi_{v}=\Pi_{n, v}\boxtimes \Pi_{n+1, v}$ such that for $\nu=n, n+1$, all places $w\vert v$, and any $\iota\colon L\into \C^{\times}$ the representation $\Pi_{\nu, w}\coloneqq \Pi_{\nu,v|\GL_{\nu}(F_{w})}$ 
is unitarily induced from a regular character of the diagonal torus.

\begin{proposition}\lb{test fin} Let $\Pi_{v}$  be a hermitian (\S~\ref{loc BC}) tempered representation of $G'_{v}$ over $L$, and let $\chi_{v}$ be a smooth character  of $F_{0,v}^{\ts}$ with values in some finitely generated extension $L'$ of $L$. For  $?\in \{\emptyset, \reg^{\pm}\}$,  denote by 
$$\sH(G'_{v}, L)_{ ?, \Pi_{v}, \chi_{v}}$$ 
the set of those Hecke measures $f_{v}'\in \sH(G_{v}', L)$ that are  supported in $G'_{?, v}$,  and satisfy $I_{\Pi_{v}}(f'_{v}, \chi_{v})\neq 0$. 

\begin{enumerate}
\item We have
$\sH(G'_{v}, L)_{ \Pi_{v}, \chi_{v}}\neq\emptyset$.
\item If $\Pi_{v}$ and $\chi_{v}$ are unramified, then $f_{v}'^{\circ}\in\sH(G'_{v}, L)_{\Pi_{v}, \chi_{v}}$.
\item \lb{test fin 3} If $v$ splits in $F$ and  $\Pi_{v}$ is a regular principal series, for every choice of sign $\pm$ there exists   a Hecke measure
$$f'_{\pm}\in \sH(G'_{v}, L)_{ \reg^{\pm}, \Pi_{v}, \chi_{v}}.$$
whose matching  measure $f_{\pm}\in \sH(G_{v}, L)$ is invariant under a deeper Iwahori subgroup.
If moreover $\Pi_{v}$ is unramified, we can take $f'_{\pm}$ to match a measure $f_{\pm}\in \sH(G_{v}, L)$ that is bi-invariant under an  Iwahori subgroup.
\end{enumerate}
\end{proposition}
The proof of part (3) relies on some  explicit results from later sections. (In fact, see \S~\ref{iw conjsy} for a definition of Iwahori   subgroups.)
\begin{proof} We omit all subscripts $v$. 
\begin{enumerate}
\item Let $K\subset G'$ be an open compact subgroup. The restriction of $I_{\Pi}(\cdot, \chi)$ to $\sH(G')_{K}$ is the inner product, for the natural pairing, of the elements 
$$P_{1, \Pi,\chi|\Pi^{K}}\circ \Pi(\cdot ) \in \Pi^{K, \vee}\ot_{L}L' , \qquad P_{2|\Pi^{\vee, K}}\in (\Pi^{\vee, K})^{ \vee} \cong \Pi^{K}.$$
 Now if $K$ is sufficiently small, both $P_{1, \Pi,\chi|\Pi^{K}}$ and $P_{2|\Pi^{\vee, K}}$ are nonzero -- the former by the theory of \cite{JPSS}, the latter  because $\Pi$, hence $\Pi^{\vee}$, is hermitian. Since $\Pi^{K}$ is irreducible as an $\sH(G', L)_{K}$-module, there exists an element $f_{L'}'\in \sH(G', L')_{K}$ such that 
$P_{1, \Pi,\chi|\Pi^{K}}\circ \Pi(f'_{L'})$ and $P_{2|\Pi^{\vee, K}}$ do not pair to zero.  Fix an embedding $\iota\colon L'\into \C$. If $\iota(L)$ is not contained in $\R$, then it is dense in $\C$, and  any  $f'\in \sH(G', L)_{K}$ that is sufficiently close to $f'_{L'}$ in the  topology induced from $\C$ by $\iota$ will also have the desired non-vanishing property. If $\iota(L)$ is contained in $\R$, note that one of $\Re\, \iota f'_{L'}$, $\Im \iota f'_{L'}$ has the non-vanishing property, and  then so does any sufficiently close  $f'\in \sH(G', L)_{K}$ (for the topology induced from $\R$ by $\iota$).
\item This follows from Remark \ref{aa1}.
\item  This will be proved  at the end of \S~\ref{sec prove test fin} based on an explicit construction from  \S~\ref{sec suff pos}.
\end{enumerate}
\end{proof}

\subsubsection{Test Gaussians}\lb{sec E gauss}

  For a pure tensor $f_{S\infty}'=f_{S}' f_{\infty}'\in \sH(G'_{S\infty}, L)^{\circ}$ and an embedding $\iota\colon L\into \C$, we define $f'^{\iota}_{S\infty}\coloneqq  \iota f_{S}'  f_{\infty}'^{\iota}$, and extend this definition to all of $\sH(G'_{S\infty}, L)^{\circ}$ by linearity. 

\begin{proposition} \lb{test gauss}
Let $\Pi$ be a trivial-weight hermitian cuspidal automorphic representation of $\G'(\A)$ over a coefficient field $L$, let $K=\prod_{v\nmid\infty}K_{v}\subset \G'(\A^{\infty})$ be an open compact subgroup such that $\Pi^{K}\neq 0$, and let $P$ be a finite set of non-archimedean places of $F_{0}$.

 There exist a finite set $S$ of split non-archimedean places of $F_{0}$, disjoint from $P$ and the set of places where at which $K_{v}$ is not maximal, and Gaussians
 $$(f_{S\infty}'^{\iota})_{\iota} \in \prod_{\iota\colon L\into \C} \sH(G'_{S\infty}, \iota L)_{K_{S}}^{\bullet}, \qquad 
 f_{S\infty}' \in \sH(G'_{S\infty},  L)_{K_{S}}^{\circ}
 $$
  such  that     for every $\iota\colon  L\into \C$:
   \begin{enumerate}
   \item the image of $f_{S\infty}'^{\iota}$ in $\sH(G'_{S\infty}, \iota L)_{K_{S}}^{\circ}$ equals $\iota f'_{S\infty}$;
 \item 
$I_{\Pi^{\iota}_{S\infty}}(f'^{\iota}_{S\infty}, \chi_{S\infty})\neq 0$ for every unramified character $\chi_{S\infty}\colon F_{0,S}^{\ts}F_{0, \infty}^{\ts}/F_{0, \infty}^{\ts}\to \C^{\ts}$;
 \item  $R(f'^{\iota}_{S\infty})$ maps $\sA(\G')^{K}$ into $(\Pi^{\iota})^{K}$. (In particular, for any $f'^{S\infty}\in \sH(G',\iota L)_{K^{S}}$, the Hecke measure $f'^{S\infty}f'^{\iota}_{S\infty}$ is quasicuspidal.)
 \end{enumerate}
\end{proposition}
The proof will be given in \S~\ref{sec: prove Gauss}.

\begin{lemma} \lb{cebo} Let $\Pi$ be a representation in $\sC(\G')$. There exist infinitely many places $v$ of $F_{0}$ that are split in $F$ such that $\Pi_{v}$ is an unramified regular principal series.
\end{lemma}
\begin{proof} This follows from the similar observation about $\Pi_{\nu}$ made in the proof of  \cite[Proposition 3.2.5]{CH}.
\end{proof}

\begin{corollary} \lb{cor test} Let $\Pi$ be a trivial-weight hermitian cuspidal automorphic representation of $\G'(\A)$ over a  coefficient field $L$,
and let $\chi\in Y({L})$. Let $P$ be a finite set of nonarchimedean places of $F_{0}$ and let $K_{P}\subset G'_{P} $ be a compact open subgroup such that $\Pi_{P}^{K_{P}}\neq 0$.

For  $?\in \{\reg^{+}, \reg^{-} , \rs\}$, there exist $L$-rational Gaussians
$$f'^{P}_{?}\in\sH(G'(\A^{P}), L)^{\circ}_{K_{P}, ?,\Pi, \chi}$$
with globally $?$-support that are adapted to $(\Pi , \chi, K_{P})$ (in the sense of \S~\ref{rha}). 
\end{corollary}

\begin{proof} In fact we construct a Gaussian $f'^{P}$ that has at the same time plus-regular support (at one place) and minus-regular support (at another place), hence globally regular semisimple  support since $\G'_{\rs}=\G'_{\reg^{+}}\cap \G'_{\reg^{-}}$. (The  construction can of course be simplified if only one of those two properties is desired.)  Let $R$ be the set of all finite places of $F_{0}$ at which $\Pi$ or $\chi$ is ramified. Let $v_+$, $v_{-}$ be  two distinct finite places of $F_{0}$, split in $F$, not in $P\cup R$, such that $\Pi_{v_{\pm}}$ is a regular  principal series and $\chi_{v_{\pm}}$ is unramified. Let $f'_{\pm, v_{\pm}}$ be as in Proposition  \ref{test fin} \eqref{test fin 3}, and for $v\in R$ let $f_{v}'$ be as in Lemma \ref{test fin} (1). Let $f'_{S\infty}$ be as given by  Proposition \ref{test gauss} for the set of places  $P'=P\cup R\cup \{v_{+}, v_{-}\}$, and any level $K$ that is maximal away from $P'$ and sufficiently small at the places in $R\cup \{v_{+}, v_{-}\}$. Then 
$$f'^{P}=f'_{+, v_{+}}f'_{-, v_{-}}f'_{S\infty}\ot\bigotimes_{v\nmid PS\infty}f_{v}'^{\circ}$$
 is as desired.
\end{proof}

\subsection{Isolation of cuspidal representations via Gaussians}
\lb{sec: prove Gauss}
In this subsection, we prove  Proposition \ref{test gauss}.

We will refine the arguments of \cite{isolation}, of which the reader is invited to open a copy. Briefly, in order to construct the desired $f_{S\infty}'$ we will start from a Gaussian $f_{1, S'\infty}'$ constructed in a simple way as a pure tensor, and then correct $f_{1, S'\infty}'$ by 
acting on it by a carefully chosen  \emph{multiplier}   of the Hecke algebra for $\G'(\A)$. 

The substance of this subsection was generously provided to us by Yifeng Liu. Of course, any defects in the following pages are to be attributed to the authors only.

\subsubsection{Archimedean multipliers annihilating  non-strongly typical cuspidal data} We momentarily consider the  more general situation of \cite[\S~3.2]{isolation}.
Consider a connected reductive algebraic group $\G$ over a number field $F_{0}$. We freely adopt some of the notation from \cite{isolation}*{\S~3}, up to cosmetic modifications to adapt to our conventions (for instance, in \emph{loc. cit.} the algebraic group is denoted by $G$ rather than $\G$). Take a unitary automorphic character $\omega\colon \RZ(\A)\to{\C}^\times$. We fix
a  character 
$$\xi_\infty\colon \mathcal{Z}({\frak g}) \to \C.$$

The following definition is modified from \cite{isolation}*{Definition~3.11}; the set  $\mathfrak{C}(\RM,\omega)^\heartsuit$, consisting of classes of cuspidal automorphic representations of $\RM(\A)$, is defined \emph{ibid.} p. 550.

\begin{definition}\label{de:typical} Let $\RM\subset \G$ be a standard Levi subgroup.
We say that an element  $\sigma\in\mathfrak{C}(\RM,\omega)^\heartsuit$ is \emph{strongly $\xi_\infty$-typical} if $\gamma_\RM(\xi_{\sigma_\infty})\subseteq\gamma_\RM(\xi_\infty)$, where $\gamma_\RM\colon \frak h^*\to  \frak h^*/\frak a^*_\RM$ is the quotient map of \cite[p. 552]{isolation}. We denote by $\mathfrak{C}(\RM,\omega)^\heartsuit_{\xi_\infty !}$ the subset of $\mathfrak{C}(\RM,\omega)^\heartsuit$ consisting of strongly $\xi_\infty$-typical elements.

\end{definition}

It is clear that the set  $\mathfrak{C}(\RM,\omega)^\heartsuit_{\xi_\infty !}$ of  strongly $\xi_\infty$-typical elements is a subset of  $\mathfrak{C}(\RM,\omega)^\heartsuit_{\xi_\infty}$, the set of   $\xi_\infty$-typical elements defined in \emph{loc. cit.} 
The following lemma slightly strengthens \cite{isolation}*{Lemma~3.14}, whose notation we simplify by putting   
$$\CM_{\infty}\coloneqq \mathcal{M}^\sharp_\theta(\mathfrak{h}_{\bC}^{\*})^{\mathsf{W}}$$
for the Weyl-fixed elements of the space of holomorphic functions from  \cite[Definition 2.8]{isolation}.

\begin{lemma}\label{le:finite3}
Fix an element $\mu_{\infty}^{0}\in \CM_\infty$ and a finite set $\frak T$ of $K_{0, \infty}^{\RG}$-types as described after \cite[Definition 3.11]{isolation}. 
For every standard Levi subgroup $\RM\subset \G$, every open compact  subgroup $K_{\RM}\subset \RM(\A^{\infty})$ and every finite set $\frak T_{\RM}$ of $K_{0, \infty}^{\RM}$-types satisfying the conditions following \cite[Definition 3.11]{isolation} with respect to $(\mu_{\infty}^{0}, \frak T)$, there exists an element
 $$\mu_\infty^\RM\in\CM_{\infty}
$$
  satisfying:
  \begin{itemize}
  \item $\mu_\infty^\RM(\xi_\infty)\neq 0$,
  \item  for every non-strongly typical
  $$\sigma\in\mathfrak{C}(\RM,\omega;K_\RM,\mathfrak{T}_\RM)^\heartsuit-\mathfrak{C}(\RM,\omega)^\heartsuit_{\xi_\infty !}$$ (see \cite[p. 553]{isolation} for the definition) and  every $s\in\mathfrak{a}_{\RM,{\bC}}^*$,  we have
$$\mu_\infty^\RM(\xi_{\sigma_{s,\infty}}^\G)=0.$$
Here, $\xi_{\sigma_{s,\infty}}^\G$ is  the infinitesimal character of $\Ind^\G_{\RP_\RM}(\sigma_{s,\infty})$, where $\sigma_{s,\infty}$ is the unramified twist defined \emph{ibid.} p. 550.
  \end{itemize}
\end{lemma}

\begin{proof}
By Definition \ref{de:typical}, it is easy to see that for each element $\sigma\in\mathfrak{C}(\RM,\omega;K_\RM,\mathfrak{T}_\RM)^\heartsuit- \mathfrak{C}(\RM,\omega)^\heartsuit_{\xi_\infty!}$, there exists a $\mathsf{W}$-invariant polynomial function $\nu_\sigma$ on $\mathfrak{h}_{\bC}^*$ satisfying $\nu_\sigma(\xi_\infty)\neq 0$ and $\nu_\sigma(\xi_{\sigma_{s,\infty}}^\G)=0$ for every $s\in\mathfrak{a}_{\RM,{\bC}}^*$. By \cite{isolation}*{Lemma~3.14}, we have an element $\nu_\infty^\RM\in \CM_\infty$
satisfying the same property but with $\mathfrak{C}(\M,\omega)^\heartsuit_{\xi_\infty !}$ replaced by $\mathfrak{C}(\M,\omega)^\heartsuit_{\xi_\infty}$. Now by \cite{isolation}*{Lemma~3.13}, the set
$$\mathfrak{C}'
 \coloneqq \mathfrak{C}(\RM,\omega;K_\RM,\mathfrak{T}_\RM)^\heartsuit
\cap\left(\mathfrak{C}(\RM,\omega)^\heartsuit_{\xi_\infty}-\mathfrak{C}(\RM,\omega)^\heartsuit_{\xi_\infty !}\right)$$ is finite. Thus, we may take
\[
\mu_\infty^\RM\coloneqq\nu_\infty^\RM\cdot\prod_{\sigma\in \mathfrak{C}'}
\nu_\sigma.
\]
\end{proof}

The background for the following proposition is the Langlands decomposition
$$
L^{2}(\G(F_{0})\bs \G(\A) , \omega)=\widehat{\bigoplus}_{(\RM, \sigma) \in {\frak D}(\RG,\omega)^{\heartsuit}} 
L^{2}_{(\RM, \sigma)}(\G(F_{0})\bs \G(\A), \omega)
$$
of \cite{isolation}*{(3.1)}
in terms of a set ${\frak D}(\RG,\omega)^{\heartsuit}$ of classes of cuspidal data. We denote by $\star$ the action of $\CM_\infty$ on $\sH(G_\infty, \C)$ from \cite[Theorem 2.18]{isolation}.
\begin{proposition}\lb{kill not st} There exists $\mu_{\infty}\in \CM_{\infty}$ and a finite set $\frak T$ of $K_\infty^\G$-types, such that: 
\begin{itemize}
\item 
$\mu_{\infty}(\xi_{\infty})=1$;
\item for every cuspidal datum  $(\RM, \sigma)$ for $\G'$ that does not belong to the finite set  ${\frak D} (\G, \omega, K, \frak T)_{\xi_{\infty}!}^{\heartsuit}$ defined in \cite[p. 555]{isolation}, and for every
$f\in \sH(\G(\A), \C)_{K}$, the endomorphism $R(\mu_{\infty} \star f)$ of $L^{2}(\G(F_{0})\bs \G(\A)/K, \omega)$ annihilates the subspace $L^{2}_{(\RM, \sigma)}(\G(F_{0})\bs \G(\A)/K, \omega)$.
\end{itemize}
\end{proposition}
\begin{proof}
    This is a direct analogue of \cite{isolation}*{Proposition~3.15}, in which we start with $(\mu_\infty^0, \frak T)$ as given at the end of p. 552 \emph{ibid.} (depending only on $\xi_\infty$), and we
    take $\mu_\infty=\mu_\infty^0\prod_\RM\mu_\infty^\RM$ where the $\mu_\infty^M$ are as provided by Lemma \ref{le:finite3} instead of \cite{isolation}*{Lemma 3.14}. 
\end{proof}

\subsubsection{Multipliers annihilating strongly typical cuspidal data for a proper Levi subgroup} 
We now specialize  back to the setup of Proposition \ref{test gauss}.
We denote by $\xi_{\infty}^{\circ}$ the infinitesimal character of~$\Pi_{\infty}^{\circ}$.  
 We still freely use  terminology and notation from \cite[\S~3]{isolation} where not in conflict with ours.

\begin{lemma}\lb{lem clozel}
Let $\RM\subset \G'$ be a standard Levi subgroup, and
let $(\RM, \sigma) \in {\frak C} (\RM, 1)^{ }$ satisfy
$$\xi^{\G'}_{\sigma, \infty}=\xi_{\infty}^{\circ}.$$
Then for every finite place $v$ at which $\sigma_{v}$ is unramified, the Satake parameters of $\sg_{v}$ are algebraic.
\end{lemma}
\begin{proof}
We start with some preliminaries. If $\RM_{\nu}$ is a Levi subgroup of $\GL_{\nu/F}$ of type $(a_{1}, \ldots, a_{r})$, put $\psi_{\RM_{\nu}}\coloneqq \boxtimes_{i=1}^{r}\vert \det\vert^{\nu-a_{i}\over 2}$. 
 If $\sg_{\nu}$ is a representation of $\RM_{\nu}(\A)$, denote $\sg_{\nu}^{\natural}\coloneqq\sg_{\nu}\ot \psi_{\RM_{\nu}}$. 
    We  extend the definitions to the case of Levi subgroups of $\G'$ in the obvious way. 
    For $\RM $ a Levi subgroup of $\G'$, let $ \boxplus$ be the isobaric sum operation  introduced (for general linear groups) in \cite[p. 85]{clozel}, which sends  cuspidal automorphic representations of $\M(\A)$ to automorphic representations of $\G'(\A)$.  Let $ \boxplus^{T}$ be the twisted version 
 $ \boxplus^{T} \sg^\nat\coloneqq (\boxplus\, \sg)^{\natural}$ of \cite[D\'efinition 1.9]{clozel}.   The operation  $ \boxplus^{T}$  preserves fields of rationality and induces the direct sum operation on infinity types  \cite[Lemme 3.9 (ii)]{clozel}.
 
Let $a_{\infty}^{\circ}$ be the infinity type associated with $\xi_{\infty}^{\circ}$, which is regular (\cite[D\'efinition 3.12]{clozel}).
Since any  direct summand of $a_{\infty}^{\circ}$
 is also regular, it follows that $\sg^{\natural}$ is regular algebraic. Thus by  \cite[Th\'eor\`eme 3.13]{clozel}, it is defined over a number field. It follows that its algebraic twist $\sg$ has algebraic  Satake parameters.
\end{proof}

 For a characteristic-zero field $L$ and a finite set $P$ of non-archimedean places of $F_0$, define 
$\BT^{{\rm spl}, P}_{L}\subset \sH(\G'(\A^{P\infty}), L)_{\prod_{v\nmid P\infty}K_{v}^{\circ}}$ to be the spherical Hecke algebra of elements supported at a set of places of $F_{0}$ split in $F$ and disjoint from $P$.
If  $L$ is a subfield of  $\C$, define  $\mathcal{M}_{\infty, L}$
to be the $L$-linear subspace of $\CM_{\infty}$ consisting of those $\mu$ such that $\mu(\xi_{\infty}^{\circ})\in L$. We put 
\[
\mathcal{M}^{{\rm spl}, P}_{L}\coloneqq  \BT^{{\rm spl}, P}_{L} \ot_{L}\mathcal{M}_{\infty,L} ,
\]
which is stable under multiplication. By \cite[Theorem 2.18]{isolation} and the spectral characterization of matching, we see  that the $\star$-action of
$\mathcal{M}^{{\rm spl}, P}_{L}$ preserves $\sH(\G'({\A}),L)^\bullet_{K}$ and descends to an action on $\sH(\G'({\A}),L)^\circ_{K}$.
We have a surjective map 
\beqq
\ [ - ]^{\circ}\colon \mathcal{M}^{{\rm spl}, P}_{L}&\to\BT'^{{\rm spl}, P}_{L}\\
\mu &\mapsto [\mu]^{\circ}
\eeqq
given by the evaluation at $\xi_{\infty}^{\circ}$. It is clear that the  action of $ \mathcal{M}^{{\rm spl}, P}_{L}$ on $\sH(\G'(\A), L)^{\circ}$ factors through $[-]^{\circ}$.

For $\tau\in {\rm Aut}(\C/\Q)$ and $\mu\in  \mathcal{M}^{{\rm spl}, P}_{L}$, we denote by $\tau.\mu\in \CM_{\tau L}$ a chosen lift of $\tau([\mu]^{\circ})$. We also define  an action of ${\rm Aut}(\C/\Q)$ on the set  of strongly $\xi_\infty^\circ$-typical cuspidal data $(\RM, \sg)$ by the requirement that $\tau.(\RM, \sg) = (\RM, \tau.\sg) $ for $\tau.\sg $  that satisfies ${\rm Ind}_{\rm P_\RM}^{\G'}(\tau.\sg)^\infty \cong {\rm Ind}_{\rm P_\RM}(\A^\infty)^{\G'(\A^\infty)}(\sg^\infty)^\tau$ (where in general we denote by $\Sigma^\tau$ the representation of $\G'(\A^\infty)$ on the space of $\Sigma\ot_{\C, \tau} \C$) . Note that 
by \cite[Théorème 3.13]{clozel}, each  orbit under this action is finite.

We denote by $\sC(\G')_{K}\subset \sC(\G')$ the subset consisting of those $\Pi'$ with $\Pi'^{K}\neq 0$.   

\begin{lemma}\lb{mu sg} Let  $L'\subset \C$ be a subfield and let $\Pi'\in \sC(\G')_{K}(L')$. Let $(\RM, \sigma)$ be a strongly $\xi_{\infty}^{\circ}$-typical cuspidal datum for $\G'$ with $$\RM\neq \G', $$ 
and 
let  $P$ be a finite set  of non-archimedean places of $F_0$.

There exists an element  $\mu'\in\mathcal{M}^{{\rm spl}, P}_{L'}$   satisfying:
 \begin{itemize}
 \item  for every $f'\in\mathcal{S}(\G'({\A}), \C)_{K}$ and every $\tau\in {\rm Aut}(\C/\Q)$,
the endomorphism $R(\tau.\mu'\star \tau f')$ of $L^2(\G'(F_{0})\backslash \G'(\A)/K)$ annihilates the subspace $L^2_{\tau.(M,\sigma)}(\G'(F_{0})\backslash \G'(\A)/K)$;
 \item $\tau.\mu'(\xi_{(\Pi')^\tau}^{P})=1$.
 \end{itemize}
\end{lemma}
\begin{proof} 
We refine the argument in the proof of \cite{isolation}*{Proposition~3.17}.\footnote{With respect to the notation of \emph{loc. cit.}, we omit the central character $\omega$, which in our setup is necessarily trivial.}  We may assume that $L'$ is a finite extension of $\Q$ inside $\C$.
Note that the subspace $\mathfrak{a}_\RM^*\subseteq\mathfrak{h}'^*$ has a natural model $\mathfrak{a}_{\RM,{\Q}}^*\subseteq\mathfrak{h}'^*_{\Q}$ over~${\Q}$. We fix a rational splitting map $\ell\colon\mathfrak{h}'^*_{\Q}\to\mathfrak{a}_{M,{\Q}}^*$ and an element $\alpha\in\xi_\infty^{\circ}$. Denote by $\ot_v\xi_v\colon \BT^{{\rm spl}, P}_{L} \to L'$ the character attached to $\Pi'$.
By  Ramakrishnan's automorphic  Tchebotarev theorem  (Proposition \ref{Ram-prop}),
for every $w\in\mathsf{W}'$, there is a finite place $v[w]\notin P$ of $F_{0}$, split in $F$, 
such that $\xi^{G'}_{\sigma_{s_w,v[w]}}\neq\xi_{v[w]}$ where $s_w\coloneqq\ell(w\alpha)-\ell(\alpha)\in\mathfrak{a}_{M,{\Q}}^*$.  

This allows us to choose an element $\nu_w\in  \sH(G'_{v[w]}, L')_{K'_{v[w]}}$
such that
\[
\nu_w(\xi_{v[w]})\neq\nu_w(\xi^{G'}_{\sigma_{s_w,v[w]}}).
\]
First we show the existence of a multiplier $\mu''\in\mathcal{M}^{{\rm spl}, P}_{L''}$ with the required properties, for some finite Galois extension $L'\subset L''\subset \C$. 
By the process in the proof of \cite{isolation}*{Proposition~3.17}, it suffices to show that for every $w'\in\mathsf{W}'$, the value $\nu_w(\xi^{G'}_{\sigma_{s_{w'},v[w]}})$ is algebraic. By Lemma \ref{lem clozel}, the Satake parameters of $\sg_{w',v[w]}$ are algebraic numbers. Since $s_{w'}\in\mathfrak{a}_{\RM,{\Q}}^*$, it follows that the Satake parameters of $\sigma_{s_{w'},v[w]}$ are all algebraic numbers as well, which implies that $\nu_w(\xi^{G'}_{\sigma_{s_{w'},v[w]}})$ is algebraic. 

Now it is clear from the construction that $\mu'\coloneqq\prod_{\tau\in \Gal(L''/L')} \tau.\mu''$ satisfies the desired properties.
\end{proof}

We now extend the result to a finite set of cuspidal data and descend it to abstract coefficient fields.

\begin{proposition} \lb{kill notG}
Let $L$ be a coefficient field and let $\Pi\in \sC(\G')_K(L)$. 
Let $\frak D$ be a finite set of  strongly $\xi_{\infty}^{\circ}$-typical cuspidal data for $\G'$ such that $$\RM\neq \G'$$ for every $(M, \sigma)\in \frak D$. 
Let  $P$ be a finite set  of non-archimedean places of $F_0$.  There exists a collection 
$$(\mu_{\frak D}^{\iota})_{\iota}\in\prod_{\iota \colon L\into \C} \CM_{\iota L}^{{\rm spl}, P}$$
satisfying: 
 \begin{enumerate}
 \item  for every $f'\in\mathcal{S}(\G'({\A}), \C)_{K}$, every $\iota\colon L\into \C$, and every $(M, \sigma)\in \frak D$, the endomorphism $R(\mu_{\frak D}^{\iota}\star f'^{\iota})$ of $L^2(\G'(F_{0})\backslash \G'(\A)/K)$ annihilates the subspace $L^2_{(M,\sigma)}(\G'(F_{0})\backslash \G'(\A)/K)$;
\item    $\mu_{\frak D}^{\iota}(\xi_{\Pi^{\iota}}^{P})=1$ for every $\iota \colon L\into \C$.
 \item there exists an element $[\mu_{\frak D}]^{\circ}\in \BT'^{{\rm spl}, P}_{L}$ such that $[\mu_{\frak D}^{\iota}]^{\circ}=\iota[\mu_{\frak D}]^{\circ} $ for every $\iota\colon L\into \C$;
 \end{enumerate}
\end{proposition} 
\begin{proof} 
We denote the elements of $\frak D$ simply by $\sg$ in order to lighten the notation. 

For every  embedding $\iota\colon L \into \C$, let
$\mu_{\sg, \iota}\in \CM_{\iota L}^{{\rm spl}, P}$ be as provided  by Lemma \ref{mu sg} applied to $\sg$ and $\Pi^{\iota}$.
Let $L'$ be a Galois extension of $\Q$ in $\C$ containing $\iota L$ and the fields of definition of $\mu_{\sg, \iota}$ for every $\sg\in \frak D$ and every $\iota\colon L\into \C$.  Now take the collection
 $$\mu_{\frak D}^{\iota}\coloneqq \prod_{\tau\in \Gal(L'/\Q) } \prod_{\sigma \in {\frak D}} \tau. \mu_{\sigma, \tau^{-1}\iota}.$$
 
By Lemma \ref{mu sg}, it satisfies the first two desiderata. For the third  one, by  Galois theory we need to check that for each $\tau'\in \Gal(L'/\Q)$, we have $\tau'. \mu_{\frak D}^{\iota} = \mu_{\frak D}^{\tau'\iota}$: indeed, by a change of variables 
\beqq
\tau'([ \mu_{\frak D}^{\iota}]^{\circ})
&= 
 \prod_{\tau\in \Gal(L'/\Q) } \prod_{\sigma \in {\frak D}} \tau'\tau( [ \mu_{\sigma, \tau^{-1}\iota}]^{\circ})=
   \prod_{\tau\in \Gal(L'/\Q) } \prod_{\sigma \in {\frak D}} \tau( [ \mu_{\sigma, \tau^{-1}\tau'\iota}]^{\circ}) =[ \mu_{\frak D}^{\tau'\iota}]^{\circ}.
\eeqq
\end{proof}

\subsubsection{Proof of Proposition \ref{test gauss}}
Let $f_{1,\infty}'\in \sH(G'_{\infty}, \Q)^{\bullet}$ be a nontrivial rational Gaussian, which exists by Proposition \ref{gauss E}. By Proposition \ref{Ram-prop}, we can find a finite set  $S_{1}$ of split places of $F_{0}$, disjoint from $P$ and the ramification set of $\Pi$, and an  element $$f_{1,S_{1}}'\in \sH(G_{S_{1}}, L)_{K_{S_{1}}},$$
 such that $\Pi(f_{1, S_{1}})\neq 0$ and $\Pi'(f_{1, S_{1}})=0$ for every $\Pi'\in \sC(\G')_{K, L} - \{\Pi\}$.  For each $\iota\colon L \into \C$, let $f_{1}'^{\iota}= \iota f'_{1, S_{1}}\ot  f'_{1, \infty}\ot\ot_{v\notin S_{1}\infty}f_{v}^{\circ}$.
   
 Let $\mu_{\infty}$ and $\frak D= {\frak D} (\G, \omega, K, \frak T)_{\xi_{\infty}!}^{\heartsuit}$ 
 be as provided by Proposition \ref{kill not st} applied to $\xi_\infty$; the set  $\frak D$ is finite and consists of  strongly $\xi_{\infty}^{\circ}$-typical cuspidal data for $\G'$. Let $(\mu_{\frak D}^{\iota})_{\iota}$, $[\mu_{\frak D}]^{\circ}$ be as provided by Proposition \ref{kill notG} for $\frak D$. Let 
  $$ (f'^{\iota})\coloneqq  \mu_{\frak D}^{\iota}\star \mu_{\infty}\star f_{1}'^{\iota}, 
  \qquad f' \coloneqq  [\mu_{\frak D}]^{\circ} \star [\mu_{\infty}]^{\circ}\star f_{1}'$$
By construction, there is a set of split places $S\supset S_{1}$ disjoint from $P$ such that for $?=\iota, \emptyset$, we have $f'^{?}=f'^{?}_{S\infty}\ot \ot_{v\notin S} f_{v}^{\circ}$ for some 
$$f'^{\iota}_{S\infty},\qquad f'_{S\infty}$$
  that satisfy the desired properties. The proof is complete.

\subsection{Proofs of the rationality statements}\lb{4pf} \lb{last 422 pf}
 We will prove  Proposition \ref{rat rtf} \eqref{rat rtf part} (recall that the other parts were proved at the end of \S~\ref{sec: rat st}) 
 and, as an interlude, Theorem \ref{rat L}. 

\subsubsection{Global distribution}
The global orbital-integral distributions $I_{\gamma}$  of part \eqref{rat orb int 2 reg}  are well-defined and we  may define the  distribution $I$ of part \eqref{rat rtf part} by its asserted geometric expansion:
$$I\coloneqq  \sum_{\gamma\in \RB'(F_{0})} I_{\gamma}.$$
We  show the sum is locally finite. We may assume that $f' \in   \sH(\G'(\A), L)^{\circ}_{\reg^+, \rm qc}$ factors as  $f'=f'^{\infty}\ot f_{\infty}'$. By definition, the sum is supported in $ \RB'(F_{0})\cap B_{\infty}'^{\circ}$. The invariant map \eqref{eq inv} sends $\RB'$ isomorphically to a closed subvariety of the affine space ${\rm Res}_{F/F_{0}}\BA^{2n+1}$. Let $\Omega^{\infty}\subset (\A_{F}^{\infty})^{2n+1}$ be the image of the support of  $f'^{\infty}\in \sH(\G'(\A^{\infty}))$, which is compact. Let $\Omega_{\infty}\subset F^{2n+1}_{\infty}$ be the image of $B'^{\circ}_{\infty}$. By definition, this is contained in the image of the positive-definite unitary group $G_{\infty}^{\circ}$ under the invariant map, which is compact. Therefore the support of the sum is in bijection with a subset of the set $F^{2n+1}\cap \Omega^{\infty}\Omega_{\infty}$; as the first intersecting set is discrete and the second one is compact, the intersection is finite.

By construction, $I$ has the geometric expansion asserted in part \eqref{rat rtf part}; 
by Prop. \ref{exps rtf}, it satisfies
\beq \lb{I rat kappa} 
I(f', \chi)=
 \kappa(\one_{\infty})^{-1} 
I_{}^{\C}(f'^{\iota}, \chi) 
\eeq
for any $\chi\in Y_{L}(\C)$ with underlying embedding $\iota\colon L\into \C$, and any $f'^{\iota}\in \sH(G'(\A), \C)^{\bullet}_{\reg^+, {\rm qc}}$ mapping to $\iota f'$.

\begin{remark} \lb{gl int ell}
By linearity, we may extend the distributions $I$, $I_{\Pi}$, $I_{\gamma}$ to distributions
(defined over $\Q$ or, for $I_{\Pi}$, over the field of definition of $\Pi$)
on the space of locally constant functions  $\ell\colon F_{0}^{\ts}\bs\A^{\ts}/F_{0, \infty}^{\ts}\to L $,  in such a way that for 
every  $\gamma'\in \RG'_{\rm rs}(F_{0})$ with image $\gamma\in \RB_{\rs}'(F_{0})$ and every $f'^{\infty}\ot f'_{\infty}\in \sH(\G(\A), L)^{\circ}$, we have
$$
I_{\gamma}(f'^{}, \ell)= {I_{\gamma}(f'_{\infty})
\over  \kappa(\one_{\infty}) \kappa_{\infty}(\gamma', \one)} 
 \int_{\RH_{1}(\A^{\infty})}\int_{\RH_{2}(\A^{\infty})} \ff'^{\infty}(h_{1}^{-1}\gamma'  h_{2}) \ell(h_{1}) \eta(h_{2}) \, {d^{\natural}{h_{1}}d^{\natural}h_{2}\over d^{\natural}g},
 $$
where $d^{\natural}x\coloneqq \prod_{v\nmid \infty}d^{\natural}x_{v}$, and the integral reduces to a finite sum.
\end{remark}
\subsubsection{$L$-function}
We are now  ready to prove the rationality of $\sL$.

\begin{proof}[Proof of Theorem \ref{rat L} (= Theorem \ref{thm A})]
For any  $\chi\in Y_L$, consider the set 
$\sH(\G'(\A^{\infty}) , L(\chi))^{\circ}_{\reg^{+}, \Pi, \chi}$  of  Gaussians with globally plus-regular support that are adapted to $(\Pi, \chi)$ in the sense of \S~\ref{rha}.  It is non-empty by Corollary \ref{cor test}.
 For    $f'\in \sH(\G'(\A^{\infty}) , L)^{\circ}_{\rs, \Pi, \chi}$, 
 we define
 $$
 \sL(\RM_{\Pi}, \cdot )_{f'}\coloneqq  {4 \cdot I(f', \cdot) \over
  (\ot_{v} I_{\Pi_{v}})(f'_{v}, \cdot)},
  $$ as a function on the open subset of $Y_{L(\chi)}$ where the denominator does not vanish.
Then  for any $\chi\in Y_{L}$, any  $\iota\colon L(\chi)\into \C$, any $f'\in \sH(\G'(\A^{\infty}) , L(\chi))^{\circ}_{\reg^{+}, \Pi, \chi}$, and any $f'^{\iota}$ as in \S~\ref{rha}, we have
$$\sL(\RM_{\Pi}, \chi)_{f'}  = 
{4 \cdot I^{\C}_{\Pi^{\iota}}(f'^{\iota},  \chi)  \over 
 \kappa(\one_{\infty})
( \ot_{v\nmid \infty} I^{\C}_{\Pi^{\iota}_{v}}\ot I^{\circ, \C}_{\Pi^{\circ}_{\infty}})(f'^{\iota}, \chi)}
= { \sL^{\infty}(1/2, \Pi^{\iota},\chi) 
\over  \ep(\textstyle{1\over 2} , \chi^{2})^{{n+1\choose 2}}}
$$
where the first equality is \eqref{I rat kappa},  and the second one is  \eqref{factor I rat}.  It is also clear from the definition that for any $\chi\in Y_L$ and any isomorphism of fields $\sg\colon L(\chi)\to L'$ fixing $L$, and any $f' \in \sH(\G'(\A^{\infty}) , L(\chi))^{\circ}_{\rs, \Pi, \chi}$, we have $\sg f' \in \sH(\G'(\A^{\infty}) , L')^{\circ}_{\rs, \Pi, \sg \chi}$
and
$\sg\sL(\RM_{\Pi}, \cdot)_{f'} = \sL(\RM_{\Pi}, \cdot)_{\sg f'}$.
Thus the functions $\sL(\RM_{\Pi}, \cdot)_{f'}$, for  $ f'\in {\bigcup_{\chi\in Y_L} \sH(\G'(\A^{\infty}) , L(\chi))^{\circ}_{\rs, \Pi, \chi}}$, descend to the desired $\sL(\RM_{\Pi}, \cdot)$. 
\end{proof}

\subsubsection{Spectral expansion} We define 
$$I_{\Pi}\coloneqq  {1\over 4} \sL(\RM_{\Pi}) \cdot \prod_{v} I_{\Pi_{v}}$$
Then the spectral expansion of part \eqref{rat rtf part} of Proposition \ref{rat rtf} follows from the definition, Proposition \ref{exps rtf}, and \eqref{I rat kappa}. {It is locally finite since for  $f' \in \sH(\G'(\A), L)^{\circ}_{\reg^+, \rm qc}$, for any compact open subgroup $K\subset \G'(\A^\infty)$ such that $f'=f'\star e_K$,  we  have $I_\Pi(f')=0$ unless $\Pi$ is  a point of the finite scheme $\sC(\G')_K$.}  This completes the proof of the proposition.

\subsection{On the Ichino--Ikeda conjecture} \lb{sec: II}
 For expository purposes, we recall an outline of the proof of the following special case of the  Ichino--Ikeda--Harris conjecture  (in its most general form, the conjecture is now   \cite[Theorem 1.1.6.1]{BPCZ}),  paying special attention to the rationality.  The basic architecture of the proof of Theorem \ref{main thm}  in \S~\ref{sec:11} will be similar.

 Let $V\in \sV^{\circ, +}$ be a coherent pair, let $\RH=\RH^V\subset \G=\G^{V}$, and let $\sA(\G)^{\circ}\coloneqq \Q[\G(F_{0}\backslash \G(\A)/ \G(F_{0, \infty})]$, which is equipped with the Petersson product with respect to the measure $dg$.    Let $\pi$ be a cuspidal automorphic representation of $\G^{V}(\A)$, trivial at infinity,  over a number field $L$.   Upon choosing an embedding in $\Hom(\pi, \sA(\G)^{\circ}_{L})$ (which is an $L$-line by \cite{KMSW}, see \cite[Proposition C.3.1 (2)]{LTXZZ}) we have an $\RH$-period 
 \beq\lb{H per pi}
 P_{\pi}\colon \pi\to L
 \eeq 
  defined as in \eqref{H per} (where the integration reduces to a finite sum). The unique embedding $\pi^{\vee}\into \sA(\G)^{\circ}_{L}$ that intertwines the natural duality $\pi \times\pi^{\vee}\to L$ with the Petersson product gives rise to the analogous period $P_{\pi^{\vee}}\colon \pi^{\vee}\to L$.

 \begin{theorem} \lb{thm: ii} 
 Assume that $\pi$ is 
  stable and cuspidal, and let $\Pi\coloneqq {\rm BC}(\pi)$.
   Then  for all $\phi\in \pi$, $\phi'\in \pi^{\vee},$ we have
 $$ 
P_{\pi}(\phi)P_{\pi^{\vee}}(\phi')
=  {1\over 4}  \sL  (\RM_{\Pi}, 0) 
\cdot \alpha(\phi, \phi').$$
 \end{theorem}
 
 We need a lemma to isolate $\pi$ within the discrete automorphic spectrum.
 
 \begin{lemma} \lb{isolate pi} Let $L$ be a coefficient field and let $V\in \sV^{\circ}$. 
 Let $\Sigma$ be a finite set of  isomorphism classes of discrete irreducible automorphic representations of $\G^{V}(\A)$ over $\ol{L}$, trivial at infinity.
 Let $\ol{\pi}\in \Sigma$ be cuspidal and tempered at every finite place, and assume that $\ol{\pi}=\pi\ot_{L}\ol{L}$ for some representation $\pi$ over $L$.  Let $P$ be a finite set of places of $F_{0}$.

 Then there is a finite set $S$ of split places of $F_{0}$, disjoint from $P$, a relative selfdual hyperspecial subgroup $K_{S}\subset G^{V}_{S}$ (\S~\ref{sss rel sd}),  and an element  $f_{S}\in \sH(G_{S}^{V}, L)_{K_{S}}$, such that
\beqq
\begin{cases}
\pi(f_{S})={\pi}(e_{K_S}),& \\
  \ol\pi'(f_{S})=0 &\text{ for all }\ol\pi'\in \Sigma \text{\ with ${\rm BC}(\ol{\pi}')^\iota\neq {\rm BC}(\ol{\pi})^\iota$ for some $\iota \colon {\ol L}\into \C$}.\end{cases} \eeqq
 \end{lemma}

 \begin{proof} By induction, we may assume $\Sigma=\{\ol{\pi}, \ol{\pi}'\}$ contains exactly two elements, and view ${\ol L}\into \C$ via an embedding such that   $\Pi\coloneqq{\rm BC}(\ol{\pi})^\iota\neq  \Pi'\coloneqq {\rm BC}(\ol{\pi}')^\iota$. 
 
{Since base-change  preserves temperedness, the representation $\Pi$ is isobaric, hence tempered at every finite place by Proposition \ref{ramanujan}. We consider the alternative of that proposition for $\Pi'$.

If $\Pi'$ is not isobaric, then it is not tempered at any finite place. 
Thus it is trivial to find a finite set (a singleton suffices) of split places $S$ disjoint from $P$, and an element $f'_{S}\in \sH(\G'_{S}, L)_{K'_{S}}$ (for $K'_{S}\subset G'_S$ maximal hyperspecial), satisfying $\Pi(f'_{S})=\Pi(e_{K'^S})$ and $\Pi'(f'_{S})=0$.
    (We can take $f'_S$ with coefficients in $L$ since $\Pi$ is defined over $L$.) Then the Hecke measure $f_{S}\in  \sH(G_{S}^{V}, L)$ matching $f'_{S}$ satisfies the desiderata. 

If $\Pi'$ is isobaric, then 
    the same picture as in the previous case holds by Proposition \ref{Ram-prop}.}
\end{proof}

{
\begin{remark}[Matching over arbitrary coefficient fields.]\lb{rat matching} 
If $v$ is a non-archimedean place, since regular-semisimple orbital integrals reduce to finite sums, we may extend to arbitrary coefficient fields the notion of (geometric) matching of local-at-$v$ Hecke measures from \S~\ref{sec:35}. If $v=\infty$, then we trivially have a notion of pure matching when restricting to Gaussians on $G'_\infty$, and to multiples of the unit measure on definite groups on the unitary side. For $\pi_v\in {\rm Temp}(G^{V_v})(L)$ (for any coefficient field $L$) and $f_v\in \sH(G^{V_v} , L)$, we may then define $J_{\pi_v}(f_v)\coloneqq I_{{\rm BC}(\pi_v)}(f'_v)$ for any matching $f'_v$. This makes sense by the equivalence of geometric and spectral matching (Proposition \ref{equiv match}).
\end{remark}
}

 \begin{proof}[Proof of Theorem \ref{thm: ii}] (For more details on the argument, see the proof of Theorem \ref{main thm}  in \S~\ref{sec:11}.)
The formula extends by bilinearity to any $\tau\in \pi\ot \pi^{\vee}$, and by multiplicity one if suffices to prove it for any $\tau$ not annihilated by $\alpha$.

 By Corollary \ref{cor test} (with $\chi=\one$ and  $P=\emptyset$) and Lemma \ref{isolate pi} (with $P$ a set of places such that the Gaussian produced by Corollary \ref{cor test} is spherical away from $P\infty$), together with the explicit matching at split places of Lemma \ref{match split}, we may construct purely matching Gaussians $f'\in \sH(\G'(\A), L)^{\circ}$ and $f\in  \sH(\G(\A), L)^{\circ}$ with globally regular semisimple  support that are adapted to $\Pi=\mathrm{BC}(\pi)$ and, respectively,~$\pi$. Then for all $v$ and matching $\gamma\in \RB_{\rm rs}'(F_{0,v})$,  $\delta\in \RB_{\rm rs}(F_{0,v})$, 
$$I_{\gamma}^{}(f'_{v}) = J_{\delta}^{}(f_{v})$$ 
so that by  \eqref{prod orb int}, and \eqref{J delta v},
 $$J_{\pi}(f) = J(f)= \sum_{\delta\in \RB^V_{\rm rs}(F_{0})}  J_{\delta}(f) 
=\sum_{\gamma\in \RB'_{\rm rs}(F_{0})^{\circ}}  I_{\gamma}(f')  =
 I(f')=
 I_{\Pi}(f'),$$
where $ \RB'_{\rm rs}(F_{0})^{\circ} \coloneqq  \RB'_{\rm rs}(F_{0}) \cap  B_{\rs, \infty}'^{\circ}$, and we have used the RTFs  for  $I$  in Proposition \ref{rat rtf} and for $J$ in \eqref{J rtf}, together with the automorphic multiplicity-one result in \cite[Theorem 1.7.1]{KMSW}.
By the factorization of $I_{\Pi}$ in Proposition \ref{factor I} and the local spectral matching  (together with the fact that $\prod_{v}\kappa(\pi^{V}_{v})=1$), we have
$$J_{\pi}(f)=I_{\Pi}(f') 
= {1\over 4\Delta_{\RH}} \sL^{}(1/2,{\Pi}) \cdot
 \otimes_{v} J^{\C}_{\pi_{v}}(f)
 =
  {1\over 4} \sL^{}(\RM_{\Pi}, \one) 
\cdot\otimes_{v\nmid \infty} J^{}_{\pi_{v}} J^{\circ}_{\pi_{\infty}}(f),
 $$
 where we have used the definitions of $J^{\circ}_{\pi_{\infty}}$ in \eqref{def J circ} and of $\sL^{}(\RM_{\Pi})$ in Theorem \ref{rat L}. This is equivalent to the desired formula for $\tau=\pi(f)$.
 \end{proof}

\section{$p$-adic relative characters}  \lb{sec: loc rtf} \lb{sec 5}
This section and the next one contain the local results needed, at $p$-adic places,  in order to develop the $p$-adic relative-trace formula; in particular, the construction of a suitable family of  Hecke measures. Remarkably, suitable members of these families can be used at any (split) place as the regular local test measures needed to prove the results of \S~\ref{sec: test}.

Throughout this section, we fix a non-archimedean  place $v $ of $F_{0}$ and work in a local situation, dropping all subscripts $v$.  We denote by $\sO$ the ring of integers of the \'etale $F_{0}$-algebra  $F$, by $\sO_{{0}}$ the ring of integers of $F_{0}$, by $\vpi\in \sO_{0}$ a chosen uniformizer, and we let $q_{0}\coloneqq  |\sO_{0}/\vpi\sO_{0}|$, $q\coloneqq  |\sO/\vpi\sO|$.

After some preliminaries in  \S~\ref{sec:51}, in  \S~\ref{sec:52} we define and study  $p$-adic modifications of the local periods of \S~\ref{sec:31}; in  \S~\ref{sec:53} we turn to  the associated relative characters. The main result on $p$-adic local periods is  Liu--Sun's evaluation in terms of explicit Coates--Perrin-Riou local factors (stated as Proposition~\ref{hyp ep}); the resulting non-vanishing properties  allow for the application to test measures alluded to above,   in Proposition~\ref{nv iw ps}.

\subsection{Group-theoretic preliminaries} \lb{sec not G0}  \lb{sec:51}
We introduce some notation and the group-theoretic foundations for  the construction of the $p$-adic distribution. 
\subsubsection{Notation}  
We expand our list of groups over $F_0$ to include 
$${G}'_{0}\coloneqq  G_{n,0}'\ts G_{n+1,0}', \qquad H_{1, 0}'\coloneqq  G'_{n, 0},$$
 so that if $F\cong F_{0}\ts F_{0}$, we have $\wt{G}'={G}'_{0}\ts {G}'_{0}$ and
$H'_{1}=H'_{1, 0}\ts H'_{1, 0}$. 
In this case, we  may  write elements of $G' = \wt{G}'/(F_{0}^{\ts})^{2}$ as $[g_{1}; g_{2}]$ with $g_{i}\in {G}'_{0}$.  
  
We will  denote all conjugation actions  by
$$x^{g}\coloneqq g^{-1}xg.$$
  
\subsubsection*{Convention}  Throughout this section, for $\nu\in \{n, n+1, \emptyset\}$ and $*\in \{\emptyset, 0\}$ we will define various subgroups and elements $\Box_{\nu,*}$ of $G_{\nu, *}'$. 
 Unless otherwise specified, \emph{we will define $\Box_{\nu,*}$ in a way that makes sense for $\nu = n, n+1$, and tacitly stipulate that $\Box_{*}$ is the product of $\Box_{n, *}$ and $\Box_{n+1, *}$}, if $*=0$, or its image via $\widetilde{G}'\to G'$ if $*=\emptyset$.

\subsubsection{Some subgroups}
The lattice $\sO_{*}^{\nu}\subset F_{*}^{n}$ induces an integral model for $\G'_{\nu, *}$ over $\sO_{0}$, still denoted by $\G'_{\nu, *}$. Let $T_{\nu , *}\subset G_{\nu, *}'$ denote the diagonal torus, and let $W_{\nu, *}$ be the associated Weyl group, identified with the permutation matrices in $G_{\nu, *}'$. We denote by  $$w_{\nu, *}\in W_{\nu, *}$$ the antidiagonal matrix $(w_{\nu, *})_{ij}=\delta_{i,\nu+1-j}$.

\subsubsection{On the torus in $G'_{\nu, *}$}
We denote by  $N_{\nu, *}\subset G'_{\nu, *}$ the set of upper-triangular unipotent matrices and by
$$
N_{\nu, *}^{\circ}\coloneqq  N_{\nu, *}\cap   \RG_{\nu, *}'(\sO_{0}).
$$
Let $T_{\nu, *}^{+}\subset T_{\nu, *}$ be the sub-monoid consisting  of those $t$ such that $N_{\nu, *}^{\circ, t} \coloneqq  (N_{\nu, *}^{\circ})^{t}   \subset N_{\nu, *}^{\circ}$, and  $T_{\nu, *}^{++}\subset T_{\nu, *}^{+}$ the multiplicative subset ot those $t$ such that 
$$
\bigcap_{r\geq 1} N_{\nu, *}^{\circ , t^{r}}=\{1\}.
$$
Concretely, ${T}_{\nu, *}^{+}$ (respectively ${T}_{\nu, *}^{++}$)  consists of matrices $\diag(t_{1} , \ldots, t_{\nu})$ with $t_{i}\in F_{*}^{\ts}$ and $v(t_{i}/t_{i+1} )\geq 0$  (respectively $>0$) for all $1\leq i\leq \nu-1$.

The group $T_{\nu, *}$ is equipped with
the involution
$$\iota \colon t\mapsto w_{\nu, *}^{-1}t^{ -1}w_{\nu, *}, $$ which preserves ${T}_{\nu, *}^{+}$ and ${T}_{\nu,*}^{++}$. We still denote by $\iota $ the resulting involution on $\bQ_{p}[T_{\nu, *}]$.

 We identify $\Z^{\nu}$ with the space of cocharacters of $T_{\nu, *}$ via $$\lm\mapsto[ x\mapsto x^{\lm}\coloneqq \diag (x^{\lm_{1}}, \ldots , x^{\lm_{\nu}})]\in T_{\nu, 0}\subset T_{\nu, *},$$
  where the inclusion is diagonal. 

We fix the elements
 \beq\lb{def t+}
{t}_{\nu, *}\coloneqq   \vpi^{(\nu-1, \ldots, 0)}\in {T}_{\nu, *}^{++}, 
\qquad     z_{\nu,*}= \vpi^{\nu-1}1_{\nu}\in G'_{ \nu, *}.
\eeq
Then
   $$t_{\nu, *}^{\iota}= z_{\nu, *}^{-1} t_ {\nu,*}, \qquad  t_{\nu, *} t_{\nu, *}^{\iota}=\vpi^{2 \rho_{\nu}}$$ 
   where $\rho_{\nu}\in \Z^{\nu}$ denotes half the sum of positive roots (with respect to $N_{\nu, *}$); 
   concretely,
   $${\rho}_{\nu}\coloneqq  {1\over 2}(\nu-1, \nu-3, \ldots, 1-\nu) \in {1\over 2}\Z^{\nu}.$$  

\subsubsection{Iwahori and deeper Iwahori subgroups}  \lb{iw conjsy}
 The standard \emph{Iwahori subgroup}
 $${\rm Iw}_{\nu, *}\subset G'_{\nu, *}$$ 
 is the set of matrices in $\G'_{\nu,*}(\sO_{0})$ whose reduction modulo $\vpi$ belongs to  the image of the upper-triangular matrices in  $\G'_{\nu,*}(\sO_{0})$. An {Iwahori subgroup}  of $ G'_{\nu, *}$ is one of the form ${\rm Iw}_{\nu, *}^{g}$ for some $ g\in G'_{\nu, *}$.
 A \emph{deeper Iwahori subgroup}  of $G'_{\nu, *}$ 
  is an open  subgroup $K\subset G'_{\nu, *}$  satisfying $ K\subset {\rm Iw}_{\nu, *}^{g}$ for some $g\in  G'_{\nu, *}$. It is  said to be \emph{semistandard} if   $N^{\circ}_{\nu, *}\subset K\subset {\rm Iw}_{\nu, *}^{g}$ for some
  $g\in N_{G'_{\nu, *}}(T_{\nu,*})$, the normalizer of  $T_{\nu,*}$ in ${G'_{\nu, *}}$; it is said to be \emph{standard} if $ K\subset {\rm Iw}_{\nu, *}$ and $K\cap N_{\nu, *} =N^{\circ}_{\nu, *}$.

For $r\in \bZ-\{0\}$, we define three families of subgroups 
\beq
 \lb{Krdef} K_{\nu, *}^{[r]}\subset K_{\nu, *}^{\la r\ra} \subset K_{\nu, *}^{(r)}
\eeq
of $\RG'_{\nu,*}(\sO_{{0}})$
 by 
\beqq
K_{\nu, *}^{(r)} & \coloneqq  \RG'_{\nu,*}(\sO_{{0}})\cap  t_{\nu, *}^{-r} \RG'_{\nu,*}(\sO_{{0}})t_{\nu, *}^{r},\\
K_{\nu,*}^{\la r\ra} & \coloneqq  \{ g\in K_{\nu, *}^{(r)} \ |\  g_{ii}\in 
 1+\vpi^{|r|-1} \sO_{{*}}, \ 1\leq i\leq \nu\} \\
K_{\nu,*}^{[r]} &\coloneqq  \{ g\in K_{\nu, *}^{(r)} \  |\  g_{ii}\in 
 1+\vpi^{|r|} \sO_{{*}}, \ 1\leq i\leq \nu\}.
\eeqq
They are standard deeper Iwahori  subgroups whenever $r\geq1$.

For $r\geq 1$, we say that a standard deeper Iwahori   $K_{\nu, *}$ has  \emph{level $\leq {r}$} if 
$$
 K_{\nu, *} \supset K_{\nu, *}^{\la r\ra}.
 $$

 \subsubsection{Iwahori--Weyl symmetries} \lb{IW sym}
 For $c\geq 1$, define
  $$
  w_{ \nu,*,  c}\coloneqq  w_{\nu, *}t_{\nu, *}^{c}\in  N_{G'_{\nu, *}}(T_{\nu,*})\subset G_{\nu,*}'.
  $$
Let  $K\subset G'_{\nu, *}$ be a semistandard deeper Iwahori subgroup. We say that $K$ is \emph{symmetric} if $K^{w_{\nu, *, c}}=K$ for some $c\geq 1$ such that $K^{\la c\ra}_{\nu, *}\subset K$.
 If $F=F_{0}\ts F_{0}$ and $*=\emptyset$, we say that $K$ is \emph{conjugate-symmetric
of  {depth} $c=c(K)\geq1$}
 if $K=K_{0}\ts K_{0}^{w_{\nu, 0,c}}$ for  some standard
  deeper Iwahori subgroup $K_{0} \subset G'_{\nu, 0}$ containing $K_{\nu, 0}^{\la c\ra}$.

\begin{remark}
 For $r\geq1$, the subgroups  $ K_{\nu, *}^{[r]}\subset K_{\nu, *}^{\la r\ra} \subset K_{\nu, *}^{(r)} \subset G'_{\nu, *}$ are all symmetric,
 whereas for $\nu\geq 3$ Iwahori subgroups are not symmetric. On the other hand, conjugate-symmetric deeper Iwahori subgroups of $G_{\nu}'$ are obviously  abundant.
\end{remark}

\subsubsection{Iwahori--Hecke algebras and the operators $U_{t}$} \lb{H dag def}
Let $K\subset G_{\nu, *}'$ be  a semistandard deeper  Iwahori subgroup. For a ring $R$, define
$$\sH^{\dag, +}_{K, *}(R)\coloneqq  C^{\infty}_{c}(K\bs KT_{\nu, *}^{+}K / K, R)\, dg \quad
 \subset \quad 
\sH_{K, *}\coloneqq C^{\infty}_{c}(K_{}\bs G_{\nu, *}'/ K_{}, R)\, dg.$$
The involution $\iota$ extends to $\sH^{\dag, +}_{K, *}(R)$ by linearity. For $x\in G_{\nu,*}'$, we define 
{\beqq
[K xK] &\coloneqq  \vol(K, dg)^{-1} \one_{Kx K} \,dg
\eeqq
in $\sH_{K, *}(\Q)$. The map
\beq \lb{IH iso}
R[T_{*}^{+}/T_{*}\cap K] & \to \sH_{K, *}^{\dag, +}(R)\\
[t]&\mapsto U_{t , K}\coloneqq [K tK].
\eeq 
is an $R$-algebra isomorphism. We define
$$\sH^{\dag}_{K, *}(R)\coloneqq  \sH_{K, *}^{\dag,+}[(U_{t,K}^{-1})_{t\in T^{+}} ]\cong R[T_{*}/T_{*}\cap K].$$

For $?=+, \emptyset$ and $R$ a $\Q$-algebra,  we define $\sH^{\dag,?}_{\nu, *}(R) \coloneqq \varprojlim_{K} \sH_{K, *}^{\dag,?}(R)$, where the limit runs over the standard deeper Iwahori subgroups and the transition maps  are $\star e_{K}\colon \sH_{K', *}^{\dag, +}\to  \sH_{K, *}^{\dag, +}$. 
By Lemma \ref{iw mult} below, the limit  $U_{t}\coloneqq \lim U_{t, K} \in \sH^{\dag,?}_{\nu, *} (R)$ is well-defined. Concretely, if we denote 
$$N_{\nu*}^{\circ,(r)}\coloneqq  t_{\nu, *}^{r} N_{\nu, *}^{\circ}t_{\nu, *}^{-r}$$
 we have 
$$U_{t_{\nu, *}} = \sum_{ x\in N_{\nu, *}^{\circ}/ N_{\nu, *}^{\circ , (1)}}  xt_{\nu, *}$$
as operators on the $N^{\circ }_{\nu,*}$-fixed points of any smooth $G'_{\nu, *}$-module.

\subsubsection{Multiplication rules  in  Iwahori--Hecke algebras} We have the following basic result.
\begin{lemma} \lb{iw mult}  Let $K\subset G_{\nu, *}'$ be a deeper
 Iwahori subgroup, and define $\ell_{K}\colon K\bs G'_{\nu, *} / K \to \bN$  by $q_{*}^{\ell_{K}(g)}\coloneqq  |K_{}gK_{}/K |=  | K_{} / K_{} \cap gK_{}g^{-1}|$. Then:
\begin{enumerate}
\item 
We have $K_{}gK_{}g'K_{} =K_{}gg'K_{}$ 
if and only if $\ell_{K}(gg')=\ell_{K}(g)+\ell_{K}(g')$.
\item Assume that $K$ is standard. Then  for all  $t'\in T_{\nu, *}^{+}$,
 $$\ell_{K}(t'w_{\nu, *})=\ell_{K}(t')+ \ell_{K}(w_{\nu, *}), \qquad \ell_{K}(w_{\nu, *}t'^{-1})=\ell_{K}(w_{\nu, *})+\ell_{K}(t'^{-1}).$$
 \item Assume that $K$ is standard, and let $K'\subset K$ be a standard deeper Iwahori subgroup. Then for all $t'\in T_{\nu, *}^{+}$,
 $$ Kt'wK'=Kt'wK, \quad e_{K}\star[K't'wK']=[Kt'wK].$$
If moreover $K$ is of level $\leq c$ and $t't_{\nu, *}^{-c}\in T_{\nu,*}^{+}$, then 
\beqq
K't'K=Kt'K, &\quad [K't'K']\star e_{K}=[Kt'K]. 
 \eeqq

\item \lb{iw mult 4} For all $g\in G'_{\nu, *}$, we have 
$$e_{K}\star ge_{K} = q_{*}^{-\ell_{K} (g)} [KgK].$$
\end{enumerate}
\end{lemma}
\begin{proof} Part (1) is well-known, see \cite[Corollary 1.1]{howe}. Consider the first equality of part (2), and drop all subscripts. By part (1), it is equivalent to proving $Kt'KwK=Kt'wK$. Since the quotient $K\backslash Kt'K$ is represented by lower-triangular matrices in $K$, it suffices to show that for such a matrix $k$, we have $t'kw\in Kt'wK$; in fact, since $K\supset N^{\circ}$ we even have $kw\in wK$. The second equation follows from taking inverses in  the identity $Kt'KwK=Kt'wK$. 

Consider now  part (3); we only prove the equalities as sets, from which the ones in Hecke algebras can be easily obtained.  For the first equality, it suffices to prove that for any lower-triangular $k\in K$ we have $t'wk \in Kt'w$, which is clear since $t'wkw^{-1}t'^{-1}\in N^{\circ}\subset K$.  For the second one, it suffices to prove that for any lower-triangular $kt'\in Kt'$ we have $kt'\in t'K$. In fact, by the assumptions we have $t'^{-1}tk\in K^{(c)}\cap K\subset K$.  
Part~\eqref{iw mult 4} follows from the definitions.
\end{proof}

\subsubsection{Twisting matrices}
Let $u \in  (\sO_{F}^{\ts})^{n}$;  we will take $u=(1, \ldots, 1)^{\rm t}$ to fix ideas in computations. Then we define the \emph{twisting matrices}\footnote{For their history, see \cite{janu-nonab} and references therein.}
\beq\lb{m def}
{m}_{n, *}\coloneqq 1_{n}, \qquad {m}_{n+1, *}\coloneqq \twomat  {w_{n}} u {} 1 , 
\eeq
 and for $r\geq 1$  we let
\beqq
 & m_{\nu,*, r}\coloneqq  m_{\nu, *}  t^r_{\nu, *}
\eeqq

\subsubsection{Subgroups of $H_{1}'$}  \lb{518} Recall that by the convention introduced at the beginning ot this subsection, $\Box_{*}$ denotes the (image of the) product of $\Box_{n, *}$ and $\Box_{n+1, *}$  in ${G}_{*}'$. 
For $r\in \bZ_{> 0}$,
let
\beq\lb{KH def}
K_{H, *}^{(r)} &\coloneqq m_{*} K_{*}^{(-r)} m_{*}^{-1}\cap \RG_{n, *}'(\sO_{0})\\
&=m_{*} K_{*}^{[-r]} m_{*}^{-1}\cap \RG_{n, *}'(\sO_{0})
\qquad \subset \RG_{n , *}'
(\sO_{0}) \subset H'_{1, *},
\eeq
where  the intersections are with respect to the usual diagonal embedding $H'_{1, *}\into  G_{*}'$.

\begin{remark}  \lb{vol0}
A simple computation shows that   ${K}_{H, *}^{(r)} $ consists of the matrices $h$ satisfying 
\beq \lb{KHdet}
\begin{cases}
h_{ij} \in \vpi^{r|i-j|}\sO_{*}\\
\sum_{j=1}^{n} h_{ij} \in 1+ \vpi^{ir}\sO_{*}
\end{cases}
\eeq
 for all $1\leq i, j\leq n$. This description also shows the equality in \eqref{KH def}. We may then compute  that  
\beq\lb{vol circ} \vol^{\circ}(K_{H, *}) \coloneqq  q_{*}^{d(n)s}
 \vol(K_{H, *}^{(s)}) = \prod_{ i=1}^{n} {1-q_{*}^{-1}\over 1-q_{*}^{-i}}
  \eeq
 is a rational number independent of $s\geq 1$, and a $p$-unit. 
\end{remark}

We record the following easily checked property, for a later use: for all $r\geq 1$, we have
\beq\lb{mr KH}
m_{*, r}^{-1} K_{H, *}^{(r)} m_{*, r}\subset K_{*}^{(2r)}\cap K_{*}^{[r]}\subset K_{*}^{\la r+1\ra}.\eeq

\subsubsection{Twisting identity} We come to the key result of this subsection, which refines  \cite[Lemma 5.2]{janu-nonab} in the spirit of \cite[Lemma 4.4.1]{loeffler}.

\begin{lemma}[Twisting identity] \lb{mat id} Let $r\geq 1$ and let $K\subset {G}'_{*}$ be a subgroup containing $K_{*}^{\la r+1\ra}$.
 For all $x \in N_{*}^{\circ}$, there exists $h_{x}\in K_{H,* }^{(r)}$ such that 
\beq \lb{eq mat id}
m_{*,r} x t K
= h_{x} m_{*,r+1} K \eeq
Moreover, the map
\beqq
N_{*}^{\circ, (1)} \bs N_{*}^{\circ}
 &\to K_{H,*}^{(r+1)} \bs K_{H,*}^{(r)} \\
[x] & \mapsto [ h_{x}]
\eeqq 
is well-defined and a group isomorphism.
\end{lemma}
\begin{proof} 
We omit the subscript `$*$' from the notation. It suffices to take $K=K^{\la r+1\ra}$.
Consider the diagram
$$ K_{H}^{(r+1)} \bs K_{H}^{(r)} \stackrel{\alpha}{\to}   K^{\la -r-1\ra }\bs  K^{[-r]} \stackrel{\beta}{\longleftarrow} N^{\circ, (r+1)} \bs N^{\circ, (r)}
\stackrel{\gamma}{\longleftarrow} N^{\circ, (1)} \bs N^{\circ}$$
where $\alpha\colon h\mapsto m^{-1} h m$,  $\beta$ is  induced by the inclusion $N^{\circ, (r)}\subset K^{[-r]}$, and $\gamma$ is the isomorphism $ x\mapsto t^{r}x t^{-r}$. 
All four  quotients  have cardinality $q^{d(n)}$ where $d(n)=\eqref{dn}$, and by \eqref{KH def}, $\alpha$ is well-defined and injective. 
Hence all three maps are isomorphisms, and the second statement of the lemma is proved with $[h_{x}]=\alpha^{-1}\circ\beta\circ \gamma ([x])$. The first statement is then easily verified using  $t^{-r-1}K^{\la  -r-1\ra} t^{r+1} = K^{\la r+1\ra}$.
\end{proof}

\begin{corollary} \lb{twisting cor}  Let $r\geq 1$,  and let $K_{*}= K_{*}^{\langle r+1\rangle}\subset {G}_{*}'$. For all $s\geq r$, we have the identities
\beq\lb{mat id integral}
m_{*,s} U_{t_{*}, K_{*}} 
&= \sum_{h\in K_{H,*}^{(s+1)} \bs K_{H,*}^{(s)}}  hm_{*,s+1} e_{K_{*}} ,\\
q_{*}^{sd(n)}\cdot m_{*,s}U_{t_{*}, K_{*}}^{-s}
&= q_{*}^{(s+1)d(n)} \cdot \sum^{\rm avg}_{h\in K_{H,*}^{(s+1)} \bs K_{H,*}^{(s)}} hm_{*,s+1} U_{t_{*}, K_{*}}^{-(s+1)},
  \eeq
where $ \sum^{\rm avg}$ denotes average. They are valid in 
 $$C_{c}( G_{*}'/K_{*}, \Q), \qquad C_{c}( G_{*}'/K_{*}, \Q)\ot_{\sH_{K_{*}}^{\dag, +}(\Q)} \sH_{K_{*}}^{\dag}(\Q)$$
respectively.
\end{corollary}

\subsubsection{Volumes} The volumes of $K_{\nu, *}^{(r)}$ and $K_{H, *}^{(r)}$ are constant multiples of $q_{*}^{-c(\nu)r}$, respectively  $q_{*}^{-d(n)r}$, where
\beq \lb{dn}
c(\nu)&\coloneqq  {1\over 6} (\nu-1)\nu(\nu+1),\\
d(n)&\coloneqq  \sum_{k=1}^{n } k^{2} = {{1\over 6} n(n+1)(2n+1)} =c(n)+c(n+1).\eeq

\subsection{$p$-adic  periods} \lb{sec:52}
 Let $\Pi$ be a tempered irreducible admissible 
 representation of $G'$ over a field $L$ of characteristic zero. Denote by $B\subset G'$ the upper-triangular Borel and by $\delta_{B}\colon T\to \Q^{\ts}$ its modulus character.

 \subsubsection{Finite-slope subspace}\lb{fss}
Let
$$\Pi^{\dag}\subset \Pi^{N^{\circ}},$$ 
be the subspace where $T^{+}$ acts invertibly. It has a structure of $\sH^{\dag}(L)$-module, and it is isomorphic as an $L[T]$-module to the normalized Jacquet module $\delta_{B}\otimes \Pi_{N}$ of $\Pi$ (see e.g. \cite[Proposition 4.3.4]{emerton}). 

We define $c(\Pi^{\dag})$ to be the minimal $c\in \Z_{\geq 1}$ such that $\Pi^{\dag}\subset \Pi^{K^{\langle c\rangle}}$.

 Denote by $\widehat{T}$ {the $L$-ind-group-scheme of smooth characters of $T$}. For a subgroup $K\subset G'$
 containing $N^{\circ}$,  let $\Pi^{K , \dag} \subset \Pi^{\dag}$  be the image of $\Pi^{K}$ under the $U_{t}$-eigen-projection from $\Pi^{N^{\circ}} $ to $ \Pi^{\dag}$, for any sufficiently positive $t$. Then there are decompositions into generalized $\sH^{\dag}_{K}$-eigenspaces  
 $$\Pi^{K, \dag}_{\ol{L}}=\bigoplus_{\xi\in \widehat{T}_{L}(\ol{L})} \Pi_{\ol{L}}^{K, \dag}[\xi],$$
and similarly $\Pi^{\dag}_{\ol{L}}=\bigoplus \Pi^{\dag}_{\ol{L}}[\xi]$. 

If $\Pi$ is a subquotient of a regular principal series (as defined in \S~\ref{def reg}) and  $\xi \in \widehat{T}$ is a character occurring in $\Pi^{\dag}_{L(\xi)}$, then  by \cite[Proposition 1.3 (ii)]{janu-nonab} (or its proof, applied to $\Pi_{n}$, $\Pi_{n+1}$), any Whittaker model  of $\Pi_{L(\xi)}$ contains  a unique  vector 
$W_{\xi}$ satisfying $W_{\xi}(1)=1$ and $U_{t}W=\xi(t) W$ for all $t\in T^{+}$.

\subsubsection{Ordinary representations}
Suppose for this paragraph only that $L$ is a finite extension of $\Q_{p}$, and denote by   $\ol{\Q}_{p}$ an algebraic closure of $L$. 
\begin{definition}\lb{ord v} Let $N^{\circ}\subset K\subset  G'$. We say that  the tempered representation $\Pi$ is \emph{$K$-ordinary} (with respect to $\Pi_{\infty}^{\circ}$) if there is a character $\xi^{\circ}\in \widehat{T}(\ol{\Q}_{p})$ occurring in $\Pi^{K, \dag}_{\ol{\Q}_{p}}$ (that is, such that  $\Pi^{K,\dag}_{\ol{\Q}_{p}}[\xi^{\circ}]\neq 0$)
satisfying  $$|\xi^{\circ}(t')|=1$$ for all $t'\in T^{+}$ and the absolute value on $\ol{\Q}_{p}$.\footnote{This definition is adapted to the local components of automorphic representations of trivial weight at infinity; in general it would need to be modified, see \cite{hida-aut}.} 
We say that  $\Pi$ is ordinary if it is $K$-ordinary for sufficiently small $K\supset N^{\circ}$.
\end{definition}

We call a character $\xi^{\circ}$ as above an \emph{ordinary refinement} of $\Pi$. By the following proposition, an ordinary refinement is unique and defined over the field of definition of $\Pi$. We will then denote
$$\Pi^{\rm ord}\coloneqq \Pi_{\ol{\Q}_{p}}[\xi^{\circ}]\cap \Pi.$$

\begin{proposition}\lb{hida m1} Let $\Pi$ be an ordinary tempered representation of $G'$ over $L$. Then
$\Pi$ is a subquotient of a regular principal series, the space
 $\Pi^{\dag}$ is $T_{L}$-semisimple, and  every $\xi\in \widehat{T}_{{L}}(\ol{\Q}_{p})$ occurring in $\Pi_{\ol{\Q}_{p}}^{\dag}$ satisfies
 $\dim_{\ol{\Q}_{p}} \Pi^{\dag}_{\ol{\Q}_{p}} [\xi]= 1$ and    is defined over $L$. Moreover the ordinary refinement $\xi^{\circ}$ is unique.  
\end{proposition}

\begin{proof} 
 This is essentially \cite[Corollary 8.3]{hida-aut}. We recall the argument, working over $\ol{\Q}_{p}$ without signalling this in the notation. Let $W$ be the Weyl group of $G'$. Recall from \S~\ref{fss} that $\Pi^{\dag}\cong \delta_{B}\ot \Pi_{N}$, the normalized Jacquet module of $\Pi$.   By Frobenius reciprocity,  $\xi$ occurs in $\Pi^{\dag}$  if and only if $\Pi$ embeds into the normalized induction ${\rm Ind}_{B}^{G'}(\wt\xi)$
where $\wt{\xi}\coloneqq \delta_{B}^{-1/2}\xi$.
Now ${\rm Ind}_{B}^{G'}  (\wt{\xi})\cong{\rm Ind}_{B}^{G'}(\wt{\xi}^{w})$ for all $w\in W$.
 If $\xi^{\circ}_{|T^{+}}$ is valued in units, then the stabilizer of $\wt{\xi}^{\circ}$ in $W$ is trivial, therefore its orbit consists of $|W|$ distinct characters $\wt{\xi}$, and ${\rm Ind}_{B}^{G'}  (\wt{\xi})$ is regular. By \cite[Theorem 5.21]{BZ}, we have $\dim \Pi_{N}\leq |W|$, hence all the characters $\xi$ occur with multiplicity one. The rationality assertion follows from the fact that the $\Gal(\ol{\Q}_{p}/L)$-action on the set of occurring $\xi$ preserves valuations. 
\end{proof}

Denote by 
$${\rm e}^{\rm ord}
\colon \Pi^{N^{\circ}}\to \Pi^{\ord}$$
the $\sH^{\dag}(L)$-eigenprojector,
and let $e^{\rm ord}_{K}\coloneqq  e^{\rm ord}e_{K}$.  Thus 
$\Pi$ is $K$-ordinary if $e_{K}^{\ord} \Pi =\Pi^{\ord}$.

\begin{lemma}\lb{nv ord unr} Suppose that $\Pi$ is ordinary and unramified, and let $K\coloneqq  \RG'(\sO_{F_{0}}) \subset G'$. Then $e^{\rm ord}_{K}\Pi=\Pi^{\rm ord}$. 
\end{lemma}
\begin{proof}  
With notation as in the proof of Proposition \ref{hida m1}, let 
$\phi_{w}$ be a generator of the line $\Pi^{\dagger}[\xi^{\circ}.w]$, where we define $\xi.w$ by $\wt{\xi.w}=\wt{\xi}^{w}$.
Write a nonzero spherical vector $\phi_{K}\in \Pi^{K}$ as 
\beq\lb{sph eb}
\phi_{K}=\sum_{w\in W_{G}}  c_{w} \phi_{w}\eeq
 with $c_{w}\in L$.
  Then we need to show $c_{1}\neq 0$. Now by \cite[Lemma 3.9]{cass1}, the expansion of \cite[Lemma 3.8]{cass1} (where $\chi=\wt\xi^{\circ}$) is of the form \eqref{sph eb}, and there one has (see Theorem 3.1 \emph{ibid.}) that $c_{1}=1$.
\end{proof}

\subsubsection{$p$-adic Rankin--Selberg period} Let $\chi\in Y_{L}$. We define a functional on $\Pi^{\dag}$ by 
\beq\lb{P1 dag def} P_{1, \Pi,\chi}^{\dag}\coloneqq  \lim_{s\rightarrow \infty} P_{1, \Pi,\chi,s }^{\dag} , \qquad 
P_{1, \Pi,\chi,s}^{\dag} \coloneqq  q^{d(n)s} P_{1, \Pi, \chi} \circ m_{s} U_{t}^{-s}\colon \Pi^{\dag}
\to L( \chi),
\eeq
where $P_{1, \Pi, \chi}\colon \Pi\to L(\chi)$ was implicitly defined in the proof of Proposition  \ref{rat rtf} \eqref{rat sph ch}.

Let $c(\chi)$ be the  conductor of $\chi$ in the usual sense: $c(\chi)=0$ if $\chi $ is unramified and otherwise $c(\chi)$ is the minimal $c\in \Z_{\geq 1}$ such that $\chi_{|1+\vpi^{c}\sO_{0}}=1$.
\begin{lemma}
The sequence in the limit \eqref{P1 dag def} stabilizes as soon as $s\geq s_{0}\coloneqq  \max\{1, c(\Pi^{\dag})-1, c(\chi)\}$.
\end{lemma}
\begin{proof} 
In the definition  of $P_{1, \Pi, \chi, s}^{\dag}(W)$ in \S~\ref{sec fact}, we may first integrate over $K_{H}^{(s_{0})}$; observing that $\chi$ is $(\det K_{H}^{(s_{0})})$-invariant, the lemma results from \eqref{mat id integral}.
\end{proof}

\subsubsection{$p$-adic pairing}
We define a (non-degenerate) pairing 
$$\vth_{\Pi}^{\dag}\coloneqq  \lim_{r\rightarrow \infty} \vth_{\Pi, r}^{\dag}  , \qquad 
\vth_{\Pi, r}^{\dag}( \cdot, \cdot) \coloneqq  q^{d(n)} \vth_{\Pi} (w_{r} U_{t}^{-r} \cdot, \cdot)\colon \Pi^{\dag}\ts \Pi^{\vee, \dag} \to L.$$
It is easy to show, using the symmetry $(K^{\la c\ra})^{w_{c}}= K^{\la c\ra}$, that the sequence in the limit stabilizes as soon as $r\geq c(\Pi^{\dag})$. 
\begin{remark} \lb{perfect dag}
 For all $t'\in T^{+}$ the $\vth_{\Pi}^{\dag}$-adjoint of  $U_{t'}$ is $U_{t'^{\iota}}$. Thus for every character $\xi$, the pairing $\vth_{\Pi}^{\dag}$ yields a perfect pairing on $\Pi^{\dag}[\xi] \times \Pi^{\vee, \dag}[\xi^{\iota}] $ and moreover, for all $r\geq c(\Pi^{\dag})$ and every semistandard deeper Iwahori subgroup $K\subset G'$, a perfect pairing on  $\Pi^{K, \dag}[\xi] \times \Pi^{\vee,  K^{w_{r}}, \dag}[\xi^{\iota}] $.\end{remark}

\subsubsection{$p$-adic Flicker--Rallis period} Suppose  that $F=F_{0}\ts F_{0}$ and that  $\Pi$ is  in the image of the local base change map \eqref{local BC}; in other words, we may write\footnote{We abusively still write $\Pi$ for the inflation via $\wt\G'_{}\to \G'_{}$, for the sake of lightness.}
 $\Pi_{}\cong \Pi_{0}\boxtimes \Pi_{0}^{\vee}$ for some representation $\Pi_{0}$ of ${G}'_{0}$.  We define
$$P_{2, \Pi}^{\dag}: = \lim_{r\rightarrow \infty}
 P_{2, \Pi, r} , \qquad  
P_{2, \Pi, r} \coloneqq  q_{0}^{d(n)r} P_{2}\circ [1; w_{0, r}]  U_{[1; t_{0}]}^{-r} 
\colon \Pi^{\dag} \to L.$$
The sequence in the limit stabilizes as soon as $r\geq c(\Pi^{\dag})$.

\subsubsection{$p$-adic Rankin--Selberg periods at $U_{t}$-eigenvectors} \lb{dec xi}
 Identify $\Pi_{n+1}$ (respectively $\Pi_{n}$) with its $\psi$- (respectively $\ol{\psi}$-) Whittaker model, and $\Pi$ with their product.  Suppose that $\Pi$ is a subquotient of a regular principal series.

 Let $\xi\in\widehat{T}_{L}$ be a character that  occurs  in $\Pi^{\dag}_{L(\chi)}$; by the argument in the  proof of Proposition \ref{hida m1}, we have $\dim_{L(\xi)} \Pi^{\dag}_{L(\xi)}[\xi]=1$.  We denote by $W_{\xi}\in \Pi_{L( \xi)}$ the vector introduced at the end of \S\ref{fss}.

For $\chi\in Y_{L(\xi)}$, define 
\beq\lb{ev def}
e(\Pi, \xi, \chi) \coloneqq  {P_{1, \Pi,\chi}^{\dag} (W_{\xi})}
\in L(\xi, \chi) \coloneqq L(\xi)(\chi).\eeq
Liu and Sun have recently proved a result that is equivalent to an explicit formula for this term.  Write $\wt{\xi}=\wt{\xi}_{n}\boxtimes \wt{\xi}_{n+1}$, and for $1\leq i\leq \nu$,  let $\wt{\xi}_{\nu, i}\colon F^{\ts}\to L(\xi, \chi)^{\ts}$ be the restriction of $\wt{\xi}_{\nu}$ to the $i^{\rm th}$ component of $T_{\nu}=(F^{\ts})^{\nu}/F_{0}^{\ts}$. For any character  $\xi'$ of $F^{\ts}$ and any factor $F_w$ in the decomposition of $F$ as a product of (one or two) fields. Denote by $\xi'_{w}\coloneqq \xi_{|F_{w}^{\ts}}$; denote by $N_{w}\colon F_{w}^{\ts}\to F_{0}^{\ts}$ the norm map.
Finally, we denote by 
$$\gamma(s,\xi_{F,w}', \psi_{F,w})^{-1}\coloneqq 
{L(s, \xi_{w}')/ \ep(s, \xi_{w}', \psi_{F,w}) L(1-s, \xi_{w}'^{-1})}$$
the inverse Deligne--Langlands $\gamma$-factor of a character of $\xi_{w}'\colon F_{w}^{\ts}\to \C^{\ts}$. If 
$$|\cdot|^{1/2}\xi'_{k}, \quad |\cdot|^{1/2} \xi''_{k}\colon F_{w}^{\ts}\into L'^{\ts}\subset \C$$ (for $1\leq k\leq M$) are characters with $\prod_{k=1}^{M}|\cdot |^{1/2}\xi'_{k} =\prod_{k=1}^{M} |\cdot|^{1/2} \xi''_{k}$, then it is easy to see that $\prod_{k=1}^{M}\gamma(1/2, \xi'_{k} , \psi_{F,w}) /\gamma(1/2, \xi''_{k} , \psi_{F,w})$ belongs to $L'$. Thus the following expression  gives an element  of   $L(\xi, \chi)$ (unless some division by zero has occurred). 

Define
 $$\hat{e}(\Pi,\xi,  \chi)\coloneqq  {\ep({1\over 2}, \chi^{2}, \psi)^{{n+1\choose 2}}\over L({1\over 2}, \Pi\ot\chi)}
\prod_{i+j\leq n} \prod_{w}\gamma({1\over 2}, \chi\circ N_{w}\cdot \wt{\xi}_{n, i,w} \wt{\xi}_{n+1,j,w}, \psi_{F,{w}})^{-1},$$
where as recalled above $w$ indexes over the field factors of $F$.

\begin{proposition}\lb{hyp ep}
There exists a constant $c_{\Pi, \xi}\in L(\xi)^\ts$ such that for every character $\chi$,
$$e(\Pi,\xi,  \chi)= c_{\Pi, \xi}\cdot \hat{e}(\Pi,\xi,  \chi).$$
\end{proposition}
\begin{proof} {We reduce this to a calculation of Liu--Sun.

For every Borel subgroup $B'\subset G'$ with unipotent radical $N'$, we denote by $\frak b'$ its Lie algebra, and  we add a bar to denote the opposite Borel and attached objects. We then let 
$\left(\widehat\Pi_{\mathfrak{b}'\textup{-smooth}}\right)^{{N}'}$ be the set of vectors in the formal completion of $\Pi$ that are fixed both by  ${N}'$ and by some closed subgroup of $G'$ with Lie algebra~$\frak b'$,  see \cite[paragraph following (1.6)]{LS} for more details.

We have a commutative diagram of isomorphisms,
equivariant for the $T$-actions up to twists,
\begin{equation*}
\xymatrix{
\Pi^\dag\ar[r]\ar[d]^{\Pi(m)} &\Pi_N\ar[r] &\left(\widehat\Pi_{\frak{\ol{b}}\textup{-smooth}}\right)^{\ol{N}}\ar[d]^{\widehat \Pi(m)}\\
\Pi_{mBm^{-1}}^\dag\ar[r] & \Pi_{mNm^{-1}}\ar[r] &\left(\widehat\Pi_{m\frak{\ol{b}}m^{-1}\textup{-smooth}}\right)^{m\ol{N}m^{-1}}
}
\end{equation*}
where in each row, the  first map is the natural projection and the second (also canonical) map is a form of Bernstein's second adjointness as recalled in \emph{loc. cit.}, with $\ol{N}$ denoting the  unipotent subgroup opposite to $N$.
The composition in the upper row is given explicitly  by $\phi\mapsto \lim_{s\to \infty} \Pi(t^s) U_{t}^{-s}\phi$.   Similarly for the lower row, where  $\Pi_{mBm^{-1}}^\dag\coloneqq\Pi(m)\Pi^\dag$, which is also the finite-slope space with respect to the Borel $mBm^{-1}$. Note that the latter is the unique Borel subgroup of $G'$ that is transversal to $H_1'$ in the sense that $mBm^{-1} H_1'$ is open in $G'$.

Starting from $W_\xi$ in the upper-left corner and going to the lower-right corner, if we move right then down, we get the formal vector $\widehat{W}_\xi$ such that $P_{1, \Pi, \chi}^\dag(W_\xi)= P_{1, \Pi, \chi}(\widehat{W}_\xi)$. On the other hand, Proposition 8.7 \emph{ibid.} computes the value of $P_{1, \Pi, \chi}$ at the image of $\Pi(m)W_\xi$ under the lower horizontal map, up to a multiple in $L(\Pi, \xi)^\ts$ depending on the choices of normalizations. (In \emph{loc. cit.}, it is further assumed that the vector $\Pi(m)W_\xi$ is `ordinary' for $mBm^{-1}$, but this assumption plays no role.) That is, it computes a nonzero constant multiple of $P_{1, \Pi, \chi}(\widehat{W}_\xi) = e(\Pi, \xi, \chi)$, with the asserted result.

The constant $c_{\Pi, \xi}$  can be made explicit by comparing Liu--Sun's formula  in the case of sufficiently ramified $\chi$ with the Local Birch Lemma of \cite[Theorem 2.8]{janu-nonab}. We leave this to the interested reader. A direct evaluation can also be obtained analogously to \cite[Theorem 5.10]{zheng}, as noted in Remark 5.11 \emph{ibid}.}
\end{proof}

The key consequences for us will be Propositions \ref{nv iw ps} and Remark \ref{unr nexc} below, both derived from the following lemma. 
We temporarily restore the  notation of the rest of the paper.

 \begin{lemma}  \lb{nv ep}
Suppose that $\Pi_{v}$ is a regular irreducible principal series that is  the local component of a representation $\Pi$ in $\sC(\G')^{\rm her}_{L}$. For every  smooth character  $\xi_{v}\in \widehat{T}_{L(\Pi)}$ that occurs in $\Pi_{v, {L(\Pi, \xi_v ) }}^{\dag}$ and every  finite-order character $\chi_{v}\in Y_{v, L(\Pi, \xi_v)}$, we have
  $$\hat{e}(\Pi_{v}, \xi_{v}, \chi_{v})\in L(\Pi,\xi_v, \chi_{v})^{\times}.$$
  \end{lemma}
\begin{proof} By \cite[Theorem 1.1]{Ca2}, for each finite place $w$ of $F$,  the semisimple Weil--Deligne  representation  attached to  $\rho_{\Pi|G_{F_{w}}}$ (cf. \eqref{rho Pi intro}) is 
 $$r_{\Pi, w}=\bigoplus_{1\leq i\leq  n, 1\leq j\leq n+1}  |\cdot|^{1/2} \wt{\xi^{}}_{n, i,w} \wt{\xi^{}}_{n+1,j,w},$$
 and it is strictly pure of some weight that is independent of $w$ (here we identify a character of $F_{w}^{\ts}$ with its correspondent on the Weil group of $F_{w}$ via class field theory). By considering $\det r_{\Pi, w}$ at an inert place $w$ we then see that the weight must be $-1$. Thus for each $(i, j)$, the character  $|\cdot|^{1/2} \wt{\xi^{}}_{n, i,w} \wt{\xi^{}}_{n+1,j,w}$ is either ramified (so that its $\gamma$-factor is an $\ep$-factor, hence  nonzero), or it is an unramified character whose value at a uniformizer of $F_{w}$ is a Weil $q_{w}$-number of weight~$-1$, which again implies the non-vanishing of each term in the $\gamma$-factors of Proposition \ref{hyp ep}. 
\end{proof}

\subsection{$p$-adic  relative characters}  \lb{sec:53}
We go back to dropping the subscripts as in  the rest of this section. We say that a subgroup $K\subset G'$ is \emph{convenient} if either $K=\G'(\sO_{0})$, or $F=F_{0}\ts F_{0}$ and $K$ is a conjugate-symmetric deeper Iwahori as defined in \S~\ref{IW sym} (henceforth: a CSDI).

\subsubsection{Finite-slope relative character}\lb{def sph ch} 
Let $K\subset G'$ be a convenient subgroup. 
We define a distribution 
$$I_{\Pi, K}^{\dag}\colon \sH^{\dag}({L}) \to\sO( Y_{v} , L)$$
by 
$$I_{\Pi,K}^{\dag}(f_{}^{\dag}, \chi) \coloneqq  
\begin{cases}
\Tr_{\vth_{\Pi}^{}}^{P_{1, \Pi, \chi}^{\dag}  \ot P_{2, \Pi}} (\Pi(f^{\dag}e_{K}))
& \text{if $K=G'(\sO_{0})$},\\
\ \\
\Tr_{\vth_{\Pi}^{\dag}}^{P_{1, \Pi, \chi}^{\dag} \ot P_{2, \Pi^{}}^{\dag}} (\Pi(f^{\dag}e_{K}))
& \text{if $F=F_{0}\ts F_{0}$ and  $K$ is a CSDI}.
\end{cases}
$$
\begin{remark} The second definition is the `correct' one from the $p$-adic point of view. The first one is made because, first, in the arithmetic side the geometry will compel us to work at spherical level; and second, we have not investigated the analogue of the notion of `conjugate-symmetric' in the nonsplit case. 
\end{remark}

\subsubsection{Eigen-decomposition} 
    Suppose that $\Pi$ is a subquotient of a regular principal series, and denote by $\Xi_{K}(\Pi)\subset \widehat{T}_L$  the  scheme (finite over $L$)  of characters of $T$ occurring in $\Pi^{ K, \dag}$.
\begin{itemize}
\item
If $K=\G'(\sO_{{0}})$ and $\Pi$ is an unramified  principal series, let $W_{0}\in \Pi^{K}$, $W_{0}^{\vee}\in \Pi^{\vee,K}$ be  generators normalized by $W_{0}(1)= W_{0}^{\vee}(1)=1$, write $W_{0}^{} = \sum_{\xi} \lm_{\xi} W_{\xi}$, and let 
\beq\lb{cKun}
  c_{K}(\Pi,\xi) \coloneqq  \lm_{\xi} P_{2, \Pi}^{}(W_{0}^{\vee})/ \vth_{\Pi}^{}(W_{0}, W_{0}^{\vee}) =  \lm_{\xi}.
\eeq 
where the second equality follows from Remark \ref{aa1}. By the same proof as for Lemma \ref{nv ord unr}, we have $\lm_{\xi}\neq 0$ for all $\xi \in \Xi_{\Iw}(\Pi)$. (An explicit formula for $\lm_{\xi}$ could be obtained from combining the formulas cited in that proof with the Casselman--Shalika formula \cite{CasSha} and the formulas of \cite[Proposition 3.1]{reeder} for Whittaker $U_{t}$-eigenfunctions.)
\item
If  $F=F_{0}\ts F_{0}$ (so that we may write $\Pi\cong \Pi_{0}\ts \Pi_{0}^{\vee}$) and $K$ is a conjugate-symmetric deeper Iwahori,  define 
\beq\lb{cdag}
c_{K}(\Pi, \xi)\coloneqq c(\Pi, \xi)\coloneqq {P_{2, \Pi}^{\dag}(W_{\xi^{\iota}})\over \vth_{\Pi}^{\dag}(W_{\xi}, W_{\xi^{\iota}})} \in L(\xi).
\eeq
Here, the denominator is non-vanishing since $U_{t}$ is $\vth_{\Pi}^{\dag}$-adjoint to $U_{t^{\iota}}$ (Remark \ref{perfect dag}). Similarly, the numerator is non-vanishing  if and only if $\Pi$ is hermitian. More precisely, when $\Pi$ is hermitian,  if we view $P_{2, \Pi}^{\dag}$ as a bilinear form on $\Pi_{0}^{\dag}\ts\Pi_{0}^{\vee, \dag}=\Pi^{\dag}$, then the operator $U_{(t_{0})}$  on $\Pi_{0}^{\dag}$ (which corresponds to $U_{[t_{0}; 1]}$ on $\Pi^{\dag}$)  is $P_{2, \Pi}^{\dag}$-adjoint to  $U_{ t_{0}^{\iota}}$ on $\Pi_{0}^{\vee, \dag}$  (which corresponds to $U_{[1;t_{0}]}$ on $\Pi^{\dag}$). In this case, the explicit value of $c_{K}(\Pi, \xi)$ has been recently computed by Zheng Wang, \cite{zheng}.
\end{itemize}

Then, in either case, by the definitions we have a decomposition 
\beq\lb{dec dec}I_{\Pi, K}^{\dag}(f^{\dag}, \chi) = \sum_{\xi \in \Xi_{K}(\Pi)}
 I_{\Pi,K, \xi}^{\dag}(f^{\dag}, \chi) \eeq
where
\beq\lb{Idag ep} I_{\Pi, K , \xi}^{\dag}(f^{\dag}, \chi) \coloneqq  
\xi(f^{\dag}) c_{K}(\Pi, \xi) e(\Pi, \xi, \chi).
  \eeq

\subsubsection{Ordinary relative character} \lb{sec: ord sph} 
Suppose for this paragraph only that  $L$ is a finite extension of $\Q_{p}$ and that there is an $\sO_{L}$-lattice $\Pi_{\sO_{L}}\subset \Pi$ that is stable under $\sH^{\dag}$. 
 Then  we have Hida's description 
$${\rm e}^{\rm ord}= \lim_{N \rightarrow \infty} U_{t}^{N!} $$
for the action of the ordinary projector on $\Pi$. 
\begin{remark} The above assumption holds whenever $\Pi$ is a local component of a global representation in $\sC(\G')_{L}$. Indeed,  representations in $\sC(\G')_{L}$ can be realised in the Betti cohomology of the locally symmetric space attached to $\G'$, and the cohomology with coefficients in $\sO_{L}$ gives a natural integral structure stable under the Hecke operators; see \cite{hida-aut} for more details.
\end{remark}

For any convenient  $K\subset G'$, we then define
$$I_{\Pi, K}^{\ord}(\chi): = \lim_{N \rightarrow \infty} I_{\Pi, K}^{\dag}(U_{t}^{N!}, \chi).$$
If $\Pi $ is ordinary, we denote
\beq\lb{e v ord} 
e(\Pi,  \chi)&\coloneqq  e(\Pi, \xi^{\circ}, \chi) \in L(\chi),\\
c_{K}(\Pi)&\coloneqq  c_{K}(\Pi, \xi^{\circ}) \in L^{\ts}
\eeq
where the right-hand sides are defined in   \eqref{ev def}, \eqref{cKun}, \eqref{cdag}.

\begin{remark}\lb{unr nexc}
 If $\Pi_{v}$, $\chi_{v}$ are  as in Lemma \ref{nv ep} and moreover $\Pi_{v}$ is  ordinary, it follows from that lemma and Proposition \ref{hyp ep} that $e(\Pi_{v}, \chi_{v})$  and $c_{K_{v}}(\Pi)$ are nonzero. 
\end{remark}

\begin{corollary} \lb{Iord ep}  Suppose that $\Pi$ admits an $\sO_{L}$-stable lattice. Then for every  $\chi\in Y_{L}$ and every convenient   $K\subset G'$,
we have 
$$I_{\Pi, K}^{\ord}( \chi) = \begin{cases}
c_{K}(\Pi) e(\Pi, \chi) & \text{if $\Pi$ is $K$-ordinary}\\
0 & \text{otherwise}.
\end{cases}
$$
\end{corollary}
\begin{proof} This  follows from  \eqref{Idag ep}.
 \end{proof}
 
\subsubsection{Relation to the character $I_{\Pi}$}\lb{sec suff pos}
 Let $\Pi$ be a tempered irreducible representation of $G'$,   let $\chi\colon F_{0}^{\ts}\to L^{\ts}$ be a smooth character,  let $K\subset G'$ be a convenient subgroup, and let $s\geq 1$. We say that $s$ is sufficiently positive for $\chi$ (respectively for $K$) if $s\geq \max\{1, c(\chi)\}$ (respectively $K$ contains a deeper Iwahori of level $c$ with\footnote{In fact, at least if $K$ is an  Iwahori subgroup or one of the subgroups \eqref{Krdef} with $r=c\geq1$, the weaker condition  $s\geq c$ will suffice; this is only used in the application of Lemma \ref{iw mult} (3) in the proof of Lemma \ref{f*}.}
 $s\geq 2c$).
 We say that $f^{\dag}\in \sH^{\dag}(L)$ is \emph{sufficiently positive for $\Pi$ (respectively for  $s_{0}$,  for $\chi$, for $K$)} if   $f^{\dag} \Pi \subset \Pi^{\dag}$ (respectively  if $U_{t}^{-s}f^{\dag}$ belongs to $\sH^{\dag, +}(L)$ for $s=s_{0}$,  for some $s$ that is sufficiently positive for $\chi$,  for some $s$ that is sufficiently positive for $K$).

 It is clear that if $f^{\dag}$ is in the span of $\{U_{t}\  |\, t\in T^{++}\}$ and $s$ and $\Pi$ are given, then some power of $f^{\dag}$ is sufficiently positive for both $s$ and $\Pi$.

\begin{lemma}\lb{I Pi dag}
For every $s$ that is sufficiently positive for $K$ and $\chi$ and 
every $f^{\dag}$ that is
 sufficiently positive for $s$ and $\Pi$,  we have
$$I_{\Pi, K}^{\dag}(f^{\dag}, \chi) = I_{\Pi}(f', \chi)$$
where  
\begin{subnumcases}{f'=f'_{K, s} : =  }
q^{d(n)s} \cdot m_{s} U_{t}^{-s}f^{\dag}e_{K}  & \text{if $K=\G(\sO_{0})$},  \lb{fp def a} 
\\
  q_{0}^{d(n)(2s-c)} \cdot m_{s} U_{t}^{-s}f^{\dag} e_{K}  U_{[t_{0}; 1]}^{c} [ w_{0, c}^{-1}; 1] 
 & \text{\parbox[t]{5cm}{if $F=F_{0}\times F_{0}$ \\ and $K$ is a CSDI of depth $c$.}} \lb{fp def b} 
\end{subnumcases}

\end{lemma}

\begin{proof} The first case is clear. Consider the second case,  dropping the subscripts $\Pi$ and  $K$ from the notation.  Let $\Pi_{\dag, K} \coloneqq  w_{c}\Pi^{\dag, K}$, and let $\vth_{|}\colon \Pi_{\dag, K}\ot \Pi^{\vee, \dag, K}\to L$ be the restriction of $\vth\colon \Pi\ot \Pi^{\vee}\to L$, which is still a perfect pairing. 
  By Lemma \ref{different RT} (using, in order,  part \eqref{drt0}, part \eqref{drt1},  and part \eqref{drt0} together with part \eqref{drt2}), 
\beqq
I^{\dag}(f^{\dag}, \chi) &
 = q_{0}^{d(n)c}  \Tr_{\vth^{\dag}}^{P_{1, \chi}^{\dag}   \ot P_{2}[1;w_{0, c}]} (\Pi( f^{\dag} e_{K} U_{[1;t_{0}^{\iota}]}^{-c}))
 \\
&
=q_{0}^{-d(n)c}\Tr_{\vth_{|}}^{P_{1, \chi}^{\dag} U_t^c w_c^{-1} \ot P_{2}[1;w_{0, c}]}
(\Pi (w_c U_t^{-c} f^{\dag}   e_{K} U_{[t_{0};t_{0}/t_{0}^{\iota}]}^{c} w_{c}^{-1} ))
\\
&=q_{0}^{d(n)(2s-c)}\Tr_{\vth}^{P_{1, \chi}\ot P_{2}}
(\Pi(m_{s} U_{t}^{-s}f^{\dag} e_{K} U_{[t_{0}; z_{0}]}^{c}   w_{c}^{-1} [1;w_{0, c}^{-1}] ) )
=I (f', \chi),
\eeqq
where $f'$ is as asserted.
\end{proof}

\subsubsection{A non-vanishing result} 
Unlike the rest of this section, the following result is not used for the $p$-local theory of the $p$-adic relative-trace formula, but rather as an input to Proposition \ref{test fin} \eqref{test fin 3}.
\begin{proposition}\lb{nv iw ps}
 Let $\Pi$, $\chi$, $K$ be as in \S~\ref{sec suff pos}. Suppose that $F=F_{0}\ts F_{0}$, $\Pi$ is a regular principal series, and  $K$ is a conjugate-symmetric deeper Iwahori such that $\Pi^{K, \dag}\neq 0$.
Then there exists an element $f^{\dag}\in \sH^{\dag}$ that is sufficiently positive for $\Pi$, $\chi$, $K$, such that  the Hecke measure $f'\coloneqq f'_{K, s}=\eqref{fp def b}$ satisfies
$$I_{\Pi}(f'_{K, s}, \chi)\neq 0.$$
\end{proposition}
\begin{proof}  Order the characters $\xi$ occurring in $\Pi^{\dag}$ as $\xi_{1}, \ldots, \xi_{r}$. Pick some sufficiently positive  $t'\in T^{++}$ such that $\xi_i(t')\neq \xi_j(t')$ for $i\neq j$(the existence of such $t'$ is easy to verify). For each integer $N\geq 1$ let $f'_{ N}$ correspond to $f^{\dag}=U_{N}$ 
 We may  and do extend scalars from $L$ to $\C$; we do not   alter the notation. 
By \eqref{dec dec}, we have
\beqq 
I_{\Pi}(f'_N, \chi)=\sum_{\xi \in \Xi_{K}(\Pi)}
  \xi(t')^N  {c}_{K}(\Pi, \xi) {e}(\Pi,\xi, \chi).
  \eeqq
We may write $I_{\Pi}(f'_N, \chi)=a_{N} x$ where $x=(c_{K}(\Pi,\xi_{i})e(\Pi, \xi_{i}, \chi))_{i}\in \C^{r}$ and for any positive integer $M$ we put  $a_{M}^{\rm t} \coloneqq (\xi_{i}(t')^{M})_{i}\in \C^{r}$.  Now all entries of the vector $x$ are nonzero by Proposition \ref{hyp ep} and Lemma \ref{nv ep}, and the (almost) Vandermonde matrix  $A$ with rows $a_{N}, \ldots, a_{2N}, \ldots, a_{rN}$ is invertible. Hence there is some $1\leq i\leq r$ such that $0\neq a_{iN}x = I_{\Pi}(f'_{iN}, \chi)$, as desired.
\end{proof}

\section{$p$-adic orbital integrals} \lb{app orb p}
We define and study certain local orbital integrals matching the relative characters just defined. After establishing their $p$-adic boundedness (as the character $\chi$ varies, in a suitable sense), the main result of this section, Proposition \ref{chi1 along}, says that in case $F=F_{0}\ts F_{0}$ and $K$ is a CSDI, our orbital integrals have plus-regular support, and it explicitly computes the values at all orbits. 

  In \S~\ref{sec prove test fin} we define the $p$-adic local orbital integrals and state their main properties, along with an application  to the existence of   local test measures with plus-regular support, which is needed (at  split, typically non-$p$-adic,  places) to complete the results of \S~\ref{sec: test}..
  The rest of the section is dedicated to proving Proposition \ref{chi1 along}.   In \S\S~\ref{sec:62}-\ref{sec:63}, we reduce the  explicit   evaluation of $p$-adic orbital integrals to a  $p$-adic linear algebra statement,  which  is proved by a somewhat elaborate induction that occupies  \S\S~\ref{sec:64}, \ref{sec:66}, \ref{sec:67}. 
We pause along the way in \S~\ref{sec:65} in order to reap  the  plus-regularity of supports. 

We continue with the notation of the previous section.

\subsection{Definition and statement of the main result} \lb{sec prove test fin}
Let $K\subset G'$ be a convenient subgroup and let $L$ be a coefficient field. 
\subsubsection{Definition of the orbital integrals}
 For $f^{\dag}\in \sH^\dag(L)$ sufficiently positive (depending on $\chi$) and $\gamma\in B'$,  let  $s$ and $f'_{K,s}$ be  as in Lemma \ref{I Pi dag}, and define
\beq\lb{def orb dag}
I_{\gamma,K}^{\dag}(f^{\dag}, \chi) &\coloneqq  L_{\gamma}(\chi)  I_{\gamma} (f'_{K, s}, \chi) 
\eeq
where the terms in the right-hand side are as given in Proposition \ref{rat rtf}.

\begin{lemma}\lb{mu bded} Let $\gamma \in B'$ and let $f\in \sH^{\dag}(L)$.
\begin{enumerate}
\item
The right hand side of \eqref{def orb dag}  is independent of the choices of an integer  $s$ that is sufficiently positive for $K$ and $\chi$, 
 so long as $f^{\dag}$ is sufficiently positive for $s$. 
\item \lb{mub 2}
Suppose that $f^{\dag}\in\sH^{\dag}(\OO_{L})$. For any $s_{0}$ that is sufficiently positive for $K$ such that $f^{\dag}$  is sufficiently positive for $s_{0}$,  the map 
$$\chi\longmapsto I_{\gamma, K}^{\dag} (f^{\dagger}, \chi)$$
extends by linearity to a functional $C^{\infty}(F_{0}^{\ts}/ (1+ \vpi^{s_{0}}\sO_{0}), \sO_{L})\to k\sO_{L}$
for the constant
\beq\lb{const k} k=q_{0}^{-d(n)c(K)}\vol(K\cap H_{2}')\vol^{\circ}(K_{H}) \qquad  \in \Q^{\ts}\eeq
depending only on $K$. (Here, $\vol^{\circ}(K_{H})=\eqref{vol circ}$.)
\end{enumerate}
\end{lemma}

\begin{proof}
If $F=F_{0}\ts F_{0}$ and $K$ is a CSDI, let $c=c(K)$ be the depth of $K$, and let $A_c=U_{[t_{0}; 1]}^{c}[w_{0 ,c}^{-1};1]$; if $K$ is unramified, let $c=0$ and $A_0={\rm id}$. For a Hecke operator $A' =\sum_{i} \lm_{i}\vol(\Omega_{i})^{-1}\one_{\Omega_{i}} dg \in \sH_{K}$, denote $\one[\gamma'\in A']\coloneqq \sum \lm_{i} \one_{\Omega_{i}}(\gamma').$

   Let $s_{0}$ be sufficiently positive for $K$ and $\chi$.  The integrand in  the explicit expression for  $I_{\gamma, K}^{\dag}(f^{\dag}, \chi)$ equals 
\beq\lb{integrand}
q_{0}^{d(n)(2s-c)}\chi(h_{1}) \one[\gamma  h_{2}\in  h_{1} m_{s} U_{t}^{-s}  f^{\dag}e_{K} A_c]
\eeq
and it is $K_{H}^{(s_{0})}\ts (K\cap H_{2}')$-invariant  by \eqref{mr KH}. 
Integrating first  over $K_{H}^{(s_{0})} \subset H_{1}'$, the relation \eqref{mat id integral} shows that \eqref{integrand} is independent of $s\geq s_{0}$. We also see that 
    the functional $I^{\dag}_{\gamma, K}(f^{\dag}, -)$ sends  $C^{\infty}(F_{0}^{\ts}/ (1+\vpi^{s_{0}}\sO_{0}), \sO_{L})$ to $k\sO_{L}$.
\end{proof}

\subsubsection{Main result and application to regular test Hecke measures}
When $F=F_{0}\ts F_{0}$ and $K$ is a conjugate-symmetric deeper Iwahori and $f^{\dag}$ is sufficiently positive, the following key result asserts that the associated $f'_{K,s}$ has plus-regular support and the  $p$-adic orbital integral may be  explicitly computed. A remarkable fact is that its value is  independent of $\chi$. 

By linearity, it suffices to study the case $f^{\dag}=U_{t'}$ for some $t'\in T^{++}$.   
\begin{proposition}  \lb{chi1 along}  Suppose that $F=F_{0}\ts F_{0}$ and let $K=K_{0}\times K_{0}^{w_{c}}\subset G'$ be a conjugate-symmetric deeper Iwahori. Assume that $f^{\dag}=U_{t'}\in \sH^{\dag}$  for some $t' \in T^{++}$. Then:
\begin{enumerate}
\item \lb{+reg spt}
for every $s$ that is sufficiently positive for $K$ such that   $f^{\dag}$ is sufficiently positive for $s$, the support of 
$$f'_{K, s}=\eqref{fp def b}$$
is contained in $G'_{\reg^{+}}$; moreover, the Hecke measure $f'_{K ,s}$ matches a measure $f_{K,s}\in \sH(G' ,L)$ that is  regularly supported and bi-invariant under a subgroup conjugate to $K_{0}$;
\item  \lb{chi1 along!}
there exists  a compact subset  
$$B_{K}^{\dag}(f^{\dag}) \subset B'$$
with the following property:
 for every  smooth character $\chi$  of $F_{0}^{\ts}$ such that  $f^{\dag}$ is sufficiently positive for $\chi$ and $K$,
we have
\beqq I_{\gamma, K}^{\dag}(f^{\dag}, \chi) = 
\begin{dcases} k' & \text{if } \gamma \in B_{K}^{\dag}(f^{\dag}) \\
0& \text{if } \gamma \notin  B_{K}^{\dag}(f^{\dag}) ,
\end{dcases}
\eeqq
where $k'= kq_{0}^{-\ell_{K_{0}}(w_{0})}$, with $k=\eqref{const k}$.
\end{enumerate}
\end{proposition}

The proof of  Proposition \ref{chi1 along} will occupy  the rest of this section, which may be skipped on a first reading.

The first part allows to complete the proof of Proposition \ref{test fin}.
\begin{proof}[Proof of Proposition \ref{test fin} \eqref{test fin 3}] 
We drop the subscript $v$ from the notation of the  statement of the proposition. Recall that we need to find a Hecke measure $f'_{\pm}\in \sH(G', L)$ that is supported in $G'_{\reg^{\pm}}$ and  adapted to a given pair $(\Pi, \chi)$. 
 
 We may take $f'_{+}$ to be the element $f'_{K,s}=\eqref{fp def b}$ associated to the data of: a  conjugate-symmetric deeper Iwahori $K$ such that $\Pi^{K}\neq 0$; an integer $s\geq 1$ that is sufficiently positive  for $\chi$ and $K$; and an element $f^{\dag}$ that is sufficiently positive for $\Pi$, $\chi$, $K$. 
Then the measure $f'_{+}$  is adapted to $(\Pi, \chi)$ by   Proposition \ref{nv iw ps}, and it has plus-regular support by Proposition    \ref{chi1 along} \eqref{+reg spt}.  If $\Pi$ is unramified we can take $K$ to be an Iwahori subgroup, hence (again by Proposition    \ref{chi1 along} \eqref{+reg spt}) we have that $f'_{+}$ matches a measure $f_{+}$ that is bi-invariant under an Iwahori subgroup.

We may take $f_{-}'\coloneqq  f_{+}'^{\diamond}$ for the involution $g^{\diamond}=g^{\rm c , -1, \rm t}$ of Remark \ref{invol +-}.  By that remark, $f_{-}'$ is minus-regular, and it is clear that  its matching $f_{-}$ is bi-invariant under an Iwahori subgroup if $f_{+}$ is. Moreover, $I_{\Pi}(f_{-}, \chi)=I_{\Pi^{\diamond}}(f_{+}, \chi^{-1})$ for  $\Pi^{\diamond}(g)\coloneqq \Pi(g^{\diamond})$; since $\Pi^{\diamond}\cong\Pi^{\rm c, \vee}\cong \Pi$, this expression is non-vanishing too.
\end{proof}

\subsection{Reduction to $p$-adic linear algebra} \lb{sec 71} \lb{sec:62}
We start working towards the proof of  Proposition \ref{chi1 along}, of which we retain all the assumptions. 
The proof of part \eqref{chi1 along!} relies  on some reductions in the present subsection and  \S~\ref{sec 72}, and on two auxiliary inductive lemmas in \S\S~\ref{sec 73},~\ref{sec 75}, and it is completed in \S~\ref{sec 76}. The proof of part \eqref{+reg spt} relies on the first auxiliary lemma, and is given in \S~\ref{sec 74}.

We keep using the notation of \S~\ref{sec not G0}; however, at various steps of our descent  into the argument, we will lighten  (and sometimes recycle) the notation  for the sake of readability. We start by dropping all apices from the notation,  writing for instance $f$ and $G$ in place of $f'$ and $G'$.

We define involutions $w$ and $\iota$ on $\Z^{\nu}$
by  
$$(\lm^{w})_{i}\coloneqq  \lm_{\nu+1-i}, \qquad \lm^{\iota}\coloneqq -\lm^{w},$$
 and a notion of positivity by declaring $ \lm\in \Z^{\nu, +}$ if $\lm_{i}\geq \lm_{i+1}$ for all $1\leq i\leq \nu-1$; thus $\iota$ preserves $\Z^{\nu, +}$. We also write $\lm\succeq \lm'$ if $\lm-\lm'\in \Z^{\nu, +}$. Then $\vpi^{\lm}\in 
\ul{T}_{\nu, *}^{+}$ if and only if $\lm\in\Z^{\nu, +}$, and $(\vpi^{\lm})^{\iota}=\vpi^{\lm^{\iota}}$. 

Extending the notation from \eqref{pr Hecke mu}, let  ${\rm p}_{\nu}\colon G_{\nu}\to G_{\nu, 0}\ts G_{\nu, 0}/F_{0}^{\ts}$ be the projection, and let ${\rm p}_{\nu, *}
\colon \sH(G_{\nu})\to  \sH(G_{\nu, 0}\ts G_{\nu, 0}/F_{0}^{\ts})$ be the pushforward map. Thus ${\rm p}\coloneqq {\rm p}_{n}\ts {\rm p}_{n+1}\colon \wt{G}\to G$ and ${\rm p}_{*}={\rm p}_{n, *}\ot{\rm p}_{n+1, *}$.

Let $c$ be the depth of $K$. By the positivity condition on $f^{\dag}$ and linearity, we may assume that
\beq \lb{fdag ass}
f_{\nu}^{\dag} &=U_{t',K}={\rm p}_{\nu, *} (f_{\nu, 1}^{\dag} \ot f_{\nu, 2}^{\dag}),\\
  f_{\nu, 1}^{\dag}&=[{K}_{\nu,0}\vpi^{\lm_{\nu, 1}} {K}_{\nu , 0}],\qquad  
f_{\nu, 2}^{\dag}=[{K}_{\nu,0}^{w_{c}}\vpi^{\lm^{w}_{\nu, 2}} {K}_{\nu , 0}^{w_{c}} ] \eeq
 for some $\lm_{\nu,i} \in \Z^{\nu}$ with $\lm_{\nu, i}, \lm_{\nu, i}^{\iota} \succeq  {(s+c)\rho_{\nu}}$.

We decompose 
$$f=f_{K, s} =\eqref{fp def b}=    q_{0}^{(2s-c)d(n)} \cdot m_{s} U_{t}^{-s} f^{\dag}  e_{K} U_{[ t_{0}; 1]}^{c} [w_{0, c}^{-1}; 1] 
 = f_{n}\ot f_{n+1} $$
where each $f_{\nu}$ is a Hecke measure  on $G_{\nu} /F_{0}^{\ts}$, and further decompose
\beqq
 f_{\nu}&= 
 q_{0}^{(2s-c)c(\nu)} m_{ \nu, s} f_{\nu}^{\dag}  U_{t_{\nu}}^{-s}e_{K_{\nu}} U_{[t_{\nu, 0}; 1]}^{-c}  [w_{0, \nu, c}^{-1}; 1] =
   {\rm p}_{\nu, *}(f_{\nu,1}\ot f_{\nu,2} ),
\eeqq
where $f_{\nu, i}\in \sH(G_{\nu ,0})$ are defined up to some scalar ambiguities that we do not need to resolve. We also denote $f_{i}=f_{n, i}\ot f_{n+1, i}$ for $i=1, 2$.\footnote{The context should prevent any possible confusion from the clash of notation with $f_{n}\in \sH(G_{n}/F_{0}^{\ts})$, since the integer in  this $f_{n}$ will never be specialized.}

Fix a representative $\gamma=[\gamma_{0}; 1]=[(\gamma_{n, 0},\gamma_{n+1, 0}); (1_{n} , 1_{n+1})]\in G$ under the decomposition $G=(G_{n,0}\ts G_{n+1, 0})^{2}/F_{0}^{\ts, 2}$.

Decompose $H_{1}=H_{1,0}^{2}$, and write $h_{1}\in H_{1}$ as $h_{1}=(h_{1, 0}, h_{1, 0}')$. In the orbital integral \eqref{orb int def}, we first integrate over
$H_{2}$ (noting that  $\eta=1$),  to find
\beq \lb{pO0}
 \chi(\gamma_{n})
I^{\dag}_{\gamma}(f^{\dag} , \chi) &=
  \int_{H_{1}}
\int_{H_{2}}
f([h_{1,0}^{ -1}\gamma_{0} h_{2}; h_{1,0}'^{ -1}  h_{2}])
\, dh_{2} \, \chi(h_{1}) \, dh_{1}
\\
 &= \int_{(H_{1, 0})^{2}} f^{\star}(h_{1, 0}^{-1}\gamma_{0} h_{1 ,0}')\, \chi((h_{1, 0}, h_{1,0}')) \, dh_{1, 0} dh_{1, 0}'
\eeq
where $$f^{\star}= f_{1}\star f_{2}^{\vee} \in \sH(G_{n,0}\ts G_{n+1, 0}).$$
(As part of the proof, we will show that the above integral always reduces to a finite sum.)

\begin{lemma}\lb{f*} Assume that $f^{\dag}=\eqref{fdag ass}$, and let   $\lm_{\nu}\coloneqq  \lm_{\nu, 1}+\lm_{\nu, 2}^{\iota}\succeq 2(s+c)\rho_{\nu} $. Then $f^{\star}=f_{n}^{\star}\ot f_{n+1}^{\star}$ for 
\beqq
f_{\nu}^{\star} \coloneqq   q_{0}^{2sc(\nu)-\ell_{K_{\nu, 0}}(w_{\nu, 0})} {m}_{s}[K  \vpi^{\lm_{\nu} -2s\rho_{\nu}} w K] m_{s}^{-1}
\qquad\in \sH_{G_{\nu}} .
 \eeqq
\end{lemma}

\begin{proof}
Let ${\sg}_{\nu}\coloneqq (\nu-1, \ldots, 0)\in \Z^{\nu}$, so that $t=\vpi^{\sg_{\nu}}$ and $\sg_{\nu}+\sg_{\nu}^{\iota}=2\rho_{\nu}$. Abbreviate  $w=w_{\nu, 0}$, $w_{c}=w_{\nu, 0, c}$; $m_{s}={m}_{\nu, 0, s}$;  $t= t_{\nu, 0}$; $K={K}_{\nu, 0}$;  $K'= K^{w_{ c}}$; $K''= K_{\nu, 0}^{\la c\ra}$.
 Then 
\beqq 
 f_{\nu}^{\star}=  f_{\nu, 1}\star f_{\nu,2}^{\vee}
&= q_{0}^{(2s-c)c(\nu)} {m}_{s}  U_{{t}}^{-s} f_{\nu, 1}^{\dag} e_{K} U_{t}^{c} w_{c}^{-1} e_{K'}  (U_{{t}}^{-s} f_{\nu,2}^{ \dag})^{\vee} {m}_{s}^{-1}\\ 
 &= q_{0}^{2sc(\nu)} {m}_{s}  
[ K'' \vpi^{\lm_{\nu, 1}-s\sg_{\nu}} {K}''] e_{K} e_{K''}t^{c}e_{K''}w_{c}^{-1}    [K'' \vpi^{-\lm_{\nu, 2}+(s-c)\sg_{\nu}} 
K'']
 {m}_{s}^{-1}\\ 
&= q_{0}^{2sc(\nu)} {m}_{s}  
[ K'' \vpi^{\lm_{\nu, 1}+(c-s)\sg_{\nu}} {K}''] e_{K} we_{K''}   [K'' \vpi^{-\lm_{\nu, 2}+(s-c)\sg_{\nu}} 
K'']
 {m}_{s}^{-1}\\ 
 &= q_{0}^{2sc(\nu)-\ell_{K''}(w)} {m}_{s}  
[ K'' \vpi^{\lm_{\nu, 1}-s\sg_{\nu}} {K}''] e_{K} [{K''}w{K''}]    [K'' \vpi^{-\lm_{\nu, 2}+s\sg_{\nu}} 
K'']
 {m}_{s}^{-1}\\ 
 &= q_{0}^{2sc(\nu)-\ell_{K''}(w)} {m}_{s}  
[ K'' \vpi^{\lm_{\nu, 1}-s\sg_{\nu}} {K}''] e_{K}   [K'' \vpi^{\lm_{\nu, 2}^{\iota}+s\sg_{\nu}^{\iota}} w 
K'']
 {m}_{s}^{-1}\\ 
  &= q_{0}^{2sc(\nu)-\ell_{K}(w)} {m}_{s}  
[ K \vpi^{\lm_{\nu, 1}-s\sg_{\nu}} {K}]   [K \vpi^{\lm_{\nu, 2}^{\iota}+s\sg_{\nu}^{\iota}} w K]
 {m}_{s}^{-1}\\ 
     &= q_{0}^{2sc(\nu)-\ell_{K}(w)} {m}_{s}  [K  \vpi^{\lm_{\nu} -2s\rho_{\nu}} w K] m_{s}^{-1},
\eeqq
where we have used the symmetry of $K''$ and the algebra rules  of Lemma \ref{iw mult}.\end{proof}

Let 
\beq \lb{Xnu0}
X_{\nu}^{\circ} \coloneqq  \vpi^{\lm_{\nu} -2s\rho_{\nu}}  w_{\nu,0} \in G_{\nu , 0}.
\eeq
By Lemma \ref{f*}, the integrand in \eqref{pO0} is non-vanishing at $h_{1}$ if and only if
\beq\lb{pO}
h_{1, 0}^{-1}\gamma_{\nu, 0}h_{1, 0}' &\in m_{\nu, 0,s} K_{\nu, 0}  X_{\nu}^{\circ}  K_{\nu, 0} m_{\nu,0, s}^{-1}
\eeq
for  $\nu=n, n+1$. Therefore, if the orbital integral $I_{\gamma}(f^{\dag} , \chi)$ is non-vanishing, up to changing the representative $\gamma_{0}$ in its $H_{1, 0}$-orbit we may and will assume that
 \beq \lb{where is gm}
\gamma_{\nu, 0} &\in m_{\nu,0, s} K_{\nu, 0}  X_{\nu}^{\circ}  K_{\nu, 0} m_{\nu,0,s}^{-1}.
\eeq

We introduce the convenient variables 
\beq\lb{Xnu def}
X_{\nu} &\coloneqq  m_{\nu, 0,s}^{-1} \gamma_{\nu, 0 } m_{\nu,0,s}.
\eeq
Then \eqref{where is gm} is equivalent to 
\beq\lb{where is Y}
X_{\nu}&\in K_{\nu,0}X_{ \nu}^{\circ} K_{\nu, 0}
\eeq
and  \eqref{pO} is equivalent to 
\beq \lb{pOY}
m_{\nu,0,s}^{-1} h_{1, 0}^{-1}m_{\nu,0,s} \cdot    X_{\nu} \cdot  m_{\nu, 0,s}^{-1}h_{1, 0}' m_{\nu, 0,s} \in K_{\nu, 0} X_{\nu}^{\circ} K_{\nu, 0}.
\eeq

We will reduce Proposition \ref{chi1 along} \eqref{chi1 along!}  to the following. 

\begin{proposition}\lb{Yhk}  Let $X_{\nu}^{\circ}\coloneqq  \vpi^{\lm_{\nu}'}w_{\nu, 0}$ for some $\lm'_{\nu}\in \Z^{\nu, +}$, and let  $(X_{n}, X_{n+1}, h_{1}) \in G_{n,0}\ts G_{n+1, 0}\ts H_{1}$
 satisfy \eqref{where is Y}, \eqref{pOY} for $\nu=n, n+1$.  Then $h_{1} =(h_{1, 0}, h_{1, 0}')\in K_{H}^{(s)}$.
\end{proposition}

\begin{lemma}\lb{suff1}  Proposition \ref{Yhk} implies Proposition \ref{chi1 along} \eqref{chi1 along!}.
\end{lemma}
\begin{proof} 
By linearity we may assume that $f^{\dag}$ is of the form \eqref{fdag ass}. 
 Let $X^{\circ}\coloneqq (X_{n}^{\circ}, X_{n+1}^{\circ})\in G_{n, 0}\ts G_{n+1, 0}$ (which depends on $f^{\dag}$) be as in \eqref{Xnu0}, and let $B_{K}^{\dag}=B_{K}^{\dag}(f^{\dag}) \subset B'$ be the image of 
$$ m_{0, s}^{-1} K_{0} X^{\circ} K_{ 0} m_{ 0, s}^{-1} \ts \{1\}
\subset G.$$ 
We have already noted that if $\gamma\notin  B_{K}^{\dag}$, then $I_{\gamma}(f^{\dag}, \chi)=0$.  Assume thus that $\gamma\in B_{K}^{\dag}$, and pick a representative of the form $[\gamma_{0};1]$. 
Proposition \ref{Yhk} and the  discussion preceding it, applied to $X_{\nu}=\eqref{Xnu def}$ and $\lm_{\nu}'=\lm_{\nu}-2s\rho_{\nu}$, show  that the integrand 
$$f^{\star}_{H, \gamma_{0}, \chi} \colon h_{1}\mapsto \chi(h_{1})f^{\star}(h_{1 , 0}^{-1}\gamma_{0} h_{1 ,0}')$$ in \eqref{pO0}  has support contained in $K_{H}^{(s)}$. 
 Thus in order to prove  Proposition \ref{chi1 along} \eqref{chi1 along!} we need to show 
\beq\lb{lolo}
f^{\star}_{H, \gamma_{0}, \chi|K_{H}^{(s)}}=q_{0}^{2sd(n)-\ell_{K_{0}}(w_{0})}.
\eeq

Recall the observation from \eqref{mr KH} that if $h_{1, 0}\in K_{H,0}^{(s)}$, then $m_{0,s}^{-1} h_{1,0}^{-1}m_{ 0,s}\in K_{0}^{\la s+1\ra}\subset K$, and similarly for $h_{1, 0}'^{-1}$. Therefore, the equivalent form \eqref{pOY} of \eqref{pO} and the fact that $\chi_{|K_{H}^{(s)}}=1$ imply \eqref{lolo}.
\end{proof}

In \S\ref{sec: contr} we reduce  Proposition \ref{Yhk} to a simpler statement, to be proved in the remainder of this section. 

\subsection{Contraction} \lb{sec: contr} \lb{sec 72}  \lb{sec:63}
From now until the end of the section, we  lighten the notation by: dropping   all subscripts $`0$';  writing $h$ in place of $h_{1, 0}$, and $h'$ in place of $h_{1, 0}'$; and  writing $m_{s}\in \GL_{n+1}(F) $ in place of $m_{n+1, s}$, whereas we recall that $m_{n, s}=t_{n}^{s}$. 

We extract, from the pair of conditions on $h$, $h'$  in \eqref{pOY}, a single condition on $h$.

Let $e_{n+1, n}={1_{n}\choose 0} \in M_{n+1,n}(F)$ be the matrix with rows $(e_{1}, \ldots, e_{n}, 0)$.
   Denote $\underline{s}\coloneqq  (s,\ldots, s)\in\Z^{n}$ and $\vpi_{n}\coloneqq \vpi^{\underline{1}} =\vpi 1_{n}\in \GL_{n}(F)$, and define the $(n+1)\ts n $  matrices 
\beqq \lb{Y*}
X&\coloneqq  X_{n+1} m_{s}^{-1}  e_{n+1, n} t_{n}^{s} X_{n}^{-1},\\
X^{\circ} &: = X_{n+1}^{\circ}  m_{s}^{-1} e_{n+1, n}t_{n}^{s}  X_{n}^{\circ, -1}
=  \vpi^{\lm_{n+1}' } w_{n+1}
 { \vpi^{-\lm_{n}' -2s\rho_{n} -\underline{s}}   \choose 0}\\
&  = \begin{pmatrix} &&&0\\&&&\vpi^{\lambda_{n}} \\ & 0&\cdots &\\ 0 &\vpi^{\lambda_{2}}&&\\ \vpi^{\lambda_{1}}&&\end{pmatrix},
     \eeqq
   where in the second-last matrix $0\in (F^{n})^{\rm t}$, and 
 $\lm_{i}\coloneqq  (\lm'_{n+1})_{n+2-i} -(\lm'_{n})_{i} -(n+2-2i)s$. Then 
$$\lm_{i+1}-\lm_{i}\geq 2s$$
for all $1\leq i\leq n-1$.  

Let
\beq \lb{hs def}
\overline{h}_{s}&\coloneqq m_{s}^{-1}h^{-1}m_{s} = 
 \twomat {t_{n}^{-s} w_{n}^{-1} h^{-1} w_{n}t_{n}^{s}} {\vpi_{n}^{-s}t_{n}^{-1}(w^{-1}h^{-1}-1_{n}) u}{} 1
  \in {\GL}_{n+1}(F), \\
   h_{s} & \coloneqq t_{n}^{-s}h t_{n}^{s}\in {\GL}_{n}(F).
\eeq

\begin{lemma}\lb{for 735} If \eqref{where is Y} and \eqref{pOY} are satisfied for $\nu=n, n +1$, then 
\beq\lb{pO:1+2}
X &\in  K_{n+1} X^{\circ} K_{n}\\
\overline{h}_{s} X h_{s} &\in  K_{n+1} X^{\circ} K_{n}. \eeq
\end{lemma}
\begin{proof}
Denote by $Y_{\nu}$ the left-hand side of  \eqref{pOY}. Then those equations imply that
\beq\lb{bbo}
Y_{n+1} m_{s}^{-1} e_{n+1, n} t_{n}^{s}Y_{n}^{-1} = \overline{h}_{s} X h_{s}\in K_{n+1} X_{n+1}^{\circ}  K_{n+1} m_{s}^{-1} {  t_{n}^{s}K_{n} X_{n}^{\circ, -1}K_{n} \choose 0}.\eeq
We simplify the right-hand side. First, we have $$K_{n}X_{n}^{\circ,-1}K_{n}= K_{n} w_{n} \vpi^{-\lm_{n}'} K_{n}= K_{n}^{\la 2s\ra} w_{n} \vpi^{-\lm_{n}'} K_{n},$$ where the group $K_{n}^{\la 2s\ra} $ is as in \S~\ref{iw conjsy}  (the second equality can be shown by observing that the quotient $K_{n}^{\la 2s\ra} \bs K_{n}$ is represented by lower-triangular matrices). By the symmetry of $K_{n}^{\la 2s\ra}$, we have
\beqq 
K_{n+1} m_{s}^{-1} {t_{n}^{s} K_{n}   w_{n}\vpi^{-\lm_{n}'} K_{n} \choose 0}
&= K_{n+1} { \vpi_{n}^{-s}t_{n}^{-s} w_{n}t_{n}^{s} K_{n}^{\la 2s \ra}   w_{n}^{-1}\vpi^{-\lm_{n}'} K_{n} \choose 0}\\
&=K_{n+1} {\vpi_{n}^{-s}\cdot \vpi^{-2s\rho_{n} }w_{n} K_{n}^{\la s \ra}  w_{n} \vpi^{2s\rho_{n}}  \vpi^{-\lm_{n}' -2s\rho_{n}}  K_{n} \choose 0}\\
& =K_{n+1} { K_{n}^{\la s \ra} \vpi^{-\lm_{n}' -2s\rho_{n}-\underline{s}}  K_{n} \choose 0}.
\eeqq

Therefore \eqref{bbo} is equivalent to 
$$\overline{h}_{s} X h_{s}\in
K_{n+1}   \vpi^{\lm_{n+1}' } w_{n+1}
K_{n+1} {K_{n}^{\la 2 s \ra}  \vpi^{-\lm_{n}' -2s\rho_{n}-\underline{s}} \choose 0}K_{n}  = K_{n+1} X^{\circ} K_{n},$$
where the identity follows from writing   
$$K_{n+1} { K_{n}^{\la 2 s \ra} \vpi^{-\lm_{n}' -2s\rho_{n}-\underline{s}}  K_{n} \choose 0}
= \lim_{\ep\sto 0} K_{n+1} \twomat  { \vpi^{-\lm_{n}'- 2s\rho_{n}-\underline{s}}}{}{}{\ep}K_{n+1} e_{n+1, n}$$
and applying the multiplication rules of Lemma \ref{iw mult}.
   We conclude that we have
\beqq
X &\in  K_{n+1} X^{\circ} K_{n}\\
\overline{h}_{s} X h_{s} &\in  K_{n+1} X^{\circ} K_{n}, \eeqq
where the first containment follows from the above calculation and \eqref{where is Y}.
\end{proof}

We show that the following solution to the contracted problem \eqref{pO:1+2} implies Proposition \ref{Yhk}.
\begin{proposition} \lb{contrX} Let $K_{\nu}$ be a deeper Iwahori of level $\leq s$.
 Let 
\beq\lb{Xcdef}
X^{\circ}= \begin{pmatrix} &&&0\\&&&\vpi^{\lambda_{n}} \\ & 0&\cdots &\\ 0 &\vpi^{\lambda_{2}}&&\\ \vpi^{\lambda_{1}}&&\end{pmatrix}
\in M_{(n+1)\ts n }(F)\eeq
 with $\lm_{i+1}\geq \lm_{i}+2s$ for all $1\leq i\leq n-1$, and let $X\in K_{n+1} X^{\circ} K_{n}$. 

If $h\in \GL_{n}(F)$ satisfies 
$$\overline{h}_{s} X h_{s} \in  K_{n+1} X^{\circ} K_{n} $$ 
with the notation \eqref{hs def}, 
then $h\in K_{H}^{(s)}$.
\end{proposition}

\begin{lemma} \lb{suff2} Proposition \ref{contrX} implies Proposition \ref{Yhk}.
\end{lemma}
\begin{proof} We revert for a moment to the notation of Proposition \ref{Yhk}. The discussion preceding Proposition \ref{contrX} shows that this proposition implies the conclusion that $h_{1, 0}\in K_{H, 0}^{(s)}$.  Observe now that  $(X_{n}^{\circ,-1}, X_{n+1}^{\circ, -1}; X_{n}^{ -1}, X_{n+1}^{-1}; h_{1,0}', h_{1, 0})$ also satisfies the hypothesis of Proposition \ref{Yhk}. Then the previous argument applied to these data shows that $h_{1, 0}'\in K_{H, 0}^{(s)}$ as well.
\end{proof}

The proof of Proposition \ref{contrX} will occupy the rest of this section.

\subsubsection{Iwahori-invariants from minors} 
We say that a size-$r$ minor $M$ of a matrix $X\in M_{m\ts n}(F_{0})$ is
\begin{itemize}
\item
 \emph{Southwest principal} if it is obtained by deleting all but the last $r$ rows and all but the first $r$ columns of $X$; 
 \item  \emph{quasi-SW-principal} if $r\geq 2$ and $M$ contains the Southwest principal minor of size $r-1$;
 \item \emph{anchored} if   $M$ contains  part of  the last row of $X$.
\end{itemize}

\begin{definition} \lb{SWMC} 
Fix integers $\lm_{1}<\cdots<\lm_{n}$.
We say that $X\in M_{(n+1)\ts n}(F)$ satisfies the \emph{Minor Condition}
 if for every $1\leq r\leq n$,  every $r\ts r$-minor $M_{X}^{(r)}$ of $X$, and the Southwest-principal  $r\ts r$-minor $P_{X}^{(r)}$, we have
\beq \lb{condX} 
v(\det M_{X}^{(r)}) \geq \sum_{i=1}^{r} \lambda_{i}, \qquad  v(\det P_{X}^{(r)}) = \sum_{i=1}^{r} \lambda_{i}
\eeq
We say that $X\in M_{(n+1)\ts n}(F)$ satisfies the \emph{Weak Minor Condition} if \eqref{condX} holds for all 
anchored
 minors. 
\end{definition}
The first example of a matrix satisfying the Minor Condition is $X^{\circ}=\eqref{Xcdef}$. 

\begin{lemma} \lb{minor lemma}
Let  $X, X'\in M_{(n+1)\ts n}(F)$.
\begin{enumerate}
\item \lb{ml1}
If
$$X'\in \Iw_{n+1} X\, \Iw_{n},$$
then $X$ satisfies the Minor Condition if and only if $X'$ does;
\item
if
$$X'\in \twomat{\Iw_{n}}{}{}1 X\, \Iw_{n},$$
then $X$ satisfies the Weak Minor Condition if and only if $X'$ does.
\end{enumerate}
\end{lemma}
\begin{proof} This follows from the Cauchy--Binet formula for minors of products.
\end{proof}

The reader may wish to glance at  the proof of the two parts of our Proposition  in \S\S\,~\ref{sec 74},~\ref{proof 6.1} before looking at the auxiliary lemmas that occupy \S\S~\ref{sec 73},~\ref{sec 75}.

\subsection{First auxiliary lemma}\lb{sec 73}  \lb{sec:64}
We define some variants of the condition $h\in K_{H}^{(s)}$.
\begin{definition}
For $s\geq 1$, we say that a matrix $h\in {\rm GL}_{n}(F)$ is 
\begin{itemize}
\item
\emph{$s$-small} if for all $1\leq i,j\leq n$, 
\beq\lb{ssm}
v(h_{ii})=0 \quad \text{and} \quad v(h_{ij})\geq |j-i|s;\eeq
\item \emph{upper-$s$-small up to row $i$} if there is a decomposition
$$h=h_{-}^{(i)}h_{+}^{(i)}$$ 
 where $h_{+}^{(i)}$ is $s$-small, and $h_{-}^{(i)}$ admits a  block decomposition
\beq \lb{ki bi}
h_{-}^{(i)}= \twomat {\alpha_{}} {}{{*}}{{*}}\eeq
such that $\alpha_{}\in M_{i}(F)$ is lower-triangular with units on the diagonal.
\item \emph{extremely $s$-small} if it is $s$-small and $(wh-1_{n})u=0$.

\end{itemize}
\end{definition}

\begin{remark} The set of extremely $s$-small matrices is a subgroup of $K_{H}^{(s)}$, which in turn is a subgroup of the group of  $s$-small matrices. If $h$ is of the form $h_{-}^{(i-1)}$ and it satisfies \eqref{ssm} for all $j\leq i$,  then $h$ is  upper-$s$-small up to row $i$. (In fact, there is a decomposition $h=h_{-}^{(i)}h_{+}^{(i)}$ with $h_{+}^{(i)}$ differing from the identity only in row $i$.)
\end{remark}

From now until the rest of this section, we write $t$ in place of $t_{n}$. We denote $h^{-w}\coloneqq  w_{n}h^{-1}w_{n}^{-1}$ for $h\in \GL_{n}(F)$,  and we simply denote by $0$ the zero row vector of length $n$. The following remark will often be used in conjunction with Lemma \ref{minor lemma}.
\begin{remark}\lb{s iw} If $h$ is $s$-small, then $t^{-s} ht$ and $t^{-s}h^{-w}t$ belong to $\Iw_{n}$. 
\end{remark}

\begin{lemma}\lb{LS1}  
Let $1\leq i\leq n$, and consider the  equation 
\beq \lb{eqY}
X._{s}h\coloneqq  \twomat {t^{-s}h^{-w}t^{s}} {}{}1 X  t^{-s}ht^{s} = X' ,
\eeq
subject to:
\begin{itemize}
\item $X$, $X' \in M_{(n+1)\ts n}(F)$ satisfy the Weak Minor Condition  of Definition \ref{SWMC};
\item the entries of the  last $i $ rows of $X$  below the lower antidiagonal are zero, that is 
\beq
 \lb{zero antid} 
 v( X_{n+2-i', i'})=\lambda_{i'}, \qquad X_{n+2-i', j}=0 \text{ for all $j>i'\leq i$} \eeq
(where the first equation is a consequence of the second one and \eqref{condX});
\item  $h\in {\rm GL}_{n}(F)$.
\end{itemize}

We have:
\begin{enumerate}
\item \lb{Y1} for given $X$, every solution $(h, X') $ has $h$   upper-$s$-small up to row $i$; 
\item \lb{Y2} if $h$ is of the form $h_{-}^{(i)}$ as in \eqref{ki bi}, then $X'$ also satisfies \eqref{zero antid}.
\item \lb{Y3} for given $X'$, there exists a solution $(h, X)$ with $h$ extremely $s$-small (and in fact upper triangular).
\end{enumerate}
\end{lemma}
\begin{proof}
We proceed by induction on $i$. Write $$X= {A \choose c}, \quad X'={A'\choose c'}$$ with $A, A'\in M_{n\ts n} (F)$

Consider first $i=1$.   The last row of \eqref{eqY} reads
\beq\lb{cj'}
c_{j}'= c_{1}h_{1j}/\vpi^{(j-1)}\eeq
for $j\leq n$. Thus if $X$, $X'$ satisfy \eqref{condX}, then $v(h_{11})=0$ and $v(h_{1 j})\geq (j-1)s$, hence the first statement  is proved and the second one is immediate. On the other hand, substituting   $h_{11}=1-\sum_{k=2}^{n} h_{ik}$, $c_{1}=c_{1}' h_{11}^{-1} $ in \eqref{cj'} gives the integral  linear system  
$$ \sum_{k=2}^{n} (c_{1}'\delta_{jk} + \vpi^{(k-1)s}c_{j}') \, \vpi^{(1-k)s} h_{1k} = c_{j}' $$
in the variables $\vpi^{(1-k)s} h_{1k}$. As the system is invertible, the third statement is proved too.

Now let $i\geq 2$ and suppose the first two statements known for $i-1$. By Remark \ref{s iw} and Lemma \ref{minor lemma}, acting on the right by $h_{+}^{(i-1)}$ preserves the Weak Minor Condition on $X'$; hence  we may and do replace $h$ by $h_{-}^{(i-1)}$ in a decomposition $h=h_{-}^{(i-1)}h_{+}^{(i-1)}$. In other words, we may assume that for $j>i'\leq i-1$, 
$$v(h_{i'i'})=0,  \qquad h_{i',j}=0.$$
The same conditions are then satisfied by $h^{-1}$. 

For $j\geq i$, let  $$M^{n+2-i, j}$$ be the quasi-SW-principal minor of $X'$ of size $i$  whose upper-right corner is $X'_{n+2-i, j}$;
  thus by the induction hypothesis $M^{n+2-i, j}$ has zero entries below the antidiagonal, and its antidiagonal entries  
  have valuations (in order, starting from the SW corner)
  $$\lm_{1}, \ldots, \lm_{i-1}, v(X'_{n+2-i, j}).$$
In particular, 
$$v(\det M^{n+2-i, j}) =\sum_{i'=1}^{i} \lambda_{i'} -\lambda_{i} + v( X'_{n+2-i, j}).$$
Hence the Minor Condition \eqref{condX} implies 
\beq \lb{condX1} -\lambda_{i}+v( X'_{n+2-i, i})= 0, \qquad -\lambda_{i}+ v( X'_{n+2-i, j})\geq  0 \text{ for all $j>i$}. \eeq
As $A'=t^{-s}h^{-w} t^{s}A t^{-s}ht^{s}$, 
 we have for all $1\leq j\leq n$:
\beq\lb{Xij}
\vpi^{-\lambda_{i}}X'_{n+2-i, j} 
&= \vpi^{-\lambda_{i}} \sum_{k=1}^{n-1} (h^{-w})_{n+2-i,n+1-k} \vpi^{(k+1-i)s} X_{n+1-k, k-1} h_{k+1, j} \vpi^{(k+1-j)s} \\
&= \sum_{k=1}^{i-1} h^{-1}_{i-1,k} \vpi^{(k+1-i)s} \vpi^{-\lambda_{i}} X_{n+1-k, k+1}  h_{k+1, j}\vpi^{(k+1-j)s} ,\eeq
by our assumptions on $h$. Moreover, for $j\geq i$ by induction hypothesis $h_{k+1, j}=0$ for all $k<i-1$, hence 
\beq\lb{Y'+}
\vpi^{-\lambda_{i}} X'_{n+2-i, j}= h^{-1}_{i-1,i-1}  h_{i, j}\vpi^{(i-j)s} \vpi^{-\lambda_{i}}X_{n+2-i, i}.
\eeq
Since $h^{-1}_{i-1, i-1} $ and $\vpi^{-\lambda_{i}}X_{n+2-i, i} $ are   units, condition \eqref{condX1} is equivalent to 
$$v(h_{i, i})=0, \qquad v( h_{i,j})\geq (j-i)s$$
for all $j>i$, establishing the first statement.  If $h$ is of the form $h_{-}^{(i)}$, the second statement is  immediate from  \eqref{Y'+}. 

Consider now the third statement.  After replacing $X'$ by $X'._{s}(h')^{-1}$ where $h'$ is as given by this statement for $i-1$, we may and do assume that $X'$ satisfies \eqref{zero antid} for $i'<i$.  We seek  $h$ extremely $s$-small, upper-triangular and differing from the identity only in row $i$; hence in particular $h$ takes the form $h^{(i-1)}_{-}$, 
 and by the second statement for $i-1$, we only need to find a solution to \eqref{Y'+} in $h$ (with the further simplification $h_{i-1, i-1}=1$). 

We set $h_{i, i}= 1-\sum_{k\neq i} h_{ik}$, necessary for  $h$ to be extremely  $s$-small,  and substitute in 
$ X_{n+2-i, i} =h_{i, i}^{-1}X'_{n+2-i, i} $. We find the  linear system 
$$\sum_{k= i+1}^{n} (\vpi^{-\lambda_{i}} X'_{n+2-i, i} \delta_{jk} + \vpi^{(k-i)s} \vpi^{-\lambda_{i}} X'_{n+2-i, j})\, h_{i k}\vpi^{(i-k)s} = \vpi^{-\lambda_{i} } X'_{n+2-i, j} $$
in the variables $\vpi^{(i-k)s} h_{ik}$ for  $k\geq i+1$. By our reductions, $-\lambda_{i}+v(X'_{n+2-i, i})=0$ and $-\lambda_{i}+v(X'_{n+2-i, j})\geq 0$, hence  the system is  integral and invertible; its solvability  implies  the third statement.
\end{proof}

\subsection{Plus-regularity of support: proof of Proposition \ref{chi1 along} \eqref{+reg spt}} \label{cor: reg sup} 
\lb{sec 74} \lb{sec:65}
We prove part \eqref{+reg spt} of Proposition \ref{chi1 along}. It follows from Lemma \ref{f*} (via Lemma \ref{match split})
 that  the Hecke measure $\ff'= \eqref{fp def b}$ matches a measure $f\in \sH(G, L)$ that is invariant under $K_{0}^{m_{s}}$. We now turn to proving that $f'$ is supported in the plus-regular locus $G'_{\reg^+}$ (thus $f$ is also regularly supported).

Recall that we defined
in \eqref{disc rs} the function
$D^\pm$ on $\wt G'$,
 by pulling back the corresponding function $D^\pm$ on the symmetric space $S$. Now that the place $v$ is split in the quadratic extension, we identify $S$ with $G_{n+1}$. Tracking the process of contracting the test measure,   it suffices to show that the measure $f_{n}^{\star}\star f_{n+1}^{\star}$ on $G_{n+1}=S$  has plus-regular support, where $f_n^\star$ and $f_{n+1}^\star$ are as in Lemma \ref{f*}.

Now we note that,
for  $\gamma=\twomat  {A}{b}{c}{d}\in G_{n+1}=S$, we have
$$
D^+\left(\gamma\right)= \pm\det(c, cA,\cdots, cA^{n-1})^{\rm t},$$
and  the quasi-invariance property that for $h\in G_n$,
$$ D^+\left(h^{-1}\gamma h\right)= \det(h) D^+\left(\gamma\right).
$$
Then by definition, $\gamma\in G_{n+1}$ is  plus-regular  if and only if $D^+\left(\gamma\right)\neq0$, or equivalently the vectors $c, cA,\cdots, cA^{n-1}$ form a basis.

We observe that the plus-regularity depends only the first $n$ columns, so that we may talk about the plus-regularity of an element ${A \choose c}\in M_{(n+1)\times n}(F)$. Therefore by  Lemma \ref{for 735} (together with the discussion preceding Proposition \ref{Yhk}),
it suffices to show that the set 
\beq\lb{eq sup in S}
m_s  K_{n+1} X^{\circ} K_{n}t_s^{-1} \qquad \subset  M_{(n+1)\times n}(F)
\eeq
is contained in the plus-regular locus. Here, $X^{\circ} $ is as in \eqref{Xcdef}.
 
 By Lemma \ref{minor lemma}\eqref{ml1}, any element in the set $K_{n+1} X^{\circ} K_{n}$ satisfies the  Minor Condition  (Definition \ref{SWMC}). It thus suffices to show that, if $X$ satisfies the Weak Minor Condition, then the element
 $$\twomat {1} {u}{}{1}\twomat {w t^{s}} {}{}{1} X  t^{-s}\in M_{(n+1)\times n} (F) $$ is plus-regular.
 Set
  $$\wt X\coloneqq \twomat {wt^{s}} {}{}{1} X  t^{-s}\in M_{(n+1)\times n} (F) .
  $$
  Then we claim that $\twomat {1} {u}{}{1} \wt X$ is plus-regular if and only if $ \wt X$ is. To see this, we write $\wt X={\wt A \choose \wt c}$ so that  $\smalltwomat {1} {u}{}{1} \wt X={\wt A+ u \wt c \choose \wt c}$. We see inductively that the span of $\wt c, \wt c(\wt A+ u \wt c),  \ldots,\wt c(\wt A+ u \wt c)^{i-1} $
 is equal to the span of   $\wt c, \wt c\wt A , \ldots,\wt c(\wt A)^{i-1} $. The claim follows.

  It remains to show that $ \wt X$ is plus-regular. By Lemma \ref{LS1}~(3) (applied to the case $i=n$), there exists $h\in \GL_n(F)$ 
such that, if we set
$$
h^{-1} \wt X h= \wt X'
$$
where $\wt X'=\twomat {wt^{s}} {}{}1 X '  t^{-s}\in M_{(n+1)\times n}$, then  the entries of $X'$  below the lower antidiagonal are all zero. It therefore suffices to show that $\wt X'$ is plus-regular. We note that 
$$ \wt X'= \begin{pmatrix}\ast &a_2&0&0\\\vdots&\vdots&\ddots & 0 \\ \ast& \ast&\cdots &a_n\\ \ast &\ast&\cdots&\ast\\ a_1 &0&\cdots&0\end{pmatrix}={\wt A' \choose \wt c'} ,
$$
where $a_1,a_2,\ldots, a_n$ are all non-zero. For $1\leq i \leq n$, denote by  $e_i\in F_{0}^{n}$ the standard basis vector,   and by $V_{i}\subset F_{0}^{n}$  the subspace spanned by $e_1,\ldots, e_{i}$.  In particular, $\wt c'=  a_1  e_1$.  Then by induction we see that 
$$
e_1^{\rm t}  (\wt A')^{i-1}\equiv  a_i e_{i}^{\rm t} \mod V_{i-1}^{\rm t}.
$$
It follows easily that the subspace spanned by $e_1^{\rm t},\ldots, e_1^{\rm t} (\wt A')^{i-1}$ is exactly $V_i^{\rm t}$,  for all $1\leq i \leq n$. The desired assertion follows.

\subsection{Second auxiliary lemma} \lb{sec 75}  \lb{sec:66}
We continue with another lemma towards the proof of part \eqref{chi1 along!} of Proposition \ref{chi1 along}.

\begin{definition}
We say that a lower-triangular matrix $h\in \GL_{n}(F)$ with units on the diagonal  is  \emph{lower-$s$-small from  column  $j$} if 
$$v(h_{ij'})\geq (i-j')s \text{ for all $i>j'\geq j$}.$$
This is equivalent to the existence of  a decomposition
$$h={}^{(j)}h_{-}  \cdot {}^{(j)} h_{--}$$ 
where  ${}^{(j)}h_{--} $ is lower-triangular and $s$-small, and 
\beq \lb{ki bi 2}
{}^{(j)}h_{-} = \twomat {\alpha_{-}} {}{{*}}{{1_{n+1-j}}} \qquad
\eeq
with    $\alpha_{-}\in M_{j-1}(F)$  lower-triangular with units on the diagonal.

\end{definition}

\begin{remark}\lb{low inv}
For a lower-triangular matrix $h$ with units on the diagonal:
\begin{itemize}
\item $h$ is  lower-$s$-small from  column $j$  if and only if $h^{-1}$ is;
\item $h$ is  lower-$s$-small from  column $1$ if and only if it is $s$-small. 
\end{itemize}
\end{remark}

\begin{lemma}\lb{LS2}
Let $1\leq j\leq n$, and consider the  equation 
\beq \lb{eqY2}
X._{s}h= \twomat {t^{-s}h^{-w}t^{s}} {}{}1 X  {t^{-s}ht^{s} }= X' ,
\eeq
subject to:
\begin{itemize}
\item $X$, $X' \in M_{(n+1)\ts n}(F)$ satisfy the Weak Minor Condition of Definition \ref{SWMC};
\item the entries  of $X$, $X'$ below the lower antidiagonal are zero, that is 
\beqq
 v( X_{n+2-i, i})=\lambda_{i}, \qquad X_{n+2-i, j'}=0 \text{ for all $j'> i$}\eeqq
(where the first equation is a consequence of the second one and \eqref{condX}), and similarly for~$X'$;
\item the entries of the last $n-j$ columns of $X$ above  the lower antidiagonal are zero, that is, 
\beq \lb{zero above antid} 
\qquad X_{n+2-i, j'}=0 \text{ for all $i> j'\geq j+1$;} \eeq
\item $h\in {\rm GL}_{n}(F)$ is lower-triangular with units on the diagonal.
\end{itemize}
We have:
\begin{enumerate}
\item \lb{Y1} for given $X$, every solution $(h, X') $ has $h$  lower-$s$-small from column $j$; 
\item \lb{Y2} if $h$ is of the form ${}^{(j)}h_{-} $ as in \eqref{ki bi 2}, then $X'$ also satisfies \eqref{zero above antid};
\item \lb{Y3} for given $X'$, there exists a solution $(h, X)$ with $h$ extremely $s$-small.
\end{enumerate}
\end{lemma}
\begin{proof} 
We prove this by decreasing induction on $j$, the case $j=n$ being trivial. Thus let $j\leq n-1$ and assume the statements proved for $j+1$. 

After replacing $h$ by ${}^{(j+1)}h_{-}$   as in the decomposition \eqref{ki bi 2}, that is acting by $._{s}{}^{(j+1)}h_{--} $ on both sides  of \eqref{eqY2}, by the induction hypothesis we are led to a  situation that is equivalent for the purposes of the first two statements. Hence we may and do assume that $h$ has the form ${}^{(j+1)}h_{-} $. 
For $i\geq j$, let $$N^{n+1-i, j+1}$$ be the quasi-SW-principal minor of $X'$  of size $j+1$ whose upper-right corner is $X'_{n+1-i, j+1}$;
thus  the matrix $N^{n+1-i, j+1}$ has vanishing entries below the antidiagonal, 
 and its antidiagonal entries (in order, starting from the SW corner)  have valuations
$$\lm_{1}, \ldots, \lm_{j}, v(X'_{n+1-i, j+1}).$$
In particular, 
$$v(\det N^{n+1-i, j+1}) =\sum_{j'=0}^{j+1} \lambda_{j'}  - \lambda_{j+1} + v( X'_{n+1-i, j+1}).$$
Hence \eqref{condX} implies 
\beq \lb{condX2} -\lambda_{j+1} +v( X'_{n+1-i, j+1})\geq 0. \eeq
The same condition holds for $X$ by assumption.

On the other hand, we have
\beq\lb{Xij2}
\vpi^{-\lambda_{j+1}} X'_{n+1-i, j+1}
& = \vpi^{-\lambda_{j+1}}\sum_{1\leq k, l\leq n} (h^{-w})_{n+1-i,n+1-k}\vpi^{(k-i)s}X_{n+1-k, l+1} h_{l+1, j+1} \vpi^{(l-j)s} \\
& = 
\sum_{k=j}^{i} h^{-1}_{ik}\vpi^{(k-i)s} 
\vpi^{-\lambda_{j+1}}
 X_{n+1-k, j+1}\eeq
where we have used our assumptions on $h$ and $X$. All terms are integral except possibly the last one, whose valuation is $v(h^{-1}_{i j}) -(i-j)s$. That this should be non-negative, for all $i>j$,  is equivalent to $h$ being lower-$s$-small from column $j$, proving the first statement. 

If moreover $h$ has the form ${}^{(j)}h_{-}$, then in \eqref{Xij2} all terms are zero unless $i=j$, in which case we  only have the term corresponding to $i=j=k$, giving $X'_{n+1-j, j+1}=X_{n+1-j, j+1}=0$. This proves the second statement.

For the third statement, we seek an extremely $s$-small matrix $h$ that differs from the identity only in column $j$. Then $h^{-1}$ satisfies the same conditions, $h$ is of the form ${}^{(j-1)}h_{-}$, and we need it to satisfy \eqref{Xij2} (for some $X$), in whose last sum  only the terms indexed by {$k=j$}
may be nonzero. Substituting
$$h_{jj}^{-1} =  1-\sum_{k= j+1}^n h_{kj}^{-1}$$
into \eqref{Xij2} for $i=j$, we find
\beq\lb{sub 665}
X_{n+1-j, j+1}=  (1-\sum_{k = j+1}^n h^{-1}_{kj} )^{-1}X'_{n+1-j,  j+1}.
\eeq
Substituting \eqref{sub 665} into \eqref{Xij2}, 
 we find 
 $$
 (1-\sum_{k\geq j+1} h^{-1}_{kj} )
\vpi^{-\lambda_{j+1}}
X'_{n+1-i, j+1}
=
h^{-1}_{ij} \vpi^{(j-i)s} 
\vpi^{-\lambda_{j+1}}
 X'_{n+1-j, j+1} .$$
This is an invertible integral linear system of  equations
$$
\sum_{k=j+1}^{n}  
\left(
\vpi^{-\lambda_{j+1}}
 X'_{n+1-j, j+1}\delta_{ki}+
  \vpi^{(k-j)s}
  \vpi^{-\lambda_{j+1}} X'_{n+1-i, j+1} \right) \vpi^{(j-k)s}h_{kj}^{-1}=
  \vpi^{-\lambda_{j+1}}
X'_{n+1-i, j+1}
$$
indexed by $j+1\leq i\leq n$,
in the variables $\left(\vpi^{(j-k)s}h_{kj}^{-1}\right)_{ j+1\leq k\leq n}$. The   solvability of the system  implies  our third statement.
\end{proof}

\subsection{Proof of Propositions \ref{chi1 along} \eqref{chi1 along!}, \ref{Yhk}, and \ref{contrX}} \lb{proof 6.1} \lb{sec 76}
 \lb{sec:67}
By Lemmas \ref{suff1}, \ref{suff2}, it suffices to prove Proposition \ref{contrX}. Thus we need to show that for  $X, Y \in K_{n+1} X^{\circ}K_{n}$, all  the solutions in $h$ to the equation
\beq 
Y=\overline{h}_{s} X h_{s}\eeq
have $h\in K_{H}^{(s)}$. 
By Lemma \ref{minor lemma}, both $X$ and $Y$ satisfy the Minor Condition.

Write  $X={A\choose c} \in M_{(n+1)\ts n}(F)$, with $c\in F^{n, \rm t}$. Then 
\beq\lb{Xdef}
Y&=\overline{h}_{s} X h_{s} = X'+ X'', \\
X'&\coloneqq   X._{s} h
=
 {{t^{-s}h^{-w}t^{s}}  A t^{-s} h t^{s} 
\choose
c t^{-s} ht^{s}},
\\
X''&\coloneqq  X.._{s}h\coloneqq 
{\vpi^{-s}t^{-s}(wh^{-1}-1_{n})u \choose 0} c t^{-s}h t^{s},
\eeq
where the notation $X._{s}h$ is as in Lemmas \ref{LS1}, \ref{LS2}.

Note that $X''=Y-X'$ is a rank-$1$ matrix whose rows are all multiples of  row $n+1$ of $X'$ (and whose last row is zero), so that the determinants of any pair of corresponding anchored minors of $Y$, $X'$ are equal. In particular, $X'$ also satisfies the Weak Minor Condition of Definition \ref{SWMC}.

We  proceed in several steps to show  that  $h\in K_{H}^{(s)}$.
\begin{enumerate}
\item 
By applying first  Lemma \ref{LS1}\,(3) for $i=n$, then  Lemma \ref{LS2}\,(3) for $j=1$, we find an extremely $s$-small $h'$ such that, first, $X.._{s}h'=0$ (which is automatic by the extreme smallness of $h'$) and, second,  $X._{s}h' = \overline{h}'_{s} Xh_{s}'$ 
has zero entries outside of the lower antidiagonal and of column $1$. Hence, up to changing variables by the action of such an element $h'$, we may assume $X$ satisfies these vanishing conditions.
\item Apply Lemma \ref{LS1}\,(1) to the equation 
\beq\lb{xxs}
 X'=X._{s}h,\eeq to deduce that $h$ is upper-$s$-small, $h=h_{-}h_{+}$ with $h_{-}$ lower-triangular with units on the diagonal and $h_{+}$  $s$-small.
\item Act on \eqref{xxs} by $._{s}h_{+}^{-1}$; by Remark \ref{s iw} and Lemma \ref{minor lemma}\,(2)  this preserves the Weak Minor Condition. We can then apply  Lemma \ref{LS2}\,(1) (with $j=1$)  to the resulting equation, to conclude that $h_{-}$ and $h$ are $s$-small. 
\item By Remark \ref{s iw} and Lemma \ref{minor lemma}\,(1), we deduce that $X'=X._{s}h$ satisfies the full Minor Condition; in particular, all entries of $X'$ have valuation no less than $v(\lambda_{1})$. Since this also holds for  the entries of $Y$, it must  hold for the entires of $X''$ too. As $\lambda_{1}^{-1}c t^{-s}h t^{s}$ is integral with first entry a unit,  the condition on $X''$ is satisfied if and only if $\vpi^{-s}t^{-s}(wh^{-1}-1_{n})u\in \sO^{n}$; that is, $h\in K_{H}^{(s)}$. 
\end{enumerate}
The proof of Propositions \ref{contrX},  \ref{Yhk} and \ref{chi1 along} (whose part \eqref{+reg spt} was proved in \S~\ref{sec 74}) is now complete.

\section{The $p$-adic relative-trace formula and $p$-adic $L$-functions} \lb{sec: prtf}
This section is dedicated to the construction of the $p$-adic $L$-function  of Theorem \ref{thm B} and the related RTF. In \S~\ref{sec 61} we give the statements. In \S~\ref{sec 62} we give the proofs: similarly to what done in  \S~\ref{sec rat}, we construct the $p$-adic relative-trace distribution from its geometric expansion, then we extract from it the $p$-adic $L$-function and deduce the spectral expansion. In \S~\ref{sec 63}, we give an RTF for the derivative of the distribution.

Throughout this section, we fix a rational prime $p$.
\subsection{Statements} \lb{sec 61}
Recall that we denote $\Gamma=\Gamma_{F_{0}}\coloneqq  F_{0}^{\ts}\bs \A^{\infty,\ts} /\widehat{\sO}_{F_{0}}^{p,\ts}$, and  
$\sY\coloneqq \Spec \Z_{p}\llb \Gamma_{F_{0}} \rrb\ot\Q_{p}$
We say that $\Pi\in \sC(\G')^{\rm her}_{\Q_{p}}$ is \emph{ordinary} if for all $v\vert p$, the representation $\Pi_{v}$ is ordinary in the sense of Definition  \ref{ord v}. The ordinary representations form an ind-subscheme 
$$\sC(\G')^{\rm her, ord}\subset \sC(\G')_{\Q_{p}}^{\rm her}.$$
 For $K_{p}=\prod_{v}K_{v}$, we let $\sC(\G')^{\rm her, ord}_{K_{p}}$ be the subscheme of those $\Pi$ which are $K_{v}$-ordinary  for all $v\vert p$. 
 
\subsubsection{$p$-adic $L$-function}
The following is Theorem \ref{thm B} from the introduction
\begin{theorem}\lb{plf thm} Let $L$ be a finite extension of $\Q_{p}$, and let $\Pi\in  \sC(\G')^{\rm her, ord}$.

Assume that for each place $v\vert p$ of $F_{0}$,  $v$ splits in $F$ or $\Pi_{v}$ is unramified.
Then there exists a unique function 
$$\sL_{p}(\RM_{\Pi})\in \sO(\sY_{L}) $$
whose restriction to $Y(p^{\infty})_{L}$ satisfies 
\beqq\lb{eq interp Lp} 
\sL_{p}(\RM_{\Pi})(\chi) = e_{p}(\RM_{\Pi\ot \chi})\,  \sL(\RM_{\Pi} )(\chi)
\eeqq
where $\sL(\RM_{\Pi})$ is  the function in Theorem \ref{rat L}, and  $e_{p}(M_{ \Pi\ot \chi})\coloneqq  \prod_{v\vert p} e(\Pi_{v}, \chi_{v})$ for the  factors of \eqref{e v ord}.
\end{theorem}

\subsubsection{Generalized Radon measures} We make the first of two preparations which will be relevant to the  $p$-adic relative-trace formula.

Recall that a \emph{Banach ring} is a topological ring equipped with a norm $|\cdot|$ for which it is complete; the relevant examples for us are the finite extensions of $\Q_{p}$ (with the $p$-adic norm) and $\sO(\sY)$ (with the Gauss norm).
\begin{definition} Let $X$ be a set and let $R$ be a Banach ring. A generalized bounded Radon measure\footnote{When $X$ is a topological space and  $L^{1, \infty}(X, \mu)$ contains $C_{c}(X)$, the functional $\mu$ is a (bounded) Radon measure in the sense of Bourbaki. } with values in $R$ is a pair $(\mu, L^{1, \infty}(X, \mu))$, where 
\begin{itemize}
\item  $L^{1, \infty}(X, \mu)\subset L^{\infty}(X) $ is a closed subspace of the $R$-Banach space of bounded $R$-valued function on $X$; 
\item $\mu\colon  L^{1, \infty}(X, \mu) \to R$ is a bounded $R$-linear functional. 
\end{itemize}
We will usually denote such measures simply by $\mu$, and for $\Phi\in L^{1, \infty}(X, \mu)$, we will  use the notation
$$ \int_{X} \Phi(x) \, d\mu(x) \coloneqq \mu(\Phi).$$
\end{definition}
When $R'\supset R$ is an extension of Banach rings, an $R$-valued   generalized bounded Radon measure 
$\mu$ gives rise to an  $R'$-valued  generalized bounded Radon measure by extension of scalars, which we will still denote by $\mu$.
We say that a function $\Phi\in L^{\infty}(X)$ is $\mu$-integrable if it belongs to $L^{1, \infty}(X, \mu)$.

\subsubsection{Local distributions at $p$ }
 Let $K_{p}=\prod_{v\vert p} K_{v} \subset \G'(F_{0, p})$
be a compact open subgroup that is \emph{convenient} in the sense that each $K_{v}$ is (as defined in \S~\ref{def sph ch}).
 We will say that $K_{p}$ is a \emph{conjugate-symmetric deeper Iwahori} (CSDI) if  each $v\vert p$ splits in $F$ and each $K_{v}$ is a CSDI (as defined in \S~\ref{IW sym}).

 For $\chi\in Y(p^{\infty})$ and $f_{p}^{\dag}=\ot_{v\vert p} f_{v}^{\dag}\in \sH_{p}^{\dag}=\bigotimes_{v\vert p }\sH_{v}^{\dag}$ that is sufficiently positive for $\Pi_{p}$, $\chi_{p}$, and $K_{p}$ (in the obvious sense derived from \S~\ref{sec suff pos} for each $v\vert p$), and for $\Pi_{p}$ a tempered irreducible representation of $G'_{p}$ and  $\gamma \in B'_{p}$, we define
\beqq
I^{\dag}_{\Pi_{p}, K_{p}}(f_{p}^{\dag}, \chi_{p}) 
&\coloneqq  \prod_{v\vert p }I_{\Pi_{v}, K_{v}}^{\dag}(f_{v}^{\dag} , \chi_{v}) , \qquad 
 I^{\dagger}_{\gamma,p,K_{p}} (f_{p}^{\dag}, \chi_{p}) \coloneqq  \prod_{v\vert p} I^{\dagger}_{\gamma,v,K_{v}} (f_{v}^{\dag}, \chi_{v}),
\eeqq
where the last factors are as in \eqref{def orb dag}; we impose the restriction  that $\gamma \in B'_{\rs, p}$ unless $F/F_0$ splits above $p$ and $K_{p}$ is a CSDI.

 \subsubsection{$p$-adic relative-trace formula} For $K_{p}$ as above,  recall the Hecke subspace
$$  \sH(\G'(\A^{p}))^{\circ}_{K_{p}, \rs, \rm qc}\subset  \sH(\G'(\A^{p}))_{\rs}^{\circ}
$$
of \S~\ref{rha}.
 We denote $U_{t_{p}}=\ot_{v\vert p} U_{t_{v}}$.

We first define the local terms (away from $p$) of the $p$-adic RTF, and the global spectral terms.

\begin{definition}  \lb{def prtf}
\begin{enumerate}
\item
For each finite place $v\nmid p$ of $F_{0}$ and for $v=\infty$, for  each $\gamma\in B'_{ v}$, and for  each tempered irreducible  representation $\Pi_{v}$ of $G_{v}'$ over $L$, let
\beqq
\sI_{\Pi_{v}} \colon \sH(G_{v}', L)^{\circ}&\to \sO(\sY_{L}), \\ 
\sI^{}_{\gamma, v}\colon \sH(G'_{v}, L)^{\circ} &\to \sO(\sY_{L[\sqrt{-1}]})
\eeqq
be the distributions obtained by pulling back the corresponding ones  of Proposition \ref{rat rtf} \eqref{rat sph ch}, 
\eqref{rat orb int reg}
via the restriction maps  $\sY \ni \chi\mapsto \chi_{v}\in Y_{v}(1)_{\Q_{p}}$. (If $v=\infty$, $Y_{v}(1)\coloneqq \Spec {\Q}$.)
  \item \lb{part IPi p} (Assuming Theorem \ref{plf thm}.) Let $K_{p}=\prod_{v\vert p} K_{v} \subset \G'(F_{0, p})$
be a
convenient
subgroup,
 and  let $L$ be a finite extension of~$\Q_{p}$. For each representation $\Pi\in \sC(\G')^{\rm her, \ord}(L)$, define a distribution
 $$\sI_{\Pi, K_{p}}^{}\coloneqq {1\over 4}c_{K_{p}}(\Pi)\, \sL_{p}(\RM_{\Pi})\, 
 \prod_{v\nmid p} \sI_{\Pi_{v}} 
 \colon \sH(\G'(\A^{p}), L)^{\circ} \to \sO(\sY_{L}), $$
 where  the constant $c_{K_{p}}(\Pi)\coloneqq \prod_{v} c_{K_{v}}(\Pi_{v})$ for the factors of   \eqref{e v ord}.
\end{enumerate}\end{definition}

\begin{theorem}[$p$-adic analytic RTF]
\lb{p-adic rtf}
 Let $K_{p}=\prod_{v\vert p} K_{v} \subset \G'(F_{0, p})$
be a
convenient
subgroup,
 and  let $L$ be a finite extension of $\Q_{p}$.

There exist: 
\begin{enumerate}
\item \lb{DR Lgamma}
For each $\gamma \in \RB'(F_{0})$, a unique bounded-by-$1$ $p$-adic $L$-function
$$L_{p,\gamma}\in \Z_{p}\llbracket \Gamma_{F_{0}}\rrbracket \subset \sO(\sY)$$
whose restriction to $Y(p^{\infty})$ equals $L_{\gamma}^{(p)}\coloneqq  L_{\gamma}/\prod_{v\vert p}L_{\gamma, v}$ (where $L_{\gamma}$ is as in Proposition \ref{rat rtf}\eqref{rat orb int 2 reg a}).
\item \lb{part glob orb}
An orbital-integral function
\beq\lb{orb int fam}
\sI^{p}
\colon  \RB'(F_{0}) \ts  \sH(\G'(\A^{p}), L)^{\circ} &\to  \sO(\sY_{L})\\
 \eeq
 defined by 
 $$(\gamma, f'^{p})\mapsto\sI^{p}_{\gamma}(f'^{p})\coloneqq  \kappa(\one_{\infty})^{-1} L_{p,\gamma}\, \prod_{v\nmid p} \sI_{\gamma,v},$$
 which is bounded in the variable $\gamma$.

 \item \lb{part mu} 
 \begin{enumerate}
  \item \lb{part mu a} for every $\chi_{p}\in Y_{p}(p^{\infty})$,  a $\Q_{p}(\chi_p)$-valued generalized bounded Radon measure $I^{\ord}_{\gamma, p,K_{p}}(\chi_{p})$  on  $\RB'_{\rs}(F_{0})^{\circ}\coloneqq \RB'_{\rs}(F_{0}) \cap B_{\rs, \infty}'^\circ$,
  defined by the limit of weighted samplings
\beq\lb{muK def}
\int_{\RB_{\rs}'(F_{0})^{\circ}} \Phi(\gamma)\, d I^{\ord}_{\gamma, p, K_{p}}(\chi_{p}) \coloneqq  
\lim_{N\sto\infty} \sum_{\gamma \in \RB'_{\rs}(F_{0})^{\circ} }
 I^{\dagger}_{\gamma,p,K_{p}} (U_{t_{p}}^{N!}, \chi_{p}) \cdot \Phi(\gamma)
\eeq
on the space of bounded functions $\Phi\in L^{\infty}(\RB_{\rs}'(F_{0})^{\circ})$ for which the sums over $\gamma$ converge  and the limit converges; 
 \item
 if $F/F_{0}$ splits above $p$ and $K_{p}$ is a CSDI, a $\Q_{p}$-valued generalized bounded Radon measure $I^{\ord}_{\gamma, p,K_{p}}$ on  $\RB'(F_{0})$ (whose restriction to   $\RB_{\rs}'(F_{0})^{\circ}$ coincides with the measure $I^{\ord}_{\gamma, p,K_{p}}(\chi_{p})$ of  \eqref{part mu a} for every $\chi_{p}$),
defined by
\beq\lb{mu def}
\int_{\RB(F_{0})^{}} \Phi(\gamma)\, dI^{\ord}_{\gamma,p, K_{p}} =
k_{p} \cdot  \lim_{N\sto \infty} 
\sum_{\gamma \in \RB'_{ }(F_{0})\cap B_{p, N}^{\dag}}
  \Phi(\gamma)
  \eeq
  on the space of bounded functions $\Phi\in L^{\infty}(\RB'(F_{0}))$ for which the sums over $\gamma$ converge  and the limit converges; 
here, $B_{p, N}^{\dag}=\prod_{v\vert p} B_{K_{v}}^{\dag}(U_{t_{v}, K_{v}}^{N!})$ and $k_{p}=\prod_{v\vert p} k_{v}$,
with the factors as in \eqref{const k}.
 \end{enumerate}

 \item \lb{part p exps} 
  A unique distribution
$$\sI^{}_{K_{p}}\colon   \sH(\G'(\A^{p}), L)^{\circ}_{K_{p}, \rm qc} \to \sO(\sY_{L}) $$
satisfying:
\begin{enumerate}
\item \lb{exp spec}
there is a spectral expansion
$$\sI^{}_{K_{p}}= \sum_{\Pi\in \sC(\G')^{\rm her, ord}_{K_{p}}} \sI_{\Pi, K_{p}};$$
\item \lb{exp a} 
for each  $f'^{p}\in \sH(\G'(\A^{p}),  L)^{\circ}_{K_{p}, \rs, \rm qc}$ and each finite-order $\chi\in \sY_{L}$, the function $I^{p}_{-}(f'^{p}, \chi^p)$ on  $\RB'_{\rs}(F_{0})^{\circ}$ is integrable for $I^{\ord}_{\gamma, p, K_{p}}(\chi_{p})$, and we have the geometric expansion in $L(\chi)$
$$\sI^{}_{K_{p}}(f'^{p}, \chi)=
\int_{\RB'_{\rs}(F_{0})^{\circ}} I^{p}_{\gamma}(f'^{p}, \chi^p) \, dI^{\ord}_{\gamma, p, K_{p}}(\chi_{p});$$
\item  \lb{exp b}
if $F/F_{0}$  splits above $p$ and $K_{p}$ is a
CSDI,  then  for each  $f'^{p}\in \sH(\G'(\A^{p}),  L)^{\circ}_{K_{p}, \rm qc}$ the function
$ \sI^{p}_{-}(f'^{p})$  on  $\RB'(F_{0})$ is integrable for $I^{\ord}_{\gamma, p, K_{p}}$, and we have the geometric expansion in $\sO(\sY_{L})$
$$\sI^{}_{K_{p}}(f'^{p})=
\int_{\RB'(F_{0})} \sI^{p}_{\gamma}(f'^{p}) \, dI^{\ord}_{\gamma, p, K_{p}}.$$
\end{enumerate}
\end{enumerate}
\end{theorem}

\subsection{Proofs}\lb{sec 62}
We will prove Theorem \ref{p-adic rtf} and, as an interlude, Theorem \ref{plf thm}. 
\subsubsection{Boundedness of   local orbital integrals} We consider the local distributions $I_{\gamma}$ of Proposition  \ref{rat rtf} 
\eqref{rat orb int reg}.

\begin{lemma} \lb{bdd loc orb}
Let $v\nmid p\infty $ or, respectively, $v=\infty$,  and let $f'_{v}\in \sH(G'_{v}, L)^{\circ}$.
 There is a constant $c(f'_{v})\in \Q^{\ts}$ such that for every $\gamma\in B'_{v}$ (respectively, for every  $\gamma \in B_{v}^{\circ}$), the Laurent polynomial
$$I_{\gamma}(f_{v}')\in \sO(Y_{v}(1))_{L}\cong L[T^{\pm 1}]$$
(respectively the number  $I_{\gamma}(f_{\infty}', \one)\in L$) belongs to $c(f'_{v})\sO_{L}[T^{\pm 1}]$ (respectively $c(f'_{v})\sO_{L}$). Moreover, for all but finitely many $v$, if $f_{v}'$ is the unit Hecke measure then we may take $c(f'_{v})=1$. 
\end{lemma}
\begin{proof}
By the definitions in \eqref{eq orb on B} and Lemma \ref{kappa i}, it suffices to consider $\hI^{}_{\gamma'}(f_{v}', \chi_{v})$  instead of $\hI^{}_{\gamma}(f_{v}', \chi_{v})$, for any $\gamma'$  in the unique plus-regular orbit above $\gamma$. 

First consider the case of $v=\infty$, for which we may assume that $f'_{\infty} $ matches $f_{\infty}^{\circ}$: then by the proof of Lemma \ref{I infty reg}, the function $\gamma\mapsto  \hI^{}_{\gamma'}(f_{\infty}')$ takes finitely many values on {$B_{\infty}'$}, so that the boundedness is trivial.

Assume now that  $v\nmid p\infty$.
 If $\gamma'$ is regular semisimple, then the orbital integral $\hI^{}_{\gamma}(f_{v}', \chi_{v})$ is equal to $\hI^{\sharp}_{\gamma'}( \ff'_{v},  \chi_{v}) $ defined by \eqref{orb int def}. The latter is a finite sum of the values of the integrand at the cosets under the maximal compact subgroup of $H'_{1,v}\times H'_{2,v}$ under which  $f'_{v}$ is invariant. Therefore the integral is a polynomial in $\chi_{v}(\vpi_{v})$ and $\chi_{v}(\vpi_{v})^{-1}$ whose coefficients' denominators are bounded by those of 
  $f_{v}'$. If $f_{v}'$ is the unit measure or more generally spherical, then the orbital integral is an integral polynomial since the volume of the compact open subgroup $\GL_{n}(\sO_{F_{v}})\times (\GL_{n}\times \GL_{n+1})(\sO_{F_{0,v}})$ of $H'_{1,v}\times H'_{2,v}$ is equal to one by our choice of measures.

In general, for a plus-regular element $\gamma'$,  by \cite[Lemma 5.14]{Lu} the integral in \eqref{orb int def} is absolutely convergent (in the archimedean topology) when the exponent of $|\chi_v|$ is small enough. This implies that for some large integer $N$, the product  $\chi(\vpi_{v}) ^{-N}I_{\gamma'}^{\sharp}(f'_{v},\chi_v)$ is  a power series in $\chi(\vpi_{v})^{-1}$. Now the normalized orbital integral $\hI_{\gamma'}^{}( \ff'_{v}, \chi_{v})$ is by definition (cf. \eqref{n orb int def}) given by $  \frac{\hI_{\gamma'}^{\sharp}( \ff'_{v}, \chi_{v})}{L_{\gamma'}(\chi_v)}$. In particular, we have an equality of power series in $\chi(\vpi_{v})^{-1}$:
  $$
   \chi(\vpi_{v}) ^{-N} I_{\gamma'}(f'_{v},\chi_v)=\frac{\chi(\vpi_{v}) ^{-N}I_{\gamma'}^{\sharp}(f'_{v},\chi_v)}{L_{\gamma'}(\chi_v)}.
  $$ The same argument as for the regular semisimple case shows that the coefficients of the power series $\chi(\vpi_{v}) ^{-N}I_{\gamma'}^{\sharp}(f'_{v},\chi_v)$ are $p$-adically  bounded, and integral if $f_{v}'$ is spherical.
   Since $L_{\gamma'}(\chi_v)$ is an  integral polynomial   in  the variable $\chi(\vpi_{v})^{-1}$, it follows that the coefficients of the power series $\chi(\vpi_{v}) ^{-N} I_{\gamma'}(f'_{v},\chi_v)$ (which is in fact a polynomial) are also bounded  in the $p$-adic topology, as desired. 
\end{proof}

\subsubsection{Proof of Theorem \ref{p-adic rtf} / I}
The definitions in part  \eqref{part glob orb} and \eqref{part mu} of Theorem \ref{p-adic rtf} are self-explanatory.  

For part  \eqref{DR Lgamma},  it suffices to  show the existence of an integral interpolation of  the functions on $Y(p^{\infty})_{\Q_{p}}\subset \sY$ given by the abelian $L$-functions in \eqref{pre DR} (with the Euler factors at $p\infty$ removed). This is a consequence of the results of Deligne--Ribet \cite{DR}, who prove that for every totally even  (respectively odd) finite Hecke  character  $\xi$  of a totally real field $F_{0}'$ and every even (respectively odd) $k\geq 0$, there is an element $L_{p}(1-k, \xi)\in \Z_{p}(\xi)\llbracket \Gamma_{F_{0}'} \rrbracket$ interpolating $\chi'\mapsto L_{p}(1-k, \xi\chi')$.\footnote{The results of Deligne--Ribet are stated for totally even characters only and include the interpolation at varying~$k$:  it is well-known  that this allows to obtain the case of a totally odd character $\xi$   by reduction to the totally even character $\xi\omega$, where $\omega$ is the Teichm\"uller character.}

The boundedness of the measures $I^{\ord}_{\gamma, p, K_{p}}(\chi_{p})$ (which is, importantly, uniform in $\chi_{p}$) follows from Lemma \ref{mu bded}~\eqref{mub 2}. Their  explicit and uniform variant over  $\RB(F_{0})^{\circ}$ when $K_{p} $ is a CSDI follows from Proposition \ref{chi1 along} \eqref{chi1 along!}.

The boundedness of  $\gamma\mapsto \sI^{p}_{\gamma}(f'^{p})$ follows from the integrality of $L_{\gamma , p}$ and
Lemma \ref{bdd loc orb}. 

\medskip

We have thus proved parts    \eqref{DR Lgamma}, \eqref{part glob orb}, \eqref{part mu} of Theorem \ref{p-adic rtf}. After some preliminaries, we will now prove the existence of the global distribution and the geometric expansion in  part \eqref{part p exps}. The spectral expansion in  part \eqref{part p exps} 
will be proved in \S~\ref{pf prtf 3}.

\subsubsection{Finite-slope distributions}
For $\gamma\in \RB'_{\rs}(\A)$, and $\Pi$ as in Theorem \ref{plf thm}, we first  define
the following distributions on the subspace of $\sH(\G'(\A^{p}), L)^{\circ} \ot \sH_{p}^{\dag}$ consisting of elements that are sufficiently positive for all the relevant data:
\beq\lb{I dag 2}
I^{\dag}_{\Pi, K_{p}}(f'^{p} f_{p}^{\dag}, \chi) 
&\coloneqq  {1\over 4}
\sL(\RM_{\Pi}, \chi)
\cdot \prod_{v\nmid p} I_{\Pi_{v}}(f_{v}, \chi_{v})   \cdot
I_{\Pi_{p}, K_{p}}^{\dag}(f_{p}^{\dag} , \chi) = I_{\Pi}(f'^{p} f'_{p}), \\
I^{\dagger}_{\gamma,K_{p}} (f'^{p}f_{p}^{\dag} ,\chi)&\coloneqq 
\kappa(\one_{\infty})^{-1} L^{p}_{\gamma}(\chi) \cdot
\prod_{v\nmid p } I_{\gamma,v}(f_{v}', \chi_{v})\cdot I^{\dagger}_{\gamma,p} (f_{p}^{\dag} ,\chi_{v}).
\eeq
Then we may define and expand
\beq\lb{I ord exp n}
 I^{\dag}_{K_{p}}(f'^{p} f_{p}^{\dag}, \chi) 
&\coloneqq \sum_{\Pi\in \sC(\G')_{\Q_{p}}^{\rm her}}  I^{\dag}_{\Pi, K_{p}}(f'^{p}f_{p}^{\dag}, \chi) 
&=
 \sum_{\gamma \in \RB'(F_{0})} 
I^{\dagger}_{\gamma,K_{p}} (f'^{p} f_{p}^{\dag} ,\chi),
\eeq
where the geometric expansion  is valid if $f'^{p}$ has globally plus-regular support or if $F/F_0$ splits above $p$ and  $K_{p}$ is a CSDI: in those cases it  is a consequence of   Proposition \ref{rat rtf}  \eqref{rat rtf part},  the definition \eqref{def orb dag}, and Proposition \ref{chi1 along} \eqref{+reg spt}.

For an integer $s\geq 1$, we denote
$\vpi_{p}\coloneqq  \prod_{v\vert p}\vpi_{v}\in \sO_{F_{0,p}}$ and
$$\Gamma_{F_{0}, s}\coloneqq  \Gamma_{F_{0}}/ (1+ \vpi_{p}^{s}\sO_{F_{0, p}}).$$

\begin{lemma}\lb{mu bded gl} 
Let $f'^{p}\in \sH(\G'(\A^{p}), L)^{\circ} $.  
There is a constant $c(f'^{p}, K_{p})\in \Q^{\ts}$ such that the following holds. For each $s$ that is sufficiently positive for $K_{v}$ for all $v\vert p$, and each $(f_{v}^{\dag})_{v\vert p}\in \prod_{v\vert p}\sH_{v}^{\dag}(\sO_{L})$ that are  sufficiently positive for $s$, if we set $f^{\dag}=f'^{p}\ot \ot_{v\vert p}f_{v}^{\dag}$,  then: 
\begin{itemize}
\item
 for  each $\gamma \in \RB'(F_{0})$, 
the map
\beqq
 Y(\vpi_{p}^{s})\ni\chi \longmapsto I_{\gamma, K_{p}}^{\dag} (f^{\dagger}, \chi) 
 \eeqq
 extends by linearity to a  functional $C(\Gamma_{F_{0},s}, \sO_{L})\to c(f'^{p}, K_{p})\sO_{L}$;
 \item 
if $F/F_{0}$ splits above $p$ and  $K_{p}$ is a CSDI or $f^{'p}$ has plus-regular support, the map
$$  Y(\vpi_{p}^{s})\ni\chi\longmapsto I_{ K_{p}}^{\dag} (f^{\dagger}, \chi)  $$
 extends by linearity to   a functional $C(\Gamma_{F_{0},s}, \sO_{L})\to c(f'^{p}, K_{p})\sO_{L}$.
\end{itemize}
\end{lemma}
\begin{proof} The desired extension of  $I_{\gamma, K_{p}}^{\dag}(f^{\dag} )$ is the convolution of the measure on $\Gamma_{F_{0}, s}$ given by $\prod_{v\vert p}I_{\gamma, K_{v}}^{\dag}(f_{v}, -)$, which are bounded in terms of $K_{p}$ by Lemma \ref{mu bded}~\eqref{mub 2}, and 
(the restriction of) $\sI_{\gamma}^{p}(f'^{p})$, which we have seen to be bounded uniformly in $\gamma$.

 For $I_{ K_{p}}^{\dag}(f^{\dag}) $,  the extension is then  defined via the   (finite) geometric expansion.
 \end{proof}

\subsubsection{Proof of Theorem \ref{p-adic rtf} / II}   
Corollary \ref{Iord ep} shows  that in the limit
 \beq\lb{I ord exp}
I_{K_{p}}^{\rm ord}(f'^{p}, \chi)  &\coloneqq  \lim_{N\sto \infty} I^{\dag}_{K_{p}}(f'^{p} U_{t, K^{p}}^{N!}, \chi) 
\\
& =\sum_{\Pi\in \sC(\G')^{\rm her, ord}}   {1\over 4}
e_{p}(\RM_{\Pi}, \chi) \sL(\RM_{\Pi}, \chi) \cdot 
(\ot_{v\nmid p} I_{\Pi_{v}})(f'^{p}, \chi^{p}). \eeq
 The existence of the limit and \eqref{I ord exp n} prove that  the orbital-integral functions 
 $$I^{p}_{(-)}(f'^{p}, \chi)\colon 
\gamma\mapsto \kappa(\one_{\infty})^{-1}\cdot (\ot_{v\nmid p}I_{\gamma, v})(f'^{p}, \chi^{p})$$ are $I^{\ord}_{\gamma,p,K_{p}}(\chi_{p})$-integrable,
 and that 
\beq\lb{I ord exp g}
I_{K_{p}}^{\rm ord}(f'^{p}, \chi)  = 
\int_{} I^{p}_{\gamma}(f'^{p}, \chi) \,dI^{\ord}_{\gamma,p,  K_{p}}(\chi_{p}),\eeq
where the integration is over $\RB'(F_{0})$ if  $F/F_{0}$ splits above $p$ and $K_{p}$ is a CSDI,  and $\RB(F_{0})'^{\circ}$ otherwise.

Now by Lemma \ref{mu bded gl}, the map  $\chi\mapsto I_{K_{p}}^{\rm ord}(f'^{p}, \chi) $ coincides, for each $s$,  with  the evaluation  of a limit of  uniformly  (in both $f^{\dag}=U_{t}^{N!} $ and  $s$) bounded Radon measures on $\Gamma_{F_{0}, s}$, hence it extends uniquely to a generalized bounded Radon measure 
\beq \lb{sI fun}
\sI_{K_{p}}^{}(f'^{p})\colon C(\Gamma_{F_{0}}, L)\to L\eeq
corresponding  to the element $\sI^{}_{K_{p}}(f'^{p})\in \sO(\sY_{L})$ of part \eqref{part p exps}. 

The geometric expansion in part  \eqref{exp a} is \eqref{I ord exp g}.
Then if $F/F_{0}$ splits above $p$ and $K_{p}$ is a CSDI, by Lemma \ref{In exp} below applied  to \eqref{mu def}, the expansions of \eqref{I ord exp g} imply that the distribution $\sI^{}_{K_{p}}$ has the geometric  expansion described in part \eqref{exp b}.  
\begin{lemma} \lb{In exp}
Let $(\sI_{N})_{N\in {\bf N}}$, $ \sI_{\infty}\in \sO(\sY)$. Suppose that for all $\chi\in Y(p^{\infty})$ we have 
$$\lim_{N\sto \infty}\sI_{N}(\chi)= \sI_{\infty}(\chi).$$ Then $\lim_{N\sto \infty} \sI_{N}=\sI_{\infty}$. 
\end{lemma}
\begin{proof} Recall that  $Y(p^{\infty})=\varinjlim Y(p^{s})$; then  observe that  the ideals   $J_{s}\coloneqq  \Ker [\sO(\sY)\to  \sO(Y(p^{s})) \subset \prod_{\chi\in Y(p^{s})} \Q_{p}(\chi)]$ form a fundamental system of neighborhoods of~$0$ in $\sO(\sY)$.
\end{proof}

We now turn to the  $p$-adic $L$-function, then to the spectral expansion of $\sI_{K_{p}}^{}$.

\subsubsection{Proof of Theorem \ref{plf thm}  (= Theorem \ref{thm B})}
Let $K_{p}=\prod_{v}K_{v}$ be a convenient subgroup such that for every place $v\vert p$ of $F_{0}$, the representation $\Pi_{v}$ is $K_{v}$-ordinary.
Similarly to the proof of Theorem \ref{rat L}, we use Corollary \ref{cor test} asserting the existence of suitable test Gaussians. (Note that as $\chi_{|\A^{p , \times}}$ is smooth, the definition of `adapted to $(\Pi, \chi, K_{p})$' in \S~\ref{rha} still makes sense, and the proof of that corollary goes through.)

 For any $\chi\in Y(p^{\infty})_{L}$ and any   $f'^{p}\in \sH(\G'(\A^{p\infty}) , L)^{\circ}_{ K_{p}, \rs ,\Pi^{p}, \chi^{p}}$, we define
 $$\sL_{p}(\RM_{\Pi}, \cdot )_{f'^{p}}\coloneqq  {4 \, \sI_{K_{p}}(f'^{p}, \cdot) \over c_{K_{p}}(\Pi) \cdot 
  (\ot_{v\nmid p} I_{\Pi_{v}})(f'_{v}, \cdot)}$$
  away from the  zero set $\sZ(f'^{p})$ of the denominator. Note that we may assume that $f'^{p}$ is a pure tensor with factors equal to the unit Hecke measure at places $v\nmid p\infty$ where $\Pi_{v}$ is unramified; if so, each $\sZ(f'^{p})$ is the pullback of a closed subset $Z(f'^{p})\subset Y_{S}(1)_{L}\coloneqq \prod_{v\in S}Y_{v}(1)_{L}$ for some fixed set of places $S$. As $\sI_{K_{p}}$ restricts to $I_{K_{p}}^{\rm ord}$,  it follows from \eqref{I ord exp} that the functions $\sL_{p}(\RM_{\Pi}, \cdot )_{f'^{p}}$ glue to a function $\sL_{p}(\RM_{\Pi}, \cdot)$ with the desired interpolation properties, on the complement of the polar locus $\sZ\coloneqq  \bigcap_{f'\in\sH(\G'(\A^{p\infty}) , L)^{\circ}_{K_{p}, \rs ,\Pi^{p}, \chi^{p}} }\sZ(f'^{p})\subset \sY_{L}$. By Corollary \ref{cor test}, the closed subset $\sZ$ is empty. The function $\sL_{p}(\RM_{\Pi})$ is still bounded since, by the Nullstellensatz applied to a finite subproduct of $\prod_{v\nmid p\infty}Y_{v}(1)_{L}$, finitely many $f'^{p}$ suffice to construct $\sL_{p}(\RM_{\Pi})$. This completes the proof of Theorem \ref{plf thm}.

\subsubsection{Proof of Theorem \ref{p-adic rtf} / III} \lb{pf prtf 3}
The global spectral terms of Definition \ref{def prtf} \eqref{part IPi p} are now justified.
 The spectral expansion of $\sI_{K_{p}}^{}$ in part \eqref{part p exps} then follows from the definitions and \eqref{I ord exp g}. This completes the proof of Theorem \ref{p-adic rtf}.

\subsection{Derivative of the analytic RTF}\lb{sec 63}
We study the derivative of the distribution $\sI_{K_{p}}$. 
\subsubsection{Notation} \lb{63 notn} 'Denote by $\frak m\subset \sO_{\sY}$  the ideal of functions vanishing at $\chi=\one$.  For  a $\sY$-scheme $\sY'$ and a function $\Phi\in \frak m \sO(\sY')$, we say that $\Phi$ \emph{vanishes at $\chi=\one$} and we denote by  $\partial \Phi $ be the image of $\Phi $ in $\frak m/\frak m^{2}\ot_{\sO_{\sY}} \sO_{\sY'}= \sO_{\sY'}\hat{\otimes}\Gamma_{F_{0}}$.

For  $V\in \sV^{\circ}$ a coherent or incoherent pair of definite hermitian spaces  as in \S\ref{coh incoh G}, and $v$ a finite place of $F_{0}$,  we  let:
\begin{itemize}
\item  $ \sC(\G')^{{\rm her, ord}, V}\subset  \sC(\G')^{{\rm her, ord}}$ be the subset of those isomorphism classes of representations $\Pi$ such that for each finite place $u$ of $F_{0}$, the space $V_{u}$ is the one attached to $\Pi_{u}$ by  the local Gan--Gross--Prasad conjecture (Proposition \ref{local ggp}).
\item  $\sH(\G'(\A^{p}), L)_{K_{p},\rm qc}^{\circ, V}\subset \sH(\G'(\A^{p}), L)^{\circ}_{K_{p}, \rm qc}$ be the subspace of those $f'^{p}$ that purely match  (spectrally and geometrically, see Proposition \ref{equiv match}) a Hecke measure in $\sH(\G^{V}(\A^{p}), L)^{\circ}\coloneqq  \sH(\G^{V}(\A^{p}), L)\ot_{L} L f_{\infty}^{\circ}$;
\item   $\one_{V}$ be  the characteristic function of $\RB_{\rs}'(\A)_{V}\coloneqq  \prod_{u}' B'_{\rs, u , V_{u}} \subset \RB'_{\rs}(\A)$;
\item  $V(v)\in \sV^{\circ}$ be the  pair such that $V(v)_{u}\cong V_{u}$ exactly for  $u\neq v$; it is coherent if and only if $V$ is incoherent.
\end{itemize}

\begin{proposition} \lb{prop:de I}
Consider the situation of Theorem \ref{p-adic rtf}, and let $V\in \sV^{\circ, -}$ be an incoherent pair.  Suppose that all places $v\vert p$ of $F_0$ split in $F$.

For all $f'^{p}\in \sH(\G'(\A^{p}), L)_{K_{p}, \rm qc}^{\circ, V}$, the following hold. 
 \begin{enumerate}
 \item  \lb{vanishhh}
 For  all $\Pi\in  \sC(\G')^{{\rm her, ord}}$ and  all $\gamma\in \RB'(F_{0})$, 
 $$
 \sI^{}_{K_{p}}(f'^{p}, \one)=\sI^{}_{\Pi, K_{p}}(f'^{p}, \one)
= \sI^{p}_{\gamma}(f'^{p}, \one)=0.$$
\item \lb{spectral der}
There is a  spectral expansion
\beqq
\partial\sI^{}_{K_{p}}(f'^{p}) = \sum_{\Pi\in \sC(\G')^{{\rm her, ord}, V}_{K_{p}}} \partial\sI^{}_{\Pi, K_{p}} (f'^{p})
\eeqq
where
\beqq
\partial\sI^{}_{\Pi, K_{p}} (f'^{p})
 = {1\over 4} \,c_{K_{p}}(\Pi)\, \partial \sL_{p}(\RM_{\Pi})\cdot (\ot_{v\nmid p} \sI_{\Pi_{v}})(f'^{p},\one).
\eeqq
 \item \lb{geom der}
Suppose that $f'^{p}$ has globally regular semisimple support.
The function $ \partial \sI^{p}_{(-)}(f'^{p})$ is integrable for the Radon measure $I^{\ord}_{\gamma, p,K_{p}}\coloneqq I^{\ord}_{\gamma, p,K_{p}}( \one)$, and 
there is a geometric expansion
\beq\lb{dec I'}
\partial\sI^{}_{K_{p}} (f'^{p}) 
&=   \int_{\RB_{\rs}'(F_{0})^{\circ}} \partial \sI^{p}_{\gamma}(f'^{p}) \, dI^{\ord}_{\gamma, p , K_{p}}\\
&=\int_{\RB_{\rs}'(F_{0})^{\circ}}\sum_{v\nmid p\infty  {\rm \ nonsplit}}
\one_{V(v)}(\gamma) \sI^{vp}_{\gamma}(f'^{vp}, \one) \cdot \partial \sI^{}_{\gamma,v} (f_{v}')  \, dI^{\ord}_{\gamma, p, K_{p}},\eeq
where for any $\gamma$ we put 
$ \sI^{vp}_{\gamma}\coloneqq  \kappa(\one_{\infty})^{-1} \cdot L_{p, \gamma} \cdot \ot_{u\nmid v p} \sI_{\gamma,u} $.
 \end{enumerate}
\end{proposition}

\begin{proof}
Consider the geometric terms $\sI_{\gamma}(f'^{p}, \one)$. For $\gamma\in \RB'_{\rm rs}(F_{0})\cap B_{\infty}'^{\circ}$, let $V_{\gamma}\in \sV^{\circ, +}$ be the unique coherent pair such that  $\gamma$ matches an orbit in $B_{\rs}(F_{0})_{V_{\gamma}}$ as  in \eqref{dec Brs gl}; let  $\Sigma(\gamma, V)$ be the non-empty finite set of non-archimedean (and necessarily nonsplit) places of $F_{0}$ such that $V_{\gamma,v}\not \cong V_{v}$. If $v\in \Sigma(\gamma, V)$, then by the assumption on $f'^{p}$ we have $I_{\gamma, v}(f'_{v}, \one)=0$; hence  $\sI^{p}_{\gamma}(f'^{p})$ vanishes at $\one $ to order at least $|\Sigma(\gamma, V)|\geq 1$.  Moreover, if $v\in \Sigma(\gamma, V)$ then 
\beq\lb{dec der aux} \partial \sI_{\gamma}(f'^{p}) =\sI^{vp}_{\gamma}(f'^{vp}, \one) \cdot \partial \sI^{}_{\gamma,v} (f_{v}'),\eeq
which can be nonzero only if  $\Sigma(\gamma, V)=\{v\}$, equivalently $V_{\gamma}=V(v)$.

Consider now a representation $\Pi =\Pi_{n}\boxtimes \Pi_{n+1}\in  \sC(\G')^{{\rm her, ord}}$. Let $V_{\Pi}\in \sV^{\circ, \e(\Pi)}$ be the pair such that $\Pi\in  \sC(\G')^{{\rm her, ord}, V_{\Pi}}$ (cf. Proposition \ref{cor C her}). If $\e(\Pi)=-1$, then  $\sL(\RM_{\Pi}, \one)=0$ by the functional equation of Rankin--Selberg $L$-functions; this implies $\sI_{\Pi, K_{p}}(f'^{p\infty}, \one)=0$. If $\ep(\Pi)=+1$, then for any finite place $v$ such that $V_{\Pi, v}\not\cong V_{v}$, we have $I_{\Pi_{v}}(f'_{v}, \one)=0$ by the assumption on $f'^{p}$.  This completes the proof of part \eqref{vanishhh}. 
More generally, we note that the last argument  shows that 
\beq \lb{vanish aux}
\Pi\not\in \sC(\G')^{{\rm her, ord}, V}\  \Longrightarrow
\  I_{\Pi_{v}}(f'_{v}, \one)=0 \text{ \ for some $v\nmid p\infty$}.\eeq
This shows that the sum in part \eqref{spectral der} indeed runs over $\sC(\G')^{{\rm her, ord}, V}$; as above, this implies that $\ep(\Pi)=-1$ and  $\sL_{p}(\RM_{\Pi}, \one)=0$, which implies the second equality in \eqref{spectral der}.

We now consider part \eqref{geom der}. By  the definition of the measure $I_{\gamma, p,K_{p}}^{\ord}$, the second equality follows from   \eqref{dec der aux}, whose right-hand side can be nonzero only if $V_{\gamma}=V(v)$. We consider the expansion in the first equality; it will hold without the condition that $f'^{p}$ is quasicuspidal, and by linearity we may thus assume that $f'^{p}=\ot_{v}f_{v}'$ is a pure tensor.
Suppose first that $F/F_0$ splits above $p$ and $K_{p}$ is a CSDI. Since $\int dI^{\ord}_{\gamma, p, K_{p}}$ is a bounded functional, we simply differentiate under the integral sign in Theorem \ref{p-adic rtf}\,\eqref{part p exps}. We now consider the general case.  Viewing  $\sI_{K_{p}}^{}$ as in \eqref{sI fun}
 we have
$$\partial  \sI^{}_{K_{p}}(f'^{p}) = \sI^{}_{K_{p}}(f'^{p}, \ell)$$ 
(see for instance, \cite[Lemma 3.42]{DL}),
where  the `logarithm' $\ell\colon \Gamma_{F_{0}}\to \Gamma_{F_{0}}^{\rm fr}\subset \Gamma_{F_{0}}\hat{\ot} \Q_{p}$ is the projection onto the maximal $\Z_{p}$-free quotient  $\Gamma_{F_{0}}^{\rm fr}$ of $ \Gamma_{F_{0}}$.

For $s\geq 1$, let $\ell_{s}\colon \Gamma_{F_{0}}\to \Gamma_{F_{0}}^{\rm fr}/p^{s}$ be the reduction map, and let $\wt{\ell}_{s}\colon \Gamma_{F_{0}}\sto  \Gamma_{F_{0}}^{\rm fr}/p^{s}\to\Gamma_{F_{0}}^{\rm fr}$ be any lift of $\ell_{s}$, which is a linear combination of characters whose conductors at places $v\vert p$ do not exceed~$s$.
By the definition of $\sI_{K_{p}}^{}$, the expansion \eqref{I ord exp}, and linearity, we  have 
$$\sI^{}_{K_{p}}(\wt{\ell}_{s}) = \lim_{N\sto \infty} I^{\dag}_{K_{p}}(f'^{p} U_{t_{p}}^{N!}, \wt{\ell}_{s})$$
Then by Proposition \ref{rat rtf} \eqref{rat rtf part}, we have
\beq\lb{i ell exp}
\sI_{K_{p}}(\wt{\ell}_{s}) = \lim_{N\sto \infty}
\sum_{\gamma \in \RB_{\rs}'(F_{0})}
  I^{\dag}_{\gamma}(f'^{p} U_{t_{p}}^{N!}, \wt{\ell}_{s}).
\eeq

By  Lemma \ref{mu bded gl},
   up to multiplying  $f'^{p}$ by a power of~$p$ independent of~$s$ we have that all terms in \eqref{i ell exp} are $p$-integral; hence  it makes sense to  consider  the reduction of that identity modulo~$p^{s}$, 
$$
\sI^{}_{K_{p}}({\ell}_{s}) =
 \lim_{N\sto \infty}
\sum_{\gamma \in \RB_{\rs}'(F_{0})}
  I^{\dag}_{\gamma}(f'^{p} U_{t_{p}}^{N!}, {\ell}_{s})$$
in  $\Gamma_{F_{0}}^{\rm fr}/p^{s}$.  
 Now $\ell_{s}=\sum_{v\nmid\infty}\ell_{s,v}$, where $\ell_{s,v}\coloneqq \ell_{s|F_{0,v}^{\ts}}$, so from Remark \ref{gl int ell}, the $\gamma$-summand equals
\beq\lb{sum orb der}
\sum_{v\nmid \infty} 
{I_{\gamma}(f'_{\infty})
\over \kappa(\one_{\infty}) \kappa_{\infty}(\gamma') } 
\int_{\RH_{1}(\A^{\infty})}\int_{\RH_{2}(\A^{\infty})} \ff'^{\infty}(h_{1}^{-1}\gamma'  h_{2}) \ell_{s,v}(h_{1,v}) \eta(h_{2}) \, {d^{\natural}{h_{1}}d^{\natural}h_{2}\over d^{\natural}g}\eeq
in $\Gamma_{F_{0}}^{\rm fr}/p^{s}$; here
 $\gamma'\in \G'_{\rm rs}(F_{0})$ is any preimage of $\gamma$.
 (Note that only finitely many $v$-summands are nonzero, hence it is trivial to interchange sum and integration.) For $v\nmid p$, the $v$-summand is 
\beqq
\ &
{I_{\gamma}(f'_{\infty}) \over \kappa(\one_{\infty}) \kappa_{\infty}(\gamma') } 
  I^{\dag}_{\gamma, p}(U_{t_{p}}^{N!}, \one)
\sI_{\gamma}^{vp\infty}(f'^{vp\infty}, \one) \cdot  \sI_{\gamma,v}  (f'_{v}, \ell_{s,v})\\
\equiv&
{1\over \kappa(\one_{\infty}) \kappa_{\infty}(\gamma') } 
  I^{\dag}_{\gamma, p}(U_{t_{p}}^{N!}, \one)
 \sI_{\gamma}^{vp}(f'^{vp}, \one)  \cdot \partial \sI_{\gamma,v}  (f'_{v}).\eeqq 
 For $v\vert p$, the $v$-summand in \eqref{sum orb der} is a multiple of $\sI_{\gamma}^{p}(f'^{p}, \one)$, which is zero by part \eqref{vanishhh}. Therefore  $\partial \sI_{K_{p}}^{}(f'^{p})$ is congruent to 
$$
 \sI^{}_{K_{p}}(f'^{p},\ell)\equiv \lim_{N\sto \infty} 
\sum_{\gamma \in \RB_{\rs}'(F_{0})}
  I^{\dag}_{\gamma, p}(U_{t_{p}}^{N!}, \one)
\cdot \partial \sI^{p}_{\gamma}(f'^{p})
$$
in $\Gamma^{\rm fr}_{F_{0}}/p^{s}$ for all $s$. We conclude that the above congruences amount to an equality in $\Gamma_{F_{0}}^{\rm fr}$; by definition, the right-hand side is 
$$
\int_{\RB'_{\rs}(F_{0})} \partial \sI^{p}_{\gamma}(f'^{p}) \, dI^{\ord}_{\gamma, p, K_{p}},$$
as desired.
\end{proof}

\part{$p$-adic heights and the arithmetic relative-trace formula}
\lb{part 2}

We now study the $p$-adic heights of Gan--Gross--Prasad cycles. In \S\ref{sec:8}, we recall the relevant Shimura varieties, the arithmetic diagonal cycles,  and their moduli interpretations over the reflex fields. In \S\ref{sec:9},  we study various integral models and prove some   vanishing results for their cohomologies. In \S~\ref{sec ht} we collect the necessary definitions and results on cycles and $p$-adic heights. In \S\ref{sec:10}, we define the  arithmetic relative-trace distribution encoding the heights of GGP cycles,  and prove the corresponding RTF.  

In \S\S\ref{sec:8}-\ref{sec:9}, we largely follow   \cite{LTXZZ,RSZ3,RSZ4} to define the relevant Shimura varieties and moduli functors, and we use slightly different notation on unitary groups from the rest of the paper.

\section{Unitary Shimura varieties and arithmetic diagonals}\lb{sec:8}
{We define the relevant unitary Shimura varieties and the arithmetic diagonal cycle on certain product unitary Shimura varieties. 
In \S~\ref{ss:Sh U}, we define Shimura varieties for certain unitary groups and a product of two unitary groups. They are the Shimura varieties in the statement of our main conjecture in \S~\ref{intro2.5};  they do not have natural moduli description. 
In \S~\ref{incoh Sh}, we discuss incoherent Shimura varieties (this subsection is not used in the rest of the paper).
In \S~\ref{ss:RSZ Sh}, we discuss the related RSZ Shimura varieties, which admit PEL-type moduli definitions given in \S~\ref{subs: RSZ}. }

\subsection{Unitary Shimura varieties}
\label{ss:Sh U}
We keep denoting by  $F$  a CM number field with maximal totally real subfield $F_0$ and nontrivial $F/F_0$-automorphism ${\rm c}\colon a \mapsto a^{\rm c}$.
For an algebraic group $\G$ over $F_{0}$, we denote its restriction of scalars to $\Q$ by
$$\G^{\flat}\coloneqq  \Res_{F_{0}/\Q}\G.$$

\subsubsection{Unitary Shimura data and the associated varieties} \lb{sh data} 
 Let $\nu$ be a positive integer.  
 \begin{comment}
 Recall from \cite[\S2.2]{RSZ4} that a \emph{generalized CM type} (or a \emph{signature  type}) of rank $\nu$ is a function $r\colon \Hom_\Q(F, \ov\Q)\to \BZ_{\geq 0}$, denoted $\varphi \mapsto r_\varphi$, such that 
\[
   r_\varphi+r_{\varphi^{{\rm c}}}=\nu \quad \text{for all}\quad \varphi\in \Hom_\Q(F, \ov\Q);
\]
here  $\vphi^{\rm c}\coloneqq  \vphi\circ {\rm c}$. 
  When $\nu=1$, a generalized CM type is ``the same'' as a usual CM type, via 
\[
   \Phi= \bigl\{\, \varphi\in \Hom_\Q(F, \ov\Q) \bigm| r_\varphi=1 \,\bigr\}.
\]

\end{comment}
{Let $(W, \sform)$ be an $F/F_0$-hermitian vector space of dimension $\nu$. We will assume that $W$ is {\em standard indefinite} relative to a fixed archimedean place $v_0\in \Hom(F_0,\R)$, in the sense that the signature of $W$ is $(\nu-1,1)$ at $v_0$, and $(\nu,0)$ at all other archimedean places.\footnote{The terminology follows \cite[Definition 3.2.1]{LTXZZ}. Note that our signature convention  differs from that of \cite[\S27]{GGP}, where the signs were $(1,\nu-1)$ at $v_0$ and $(0,\nu)$ elsewhere.}
}

Consider the unitary group
\begin{equation}\label{G}
   \G\coloneqq
    \U(W).
\end{equation}

 Fix an element $\varphi_0\in\Hom(F,\C)$ above $v_0$. We now define a Shimura datum $(\G^{\flat}, \{h_\G^{\flat}\})$, where  $ \{h_\G^{\flat}\}$ is a $\G^{\flat}(\R)$-conjugacy class of homomorphisms $\BS\coloneqq \Res_{\C/\R}\bG_{m, \R}\to  \G^\flat(\R) $
\begin{equation*}\label{G^Q decomp}
h_\G^{\flat}\colon
\begin{tikzcd}[baseline=(x.base), row sep=-.5ex, ampersand replacement=\&]
      |[alias=x]| \BS\ar[r]  \&  \G^\flat(\R) = \U(W_{v_0})\times \prod_{v\neq v_0,\atop v\mid\infty}\U(W_{v})\\
     z \ar[r, mapsto]
        \&  \bigl( \diag(1_{\nu-1},  z^{{\rm c}}/z), 1_{\nu},\cdots,1_{\nu}\bigr)
   \end{tikzcd}
\end{equation*}
Here we view $\U(W_{v_0})$ as a subgroup of $\GL_\C(W\otimes_{F,\varphi_0}\C)$ and identify $\GL_\C(W\otimes_{F,\varphi_0}\C)$ with $\GL_\nu(\C)$  by choosing a $\C$-basis of $W\otimes_{F,\varphi_0}\C$ with respect to which the matrix of the hermitian form is given by  $\diag(1_{\nu-1},  -1)$.

The reflex field $E({\G^\flat},\{h_{\G^\flat}\})$ of this Shimura datum is $F$, as a subfield of $\C$ via the embedding $\varphi_0:F\to\C$.

\begin{comment}
the reflex field $E_{r^\natural}$ of the function $r^\natural$, characterized by
\begin{equation}\label{eq:reflex r nat}
   \Gal(\ov\Q/E_{r^\natural})= \bigl\{\, \sigma\in\Gal(\ov\Q/\Q)\bigm| \sigma^*(r^\natural)=r^\natural \,\bigr\},
\end{equation}
where we define a modified function
\begin{equation}\label{eq: r^nat}
   r^\natural \colon
   \begin{tikzcd}[baseline=(x.base), row sep=-.5ex, ampersand replacement=\&]
      |[alias=x]| \Hom_\Q(F, \ov\Q) \ar[r]  \&  \BZ_{\geq 0}\\
      \varphi \ar[r, mapsto]
        \&  \begin{cases}
              0,  &  \varphi\in \Phi;\\
              r_\varphi,  &  \varphi \in \Phi^{\rm c}.
           \end{cases}
   \end{tikzcd}
\end{equation}
\end{comment}

We then obtain a tower of  Shimura varieties $(\Sh_K({\G^\flat}, \{ h_{\G^\flat} \}))_{K \subset {\G}(\A^{\infty})}$ over $F$. 
Since this is the only Shimura datum we use, we will abbreviate $$\Sh_K(\G)\coloneqq \Sh_K({\G^\flat}, \{ h_{\G^\flat} \}),$$ 
omitting the superscript $\flat$ and suppressing the datum $\{h_{\G^{\flat}}\}$.

\begin{remark}
The  Shimura variety $\Sh_K({\G^\flat}, \{ h_{\G^\flat} \})$  is not of PEL type, i.e. it is not related to a moduli problem of abelian varieties (this can be seen already from the fact that the restriction of $\{ h_{\G^\flat} \}$ to $\bG_m\subset \BS$ is not mapped via the identity map to the center of ${\G^\flat}$). However, this Shimura variety is of abelian type.
\end{remark}

\begin{comment}
    
\subsubsection{A special signature type}\lb{ss: fake Dr}
In the rest of the paper we will only consider the following generalized CM type  $r$
\begin{equation}\label{SFDsig}
   r_\varphi=
   \begin{cases}
      \nu-1,  &  \text{$\varphi = \varphi_0$},\\
      \nu,  &  \varphi \in \Phi \ssm \{\varphi_0\},
   \end{cases}
\end{equation}
for some CM type  $\Phi$ and some $\vphi_{0}\in \Phi$. This signature type is called {\em strictly fake Drinfeld} (with respect to $(\Phi, \vphi_{0})$) in \cite[Example 2.3~(ii)]{RSZ4} (cf. {\it loc. cit.} for some other generalized CM types). 
In this case,  we have $\varphi_0 \colon F \isoarrow E_{r^\natural}$ for all $\nu \geq 1$; in other words, the reflex field is $F$ via the embedding $\varphi_0\colon  F\to \C$.  Since this is the only generalized CM type we use, we will abbreviate $$\Sh_K(\G)\coloneqq \Sh_K({\G^\flat}, \{ h_{\G^\flat} \}),$$ 
omitting the superscript $\flat$ and suppressing the datum $\{h_{\G^{\flat}}\}$.
\end{comment}

 \subsubsection{Hecke correspondences}\lb{Hecke corr}
Recall that if  $\sX$ is a scheme, a \emph{correspondence} on $\sX$ is a diagram of finite morphisms
\begin{equation*}
\begin{gathered}
   \xymatrix{
	     &  \sX' \ar[dl] \ar[dr] &\\
\sX & &  \sX.
	}
\end{gathered}
\end{equation*}
It is said to be \emph{ \'etale} if both morphisms are  \'etale. Correspondences on $\sX$ form a monoidal category under  composition.
If $L$ is a ring, we denote by  ${\text{\'{E}tCorr}}(\sX)_{L}$ the $L$-algebra generated by isomorphism classes of  \'etale correspondences on $\sX$. It acts (on the right) on cycles and  cohomology of $\sX$ by pullback and pushforward.

For each $K\subset \G(\A^{\infty})$, and each characteristic-zero field $L$, we have an $L$-algebra homomorphism 
\beq\label{hecke Sh}
\hT\colon \sH(\G(\A^{\infty}), L)_{K}  &\to \EC(\Sh_{K}(\G))_{L} \\ 
[KgK]
&\mapsto 
\left[\begin{gathered}
   \xymatrix{
	     &\Sh_{K'}(\G) \ar[dl]_-{\text{nat}_1} \ar[dr]^-{\text{nat}_g}&\\
	  \Sh_K(\G)  & &  \Sh_K(\G)
	}
\end{gathered}\right]
\eeq
where $[KgK]\coloneqq  \vol(K)^{-1} \one_{KgK}\, dg$,  $K'\coloneqq  K\cap gK g^{-1}$, and the map $\text{nat}_1$ is the natural map induced by the embedding $K'\subset K$ while $\text{nat}_g$  is induced by the composition $$\begin{gathered}
   \xymatrix{  \Sh_{K'}(\G)\ar[r]^-{. g} &  \Sh_{g^{-1}K'g}(\G)\ar[r]&   \Sh_K(\G)}\end{gathered}.$$ 
For the other Shimura varieties in this section, we also have Hecke correspondences defined in an entirely analogous way.

\subsubsection{Product Shimura varieties and the arithmetic diagonal}\lb{ss:AD}
Fix an archimedean place $\varphi_0\in\Hom(F,\C)$ of $F$ and $v_0\in \Hom(F_0,\R)$ the induced place of $F_0$. Let $V_{n}$ be a hermitian space of dimension~$n\geq1$. Let $V_{n+1}=V_n\oplus Fu$ where the special vector $u$ has totally positive norm $\nrm$. We assume that $V_n$ is standard indefinite (relative to the fixed place $v_0$ of $F_0$); then  so is $V_{n+1}$. Let $\G_{\nu}=\RU(V_{\nu})$ for $\nu=n, n+1$, and let $(\Sh_{K_{\nu}}(\G_{\nu}))_{K_{\nu}}$ be the corresponding tower of Shimura varieties as defined in \S\ref{sh data} (relative to the place $\varphi_0$ of $F$). 
We  also have a product Shimura variety $\Sh_{K}(\G) =  \Sh_K(\G, \{h_{\G}\})$ associated with $\G=\G_{n}\times \G_{n+1}$ and $h_{\G}=h_{\G_{n}^{\flat}}\times h_{\G_{n+1}^{\flat}}$.  
Denote $\RH\coloneqq \G_{n}$. The map 
$$\jmath\colon \RH\to \G
$$ 
that is the graph of the natural embedding $\G_{n}\to \G_{n+1}$ induces a corresponding map of Shimura varieties 
\begin{equation}\lb{eq:ADC}
\jmath\colon \Sh_{K_{ \RH}}(\RH)\to \Sh_{K_{ \RG}}(\RG)
\end{equation}
whenever  $K_{\RH}\subset \jmath^{-1}(K_{\RG})\subset \RH(\A^{\infty})$. 

 The target Shimura variety has dimension $2n-1$, and the image of $\jmath$ has codimension $n$, in the arithmetic middle dimension (i.e.,\ the codimension is just more than half the dimension of the ambient variety). We thus call the map \eqref{eq:ADC} the arithmetic diagonal, and the image cycle (defined in more detail in \S~\ref{ADC 10})  the arithmetic diagonal cycle in $\Sh(\G)$.

\subsection{Incoherent Shimura varieties} \lb{incoh Sh}
For our specific signature type,  we may present the above Shimura varieties more symmetrically using incoherent hermitian spaces. {Since the existing literature is not always fully detailed, at the suggestion of the referee we only include this subsection for expository purposes. It seems plausible that  the theory recently introduced in \cite{STay} should provide
a favorable framework for  a fully rigorous treatment, which would be a valuable contribution to the literature.}

\subsubsection{Incoherent unitary groups and incoherent Shimura varieties} \lb{ssec incoh sh}  Let $W_{}$ be a totally positive definite  incoherent $F/F_{0}$-hermitian space of dimension~$\nu$ (that is, $W$ is shorthand for an $\A_F/\A$-hermitian space $W_\A$ that is definite at all archimedean places). Let $\G=\U(W)$ be the corresponding incoherent unitary group analogously to \S~\ref{sec incoh}. The theory of conjugates of Shimura varieties (\cite{MilS}; see also \cite{gross-incoh, STay})
shows that there exists a unique-up-to-isomorphism tower 
$$(\Sh_K(\G))_{K \subset \G(\A^{\infty})}$$ over $\Spec F$ with the following property. For any archimedean place $v_0$ of $F_0$, let $\varphi_0$ be a place of $F$ above $v_0$. Then the nearby hermitian space $W_{}(v_0)$ is standard indefinite. Let $\G_{}^{(v_0)}=\RU(W_{}(v_0))$ be the
unitary group associated to $W_{}(v_0)$, and let   $(\Sh_K(\G^{(v_0)}))_{K}$ be the tower of  Shimura varieties defined in \S\ref{sh data} using the place $\varphi_0$ of $F$. (We warn the reader that the notation $(\Sh_K(\G^{(v_0)}))_{K}$ has suppressed the dependence on $\varphi_0$.)
Then
$$
\xymatrix{ \Sh_K(\G)\times_{\Spec F,\,\varphi_0}\Spec \vphi_{0}(F)\ar[r]^-{\sim}&  \Sh_K(\G^{(v_0)})}.
$$
Here we fix once for all an isometry $W(v_0)\otimes\A^{\infty}\simeq W_{\A^{\infty}}$; it induces an isomorphism $\G(\A^{\infty})\simeq \G^{(v_0)}(\A^{\infty})$. 
We will call $\Sh_K(\G)$  the incoherent 
Shimura varieties
attached  $W$. 

From now on we will also make the assumption that  all our unitary groups  $\G$ are {\em anisotropic}; in the incoherent case this means that $\G^{(v_0)}$ is anisotropic for any (hence every)  archimedean place  $v_0\in \Hom(F_0,\R)$. Then $\Sh_K(\G)$ is proper for any compact open subgroup $K\subset \G(\A^{\infty})$;
 this is guaranteed if $F_{0}\neq \Q$.

\subsubsection{The arithmetic diagonal for incoherent Shimura varieties}
Fix now an incoherent pair   $V=(V_{n}, V_{n+1})\in \sV^{\circ, -}$. We denote 
 $\RG=\RG^V=\G^{V}_{n}\times \G^{V}_{n+1}\coloneqq \U(V_{n})\times \U(V_{n+1})$ (an incoherent unitary group as in \S~\ref{sec incoh}), and let $\Sh_{K_\RG}(\RG)$ denote the product of the Shimura varieties constructed in \S\ref{ssec incoh sh}. 
 For every place $v_0\in \Hom(F_0,\R)$,
 let $\G^{(v_0)}\coloneqq \G^{V(v_{0})}$ be the product of unitary groups associated to the coherent nearby pair $V(v_0)$ of hermitian spaces.
   Then there exists  a projective system  of varieties $(\Sh_K(\G))_{K \subset \G(\A^{\infty})} $ over $\Spec F$ such that,  for every embedding $\varphi_0\colon  F\to \C$ extending $v_0$ and  every choice of CM type $\Phi$ such that $\varphi_0\in\Phi$ we have
\beq \lb{isom v0}
\xymatrix{ \Sh_K(\G)\times_{\Spec F}\Spec\vphi_{0}(F)\ar[r]^-{\sim}&  \Sh_K(\G^{(v_{0})})}\eeq
where   $\Sh_K(\G^{(v_{0})})=  \Sh_K(\G^{(v_{0}), \flat}, h_{\G^{(v_{0}), \flat}})$ with   $h_{\G^{(v_{0}), \flat}}=h_{\G_{n}^{(v_{0}), \flat}}\times h_{\G_{n+1}^{(v_{0}), \flat}}$ (the latter defined in \S~\ref{sh data}). 
Similarly, we have incoherent Shimura varieties $\Sh_{K_\RH}(\RH)$  for the group $\RH=\RH^V=\U(V_{n})$. 

As in \S~\ref{ss:AD}, we expect to have (finite) maps
\beq\lb{eq: AD incoh}
\jmath \colon \Sh_{K_\RH}(\RH)\to \Sh_{K_\RG}(\RG),
\eeq
which are the pullbacks of \eqref{eq:ADC} via \eqref{isom v0}.

\subsection{RSZ Shimura varieties}\label{ss:RSZ Sh}
The unitary Shimura varieties above do not admit natural moduli descriptions. Hence we will relate them to RSZ Shimura varieties, which admit a PEL type moduli definition. They will play an auxiliary role when computing local heights. We will follow \cite{RSZ4}. 

\subsubsection{Shimura varieties for unitary similitude groups} \lb{ss: wt G}
We resume the notation from \S~\ref{ss:Sh U} and start with a standard indefinite hermitian space $W$ of dimension $\nu$ relative to a fixed archimedean place of $v_0$. Moreover we fix a CM type $\Phi$, and let $\vphi_{0}$ be the unique element of $\Phi$ above $v_0$.   

We denote by $\ov\Q$ the algebraic closure of $\Q$ in $\C$; from now on we will view  a CM type $\Phi$ as a subset of $\Hom(F,\ov\Q)$.  
It is convenient to attach to $(\Phi, \vphi_0)$ the unique function
$r\colon \Hom_\Q(F, \ov\Q)\to \BZ_{\geq 0}$ (denoted $\vphi\mapsto r_\vphi$)
such that 
\begin{equation}\label{SFDsig}
   (r_\varphi, r_{\varphi^{\rm c}})=
   \begin{cases}
      (\nu-1, 1),  &  \text{if $\varphi = \varphi_0$},\\
      (\nu, 0),  &  \text{if $\varphi \in \Phi \ssm \{\varphi_0\}$}.
   \end{cases}
\end{equation}
where $\vphi^{\rm c}\coloneqq  \vphi\circ {\rm c}$. (The function $r$ is a \emph{signature  type} of rank $\nu$ in the sense of  \cite[\S2.2]{RSZ4}.)

 We first consider the group (over $\Q$)
 \begin{equation}\label{G Q}
 \G^\Q\coloneqq \Res_{F_{0}/\Q} {\rm GU}(W)\times_{\Res_{F_{0}/\Q}\bG_{m, F_{0}}}\bG_{m}
 \end{equation}
 of unitary similitudes of $(W, \sform)$ with similitude factor in $\bG_{m}$. (Recall that $\bG_{m}$ denotes the multiplicative group over $\Q$.)
 
Let $\{ h_{\G^\Q} \}$ be  the $\G^\Q(\R)$-conjugacy class of the homomorphism $h_{\G^\Q}=(h_{\G^\Q,\varphi})_{\varphi\in\Phi}$, where the components $h_{\G^\Q,\varphi}$ are defined with respect to the inclusion
\begin{equation*}\label{G^Q decomp}
   G^\Q(\R) \subset \GL_{F \otimes \R}(W \otimes \R) \xra[\undertilde]{\iota_\Phi} \prod_{\varphi \in \Phi} \GL_\C(W_\varphi),
\end{equation*}
and where each component is defined on $\C^\times$ by
\[
   h_{G^\Q,\varphi}\colon z \mapsto \diag(z \cdot 1_{r_\varphi},  z^{\rm c} \cdot 1_{r_{\varphi^{\rm c}}}).
\]
Here the isomorphism $\iota_\Phi$ is given by choosing suitable orthonormal basis as we did in \S\ref{sh data}.

We single out the special case $\nu=1$. We let  $W=W_0$ be totally definite and we write $\RZ^\Q \coloneqq  \G^\Q$ (a torus over $\Q$) and $h_{\RZ^\Q} \coloneqq  h_{\G^\Q}$.  Explicitly,
\[
   \RZ^\Q = \bigl\{\, z\in \Res_{F/\Q} \bG_m \bigm| \Nm_{F/F_0} (z)\in \bG_m \,\bigr\}.
\]
The reflex field of $(\RZ^\Q,\{h_{\RZ^\Q}\})$ is $E_\Phi$, the reflex field of the CM type $\Phi$.

\subsubsection{RSZ Shimura varieties} \lb{ss: RSZ}
The Shimura varieties of \cite{RSZ3} are attached to the group
\begin{equation}\label{wtG}
   \wt \G \coloneqq  \RZ^\Q \times_{\bG_m} \G^\Q,
\end{equation}
where the maps from the factors on the right-hand side to $\bG_{m}$ are respectively given by $\Nm_{F/F_0}$ and the similitude character.  In terms of the Shimura data already defined, we obtain a Shimura datum for $\wt \G$ by defining the Shimura homomorphism to be
\[
 h_{\wt \G}\colon \C^\times \xra{(h_{\RZ^\Q}, h_{\G^\Q})} \wt \G(\R) .	
\]
Then $(\wt \G, \{h_{\wt \G}\})$ has reflex field $E \subset \ov\Q$ characterized by
\begin{equation}\label{defE}
\begin{aligned}
   \Gal(\ov\Q/E) &= \bigl\{\, \sigma\in\Gal(\ov\Q/\Q) \bigm| \sigma \circ\Phi=\Phi , \sigma^*(r)=r \,\bigr\}\\
   &=\bigl\{\, \sigma\in\Gal(\ov\Q/\Q) \bigm| \sigma \circ\Phi=\Phi , \sigma\circ\varphi_0=\varphi_0 \,\bigr\}.
\end{aligned}
\end{equation}
In other words, 
the reflex field is the composite 
$$E \coloneqq E_\Phi F,$$ 
where $F$ is embedded via $\varphi_0$. In particular, the field $F$ always embeds via $\varphi_0$ into the reflex field $E$.

\begin{remark}
The RSZ Shimura varieties are related to the unitary Shimura varieties as follows.
The torus $\RZ^\Q$ embeds naturally as a central subgroup of $\G^\Q$, which gives rise to a product decomposition
\begin{equation}\label{proddec}
   \begin{gathered}
   \begin{tikzcd}[row sep=0ex]
      \wt \G \ar[r, "\sim"]  &  \RZ^\Q \times \G^{\flat}\\
      (z, g) \ar[r, mapsto]  &  (z, z^{-1}g),
   \end{tikzcd}
   \end{gathered}
\end{equation}
where $\G^{\flat} \subset \G^\Q$ is the restriction of scalars of the unitary group \eqref{G}.  The isomorphism \eqref{proddec} extends to a product decomposition of Shimura data,
\begin{equation}\label{reltogross}
   \bigl(\wt \G, \{h_{\wt \G}\}\bigr) \cong \bigl(\RZ^\Q,\{h_{\RZ^\Q}\}\bigr) \times \bigl(\G^{\flat},\{h_{\G^{\flat}}\}\bigr).
\end{equation}
Hence,  for a decomposable compact open subgroup $K_{\wt \G}=K_{\RZ^\Q} \times K_{ \G^{\flat}}$,
there is a product decomposition
  \[
  \Sh_{K_{\wt \G}} \bigl(\wt \G, \{h_{\wt \G}\}\bigr) 
      \cong  \Sh_{K_{\RZ^\Q}} \bigl(\RZ^\Q,\{h_{\RZ^\Q}\}\bigr)  \times \Sh_{K_{ \G^{\flat}}} \bigl(\G^\flat,\{h_{ \G^{\flat}}\}\bigr),
\]
 of Shimura varieties over $E$.
\end{remark}

\subsubsection{Product Shimura varieties and the arithmetic diagonal} \lb{wt prod} 
We continue to fix
$v_0\in\Hom(F_0,\R)$ and a CM type $\Phi$. {Recall from \S~\ref{sec incoh} the fixed $\nrm\in F_0^\ts$ and  the associated set $\sV^\circ$ of pairs of (coherent or incoherent) hermitian spaces.} Let  $V\in \sV^\circ$ be a pair of incoherent hermitian spaces. {Let $V_\nu(v_0)$ be the $v_0$-nearby hermitian space of $V_\nu$, for $\nu=n,n+1$. Let $\G=\G^{(v_0)}=\U(V_n(v_0))\times \U(V_{n+1}(v_0))$ 
  be as in \S~\ref{ss:AD}.}
Similar to \eqref{wtG} we set
\begin{equation}\label{wt G}
   \wt \G \coloneqq  \RZ^\Q \times_{\bG_m} \G_{n}^\Q\times _{\bG_m} \G_{n+1}^\Q.
\end{equation}
where $\G_{\nu}^\Q$ is the similitude unitary group attached to the hermitian space $V_\nu$ as in \eqref{G Q}. We have an analogous Shimura datum with the reflex field $E=E_\Phi F$,
and  an isomorphism induced by \eqref{proddec}
\begin{equation} \lb{isom prod}
   \begin{gathered}
   \begin{tikzcd}[row sep=0ex]
      \wt \G \ar[r, "\sim"]  &  \RZ^\Q \times \G^{\flat}.
   \end{tikzcd}
   \end{gathered}
\end{equation}
In this situation, we  will always assume that the open compact $K_{\wt\G}$ is  {\em decomposable} of the form $K_{\wt \G}=K_{\RZ^\Q}\times K_{\G}=K_{\RZ^\Q}\times K_{n}\times K_{n+1}$.
In particular, we have a finite \'etale morphism $\Sh_{K_{\wt \G}}(\wt\G)\to \Sh_{K_{ \G^{}}}(\G^{})_{E}$ over $\Spec E$.

Moreover, let $\wt \RH\coloneqq  \wt\G_{n}$. Then we have a map $\jmath \colon\wt \RH\to \wt\G$ and corresponding maps
\beq\lb{wt j} \Sh_{K_{\wt \RH}} (\wt \RH)
\to 
 \Sh_{K_{\wt \G}} (\wt \G)\eeq
 that are the pullbacks  of \eqref{eq:ADC} along the projection $ \Sh_{ K_{\RZ^\Q}\times K_{\G} }(\wt \G)\to \Sh_{K_{\G}}(\G)$ given by \eqref{isom prod}.

\subsection{Moduli functors over $E$}\label{subs: RSZ}
We formulate the PEL type moduli functor for RSZ Shimura varieties, following \cite[\S3]{RSZ4}. Denote by  $\LNSch_{/R}$ the category of locally noetherian schemes over a ring $R$, and by ${\rm Sets}$ the category of sets.

\subsubsection{The torus case}
First we consider  the torus $\RZ^\Q$. The construction of \cite[\S2.2]{RSZ4}, specialized to $\nu=1$, gives  a Kottwitz PEL moduli functor  
$\LNSch_{/E}\to {\rm Sets}$, which  is represented by a  finite \'etale  stack   $M_{0, K_{\RZ^\Q}}$ over $E_{\Phi}$. Since the precise definition of this functor plays only a minor auxiliary role in this paper, we  omit it  and refer the interested readers to \cite[\S3.1]{RSZ4}; it suffices to recall that (among other data) one needs to fix a certain $F/F_0$-traceless element $\sqD \in F^\times$ adapted to the CM type $\Phi$.  The stack $M_{0, K_{\RZ^\Q}}$ is isomorphic, over $E$, to finitely many copies of the Shimura variety $\Sh_{K_{\RZ^\Q}}(\RZ^\Q)$. For our purposes, it suffices to work with a fixed copy, which we denote by $M_{0,K_{\RZ^\Q}}^\tau$.

\subsubsection{Definition of the moduli functor}
Let now  $W$ be of dimension $\nu$ as in \S\ref{sh data} and \ref{ss: wt G}.
We now present the moduli functor $M_{K_{\wt \G}}$ represented by the Shimura variety $\Sh_{K_{\wt \G} }(\wt \G)$. For simplicity, we will always assume
$$
K_{\wt \G}=K_{\RZ^\Q} \times K_{ \G}
$$ 
where $K_\G\subset \G(\A^{\infty})$ is a compact open subgroup.
For each scheme $S$ in $\LNSch_{/E}$, $M_{K_{\wt \G}}(S)$ is by definition the groupoid of tuples $(A_0,\iota_0,\lambda_0,\ov\eta_0, A,\iota,\lambda,\ov\eta)$, where
\begin{altitemize}
\item $(A_0, \iota_0, \lambda_0, \ov\eta_0)$ is an object of $M_{0,K_{\RZ^\Q}}^\tau(S)$;
\item $A$ is an abelian scheme over $S$;
\item $\iota\colon F\to\End^0(A)\coloneqq \End(A)\otimes_\BZ\Q$ is an action of $F$ on $A$ up to isogeny satisfying the Kottwitz condition
\begin{equation}\label{kottcond}
    \charac\bigl(\iota(a) \mid \Lie A\bigr) = \prod_{\varphi \in \Hom(F,\ov\Q)} \bigl(T-\varphi(a)\bigr)^{r_\varphi}
    \quad\text{for all}\quad
    a\in F,
\end{equation}
where $r$ is the  signature type $r$ of \eqref{SFDsig};
\item $\lambda$ is a quasi-polarization on $A$ whose Rosati involution satisfies 
\begin{equation}\label{Ros}
   \Ros_\lambda \bigl(\iota(a)\bigr)=\iota(\ov a)
   \quad\text{for all}\quad
   a\in F,
\end{equation}
 and
\item $\ov\eta$ is a $K_\G$-orbit (equivalently, a ${K_{\wt \G}}$-orbit, where $K_{\wt \G}$ acts through its projection $K_{\wt \G} \to K_\G$) of isometries of $\A^{\infty}_{F}/\A^{\infty}$-hermitian modules\footnote{We note that here our convention slightly deviates  from \cite{RSZ4}, see \cite[Remark 3.4]{RSZ4}: our $W$ would correspond to $V$ there and a certain $W'$ such that $\Hom(W_0,W')\simeq W$ would correspond to $W$ there (and we are using the canonical isomorphism $\U(W)\simeq \U(W')$). }

\begin{equation}\label{eta}
   \eta\colon \wh \RV(A_0, A) \isoarrow W\otimes_F\A_{F}^{\infty}.
\end{equation}
\end{altitemize}
Here, denoting by  $\wh \RV(A')$ the adelic Tate module  of an abelian variety $A'$, 
\begin{equation}\label{V(A_0,A)}
   \wh \RV(A_0,A) \coloneqq  \Hom_{\A^{\infty}_{F}}\bigl(\wh \RV(A_0), \wh \RV(A)\bigr) ,
\end{equation}
endowed with its natural $\A_{F}^{\infty}$-valued hermitian form $h$,
\begin{equation}\label{h_A def}
   h(x, y) \coloneqq  \lambda_0\i \circ y^\vee\circ\lambda\circ x\in\End_{\A^{\infty}_{F}}\bigl(\wh \RV(A_0)\bigr)=\A^{\infty}_{F}, \quad x,y \in \wh \RV(A_0,A).
\end{equation}

Finally, there are natural functors interpreting Hecke correspondences $T(KgK)$ for $g\in \G(\A^{\infty})$. 

\begin{proposition}[{\cite[Thm.~3.5]{RSZ4}}]
 The functor
$M_{K_{\wt \G}}$ is represented by  $\Sh_{K_{\wt \G} }(\wt \G)$.
\end{proposition}

\section{Integral models} \lb{sec:9}

We define and study various integral models of the RSZ unitary Shimura varieties introduced in the last section. In \S~\ref{ss:RSZ int}, we define the integral models for  RSZ Shimura varieties with vertex parahoric levels (at places unramified in $F/F_0$). In \S~\ref{ss:int model spl}, we define integral models at split places in two cases: at deeper level we obtain regular integral models, and at Iwahori levels we obtain strictly semistable integral models. In \S~\ref{ss:int mod}, we study the integral model of the product Shimura variety and resolve the product singularity. In \S~\ref{ss: abs coh}, we prove a vanishing result of the absolute cohomology for the integral model of our product Shimura variety (after resolving the singularity) after localization under the action of the Hecke algebra.

\subsection{Integral models with vertex parahoric levels}\label{ss:RSZ int}
 We continue with the setup of \S~\ref{subs: RSZ}. In particular, $W$ is a hermitian space of dimension $\nu$ and $\G=\U(W)$. We fix a rational  prime  $\ell$, and denote by $\CV_\ell$ the set of places of $F_0$ over $\ell$.
 If $\ell = 2$, then we assume that every $v \in \CV_\ell$ is {\em split} in $F$.  We follow \cite[\S4-5]{RSZ4} and \cite[\S6]{LMZ} with a slightly different formulation to define the integral models with  parahoric level structures.

 We will assume that $K_{\RZ^\Q,\ell}\subset \RZ^\Q(\Q_{\ell})$ is maximal. Then the auxiliary moduli stack  $M_{0,K_{\RZ^\Q}}$ (respectively its substack $M^\tau_{0, K_{\RZ^\Q}}$) has a natural integral model $\CM_{0,K_{\RZ^\Q}}$ (respectively $\CM_{0,K_{\RZ^\Q}}^{\tau}$), which is finite \'etale over $\Spec\sO_{E, (\ell)}$.   
 
 For each $v \in \CV_\ell$, we fix a vertex lattice $\Lambda_v\subset W_v \coloneqq  W \otimes_F F_v$ (cf. \S~\ref{sss:vert}). We recall from \cite[\S4.1]{RSZ4} a variant of the dual lattice, which is more convenient for formulating the polarization in the moduli functor.  
We endow the $F_v/F_{0,v}$-hermitian space $W_v $ with the $\Q_\ell$-valued alternating form $\tr_{F_v/\Q_\ell} \sqD\i\sform$, 
and let $\Lambda_v^\vee \subset W_v$ be the dual lattice with respect to this alternating form. We call $\Lambda_v^\vee$ the {\em absolute} dual lattice of $\Lambda_v$ (though it depends on the fixed choice of $\sqD$). We call $\Lambda_v \subset W_v$ an \emph{ absolute vertex lattice} if we have
\begin{equation*}
   \Lambda_v\subset \Lambda_v^\vee \subset\varpi_v\i\Lambda_v .
\end{equation*}
Similarly, we define the {\em absolute type} $t_{\rm abs}(\Lambda_v)$.

We have $\Lambda_v^\vee=\Lambda_v^\ast$ if $v$ is unramified over $\ell$ and $\sqD $ is a $v$-adic unit.
In general, to pass between these two versions of dual lattices, we note that, under the unramifiedness assumption on $v$ in $F/F_0$, for any  vertex lattice  $\Lambda_v\subset W_v$, there is a unique absolute vertex lattice in the lattice chain generated by $\Lambda_v$ and $\Lambda_v^\ast$ (namely, all lattices of the form $\varpi_v^i \Lambda_v, \varpi_v^i \Lambda_v^\ast$, $i\in\Z$); we denote the associated absolute vertex lattice by $\wt\Lambda_v$. Then the absolute type $t_{\rm abs}(\wt\Lambda_v)$ is either equal to  $t(\Lambda_v) $ or $\nu-t(\Lambda_v)$.  Note that the stabilizers (under the action of $G_v=\U(W)(F_{0, v})$) of the two lattices $\Lambda_v$ and $\wt\Lambda_v$ coincide and are maximal parahoric subgroup of $G_v$. The moduli functor below will be formulated in terms of $\wt\Lambda_v$.

 We assume that $K_\G \subset \G(\A_{F_{0}}^{\infty})$ is of the form $K_\G = K_\G^\ell \times K_{\G, \ell}$, where $K_\G^\ell\subset \G(\A^{\ell\infty})$ is arbitrary, and where
\[
   K_{\G, \ell}=\prod_{v\in \CV_\ell}K_{\G, v}\subset \G(F_{0,\ell}) = \prod_{v \in \CV_\ell}  G_{v},
\]
with 
\begin{equation}\label{K_G,p decomp}
   K_{\G, v} \coloneqq   {\rm Stab}_{G_{v}} (\Lambda_v) . 
\end{equation}

We then define 
 $\CM_{K_{\wt \G}}$ as the functor that associates to each scheme $S$ in $\LNSch_{/\sO_{E,(\ell)}}$ the groupoid of tuples 
$(A_0,\iota_0,\lambda_0,A,\iota,\lambda,\ov\eta^\ell)$, where
\begin{altitemize}
\item $(A_0,\iota_0,\lambda_0)$ is an object of  $\CM_{0,K_{\RZ^\Q}}^{\tau}(S)$;
\item $A$ is an abelian scheme over $S$;
\item $\iota\colon \sO_{F,(\ell)} \to \End_{(\ell)}(A)$ is an action up to prime-to-$\ell$ isogeny satisfying the Kottwitz condition \eqref{kottcond} on $\sO_{F,(\ell)}$;
\item $\lambda \in \Hom(A,A^\vee)_{\Z_{(\ell)}}$ is a quasi-polarization on $A$ whose Rosati involution satisfies condition \eqref{Ros} on $\sO_{F,(\ell)}$; and
\item $\ov\eta^\ell$ is a $K_\G^\ell$-orbit of isometries of $\A^{\ell\infty}_{F}/\A_{F_0}^{\ell\infty}$-hermitian modules
\begin{equation}\label{levelprimetop}
   \eta^\ell\colon \wh \RV^\ell(A_0,A) \isoarrow W \otimes_F \A_{F}^{\ell\infty},
\end{equation}
where
\begin{equation}
   \wh \RV^\ell(A_0,A) \coloneqq  \Hom_{\A_{F}^{\ell\infty}}\bigl(\wh \RV^\ell(A_0), \wh \RV^\ell(A)\bigr),
\end{equation}
and where the hermitian form on $\wh \RV^\ell(A_0,A)$ is the obvious prime-to-$\ell$ analog of \eqref{h_A def}.
\end{altitemize}
We impose the following further conditions on the above tuples.
\begin{altenumerate}
\item\label{sglob cond i} Consider the decomposition of $\ell$-divisible groups
\begin{equation}\label{decofpdivgp}
   A[\ell^\infty] =  \prod_{v \in \CV_\ell} A[v^\infty]
\end{equation}
induced by the action of $\sO_{F_0}\otimes\BZ_\ell \cong \prod_{v \in \CV_\ell} \sO_{F_0,v}$.
Since $\Ros_\lambda$ is trivial on $\sO_{F_0}$, $\lambda$ induces a polarization $\lambda_v \colon A[v^\infty] \to A^\vee[v^\infty] \cong A[v^\infty]^\vee$ of $\ell$-divisible groups for each $v$. The condition we impose is that $\ker\lambda_v$ is contained in $A[\iota(\pi_v)]$ and has  order $\#(\wt\Lambda_v^\vee/\wt \Lambda_v)$ for each $v \in \CV_\ell$.
\item\label{sglob cond ii} We require that the \emph{sign condition} and the {\em Eisenstein condition} hold; we omit the definitions and refer  to \cite[\S5]{RSZ4} (for $v$ unramified over $\Q$) and \cite[\S6]{LMZ}.
\end{altenumerate}
The morphisms in the groupoid  $\CM_{K_{\wt \G}}(S)$ are the obvious ones.

We have the following result due to \cite[\S5.3]{RSZ4} and \cite{LMZ}.
\begin{proposition}\label{th:reg vert}
The  stack $\CM_{K_{\wt \G}}$ is Deligne--Mumford, and regular with strictly semistable reduction at all places $u$ of $E$ above $\ell$, provided  that the place $u$ is unramified over $F$. The generic fibre of $\CM_{K_{\wt \G}}$ is $M_{K_{\wt \G}} $. 
It is smooth over $\Spec \sO_{E,(\ell)}$ if the lattices $\Lambda_v$ have type $0$ or $\nu$ for every $v\mid \ell$.  
\end{proposition}
\begin{proof}
When $v$ is  unramified over $\Q$, the assertion follows from \cite[\S5.3]{RSZ4}. When $v$ is ramified over $\Q$, we use the local model diagram \cite[\S6.4,~(6.16),~(6.17)]{LMZ} to reduce the regularity and semistability to the corresponding statements on  the local model and then apply \cite[\S3.2, Theorem 3.14 and 3.17]{LMZ}, noting that in our case the place $u$ of the reflex field  is assumed to be unramified over $F$.
\end{proof}

Finally, there are natural functors interpreting Hecke correspondences $\hT(f^{\ell})$ for all $f^{\ell}\in \sH(\G(\A^{\ell\infty}))_{K^\ell_{\G}}$.
 The correspondences $\hT(f^{\ell})$ are all \'etale.

\subsection{More integral models at split places}
\label{ss:int model spl}
We need regular integral models for deeper levels at split places. We will consider two cases: the principal congruence subgroup case and the Iwahori case.

\subsubsection{Setup}  We continue with the setup of \S~\ref{ss:RSZ int} and we  fix a place $v \in \CV_\ell$ that splits in $F$, say $v = w \ov w$. Let $u\colon E\to  \ov\Q_\ell$ be a place of $E$ above $v$; we will assume that $E_{u}$ is unramified over $F_{0,v}$. Let $\wt u\colon  \ov\Q\incl \ov\Q_\ell$ be an embedding extending $u$. Then $\wt u$ induces a bijection $\Hom(F,\ov\Q)\simeq \Hom(F,\ov\Q_\ell)$. Let
$ \Hom_{w}(F,\ov\Q)$ be the subset of $\Hom(F,\ov\Q)$ consisting of $\varphi\in \Hom(F,\ov\Q)$ such that $\wt u\circ\varphi$ induces $w$.  The set $ \Hom_{w}(F,\ov\Q)$ depends only on $u$ but not on the choice of $\wt u$. Note that the distinguished element $\varphi_0$ belongs to $\Hom_w(F,\ov\Q)$. We will assume that  the \emph{matching condition} between the CM type $\Phi$ and the chosen place $u$ of $E$ is satisfied:
\begin{equation}\label{cond CM}
\Hom_w(F,\ov\Q)\subset \Phi,
\end{equation}
cf.  \cite[\S4.3]{RSZ3}. Note that, for our signature type \eqref{SFDsig}, this is equivalent to the condition that
the restriction $r|_{\Hom_w(F,\ov\Q)}$ of the signature function is of the form
\begin{equation}\label{drinfeld function}
   r_\varphi =
   \begin{cases}
      \nu-1,  &  \text{$\varphi = \varphi_0 \in \Hom_w(F,\ov\Q)$};\\
      \nu,  &  \varphi \in \Hom_w(F,\ov\Q) \ssm \{\varphi_0\}.
   \end{cases}
\end{equation}

\subsubsection{Principal congruence subgroups}
We now recall from  \cite[\S4.3]{RSZ3} the moduli problem in the case of principal congruence subgroups. Let $m$ be a nonnegative integer, and define $K_{\G,v}^m$ to be the principal congruence subgroup mod $\fkp_v^m$ inside $K_{\G,v}$, where $\fkp_v$ denotes the prime ideal in $\sO_{F_0}$ determined by $v$.  Let
\[
   K_{\wt \G}^m \coloneqq  K_{\RZ^\Q}\times K_\G^\ell \times K_{\G,v}^m \times \prod_{v'\in\CV_\ell\ssm \{v\}}K_{\G, v'} \subset K_{\wt \G}.
\]
Then one can extend the definition of $\CM_{K_{\wt \G}, \sO_{E,u}}$ to the case of the level subgroup $K_{\wt \G}^m$ by adding a Drinfeld level-$m$ structure at $v$.
More precisely, consider the factors occurring in the decomposition \eqref{decofpdivgp} of the $\ell$-divisible group $A[\ell^\infty]$,
\begin{equation}\label{dec pdiv}
   A[v^\infty] = A[w^\infty] \times A[\ov w^\infty] .
\end{equation}
  The condition \eqref{drinfeld function} implies that $A[\ov w^\infty]$ is a one-dimensional formal $\sO_{F,w_0}$-module.  We introduce  $T_{\ov w}(A_0,A)[w^m]\coloneqq  \uHom_{\sO_{F, \ov w}}(A_0[\ov w^m],A[\ov w^m])$ and $T_{\ov w}(A_0,A)\coloneqq \varinjlim_{m} T_{\ov w}(A_0,A)[w^m]$. Note that $T_{\ov w}(A_0,A)$ is  a $1$-dimensional formal $\sO_{F,w_0}$-module.  The datum we add to the moduli problem is an $\sO_{F,\ov w}$-linear homomorphism of finite flat group schemes, 
\begin{equation}\label{dr level}
   \phi\colon \pi_{\ov w}^{-m}\Lambda_{\ov w}/\Lambda_{\ov w} \to   T_{\ov w}(A_0,A)[\ov w^m],
\end{equation}
which is a Drinfeld $\ov w^m$-structure on the target.  Here, $\Lambda_{\ov w}$ is the summand attached to $\ov w$ in the natural decomposition 
\begin{equation}\label{dec Lam v}
\Lambda_v = \Lambda_w \oplus \Lambda_{\ov w}
\end{equation} 
with $\Lambda_v$ the vertex lattice at $v$ chosen in \S\ref{ss:RSZ int}.  See \cite[\S4.3]{RSZ3} (which we note interchanges the roles of $w$ and $\ov w$) for more details.

Then by \cite[Theorem 4.7]{RSZ3},
the moduli problem $\CM_{K^m_{\wt \G}}$  is relatively representable by a finite flat morphism to  $\CM_{K_{\wt \G}}$ and consequently it coincides with the normalization of $\CM_{K_{\wt \G}}$ in the generic fiber of $\CM_{K^m_{\wt \G}}$. It is regular and flat over $\Spec \sO_{E,(u)}$.    Furthermore, the generic fiber $\CM_{K^m_{\wt \G}} \times_{\Spec \sO_{E,(u)}} \Spec E$ is canonically isomorphic to $M_{K^m_{\wt \G}}$.

\subsubsection{Iwahori subgroups}
We will also need the Iwahori case. 
For simplicity we assume that the vertex lattice $\Lambda_v$ in \eqref{dec Lam v} is selfdual. We now choose a chain of $\sO_{F,w}$-lattices
$$
\Lambda_{\ov w}=\Lambda^{(0)}_{\ov w}\subset  \Lambda^{(1)}_{\ov w}\subset\cdots \subset \Lambda^{(\nu)}_{\ov w}=\pi_w^{-1}\Lambda_{\ov w},
$$
where each inclusion has colength one. Equivalently, we choose a full flag in the $k_v$-vector space $\Lambda_w/\pi_w \Lambda_w$.
This chain determines uniquely a chain of vertex $\sO_{F,v}=\sO_{F,w}\times \sO_{F,\ov w}$-lattices $\Lambda_v^{(i)}\coloneqq  \Lambda_w \oplus \Lambda^{(i)}_{\ov w}, 0\leq i\leq \nu$. The stabilizer of the chain  $\Lambda_v^{(i)}$ is an Iwahori subgroup ${\rm Iw}_v$ of ${\rm Stab}(\Lambda_v)$. To the moduli problem $\CM_{K_{\wt \G}, \sO_{E,u}}$, we add the datum of a chain of isogenies of $\sO_{F,\ov w}$-divisible modules 
\begin{align}\label{eq:Iw level}
\CG_0=T_{\ov w}(A_0,A) \to \CG_1\to \cdots \to\CG_\nu=\CG_0/\CG_0[\ov w]
\end{align} with equal heights $\# k_v$. An equivalent datum is an ${\rm Iw}_v$-orbit of  the  Drinfeld level structure 
$$ \phi\colon \pi_{\ov w}^{-1}\Lambda_{\ov w}/\Lambda_{\ov w} \to   T_{\ov w}(A_0,A)[\ov w].$$
The resulting moduli functor is then denoted by $\CM_{K_{\wt \G}^{\rm Iw_v}}$, where $K_{\wt \G}^{\rm Iw_v}$ denotes the compact subgroup of $K_{\wt \G}$ with the Iwahori factor at $v$. Then the moduli problem $\CM_{K_{\wt \G}^{\rm Iw_v}}$  is relatively representable by a finite flat morphism to  $\CM_{K_{\wt \G}}$ and consequently it coincides with the normalization of $\CM_{K_{\wt \G}}$ in the generic fiber of $\CM_{K^{\rm Iw_v}_{\wt \G}}$. It is regular, proper and flat over $\Spec \sO_{E,(u)}$. Moreover, by the theory of local models, the scheme $\CM_{K_{\wt \G}^{\rm Iw_v}}$ has strictly semistable reduction  over $\Spec \sO_{E,(u)}$ (namely, its generic fiber is smooth and every closed point of the special fiber admits an open neighborhood which is
 smooth over the scheme $\Spec  \sO_{E,(u)}[x_1,\cdots, x_m]/(\prod_{i=1}^m x_i- \varpi)$ for some $m\geq 1$, cf. \cite[Prop. 1.3]{Hartl}). 
 There is a stratification of the special fiber $\CM_{K_{\wt \G}^{\rm Iw_v}}\otimes k_u$, where $k_u$ denotes the residue field of $\sO_{E,(u)}$:
\begin{equation}\label{eq:str Iw}
\CM_{K_{\wt \G}^{\rm Iw_v}}\otimes k_u=\bigcup_{1\leq i\leq n}  \CM_{K_{\wt \G}^{\rm Iw_v},k_u, i},
\end{equation}
where $ \CM_{K_{\wt \G}^{\rm Iw_v},k_u, i}$ is the closed subscheme where the kernel of the isogeny $\CG_{i-1}\to \CG_i$ in \eqref{eq:Iw level} is connected, cf. \cite[\S3]{TY} for a similar case.   By \cite[Prop. 3.4]{TY} (or rather its proof), each  $\CM_{K_{\wt \G}^{\rm Iw_v},k_u, i}$ is smooth over $\Spec k_u$.

\subsubsection{Hecke correspondences}
We recall from \cite[\S4.3]{RSZ3} that, in each of the above two cases (principal and Iwahori level), 
 there are natural functors interpreting Hecke correspondences attached to functions ${\bf 1}_{K g K}$  for any  $g\in \G(\A^{\ell\infty})\times\G(F_{0,v}) $, where we simply denote $K=K_\G$: 
\begin{equation}\label{hecke G}
\begin{gathered}
   \xymatrix{
	     & \CM_{K'_{\wt \G} } \ar[dl]_-{\text{nat}_1} \ar[dr]^-{\text{nat}_g} &\\
	   \CM_{K_{\wt \G}}  & &  \CM_{K_{\wt \G}}
	}
\end{gathered}
\end{equation}
where $K'_{\wt \G}=K_{\RZ^\Q}\times K'_{\G}$, and $K'_{\G}=K_{ \G}\cap gK_{\G}g^{-1}$. We refer to \cite[\S4.3]{RSZ3} for the unexplained notation. ({Note that in \emph{loc.\ cit.}, the authors only consider the case of a principal congruence subgroup $K_{\G,v}=K_{\G,v}^m$. The Iwahori case is similar and may be reduced to the case $K_\G=K_\G^m$ as follows. We can factorize $[{\Iw_v g_v \Iw_v}]$ as $ e_{\Iw_v } \star [{K_v g_v K_v} ]\star e_{\Iw_v }$ for some $K_v=K_{\G,v}^m\subset \Iw_v$, and accordingly we define the correspondence  for $[{\Iw_v g_v \Iw_v}]$ as the composition of the three factors:  the middle one  is as in \emph{loc.\ cit.}, and the other two are given by the natural map from the principal level to the Iwahori level.})

Both maps $\text{nat}_1$ and $\text{nat}_g$ are finite flat, and \'etale if $g_\ell=1$. The Hecke correspondence \eqref{hecke G} induces an endomorphism  (by the usual pull-back and then push-forward maps) on the group of cycles (with $L$-coefficients), rather than merely cycles modulo rational equivalence. This endomorphism is independent of the choice of $K'_{\wt \G}$ in the diagram above. The resulting map $$
\begin{gathered}
   \xymatrix{
  T\colon \sH(\G(\A^{\ell\infty}),L)_{K^\ell_\G}\ar[r]&   \EC(\CM_{K_{\wt \G}})_{L}
   	}
\end{gathered}$$ 
is a ring homomorphism.\footnote{However, we do not know if the assertion remains true for the full Hecke algebra $ \sH(\G(\A^{\infty}),L)_{K_\G}$. When $m=0$, the recent work of Li--Mihatsch \cite[Proposition 3.4]{LM} shows that the assertion holds.}
Moreover, in the Iwahori case, the away-from-$\ell$ Hecke correspondences preserve the stratification \eqref{eq:str Iw}.

\subsection{Moduli functors for the product Shimura varieties}
\label{ss:int mod}
We now continue with the setup of \S~\ref{wt prod}. In particular, we have a pair $(V_n,V_{n+1})$ of standard indefinite hermitian spaces (relative to the fixed place $v_0\in\Hom(F_0,\R)$), and we have the product unitary group $\G= \G_{n}\times \G_{n+1}= \RU(V_{n})\times \RU(V_{n+1})$, and the related group $\wt \G$ defined by \eqref{wt G}.
It is now straightforward to extend the constructions in \S~\ref{subs: RSZ}, \S~\ref{ss:RSZ int} and \S~\ref{ss:int model spl} to the product unitary group. These are analogous moduli functors over $E$ and over $\sO_{E, (\ell)}$. For example, the $\ell$-integral model may be succinctly defined as 
\begin{equation}\label{eq:prod int mod}
\CM_{K_{\wt\G}}=\CM_{K_{\wt \G(V_{n})}}\times_{\CM_{0}^\tau} \CM_{K_{\wt \G(V_{n+1})}},
\end{equation}
where $K_{\wt \G(V_{\nu})}=K_{\RZ^\Q}\times K_{\nu}$ for $\nu\in\{n,n+1\}$.  Here we denote  by $\wt \G(V_{\nu})$ the group associated to $V_\nu$  defined by \eqref{wtG}.

The product $\CM_{K_{\wt\G}}$ may no longer be regular even if both factors are regular, and we may need to resolve the product singularity.
We will need to study two cases: the vertex parahoric case at an inert place, and the Iwahori case at a split place.

\subsubsection{Vertex parahoric level at an inert place}\label{sss:ver level}
We first consider the vertex parahoric case from \S\ref{ss:RSZ int}. Fix  a place $v \in \CV_\ell$ that is inert in $F$, let $w$ denote the unique place of $F$ above~$v$. Let $\e=\e_v\in \{0, 1\}$ have the same parity as $v(\nrm)$, where $\nrm\in F_0^\ts$ is the hermitian norm of the special vector $u$ as fixed in \S~\ref{sec incoh}. We consider a vertex parahoric subgroup $K_v=K_{n,v}\times K_{n+1,v}$ of type $(t,t+\epsilon)$ in the sense of \S\ref{sss:vert}.
   
Let $u:E\to  \ov\Q_\ell$ be a place of $E$ above $v$, and  further assume that $E_u$ is unramified over $F_{0,v}$. In this case, the integral models $\CM_{K_{\wt \G(V_n)}}$ and $\CM_{K_{\wt \G(V_{n+1})}}$ have strictly semistable reduction over $\Spec \sO_{E,u}$; and $\CM_{K_{\wt \G(V_n)}}$ (resp. $\CM_{K_{\wt \G(V_{n+1})}}$) is smooth over $\Spec \sO_{E,u}$ only when $t\in\{0,n\}$ (resp. $t+\epsilon\in\{0,n+1\}$); see Proposition \ref{th:reg vert}. When $\CM_{K_{\wt \G(V_n)}}$ or $\CM_{K_{\wt \G(V_{n+1})}}$ is non-smooth over $\Spec \sO_{E,u}$,  
 its special fiber admits a ``balloon--ground" stratification (\cite[\S5.2]{LTXZZ} for $t=1$,  \cite[\S7]{zhiyu} and \cite[\S7]{LMZ} for general $t$): the special fiber is a union of two Weil divisors
 \begin{equation}\label{eq:bu-circ}
 \CM_{K_{\wt \G(V_n)}, k_u}= \CM_{K_{\wt \G(V_n)},k_u}^\circ\cup \CM_{K_{\wt \G(V_n)},k_u}^\bullet
\end{equation}
where
  the first one $\CM_{K_{\wt \G(V_n)},k_u}^\circ$ is  called the balloon stratum and the second one $ \CM_{K_{\wt \G(V_n)},k_u}^\bullet$ is called the ground stratum. (When $t\in\{0, n\}$ we understand that the balloon stratum is empty.) When  $\CM_{K_{\wt \G}}$ is not regular, we let $\wt\CM_{K_{\wt \G}}$ be the blow up along the product of the balloon strata of the two factors, and denote the blow-up morphism \begin{equation}
   \begin{gathered}
   \begin{tikzcd}[row sep=0ex]
      \pi\colon \wt\CM_{K_{\wt\G}}  \ar[r]  &  \CM_{K_{\wt\G}.}
   \end{tikzcd}
   \end{gathered}
\end{equation}
For $(?_n,?_{n+1})\in \{\circ,\bullet\}^2$, we denote by $\wt\CM_{K_{\wt\G},k_u}^{(?_{n},?_{n+1})}$ the strict transform of  $\CM_{K_{\wt \G(V_n)},k_u}^{?_n}  \times \CM_{K_{\wt \G(V_{n+1})},k_u}^{?_{n+1}} $.
 For later reference, we record the following result from \cite[Lemma 5.11.3]{LTXZZ} for $t=1$,  \cite{zhiyu} for general~$t$. 

\begin{proposition}\label{prop:small in}
   The scheme $\wt\CM_{K_{\wt\G}}$ is regular 
   with strictly semistable reduction
\begin{equation}\label{eq:str in}
\wt\CM_{K_{\wt \G}}\otimes k_u=\bigcup_{(?_n,?_{n+1})\in \{\circ,\bullet\}^2 }  \wt\CM_{K_{\wt\G},k_u}^{(?_{n},?_{n+1})},
\end{equation}
where the schemes $\wt\CM_{K_{\wt\G},k_u}^{(?_{n},?_{n+1})}$ are smooth of pure dimension $2n-1$.

The map $\pi$ is small,
    i.e., a proper birational morphism with the property that 
    $${\rm codim}\{z\in  \CM_{K_{\wt \G}}\mid \dim\pi^{-1}(z)\geq i \}\geq 2i+1,
    $$ for all $i\geq 1$.

\end{proposition}

\subsubsection{Iwahori level at a split place}
 Fix as in \S\ref{ss:int model spl} a place $v \in \CV_\ell$ that splits in $F$ into $v = w \ov w$ and we let $u\colon E\to  \ov\Q_\ell$ be a place of $E$ above $v$. We further assume that $E_u$ is unramified over $F_{0,v}$. Then the integral model $\CM_{K_{\wt\G}}$ over $\Spec \sO_{E,(u)}$ is regular if one of the two compact open subgroups  $ K_{n,v}$ and $ K_{n+1,v}$ is hyperspecial. When both   $ K_{n,v}$ and $ K_{n+1,v}$  are Iwahori, $\CM_{K_{\wt\G}}$ is no longer regular and we need to resolve the product singularity. More precisely, we consider the fiber product of the stratifications from \eqref{eq:str Iw}
\begin{equation}\label{eq:Zij}
 \CM_{K_{\wt\G}, k_u, (i,j)}\coloneqq \CM_{K_{\wt \G(V_{n})},k_u,i}\times_{\CM_{0}^\tau} \CM_{K_{\wt \G(V_{n+1})},k_u,j}.
\end{equation}
We choose an ordering of the set $\{(i,j)\mid 1\leq i\leq n,1\leq j\leq n+1\}$, and rename the component $\CM_{K_{\wt\G}, k_u, (i,j)}$ as $\CM_{K_{\wt\G}, k_u, r}$ 
where $ 1\leq r\leq n(n+1)$. Let $\CM_{K_{\wt\G}}^{(0)}\coloneqq \CM_{K_{\wt\G}}$ and for $ 1\leq r\leq n(n+1)$   let $\CM_{K_{\wt\G}}^{(r)}$  be the blow-up of $\CM_{K_{\wt\G}}^{(r-1)}$  along (the strict transforms of) $\CM_{K_{\wt\G}, k_u, r}$. We write   $\wt\CM_{K_{\wt\G}}$ for  $\CM_{K_{\wt\G}}^{(n(n+1))}$, and   $\wt \CM_{K_{\wt\G}, k_u, (i,j)}$ for the strict transform of $\CM_{K_{\wt\G}, k_u, (i,j)}$. 
The composition of the natural blow-up maps is denoted as
\begin{equation}
   \begin{gathered}
   \begin{tikzcd}[row sep=0ex]
      \pi\colon \wt\CM_{K_{\wt\G}}  \ar[r]  &  \CM_{K_{\wt\G}.}
   \end{tikzcd}
   \end{gathered}
\end{equation}
(We also note that the resolution in the inert case earlier can also be view a special case of the current procedure: one simply orders the components such that the first one is the product of the balloon strata.)

\begin{proposition}\label{prop:small spl}
The scheme $\wt\CM_{K_{\wt\G}}$ is regular with strictly semistable reduction
\begin{equation}\label{eq:str Iw 2}
\wt\CM_{K_{\wt \G}, k_u}=\bigcup_{1\leq i\leq n, \atop 1\leq j\leq n+1 }  \wt\CM_{K_{\wt \G},k_u, (i,j)},
\end{equation}
where the schemes  $ \wt\CM_{K_{\wt \G},k_u, (i,j)}$ are smooth of pure dimension $2n-1$.
 The map $\pi$ is a small map.
\end{proposition}
\begin{proof}
The first part is well-known, for example  see \cite[Prop. 2.1]{Hartl} or \cite{GS}. For the smallness,  we use the explicit description as in the proof of \cite[Prop. 2.1]{Hartl}.
Consider a point $P$ on the special fiber of $\CM_{K_{\wt \G}, k_u}$ with an open neighborhood that is
 smooth over $$\Spec  \sO_{E,(u)}[x_1,\cdots, x_r, y_1,\cdots y_s]/(\prod_{i=1}^r x_i- \varpi,\prod_{j=1}^s y_j- \varpi  )$$ for some (uniquely-determined) integers $r,s\geq 1$ such that $P$ lies over a point defined by $x_i=0,y_j=0, 1\leq i\leq r,1\leq j\leq s$. Then keeping track of the steps of the blow-ups in {\it loc. cit.} shows that the dimension of the fiber  $\pi^{-1}(P)$ is $\min\{r-1,s-1\}$.
 Note that the locus of $P$ with fixed $r,s\geq1 $ is contained in the union of 
 $$( \CM_{K_{\wt \G(V_{n})},k_u,i_1}\cap\cdots\cap \CM_{K_{\wt \G(V_{n})},k_u,i_r})\times   ( \CM_{K_{\wt \G(V_{n+1})},k_u,j_1}\cap\cdots\cap \CM_{K_{\wt \G(V_{n+1})},k_u,j_s})
 $$
for all possible $1\leq i_1<\cdots< i_r\leq n, 1\leq j_1<\cdots < j_s\leq n+1$. The codimension of such loci in $\CM_{K_{\wt\G}}$  is $r+s-1\geq 2 \min\{r-1,s-1\}+1$, which proves the smallness of the map $\pi$.
\end{proof}

This procedure depends on the choice of an ordering and therefore it is not canonical. Nevertheless the smallness of $\pi$ shows that the resolution has the property that $\pi_\ast \Q_{p}\simeq {\rm IC}$, the latter being the intersection complex of the $\Q_p$-sheaf (for $p\neq \ell$). Moreover, the resulting $\wt\CM_{K_{\wt\G}}$ and each $  \wt\CM_{K_{\wt \G},k_u, (i,j)}$ still have an action of  $\sH(\G(\A^{\ell\infty}))_{K^\ell_{\G}}$ by \'etale correspondences.

\subsubsection{Integral arithmetic diagonal} \lb{ss:ADC int mod}
We have an integral model 
\beq\lb{eq: int AD} \jmath\colon
\xymatrix{\CM_{K_{\wt \RH}} \ar[r]& \CM_{K_{\wt\G}} }
\eeq
 of the morphism \eqref{wt j}.
 In the two cases discussed above, over a place $u$ of $E$, we have the small resolution $\wt \CM_{K_{\wt\G}}$ and we denote by 
 \beq \lb{eq: int AD blow}
 \wt{\jmath}
 \colon \wt\CM_{K_{\wt \RH}} \to \wt \CM_{K_{\wt\G}}\eeq
  the strict transform of $\CM_{K_{\wt \RH}}$ along the resolution morphism. For uniformity of notation, we will put $ \wt\CM_{K_{\wt \RH}}\coloneqq   \CM_{K_{\wt \RH}} $,  $  \wt\CM_{K_{\wt\G}} \coloneqq    \CM_{K_{\wt\G}}$, $\wt \jmath\coloneqq \jmath$ in the cases where those schemes 
are already regular.

\subsection{Vanishing of absolute cohomology}
\lb{ss: abs coh}
We continue with the setup of \S\ref{ss:int mod}. Our goal in this subsection is to prove the vanishing of the top-degree absolute cohomology of the scheme  $\wt\CM_{K_{\wt \G}}$, for certain  levels $K_{\wt \G}$ and after suitable localizations.

\subsubsection{Correspondences that annihilate the cohomology}
We  use some  general results from \cite{LL1,LL2}, which we now  recall. Let $L$ be a finite extension of $\Q_{p}$. Following \cite[Appendix B]{LL1}, we define a
 \emph{commutative $L$-algebra of \'etale correspondences} on a scheme $\sX$   to be  a commutative $L$-algebra $\BT$ equipped with a homomorphism $\BT\to {\EC}(\sX)_{L}$.

Let $\sX$  be a regular scheme, proper and flat   of relative dimension $2n-1$ (not necessarily strictly semistable) over the ring of integers of a non-archimedean local field, with residue field $k$;  we assume that the generic fiber $X$ is smooth. Let $\BT$ be a commutative $L$-algebra of \'etale correspondences on $\sX$ and let $\fkm\subset \BT$ be a  maximal ideal.  Let $Y$ denote the {\it reduced} special fiber of $\sX$. Assume that 
there is a stratification $Y=Y^{[0]}\supset \cdots \supset Y^{[d]}$ by closed subschemes and, for each $0\leq i\leq d$, a refinement of $Y^{(i)}\coloneqq Y^{[i]}\setminus Y^{[i+1]}$ as a disjoint union of open and closed subschemes  of $Y^{(i)}$ of pure dimension $d_i$:
$$
Y^{(i)}=\coprod_{M\in\fkS^{i}} Y^{(M)},
$$ 
over a finite set of indices $\fkS^{i}$, such that 
\begin{enumerate}
\item For every $i$ and $M\in\fkS^{i}$, denoting by $Y^{[M]}$ the Zariski closure of $Y^{(M)}$, the scheme  $Y^{[M]}$ is {\em smooth}  and it is a disjoint union $\coprod_{M'\in \fkS_M} Y^{(M')}$ where $\fkS_M$ is a subset of $\fkS\coloneqq \coprod \fkS^{(i)}$;
\item For every $i$ and $M\in\fkS^{i}$, the scheme $ Y^{(M)}$ is stable under the action of $\BT$.

\end{enumerate}

\begin{proposition}[Li--Liu]\label{th:LL}
Under the above assumptions, 
if we further suppose that either  of the following two  conditions holds:
\begin{enumerate}
\item 
$H^j(Y^{[M]}\otimes_k \ov k, L)_\fkm=0$ whenever $j\neq \dim Y^{[M]}$ for every $M\in \fkS$,
\item   $H^{2n}(X, L(n))_\fkm=0$ and $H^{j}(Y^{(i)}\otimes_k \ov k,L)_\fkm=0$   whenever $j\leq \dim Y^{(i)}-{\rm codim}_{\sX}Y^{(i)} $ for every~$i$,
\end{enumerate}
then 
$H^{2n}(\sX, L(n))_\fkm=0$.
\end{proposition}

\begin{proof}
Case (1). The vanishing assumption $H^j(Y^{[M]}\otimes_k \ov k, L)_\fkm=0$ is the assertion of \cite[Prop. 4.25]{LL2}. 
 The proof of \cite[Theorem 4.21]{LL2} applies verbatim to show that the assumptions imply the desired vanishing $H^{2n}(\sX, L(n))_\fkm=0$.
 
 Case (2). This is  \cite[Corollary B.15]{LL1}. We sketch their proof for the convenience of the reader.

By the assumption $H^{2n}(X, L(n))_\fkm=0$ and the exact sequence
\begin{equation*}
\xymatrix{    H^{2n}_{Y}(\sX)  \ar[r]&H^{2n}(\sX) \ar[r]  & H^{2n} ( X)
}
\end{equation*}
it suffices to show $ H^{2n}_{Y}(\sX)_\fkm=0.$ 
This follows from an induction using 
\begin{itemize}
\item
the exact sequences
\begin{equation*}
\xymatrix{    H^{2n}_{Y^{[j+1]}}(\sX)  \ar[r]&H^{2n}_{Y^{[j]}}(\sX) \ar[r]  & H^{2n}_{Y^{(j)}} (\sX\setminus Y^{[j+1]}),
}
\end{equation*}
\item the absolute purity theorem of Gabber  $H^{2n}_{Y^{(j)}} (\sX\setminus Y^{[j+1]})\simeq  H^{2n-2n_j}(Y^{(j)})$   for the regular local immersion $Y^{(j)}\incl  \sX\setminus Y^{[j+1]}$ of codimension $n_j$,
\item the Hochschild--Serre spectral sequence $H^{r}(k,H^s(Y^{(j)}\otimes_k \ov k)(n))\imp  H^{r+s}(Y^{(j)})(n)$.  In particular, it suffices to replace (2) by the weaker  assumption: $$H^{2n-2n_j}(Y^{(j)}\otimes_k \ov k,L)_\fkm=H^{2n-2n_j-1}(Y^{(j)}\otimes_k \ov k,L)_\fkm=0,$$ namely $H^{d_{Y^{[j]}}-c_{Y^{[j]}}-i}(Y^{(j)}\otimes_k \ov k,L)_\fkm=0$ for $i=0,1$ where $d_{Y^{[j]}}$ and $c_{Y^{[j]}}$ denote respectively the dimension of $Y^{[j]}$ and the codimension of $Y^{[j]}$ in $\sX$.
\end{itemize}
\end{proof}

\subsubsection{The vanishing result}
We consider the scheme  $\wt\CM_{K_{\wt \G}}$ over $\Spec \sO_{E,u}$ of \S~\ref{ss:int mod} where $v$ and $K_{\G,v}$ are in one of the following cases:
\begin{enumerate}
\item the split-(Drinfeld-level, hyperspecial) case: the place $v$ is split in $F$ and in the product \eqref{eq:prod int mod} one of the two factors has Drinfeld-level for some integer $m$ and the other has hyperspecial level;
\item  the split-(Iwahori, Iwahori) case: the place $v$ is split in $F$ and in the product \eqref{eq:prod int mod} both factors have Iwahori level; in this case $\wt\CM_{K_{\wt \G}}$ is the small resolution in Proposition \ref{prop:small spl};
\item the inert-vertex-parahoric case: the place $v$ is inert in $F$ and in the product \eqref{eq:prod int mod} both factors have vertex-parahoric levels (of type $(t, t+\epsilon)$); in this case $\wt\CM_{K_{\wt \G}}$ is the small resolution in Proposition \ref{prop:small in}.
\end{enumerate}

\begin{proposition}\label{thm:van}
Let $S$ be a finite set of nonarchimedean places of $F_{0}$ containing those above $\ell$ and $p$ and those  where $K_{\G}$ is not maximal hyperspecial, and let 
$$\BT=\BT^{{\rm spl}, S}\coloneqq \bigotimes_{v\notin S\atop \textup{split}} \sH(G_{v}, L)_{K_v} \subset \sH(\G(\A^{S}), L)^{\circ}_{K^{S}}.$$
Let $\fkm\subset \BT$ be the maximal ideal attached  to a representation $\pi\in \wt\sC(\G)(L)$.
 Suppose we are in one of the above three cases, and suppose moreover that the following hold:
\begin{enumerate}
\item In the split-(Iwahori, Iwahori) case, the representation $\pi_v$ is a (tempered) principal series.
\item
 In the inert-vertex-parahoric case, the type 
$(t_n, t_{n+1})$ is one of the following:
$$(0,0),\quad (0,1),\quad(1,1),\quad (n-1,n),\quad (n,n),\quad (n,n+1). $$
\end{enumerate} Then we have
\begin{equation}\label{eq:coh Iw}
H^{2n}(\wt\CM_{K_{\wt \G}},L(n))_\fkm=0.
\end{equation}
\end{proposition}

\begin{proof}
We wish to apply the vanishing theorem of Li--Liu given in   Proposition \ref{th:LL}.
For this,  we need to specify a stratification of the reduced special fiber of $\wt\CM_{K_{\wt \G}}$.
In the split-(Drinfeld-level, hyperspecial) case, for simplicity we consider the case in which the (possibly) nontrivial Drinfeld level appears in the first factor $\CM_{K_{\wt \G(V_{n})}}$. Then the special fiber, denoted by $Y_{n+1}$, of the second factor in  the product \eqref{eq:prod int mod} is smooth. In \cite[\S4.3]{LL2} the authors have defined a stratification of  the reduced special fiber, denoted by  $Y_n$, of $\CM_{K_{\wt \G(V_{n})}}$, essentially a refinement of the Newton stratification 
$$
Y_n=\coprod_{i=0}^{n-1} \coprod_{M\in \fkS_{i}}  Y_n^{(M)},
$$
where  $\fkS_{i}$ denotes the $ \fkS^{i}_m$ in {\it loc. cit.}. 
Here we simply take the stratification of $Y=Y_n\times Y_{n+1}$ as the product of the stratification of $Y_n$ with $Y_{n+1}$
$$
Y=\coprod_{i=0}^{n-1} \coprod_{M\in \fkS_{i}}  Y_n^{(M)}\times Y_{n+1}.
$$ By \cite[\S4.3]{LL2} this stratification of $Y_n$ verifies the two conditions stated before Proposition \ref{th:LL} (for $\sX=\CM_{K_{\wt \G(V_{n})}}$). It follows easily that the above  stratification of $Y$ satisfies  the two conditions stated before Proposition \ref{th:LL} (for $\sX=\CM_{K_{\wt \G }}$).

In the split-(Iwahori, Iwahori) case and the  inert-vertex-parahoric case, the scheme $\sX=\wt\CM_{K_{\wt \G }}$ has strictly semistable reduction. The special fiber $Y$ is already reduced  and we define the stratification induced by the union \eqref{eq:str Iw 2} and \eqref{eq:str in} respectively, as follows. Let $\CJ$ denote the set of indices in \eqref{eq:str Iw 2} and \eqref{eq:str in}, and denote $Y=\cup_{j\in\CJ }Y_j$. Then we define $\fkS^i$ to be the set of subsets $M$ of $\CJ$ with $\#M=\#\CJ-i$ such that $Y^{[M]}\coloneqq \cap_{j\in M} Y_j$ is non-empty (then it has  codimension $\#M$ in $\sX$). Set $Y^{[i]}=\cup_{M\in \fkS^i}Y^{[M]}$ and $Y^{(M)} =Y^{[M]}\setminus Y^{[\#M+1]}$. Then we have the resulting stratification 
 \begin{equation}\label{eq:str prod bl}
Y=\coprod_{i=0}^{\#\CJ} Y^{[i]}=  \coprod_{i=0}^{\#\CJ} \coprod_{M\in \fkS_{i}}  Y^{(M)}.
\end{equation}
 The strict semistability of $\sX$ implies that the stratification satisfies the two conditions stated before Proposition \ref{th:LL}. (Note that the scheme $Y^{[i]}$ is empty once $i>\dim Y$.)

We write $\BT=\BT_{n}\otimes\BT_{n+1}$ and $\fkm$ corresponding to $(\fkm_n,\fkm_{n+1})$  for maximal ideals $\fkm_\nu$ of $\BT_{\nu}, \nu\in\{n,n+1\}$. We will distinguish the three cases.

\smallskip

\paragraph{\em Split-(Drinfeld-level, hyperspecial) case}
By  Proposition  \ref{th:LL} (1), it suffices to verify that, for every $M\in \fkS$, we have $H^j(Y^{[M]}\otimes_k \ov k, L)_\fkm=0$ whenever $j\neq \dim Y^{[M]}$.  This follows from the K\"unneth formula, \cite[Prop. 4.25]{LL2} for $H^{j}(Y_n^{[M]}\otimes_k \ov k, L)_{\fkm_n}=0, j\neq\dim Y_n^{[M]}$, and the similar vanishing result for $H^{j}(Y_{n+1}\otimes_k \ov k, L)_{\fkm_{n+1}}=0, j\neq\dim Y_{n+1}$.

\smallskip

\paragraph{\em Split-(Iwahori, Iwahori) case}
We first define a stratification of the special fiber $Z$ of $\CM_{K_{\wt\G}}$ prior to the resolution, similar to \eqref{eq:str prod bl}:
 \begin{equation*}\label{eq:str prod bl Z}
Z=\coprod_{i=0}^{\#\CJ} Z^{[i]}=  \coprod_{i=0}^{\#\CJ} \coprod_{M\in \fkS_{i}}  Z^{(M)}.
\end{equation*}
Then, under the condition (1), it follows from \cite[(3), p.~859]{LL1} that $H^{i}_c(Z^{(M)})_\fkm=0$ for all $i$ and  $M\in \fkS$, unless $Z^{(M)}$ is maximal dimensional, in which case $H^{i}_c(Z^{(M)})_\fkm=0$ unless $i=\dim Z^{(M)}$. (In {\it loc. cit.} the authors only treated the case of Drinfeld levels; but the proof applies verbatim to the Iwahori case.)
Now we  return to the stratum $Y^{(M)}$ in \eqref{eq:str prod bl}. It is easy to see that the natural map $\pi_M\colon Y^{(M)}\to Z^{(M)}$ is smooth and the direct images $R\pi^j_{M,!} L$ are constant on $Z^{(M)}$. It follows that $H^{i}_c(Y^{(M)})_\fkm=0$ for all $i$ and $M\in \fkS$, unless $Y^{(M)}$ is maximal dimensional hence equal to $Z^{(M)}$, in which case $H^{i}_c(Y^{(M)})_\fkm=0$ unless $i=\dim Z^{(M)}$. It follows from the cohomological exact sequence associated to $Y^{[M]}=Y^{(M)}\cup (Y^{[M]}\setminus Y^{(M)})$ (see for example \eqref{eq:coh ES} below) and an induction that  $H^{i}(Y^{(M)})_\fkm=0$ for all $i$ and $M\in \fkS$, unless $Y^{(M)}$ is maximal dimensional, in which case 
 $H^{i}(Y^{[M]})_\fkm=0$ for $i>\dim Y^{[M]}$ and by Poincar\'e duality $H^{i}(Y^{[M]})_\fkm=0$ for $i\neq \dim Y^{[M]}$.  We have thus verified the condition in case (1) of Proposition  \ref{th:LL} and therefore we have proved $H^{2n}(\sX)_\fkm=0$ in this case.

\smallskip

\paragraph{\em Inert-vertex-parahoric case}  
We note that the moduli space $\CM_{K_{\wt G(V_\nu)}}$ for type $t_\nu$ (at $v$) is isomorphic to another similarly defined moduli space of type $\nu-t_{\nu}$. Therefore it suffices to  consider the cases when $t_n,t_{n+1}\in\{0,1\}$.
 We first recall from \cite[\S9, p.~868]{LL1}\footnote{Note that the unramifiedness of $v$ (over $\BQ$) in \cite{LL1} can be dropped by \cite[\S1.5]{LMZ}} that, when the type is $t_n=1$,  the cohomology of the balloon and the ground strata satisfy
 \begin{equation}\label{eq:van t=1}
H^i(Z\otimes_k \ov k, L)_\fkm= 0, \quad i\neq \dim  Z  
 \end{equation}for $Z=Y_n^\circ, Y_n^\bullet, Y_n^\dagger$ respectively,
 where we simplify the notation $Y_n^?= \CM_{K_{\wt \G(V_n)},k}^?$ in \eqref{eq:bu-circ} for $?\in \{\circ, \bullet \}$, and define $Y_n^\dagger=Y_n^\circ\cap Y_n^\bullet$. If one of $t_n,t_{n+1}$ is $0$, the proof is now similar to the split-(Drinfeld-level, hyperspecial) case,  
Proposition \ref{th:LL} (1). It remains to consider the type $(1,1)$ case. For $(?_n,?_{n+1})\in \{\circ,\bullet\}^2$, we will write $Y^{?_n,?_{n+1}}$ for the strict transform of $Y^{?_n}_n\times Y^{?_{n+1}}_{n+1}$. Then by  the formula for cohomology of blow-up, $H^i(Y^{?_n,?_{n+1}} \otimes_k \ov k, L )$ is isomorphic to
$$
H^i((Y^{?_n}_n\times Y^{?_{n+1}}_{n+1} )\otimes_k \ov k, L)\oplus\begin{cases} 0,&  (?_n,?_{n+1})= (\circ,\bullet)\text{ or }  (\bullet, \circ)\\
H^{i-2} ( (Y^{\dagger}_n\times Y^{\dagger}_{n+1}) \otimes_k \ov k, L), &  (?_n,?_{n+1})= (\circ,\circ)\text{ or }  (\bullet, \bullet).
 \end{cases}
$$
Similarly we can compute the cohomology of all of the closed strata  $Y^{[M]}$ in terms of the notation in \eqref{eq:str prod bl}  using  \eqref{eq:van t=1}  
$$
H^i(Y^{[M]}\otimes_k \ov k, L)_\fkm= 0, \quad i\neq \dim  Y^{[M]},
$$
for all $M\in\fkS$ with one exception: the stratum $Y^{[M_0]} \coloneqq Y^{\circ,\circ}\cap Y^{\bullet, \bullet}$, which is a $\BP^1$-bundle over $Y^{\dagger,\dagger}$. Nonetheless the exceptional case  has vanishing (localized at $\fkm$) cohomology at all degrees outside $i=\dim Y^{\dagger,\dagger} $ and $i=\dim Y^{\dagger,\dagger}+2$. Using \eqref{eq:van t=1} (for both $n$ and $n+1$) we can deduce that  $$
H^i_c(Y^{(M)}\otimes_k \ov k, L)_\fkm= 0, \quad i> \dim  Y^{[M]}
$$
for $M\neq M_0 \in\fkS$. To treat the exceptional case, we use the exact sequence
\begin{equation}\label{eq:coh ES}
\xymatrix{    H^{i-1}(Y^{[M_0]}\setminus Y^{(M_0)})  \ar[r]&H^{i}_c(Y^{(M_0)}) \ar[r]  & H^{i} ( Y^{[M_0]}).
}
\end{equation}
 Note that the stratum $Y^{[M_0]}$ has codimension $2$ in $\sX$, and $Y^{[M_0]}\setminus Y^{(M_0)}$ is smooth of codimension~$3$ in $\sX$. Since $  H^{i} ( Y^{[M_0]})_\fkm=0$ when $i\geq \dim Y^{[M_0]}+2$, and $ H^{i}(Y^{[M_0]}\setminus Y^{(M_0)})_\fkm=0$ when $i\neq \dim Y^{[M_0]}\setminus Y^{(M_0)}=\dim Y^{[M_0]}-1$, we conclude that 
 $H^{i}_c(Y^{(M_0)})_\fkm=0$ when $i\geq \dim Y^{[M_0]}+2$. By Poincar\'e duality we have  $H^{i}(Y^{(M_0)})_\fkm=0$ when $i\leq \dim Y^{[M_0]}-2$. Since $H^{2n}(X,L(n))_\fkm=0$, we have verified the condition in case (2) of Proposition  \ref{th:LL} and therefore we have proved $H^{2n}(\sX)_\fkm=0$.
\end{proof}

\begin{remark}
The condition (1) in Proposition \ref{thm:van} may be unnecessary if one makes a more careful study of the stratification of the special fiber of the small resolution.
\end{remark}

\section{$p$-adic Abel--Jacobi maps and $p$-adic heights}\lb{sec ht}
We summarize the definitions and results we need from the theory of $p$-adic heights. For more details or more general setups, see \nek's original paper \cite{N93} and \cite[Appendix A]{DL}; our constructions follow the sign conventions of the latter reference. Nothing in this section is new.

In  \S~\ref{sec:102} we give the definitions and some basic properties  of local and global height pairings on Galois cohomology in terms of biextensions of Galois representations. Thanks to a construction recalled in  \S~\ref{sec:101}, we obtain height pairings of algebraic cycles, which are related to arithmetic intersections by two results reviewed in \S~\ref{sec:103}. 

The notation of this section is independent of that of the rest of the paper.
We denote by $L$ a finite extension of $\Q_{p}$, and by $\Gamma$ a finite-dimensional $L$-vector space.  

\subsection{$p$-adic Abel--Jacobi maps and biextensions} \lb{sec AJ} \lb{sec:101}
Let $F$ be a field of characteristic different from~$p$, and let $X$ be a smooth projective scheme over $F$ of pure dimension $m-1\geq1$. We denote by $\RZ^{\bullet}(X)_{R}$ the module of $\bullet$-codimensional algebraic cycles with coefficients in a ring $R$ (omitted from the notation when $R=\Z$), and by $\Ch^{\bullet}(X)_{R}=\RZ^{\bullet}(X)_{R}/(\text{rational equivalence)}$ the Chow $R$-module.
We denote $H^{i}(F, -)=H^{i}_{\rm cont}(G_{F}, -)$ where $G_{F}$ is the absolute Galois group of $F$.

\subsubsection{$p$-adic Abel--Jacobi maps}\lb{sec pAJ}
Let $0\leq d\leq m$ and consider the absolute \'etale cohomology $H^{2d}(X, L(d))$.  By the Hochschild--Serre spectral sequence, it has a filtration ${\rm Fil}^{\bullet}$ with 
$$0\to  H^{1}({F},H^{2d-1}(X_{\ol{F}}, L(d))) \to H^{2d}(X, L(d))/{\rm Fil}^{2} \to H^{0}({F}, H^{2d}(X_{\ol{F}}, L(d)))\to 0.$$
We denote by $\ol{\rm cl}\colon {\RZ}^{d}(X)_{L} \to  H^{0}(G_{F}, H^{2d}(X_{\ol{F}}, L(d)))$ the geometric cycle class map, by
 ${\RZ}^{d}(X)_{L}^{0}$ its kernel, and we let 
\beqq
\wt{\rm cl} \colon {\RZ}_{d}(X)_{L} &\to  H^{2d}(X, L(d))/{\rm Fil^{2}} ,\\
{\rm cl} \colon {\RZ}_{d}(X)_{L}^{0}&\to H^{1}(F, H^{2d-1}(X_{\ol{F}}, L(d))) 
\eeqq
be the absolute cycle class map and the Abel--Jacobi map, respectively. 
The maps $\ol{\rm cl}$, $\wt{\rm cl}$ factor through the Chow group $\Ch^{d}(X)$, and the map ${\rm cl}$ factors through the image $\Ch^{d}(X)^0\subset {\Ch}^{d}(X)$  of  ${\RZ}^{d}(X)^0$.

If $M\subset H^{2d-1}(X_{\ol{F}}, L(d))$ is a $G_{F}$-stable subspace, we denote 
 by
 $$\RZ^{d}_{M}(X)^{0}_{L} , \quad \Ch^{d}_{M}(X)^{0}_{L}$$
 the preimages in $\RZ^{d}(X)^{0}_{L}$, $\Ch^{d}(X)^{0}_{L}$ of $H^{1}(F,M) \subset H^{1}(F, H^{2d-1}(X_{\ol{F}}, L(d)))$ under the Abel--Jacobi map.
 
 Suppose that  $F$ is a non-archimedean local field of residue characteristic $\ell$. We will consider subspaces $M$ satisfying the condition:
 \begin{enumerate}
 \item\lb{cond ell} if $\ell\neq p$:   $H^{1}(F, M)=0$;
 \item\lb{cond p} if $\ell=p$:  $H^{1}_{\rm st}(F, M) =  H^{1}_{f}(F, M)$.  \end{enumerate}
 \begin{remark} Since by  \cite[Theorem B]{NN} the map ${\rm cl}$  takes values in the subspace  
 $$
 H^{1}_{\rm st}(F,  H^{2d-1}(X_{\ol{F}}, L(d))),
 $$ the  conditions above imply that
 $${\rm cl}(\RZ^{d}_{M}(X)^{0}_{L} )\subset  H^{1}_{f}(F, M).$$
If $M$ is pure of weight $-1$ (as is implied for all $M\subset H^{2d-1}(X_{\ol{F}}, L(d))$ by the weight-monodromy conjecture), then the relevant one among the conditions above is satisfied. 
 \end{remark}

\subsubsection{Biextensions from algebraic cycles}\lb{biext}
Let $d_{1}, d_{2}\geq 0$ be integers with $d_{1}+d_{2}=m$, and let
 $$Z_{1}\in \RZ^{d_{1}}(X)^{0}_{L},\quad   Z_{2}\in \RZ^{d_{2}}(X)^{0}_{L}$$   
 be cycles with disjoint supports. Let $M_{i}\coloneqq  H^{2d_{i}-1}(X_{\ol{F}},L(d_{i}))$. To each $Z_{i}$ is associated an extension of $L[G_{F}]$-modules
$$0  \sto M_{i}\sto E_{i}\sto L \sto 0$$
whose class in $H^{1}(F,M_{i})$ is the $p$-adic Abel--Jacobi image ${\rm cl}(Z_{i})$. A further geometric construction yields the \emph{biextension} $E_{{1}}^{{2}}= E_{Z_{1}}^{Z_{2}}$ fitting in the following exact diagram  
\begin{equation}\lb{eq: biext}
\xymatrix{
   &   &0\ar[d]  &0\ar[d]  &  \\
0\ar[r] & L(1) \ar[r]\ar@{=}[d]&   E^{2}
    \ar[d]\ar[r] & M_{1}
  \ar[r]\ar[d]& 0\\
0\ar[r] & L(1) \ar[r] & E_{{1}}^{{2}}   \ar[d]\ar[r] &   E_{{1}}   \ar[r]\ar[d] &0\\
&&    L  \ar@{=}[r]\ar[d]&  L \ar[d]&\\
&&0&0&
}
\end{equation}
where $M_{1} = M_{2}^{*}(1)$ via Poincar\'e duality, and $E^{2}\coloneqq  E_{{2}}^{*}(1)$. We denote its class by $[E_{{1}}^{{2}}]\in H^{1}(F, E^{{2}})$.

\subsection{Height pairings}  \lb{sec:102}
We collect some definitions and properties of  local and global height pairings.

\subsubsection{Local height pairings of algebraic cycles} \lb{pair def}
Suppose  that $F$ is a non-archimedean local field  of characteristic zero and residue characteristic $\ell$.
Let $\lm\colon F^{\ts}\hat{\ot} L\to \Gamma$ be an $L$-linear map.

For $i=1, 2$ let  $M_{i}\subset H^{2d_{i}-1}(X_{\ol{F}},L(d_{i}))$  be $L[G_{F}]$-submodules,  and denote still by $\la\, , \, \ra\colon M_{1}\ot_{L} M_{2}\to L(1)$ the restriction of  the Poincar\'e pairing 
 $$
 \la\, , \,\ra\colon   H^{2d_{1}-1}(X_{\ol{F}},L(d_{1}))\ot_{L}  H^{2d_{2}-1}(X_{\ol{F}},L(d_{2}))\stackrel{\cup}{\to}  H^{2m-2}(X_{\ol{F}},L(m))\stackrel{\Tr}{\to} L(1),
 $$
where the map $\Tr$ is the sum of the trace maps for the connected components of $X$.
Assume that $M_{1}$, $M_{2}$ satisfy the following conditions:
\begin{enumerate}
\item $\la\, , \,\ra\colon M_{1}\ot_{L} M_{2}\to L(1)$ is a perfect pairing;
\item \lb{ht cond l} if $\ell\neq p$, we have  $H^{0}(F, M_{i})=0$ for $i=1, 2$; this implies condition \eqref{cond ell} for $M_{1}$, $M_{2}$ in \S~\ref{sec pAJ}, and is implied by the condition that $M_{i}$ is pure of weight $-1$;
\item \lb{ht cond p} if $\ell=p$:
 \begin{itemize}
 \item  $M_{i}$ is crystalline with $\RD_{\rm crys}(M_{i})^{\vphi=1}=0$ for $i=1, 2$; this implies condition \eqref{cond p} for $M_{1}$, $M_{2}$ in \S~\ref{sec pAJ}, and  is implied by the condition that $M_{i}$ is  crystalline and pure of weight $-1$;
 \item the \emph{Panchishkin condition}: there is a (necessarily unique) extension of crystalline representations
$$0\sto M_{i}^{+}\sto M_{i}\to M_{i}^{-}\sto 0$$
such that ${\rm Fil}^{0}\BD_{\rm dR}(M_{i}^{+})= \BD_{\rm dR}(M_{i}^{-})/ {\rm Fil}^{0}\BD_{\rm dR}(M_{i}^{-})=0$; this implies that the natural map 
\beq
\lb{split hdg}\BD_{\rm dR}(M_{i}^{+}) \oplus  {\rm Fil}^{0}\BD_{\rm dR}(M_{i})\stackrel{\cong}{\to} \BD_{\rm dR}(M_{i})\eeq
is a splitting of the Hodge filtration on $\BD_{\rm dR}(M_{i})$.
 \end{itemize}
\end{enumerate}

Assume that $Z_{1}\in \RZ^{d_{1}}_{M_{1}}(X)^{0}$, $Z_{2}\in \RZ^{d_{2}}_{M_{2}}(X)^{0}$. Then  the biextension class $[E_{{1}}^{{2}}]=[E_{Z_{1}}^{Z_{2}}]$ belongs to the preimage $H^{1}_{\text{$M_{1}$-$f$}}(F, E^{{2}})\subset H^{1}(F, E^{{2}})$ of $H^{1}_{f}(F, M_{1})$ under the natural map $H^{1}_{f}(F, E^{2})\to H^{1}(F, M_{1})$.
This group sits in the (pushout) diagram of exact sequences\footnote{This diagram should also replace an incorrect one in \cite[(4.1.4)]{D17}.}
\begin{equation}\label{ex seq me}
\xymatrix{
0\ar[r] &H^{1}(F, L(1)) \ar[r]^{\stackrel{\sg}{\longleftarrow}}  & H^{1}_{\text{$M_{1}$-$f$}}(F, E^{2})\ar[r] &H^{1}_{f}(F, M_{1})\ar[r] &0\\
0\ar[r] &H^{1}_{f}(F, L(1))\ar[r]^{\stackrel{\sg_{f}}{\longleftarrow}}\ar[u]  &H^{1}_{f}(F, E^{2})\ar[r]\ar[u] &H^{1}_{f}(F, M_{1})\ar[r]\ar@{=}[u] &0
}
\end{equation}
admitting canonical splittings $\sg$, $\sg_{f}$. These are obvious if $\ell\neq p$, as then  $H^{1}(F, M_{1})=0$;  for $\ell=p$, they are induced by \eqref{split hdg} (see \cite[\S~4]{N93}).
Moreover, the Kummer map identifies $H^{1}(F, L(1))\cong F^{\ts}\hat{\ot} L$.
\begin{definition}\lb{ht zz} Let $M_{1}$, $M_{2}$, $Z_{1}$, $Z_{2}$ be as above. 
 We define
\beq\lb{hw def}
h_{X, \lm}(Z_{1},Z_{2})\coloneqq  
\lm\circ \sg([E_{{1}}^{{2}}]) \in \Gamma.\eeq
\end{definition}
\begin{remark} Since the conditions on the pair $(M_{1}, M_{2})$ are stable under subobjects and extensions (see \cite[Lemma A.14]{DL} for extensions when $\ell =p$), there is a maximal  pair satisfying those; this justifies omitting 
 $(M_{1}, M_{2})$ from the notation.
\end{remark}

\begin{remark} \lb{rmk crys} If $\ell=p$, it follows from the previous discussion that $\sg([E_{1}^{2}])\in \sO_{F}^{\ts}\hat{\ot}L\subset F^{\ts}\hat{\ot}L$ if and only if $[E_{1}^{2}]$ is crystalline (that is,  belongs to $H^{1}_{f}(F, E^{2})$).
\end{remark}

\begin{lemma}[Base change] \lb{bc ht}  Consider the setup of Definition \ref{ht zz}.
\begin{enumerate}
\item Let $F'/F$ be a finite extension, and let $\lm'\coloneqq \lm\circ {\rm N}_{F'/F}$. Then  for any $Z_{1}\in\RZ^{d_{1}}_{M_{1}}(X)^{0}$, $Z_{2}'\in \RZ^{d_{2}}_{M_{2}}(X_{F'})^{0}$, 
$$h_{X_{F}, \lm} (Z_{1}, {\rm N}_{F'/F}Z_{2}') = h_{X_{F'}, \lm'} (Z_{1, F'}, Z'_{2}).$$
\item Let $u\colon X'\to X$ be a finite \'etale morphism, and let $Z_{1}\in\RZ^{d_{1}}_{M_{1}}(X)^{0}$, $Z_{2}\in \RZ^{d_{2}}_{M_{2}}(X)^{0}$. Denote by $Z_{i}'$ the pullback of $Z_{i}$ to $X'$. Assume $\ell \neq p$. Then 
$$h_{X, \lm}(Z_{1}, Z_{2})={1\over \deg u} h_{X', \lm}(Z_{1}', Z_{2}').$$
\end{enumerate}
\end{lemma}
\begin{proof} Part (1) is \cite[(II.1.9.1)]{N95}. Part (2) follows from \cite[Lemma B.3]{LL1} and \cite[Proposition A.7]{DL}. 
\end{proof}

\subsubsection{Global height pairings for Selmer groups}\lb{gl ht}
Let now $F$ be  a number field and $\lm\colon \Gamma_{F, L}\to \Gamma$ be an $L$-linear map. 

Let $M_{1}$, $M_{2}$ be $L$-vector spaces endowed with continuous $G_{F}$-representations that are unramified at all but finitely many places of $F$, and de Rham at all the $p$-adic places. Assume moreover that $M_{1}$, $M_{2}$ are endowed with a perfect $G_{F}$-equivariant pairing $\la\, , \, \ra\colon M_{1}\ot_{L}M_{2}\to L(1)$, and that for each $i$ and each finite place $w$ of $F$, the representation $M_{i}$ restricted to $G_{F_{w}} $ satisfies the conditions \eqref{ht cond l}, \eqref{ht cond p} of \S~\ref{pair def}. 

Under these conditions,\footnote{These are not the most general possible; for instance, the crystalline condition at $p$-adic places is not necessary.} \nek\ \cite{N93} defined a bilinear height pairing on the Bloch--Kato Selmer groups
\beq\lb{ht sel}
h_{M_{1}, \lm}\colon H^{1}_{f}(F, M_{1})\ot_{L} H^{1}_{f}(F, M_{2})\to \Gamma\eeq
as follows. For $i=1,2$ pick representatives $E_{i}$ of the extension classes $[E_{i}]\in  H^{1}_{f}(F, M_{i})$, and let $E_{1}^{2}$ be a biextension fitting in a diagram  \eqref{eq: biext} of $G_{F}$-representations. For each finite place $w$ of $F$, one can then define $h_{M_1, w}^{E_{1}^{2}}([E_{1}], [E_{2})$ by the right-hand side of \eqref{hw def} (where everything is viewed as a representation of $G_{w}$); the sum 
$$h_{M_1, \lm}([E_{1}], [E_{2}])\coloneqq  
\sum_{w} h_{M_1, w, \lm}^{E_{1}^{2}}([E_{1}], [E_{2}])$$ 
 does not depend on the choice of $E_{1}^{2}$.

\begin{lemma}[Projection formula] \lb{proj formula}
Let $(M_{1}, M_{2})$ and $(M_{1}', M_{2}')$ be as above. Let  $\phi\colon M_{1}'\to M_{1}$ be a   map of $G_{F}$-representations, and let $\phi^{*}(1)\colon M_{2}\to M_{2}'$ be the dual map. Let $[E_{1}']\in H^{1}_{f}(F, M_{1}')$, $[E_{2}]\in H^{1}_{f}(F, M_{2})$. Denote by $E_{2}'\coloneqq  \phi^{*}(1)_{*}E_{2}$, $E_{1}\coloneqq  \phi_{*}E_{1}'$ the pushouts.  Then 
$$h_{M_{1}', \lm}([E_{1}'], [E_{2}']) =h_{M_{1}, \lm}([E_{1}], [E_{2}]). $$
\end{lemma}
\begin{proof} Let $E' \in H^{1}_{f} (F, E^{2}{}')$ be a biextension (as in \eqref{eq: biext}) of $E_{1}'$ and $E^{2}{}'\coloneqq  \phi^{*}(E_{2}^{*}(1))=E_{2}'^{*}(1)$. The map $\phi\colon M_{2}'^{*}(1)\cong M_{1}'\to M_{1}\cong M_{2}^{*}(1)$ induces by pullback a map  $\phi\colon  E^{2}{}' \to E^{2}$. Then a diagram chase shows that $\phi_{*}E' \in H^{1}_{f}(F, E^{2})$ is a biextension of $E_{1}$ and $E^{2}$. 
\end{proof}
\subsubsection{Decomposition in the case of algebraic cycles} \lb{sec: ht gl cyc}
Let $X$ be a proper smooth scheme over $F$ of dimension $m-1$, and suppose that $M_{i}\subset H^{2d_{i}-1}(X_{\ol{F}}, L(d_{i}))$ are $L[G_{F}]$-submodules satisfying the above conditions with respect to a pairing $\la\, , \,\ra$ that is the restriction of the Poincar\'e pairing. We then denote $h_{X,\lm}\coloneqq h_{M_{1},  \lm}$, for which we have
\beq\lb{dec ht}
h_{X, \lm}({\rm cl}(Z_{1}), {\rm cl}(Z_{2}))=\sum_{w\nmid \infty} h_{X_{w}, \lm_{w}}
(Z_{1}, Z_{2}),\eeq
where the sum runs over all the non-archimedean places of $F$, and $X_{w}\coloneqq X_{F_{w}}$, $\lm_{w}\coloneqq  \lm_{|F_{w}^{\ts}\hat{\ot}L}.$

\subsection{Relation to arithmetic intersection theory} \lb{sec:103}
We collect  two results relating local heights away from $p$, and the crystalline property of biextensions at $p$, with arithmetic intersections. 
 We start  with some preliminaries. For more details on the background, see \cite[Appendix B]{LL1} and references therein.

\subsubsection{Extensions of algebraic cycles}
Let $\sX$ be a regular scheme; for a closed subset $\sY$ (omitted from the notation if $\sY=\sX$) we denote by  $K_{0}^{\sY}(\sX)$ the $K$-group of complexes of coherent sheaves on $\sX$ with cohomology supported in $\sY$. We denote by ${\RF}^{\bullet}$ the filtration on $K_{0}^{\sY}(\sX)$ by the codimension of support. We have an $L$-linear map
\beq\lb{kappa}
\kappa\colon \RZ^{d}(\sX)_{L}&\to \RF^{dt}K_{0}(\sX)_{L}\eeq
such that if $\sZ\subset \sX$ is an integral subscheme, then $\kappa([\sZ])=[\sO_{\sZ}]$. 

Let now $F$ be a nonarchimedean local field and denote by $k$ its residue field. Assume that the regular scheme $\sX$ is endowed with a   projective and flat map $\pi\colon \sX\to \Spec \sO_{F}$, and denote by $X$ and $\sX_{k}$ respectively the generic and special fibre of $\sX$.

\begin{definition}\label{extension}
Let  $Z\in\RZ^d(X)_L$, and denote by $|Z|\subset X$ its support. We say that an element 
$$\sZ\in\RF^d K_0^{\sX_{k}\cup |Z|}(\sX)_L \subset \RF^d K_0(\sX)_L $$ 
is an \emph{extension} of $Z$ if $\sZ_{|_X}\in\RF^d K_0(X)_L$ coincides with $\kappa(Z)$.
\end{definition}

\subsubsection{Intersection pairing}
Suppose that $X$ has dimension $m-1\geq 1$. For a pair of integers $d_{1}, d_{2}\geq 0$ with $d_{1}+d_{2}=m$,  and cycles $\sZ_{i}\in \RF^{d_{i}}K_{0}(\sX)$ with $|\sZ_{1}|\cap |\sZ_{2}|\subset  |\sX_{k}|$, we define their \emph{intersection} by
$$(\sZ_{1}\cdot \sZ_{2})\coloneqq  \chi(\pi_{*}(\sZ_{1}\cup \sZ_{2})),$$
where
 $$\cup \colon\RF^{d_{1}} K_{0}^{|\sZ_{1}|}(\sX)\ot \RF^{d_{2}}K_{0}^{|\sZ_{2}|}(\sX)\to \RF^{m}K_{0}^{\sX_{k}}(\sX)$$
 is the cup product, and $\chi\colon K_{0}(\Spec k)\to \Z$ is the Euler characteristic. The definition is extended linearly to cycles with coefficients in $L$.

\subsubsection{Arithmetic intersections and the crystalline property at $p$}\lb{sec crys} 
Consider the setup of \S~\ref{pair def} with $\ell=p$. 
\begin{proposition}   \lb{crys prop}
Assume that $p>m$ or $m=2$,  and that $X$ admits a proper smooth model  $\sX/\sO_{F}$. If the supports of the  Zariski closures $\sZ_{1}$, $\sZ_{2}$ of $Z_{1}$, $Z_{2}$ in $\sX$ are disjoint, then the biextension $[E_{{1}}^{{2}}]$ is crystalline. 
\end{proposition}
\begin{proof} If $p>m$, this is a special case of \cite[Theorem A.8]{DL}. If $m=2$, this is a special case of \cite[Proposition 4.3.1]{D17}.
\end{proof}

\subsubsection{Arithmetic intersections and local heights away from $p$} 

Consider  the setup of \S~\ref{pair def} with $\ell\neq p$.  
\begin{proposition} \lb{B13}
 Assume that $m=2n$ and $d_{1}=d_{2}=n$.
 Let $\hT_{1}, \hT_{2}\in \EC(\sX)_{L}$, and assume that $Z_{1}.\hT_{1}$ and $Z_{2}.\hT_{2}$ have disjoint supports. Let $\sX$ be a regular flat projective scheme over $\sO_{F}$ with generic fibre $X$, and let  $\sZ_{i}\in \RF^{d_{i}}K_{0}(\sX)$ be an extension of $Z_{i}$ for $i=1,2$.
 
  Suppose that  one of the following conditions holds:
\begin{enumerate}
\item $\sX$ is smooth over $\sO_{F}$,  $\sZ_{i}$ is  (the image under $\kappa=\eqref{kappa}$ of)  the Zariski closure of $Z_{i}$, and $\hT_{i}={\rm id}$;
\item $\hT_{1}$, $\hT_{2}$ annihilate $H^{2n}(\sX, L(n))$.
\end{enumerate}
Then 
$$h_{X, \lm}(\sZ_{1}.\hT_{1}, \sZ_{2}.\hT_{2})= -((\sZ_{1}.\hT_{1})\cdot (\sZ_{2}.\hT_{2}) )\, \lm(\vpi)$$
where $\vpi\in F^{\ts}$ is a uniformizer. 
 \end{proposition} 
\begin{proof} In case (1), this is a special case of  \cite[Proposition B.10]{LL1} combined with \cite[Proposition A.7, Remark A.6]{DL}. In case (2), this is  \cite[Proposition B.13]{LL1} combined with \cite[Proposition A.7, Remark A.6]{DL}.
\end{proof}

In favorable cases, correspondences satisfying condition (2) of the  proposition can be found using Proposition 
\ref{thm:van} as in \cite{LL1,LL2}.

\section{The arithmetic relative-trace formula}
\lb{sec:10}
Let $V=(V_n,V_{n+1})\in \sV^{\circ, -}$ be an incoherent pair, and let $\G=\G^{V}$, $\RH= \RH^{V}$. In this section, we define our cycles of interest, 
and a distribution $\sJ=\sJ_{K_{p}}$ on (part of) the Hecke algebra for $\G(\A^{p})$ that encodes their  $p$-adic heights. The main result of this section is the arithmetic RTF for $\sJ$ (Theorem \ref{ARTF thm}).

We  fix an element $\vphi_0\in \Hom(F,\C)$ with $v_0$ the induced place of $F_0$. Let $V_\nu(v_0)$ be the $v_0$-nearby hermitian space of $V_\nu$, for $\nu=n,n+1$. Let $\G^{(v_0)}=\U(V_n(v_0))\times \U(V_{n+1}(v_0))$ and $\RH^{(v_0)}=\U(V_{n}(v_0))$. Then by \S~\ref{ss:AD} we obtain Shimura varieties  (with the fixed choice $\vphi_0\in \Hom(F,\C)$)
$$
X_K \coloneqq \Sh_{K}(\G^{(v_0)}),\quad Y_{K_{\RH}}\coloneqq \Sh_{K_{\RH}}(\RH^{(v_0)})
$$
with the reflex field $\varphi_0:F\to \C$. We fix an isometry $V_\nu(v_0)\otimes\A^{\infty}\simeq V_{\nu,\A^{\infty}}$, for $\nu=n,n+1$. We have an induced isomorphism $\G^{(v_0)}(\A^{\infty})\simeq \G(\A^{\infty})$.

From now on, we assume that $F_{0}\neq \Q$, which implies that these Shimura varieties are proper.

 In \S~\ref{sec: coh rho} we study the \'etale cohomology of $X_{K}$ and define the Galois representation of interest. In \S~\ref{sec: ggp cycle}, we define and study the arithmetic diagonal cycles and Gan--Gross--Prasad cycles. In \S~\ref{sec: J def} we define $\sJ$ by means of height pairings of those cycles, and give its spectral expansion. In \S~\ref{sec: J dec} we prove some vanishing results to decompose $\sJ$ as a sum indexed by the nonsplit finite places of $F_{0}$. Finally, in \S~\ref{sec: J geom} we state the geometric expansion of $\sJ$.

\subsection{Cohomology and automorphic Galois representations}\lb{sec: coh rho}
 Let  $L$  be an algebraic extension of $\bQ_{p}$.

\subsubsection{Ordinary representations of $\G(\A)$} We say that $\pi\in\wt\sC(\G)(L)$ is \emph{ordinary} if for every place $v\vert p$ of $F_{0}$, the base-change ${\rm BC}(\pi_{v})$ satisfies the ordinariness conditions of \S~\ref{p interp}.
This condition defines an ind-subscheme
$$ \sC(\G)^{\ord}\subset\sC(\G)_{\Q_{p}}.$$
We also define $\wt{\sC}(\G)^{\ord}(L)$ as the corresponding set of isomorphism classes of representations such that $\sC(\G)^{\ord}(L)=\wt{\sC}(\G)^{\ord}(\ol L)/G_{L}$ (cf. \S~\ref{sec: G BC}).

\subsubsection{Duals and Hecke actions}  If $S$ is a finite set of places of $F_{0}$ and $M$ is an admissible (left) $L[\G(\A^{S})]$-module, we denote  
$$M^{*}\coloneqq \varprojlim_{K\subset\G(\A^{S\infty})} M^{K, \vee}$$ the algebraic dual of $M$ (where the transition maps are restrictions of functionals), whereas as usual we denote by $M^{\vee}=\varinjlim_{K^{S}} M^{K^{S}, \vee}$ the contragredient (where the transition maps $M^{K, \vee}\to M^{K', \vee}$ are given by $x\mapsto x\circ e_{K \vert M^{K'}}$ if $e_{K}\colon M \to M^{K}$ is the natural $K$-projection). For any compact subgroup $K\subset \G(\A^{S\infty})$, we denote by  $M^{*}_{K} $ the $K$-coinvariants (thus the natural map $M^{\vee, K}\to M^{*}_{K}$ is an isomorphism if $K$ is open).  We have a   map
\beq\lb{inv to dir}
 M^{\vee}  &{\to} \ M^{*}\\
x&\mapsto \lim_{K} x\circ e_{K}
\eeq
The left Hecke action on $M$ induces a right action
$$T\colon  \sH(\G(\A^{S}), L)\to \Hom(M^{*},  M^{\vee}).$$

\subsubsection{Hecke and Galois actions on the cohomology of unitary Shimura varieties} \lb{sec: hecke gal coh} 
For   $i\in \Z$, we put
\beq
M^{i, K}&\coloneqq  H^{i}(X_{K, \ov F},\Q_{p}(n)),\qquad &
  & M^{i}\coloneqq \varinjlim_{K} M^{i}_{K},
\eeq
where the limit is with respect to the 
  pullback maps. For $?=\emptyset, K$, we also put $M^{\oplus,?}\coloneqq \bigoplus M^{i,?}$; it has a natural (left) action by $\sH(\G(\A^{\infty}),\Q_{p})_{?}$
    and by the Galois group $G_{F}$.

Let $\diamond\in \Z\cup \{\oplus\}$. The $\sH(\G(\A^{\infty}), \Q_{p})$-action on $M^{\diamond}$ makes it into an admissible $\G(\A^{\infty})$-module, so that we may consider $M^{\diamond , *}$.    It is  helpful to  think  of $M^{\diamond,*}$ as the inverse limit of homology (and of  $M^{\diamond}$, $M^{\diamond, \vee} $ as the direct limits of cohomology, respectively homology).

For $\pi\in\wt\sC(\G)(L)$, let
\beqq
\rho[\pi]^{\diamond}&\coloneqq  \Hom_{\sH(\G(\A^{\infty}))}(\pi^{\vee}, M^{\diamond, \vee}_{L}(1)),\\
M^{\diamond, \pi}&\coloneqq  \pi^{*}\boxtimes\rho[\pi] \subset M^{\diamond, *}_{L}(1),
\eeqq
so that we have a Hecke-equivariant map 
\beq\lb{Mpi dec}
\pi \to \Hom_{G_{F}}(M^{\diamond, *}(1), \rho[\pi])
\eeq
factoring through $\Hom_{G_{F}}(M^{\diamond, \pi}, \rho[\pi])$.
In fact, it is known (see \cite[Theorem III.5.1]{BW80}) that  temperedness implies
\beq\lb{temp mid deg}
M^{\oplus, \pi} =M^{2n-1, \pi}\eeq
so that we will simply write $M^{\pi}\coloneqq  M^{2n-1, \pi}$, $\rho[\pi]\coloneqq  \rho[\pi]^{2n-1}$.

We put $M_{\pi^{\vee}}^{K}\coloneqq  (M^{\pi}_{K})^{\vee}$ and $M_{\pi^{\vee}}\coloneqq \varinjlim_{K} M_{\pi^{\vee}}^K$, so that $M^{\pi}=M_{\pi^{\vee}}^{*}(1)$.  
  For $?\in \{{\rm temp}, \tord\}$,\footnote{The abbreviation `$\tord$' stands for tempered ordinary.}  we put 
  $$
  M_{?,\ol{\Q}_{p}}\coloneqq  \bigoplus M_{\pi}\subset M^{\oplus}_{\ol{\Q}_{p}}, 
\qquad
  M_{\ol{\Q}_{p}}^{?}\coloneqq  \bigoplus M^{\pi}\subset M^{\oplus, *}_{\ol{\Q}_{p}}(1)$$
 where the sums run over $\sC(\G)(\ol{\Q}_{p})$ and $\sC(\G)^{\ord}(\ol{\Q}_{p})$ respectively.
These are base-changes of $L$-subspaces $M_{?}\subset M^{2n-1}$, $M_{L}^{?}\subset M^{2n-1, *}(1)$. 
 Poincar\'e duality gives an isomorphism $M_{K}\cong M_{K}^{*}(1)$, which induces isomorphisms 
 \beqq\lb{iso PD}
 M_{?}^{K}\cong M^{?}_{K}\eeqq
  for $?\in \{{\rm temp}, \tord\}\cup \wt{\sC}(\G)(L)$.

\subsubsection{Automorphic Galois representations and decomposition of the cohomology}
Assume from now on that the extension  $L$ of $\Q_{p}$ is finite,  and denote by $\ol{\Q}_{p}$ an algebraic closure of $L$.
Let $\pi=\pi_{n}\boxtimes \pi_{n+1}\in \wt\sC(\G)(L)$.
\begin{lemma}
For  $\nu\in \{n, n+1\}$ there is a semisimple continuous representation
$$\rho_{\pi_{\nu}, \ol{\Q}_{p}}\colon G_{F}\sto \GL_{\nu}(\ol{\Q}_{p})$$ 
characterized, up to isomorphism, by the property that for all but finitely many places $w$ of $F$ split over $F_{0}$, the restriction $\rho_{\pi_{\nu}, \ol{\Q}_{p}|G_{F_{w}}}$ is unramified, 
and a geometric Frobenius at $w$ acts with a characteristic polynomial equal to the Satake polynomial of $\pi_{w}$ viewed as a representation of $\GL_{\nu}(F_{w})$. 
If $\pi_{\nu}$ is stable, then 
\beq\lb{rhopiPi}
\rho_{\pi_{\nu}, \ol{\Q}_{p}}\cong \rho_{{\rm BC}(\pi_{\nu}), \ol{\Q}_{p}}
\eeq
 (where the latter is as in \S~\ref{intro2}).
\end{lemma}
\begin{proof}
The construction is as in  \cite[Lemma 4.10]{DL}, using  \cite[Proposition 3.2.8]{LTXZZ} (due to Shin) instead of \cite{mok}.
Property \eqref{rhopiPi} is immediate from the construction.
\end{proof}
Let
$$\rho_{\pi, \ol{\Q}_{p}}\colon G_{F}\sto \GL_{n(n+1)}(\ol{\Q}_{p})$$
be defined by
$\rho_{\pi, \ol{\Q}_{p}}(-n)\coloneqq  \rho_{\pi_{n}, \ol{\Q}_{p}}\ot 
\rho_{\pi_{n+1}, \ol{\Q}_{p}}.$   If $\rho\colon G_{F}\to \GL_{N}(L)$ is a continuous representation, denote  by $\rho_{\ol{\Q}_{p}}\coloneqq  \rho\ot_{L}\ol{\Q}_{p}$ the base-change and by $\rho_{\ol{\Q}_{p}}^{\rm ss}$ its semisimplification.

The following key hypothesis gives an explicit description of $\rho[\pi]$ (at least in the stable case).

 \begin{hypothesis}\lb{hyp coh}
Let  $\pi\in \wt\sC(\G)(L)$, and let $K\subset \G(\A^{\infty})$ be an open compact subgroup. Then $\rho[\pi]_{\ol{\Q}_{p}}^{\rm ss}$  
 is isomorphic to a direct summand of $\rho_{\pi , \ol{\Q}_{p}}$. Moreover, if $\pi$ is stable then
$\rho[\pi]_{\ol{\Q}_{p}}^{\rm ss}\cong \rho_{\pi, \ol{\Q}_{p}}$.
\end{hypothesis}

\begin{remark}\lb{construct rho}
 Let $\Pi\in \sC(\G')^{\rm her, -}_{\Q_{p}}$, and let $L=\Q_{p}(\Pi)$. 
Let $\ol{\pi}\in \sC(\G^{V})(\ol{\Q}_{p})$ be the preimage of $\Pi$ under \eqref{BC rat st}; a priori we know it is isomorphic to its $G_{L}$-conjugates but not that it arises from some $\pi\in \wt{\sC}(\G^{V})(L)$. 
 Assume that Hypothesis \ref{hyp coh} holds.   By the definitions, the  space
$ M^{\ol{\pi}} $ is isomorphic  to $\ol{\pi}^{*}\boxtimes \rho[\ol{\pi}]$ as a Hecke- and $\ol{\Q}_{p}[G_{F}]$-module, and it  is a $G_{L}$-invariant subspace of $M_{\ol{\Q}_{p}}^{2n-1,*}(1)$. Let $M^{\pi}\coloneqq  ( M^{\ol{\pi}}) ^{G_{L}} \subset M_{L}^{2n-1, *}(1)$,  and define the $L[G_{F}]$-module
$$\rho_{\Pi}\coloneqq  (M^{\pi})^{\RH(\A^{\infty})}$$
Then we have 
$$\rho_{\Pi}\ot_{L}\ol{\Q}_{p} \cong (\ol{\pi}^{*})^{\RH(\A^{\infty})} \ot_{\ol{\Q}_{p}} \rho[\ol\pi].$$
The first tensor factor is 1-dimensional, so that by Proposition \ref{cor C her} and \eqref{rhopiPi}, the representation 
 $$\rho_{\Pi}\colon G_{F}\to \GL_{n(n+1)}(L)$$
satisfies 
$$(\rho_{\Pi}\ot_{L}\ol{\Q}_{p})^{\rm ss}\cong 
 \rho_{\Pi_{n},{\ol{\Q}_{p}}}\ot   \rho_{\Pi_{n+1}, \ol{\Q}_{p}}(n).$$
(In fact, it is conjectured that $\rho_{\Pi_{\nu, \ol{\Q}_{p}}}$ is irreducible for $\nu=n, n+1$, so that the semisimplification should be superfluous.) This also implies that $\ol{\pi}$ has a model $\pi=\Hom_{L[G_{F}]}(M^{\pi}, \rho_{\Pi})$ defined over~$L$; in other words, for an incoherent $V\in \sV^{\circ, -}$  we have $\wt\sC(\RH^{V} \bs \G^{V})_{\Q_{p}}^{\rm st}=\sC(\RH^{V}\bs \G^{V})_{\Q_{p}}^{\rm st}$.\footnote{It is plausible that this kind of equality holds more generally, but we do not explore this here.} 
\end{remark}

\subsubsection{Properties of automorphic Galois representations}

\begin{proposition} \lb{wt-mon} 
Let $\pi\in \sC(\G)_{\Q_{p}}(L)$. The Galois representation $\rho\coloneqq  \rho_{\pi, \ol{\Q}_{p}}$ satisfies the following properties:
\begin{enumerate}
\item  For every nonarchimedean place $w$ of $F$, 
the representation $\rho_{|G_{F_{w}}}$ 
is pure of weight $-1$ in the sense of \cite[Definition A.11]{DL}.
 \item The representations $\rho^{\rm c}$ and $\rho_{}^{*}(1)$ are isomorphic.
 \item For every place $v\vert p$ of $F_{0}$ and every place $w\vert v$ of $F$:
  \begin{enumerate}
 \item if $\pi_{v}$  is unramified, then $\rho_{|G_{F_{w}}}$
 is crystalline;
 \item  if moreover $\pi_{v}$  is ordinary, then $\rho_{ |G_{F_{w}}}$ 
 is Panchishkin-ordinary (i.e., it satisfies the Panchishkin condition recalled in \S~\ref{pair def}).
 \end{enumerate}
  \end{enumerate}
  
If Hypothesis \ref{hyp coh} holds, then 
 the conclusions (1)-(3) above also hold for $\rho=\rho[\pi]$ and $\rho=M^{\pi}_{K}$. 
\end{proposition}
\begin{proof} 
Part (1) is a fundamental result of Caraiani \cite{Ca, Ca2} (see also \cite[Lemma 1.4 (3)]{TY}).
  Part (2) follows from the last statement in \cite[Lemma 4.10]{DL} for $\rho_{\pi, \ol{\Q}_{p}}$, and from the Galois-equivariance of the Poincar\'e pairing for $\rho[\pi]$, $M_{\pi}$. 
The proof of part (3) is as in \cite[Lemmas 4.9, 4.14]{DL}.
 (In fact, the assumption on $\pi_{\nu,v}$ in (b)  is stronger than the analogous assumption in \emph{loc. cit.}; correspondingly, each factor $\rho_{\pi_{\nu}|G_{F_{w}}}$ is also ordinary in the sense of \cite[Definition 1.29]{N93}; however, only  Panchishkin-ordinariness is stable under tensor products.)
\end{proof}

For the rest of the paper, we will assume Hypothesis \ref{hyp coh} for every\footnote{In fact, it would be enough to assume it for the representation $\pi$ in order to prove Theorem \ref{main thm} for $\pi$, at the cost of some complication in the exposition.} representation
 $\pi\in \wt\sC(\G)(\ol{\Q}_{p})$.

\subsection{Gan--Gross--Prasad cycles}
\lb{sec: ggp cycle}
We define our cycles and study an `ordinary' modification.
\subsubsection{Arithmetic diagonal cycles}\lb{ADC 10}
We have a fundamental  cycle
$$[Y]^{\circ}=([Y_{K_\RH}]^{\circ}) \in  \varprojlim_{K_\RH}  \RZ^{0}(Y_{K_\RH})_{\bQ},$$
where the transition maps on the right are pushforwards and $[Y_{K_\RH}]^{\circ} =\vol(K_\RH, dh) [Y_{K_\RH}]$. 
Let $\jmath$ be the (system of)  arithmetic diagonal maps \eqref{eq:ADC}. The  \emph{arithmetic diagonal cycle}\beq\lb{ADC def}
Z\coloneqq \jmath_{*}[Y]^{\circ} \in  \varprojlim_{K}  \RZ^{n}(X_{ K})_{\Q}\eeq
is well-defined. 
We denote by $Z_{K}$ its image in $\RZ^{n}(X_{ K})_{\Q}$.
\subsubsection{Limits of  Selmer groups}
Let $L$ be a finite extension of $\Q_{p}$, and  let $\pi\in \wt\sC({\G})(L)$.
For $?\in \{ \rm temp, \tord, \pi\}$, define 
$$ H^{1}_{f}(F, M^{?})\coloneqq \varprojlim_{K}  H^{1}_{f}(F, M^{?}_{K}), \qquad  H^{1}_{f}(F, M_{?})\coloneqq \varinjlim_{K}  H^{1}_{f}(F,M_{?}^{K}).$$

\subsubsection{GGP cycles and associated functionals}\lb{ssec: ggp cycle} 
Let 
 $$Z_{\pi, K}\in H^{1}_{f}(F, M^\pi_{ K})$$
 be the Hecke-eigencomponent of $\wt{\rm cl}(Z_{K})$ (where $\wt{\rm cl}$ is as in \S~\ref{sec pAJ}). 
   Here, by the discussion in \S~\ref{sec AJ}, the fact that $Z_{\pi,K}$ belongs to the Bloch--Kato Selmer group is a consequence of the vanishing of $M^{2n,\pi}$ and Proposition \ref{wt-mon}.

\begin{definition}
The \emph{Gan--Gross--Prasad cycle} of $\pi$ is
$$Z_{\pi}\coloneqq \varprojlim_{K}Z_{\pi, K} \quad \in H^{1}_{f}(F, M^{\pi}).$$ 
The $\RH(\A)$-invariant functional associated to it via \eqref{Mpi dec} will  still be denoted  by
\beqq
Z_{\pi}\colon \pi &\to H^{1}_{f}(F, \rho[\pi])\\
\phi&\mapsto Z_{\pi}\phi\coloneqq  \phi_{*}Z_{\pi}.
\eeqq
\end{definition}
From the $\RH(\A)$-invariance it follows that
$Z_{\pi}$ vanishes unless $\pi\in\sC(\RH \bs \G)$.
\begin{remark}
The linear functional $Z_{\pi}$ valued in the Selmer group can be viewed as an arithmetic analogue of the automorphic period functional $P_{\pi}$ of \eqref{H per pi}.
\end{remark}

\subsubsection{Ordinary cycles}\lb{fs cycle} Suppose that every place $v\vert p$ of $F_{0}$ splits in $F$. 
For each $v$, we may fix a place $w\vert v$ of $F$ and compatible bases of $V_{\nu, w}$, giving isomorphisms   $\RG_{F_{0,v}}\cong \GL_{n}\times \GL_{n+1}$, $\RH_{F_{0,v}}\cong \GL_{n}$ as algebraic groups over $F_{0,v}=F_{w}$.  Then we may and will use the notation, definitions and results of \S~\ref{sec 5}; we generally also denote $\Box_{p}\coloneqq  \prod_{v\vert p}\Box_{v}$; for instance, $t_{0, p}=\prod_{v\vert p}t_{0,v}$, $N_{0,p}^{\circ}=\prod_{v\vert p}N_{0,v}^{\circ}\subset \GL_{n}(F_{0, p})\times \GL_{n+1}(G_{F_{0, p}}) \cong \RG(F_{0,p})$.   We define an operator $e^{\rm ord}\coloneqq  \lim U_{t_{0, p}}^{N!}$;  it acts on $M_{N^{\circ}_{0,p}}$ and
 on $\pi^{N^{\circ}_{0,p}}$ for any $\pi\in \sC(\G)_{\Q_{p}}$; the representation   $\pi$ is ordinary if and only if $\pi^{N^{\circ}_{0,p}}$ is not annihilated by $e^{\ord}$.

Let $K_{p}\subset \RG(F_{0,p})$ be an open compact subgroup containing $N_{0,p}^{\circ}$, and let $c\geq1$ be such that  $K_{p}$ contains $K_{0,p}^{\la c+1\ra}$. For  positive integers $r$, $N$ with $N!\geq r\geq c$, set  $m_{0, r}=\prod_{v\vert p }m_{0,r,v}$ (cf. \eqref{m def} for the definition of twisting matrices), and define
\beqq\lb{Zvee}
Z_{K_{p}}^{\dag, N}\coloneqq   \prod_{v\vert p}q_{v}^{rd(n)}\cdot Z.T(m_{0,r} U_{t_{0,p}}^{N!-r} e_{K_{p}})_{\Q} \ \in \RZ^{n}(X_{ K_{p}}), \eeqq
which is independent of $r$ by Corollary \ref{twisting cor}.  

We define the \emph{ordinary arithmetic diagonal cycle}  by
$$Z_{K_{p}}^{\ord}\coloneqq   \lim_{N\sto \infty}\wt{\rm cl}(Z_{K_{p}}^{\dag, N}) \quad \in \varprojlim_{K^{p}}(H^{2n}(X_{K^{p}K_{p}}, \Q_{p}(n))/{\rm Fil}^{2}).e^{\ord}.$$

For any $\pi\in \sC(\RH\bs \RG)^{\ord}$,  we define the ordinary GGP cycle  
$$Z_{\pi, K_{p}}^{\ord}\in H^{1}_{f}(F, M^{\pi}_{K_{p}})$$
to  be the eigencomponent of $Z_{K_{p}}^{\ord}$.  By the definitions, for any sufficiently large $r$,
 $$Z_{\pi, K_{p}}^{\ord} = \prod_{v\vert p}q_{v}^{rd(n)} \lim_{N\sto \infty} Z_{\pi}.T(m_{0,r}U_{t_{0}, p}^{N!-r}e_{K_{p}}).$$
We have  an induced  $\RH(\A^{p})$-invariant functional   still denoted by the same name
$$Z_{\pi, K_{p}}^{\ord} \colon \pi^{K_{p}} \to H^{1}_{f}(F, \rho[\pi]).$$
It factors through $e^{\ord}e_{K_{p}}$.

\subsubsection{Norm relation}
Continue with the notation and assumptions  of \S~\ref{fs cycle}. In order to study $p$-adic heights, it will be useful to know that $Z_{ K_{p}}^{\dag, N}$ is a norm from some ring class fields of $F$ of conductors that are high powers of $p$.

Let $\RT$ be the unitary group of the $1$-dimensional hermitian space $(F, N_{F/F_{0}})$. We have a map
$${\rm rec} \colon G_{F}\to \ol{F^{\times}} \bs\A_{F}^{\infty,\times}\to \ol{\RT(F_{0})}\bs \RT(\A^{\infty}),$$
where the first map is the reciprocity law of class field theory and the second map is $x\mapsto x^{\rm c}/x$ (and the bars denote Zariski closures). For $v\vert p$ and $r\geq0$, let $K_{T, v}^{(r)}\coloneqq \RT(\sO_{F_{0,v}}) \cap  1+v^{r}\sO_{F_{0, v}}$,  let $\Gamma_{r}=\Gamma_{r}^{(v)}\coloneqq   \RT(F_{0})\bs \RT(\A^{\infty}) / \RT(\widehat{\sO}_{F_{0}}^{v}) K_{T,v}^{(r)}$, and let 
$$F_{r}=F_{r}^{(v)}/F$$ be the abelian extension such that $\Gal(F_{r}/F)\cong \Gamma_{r}$
 under the reciprocity map.  We have the norm map
$$\RN_{F_{r}/F}\colon \RZ^{n}(X_{ K, F_{r}}) \to \RZ^{n}(X_{ K}).$$ 

\begin{lemma} \lb{N rel} Fix a place $v\vert p$ of $F_{0}$. For any $f^{p}\in \sH(\G(\A^{p\infty}), \sO_{L})_{K^{p}}$ and any integer $r$ with  $ \max\{1, c(K_{v})-1\}\leq r\leq N!$,  we have
 $$Z_{K^{p}}^{\dag, N}(f^{p})\quad  \in \quad\RN_{F_{r}/F} ( \RZ^{n}(X_{ K, F_{r}})_{\sO_{L}})\quad  \subset  \RZ^{n}(X_{ K})_{\sO_{L}}.$$
\end{lemma}
\begin{proof}
We may assume $f^{p}=e_{K^{p}}$, and abbreviate $Z_{K}^{\dag, N}=Z_{K^{p}}^{\dag, N}(f^{p}) $. 
Let $K_{H}^{p}\coloneqq  \RH(\A^{p\infty})\cap K^{p}$ and let $Y_{r}\coloneqq Y_{K_{H}^{p}K_{H,0, p}^{(r)}}$. 
Then by \eqref{mr KH}, the map $Y\stackrel{\jmath}{\to} X\stackrel{m_{ 0,r}}{\to}X\to X_{K}$ factors through $Y_{r}$, and we can write
$$Z_{K}^{\dag, N}
=  
 \prod_{v\vert p}q_{v}^{r d(n)}
\cdot (\jmath_{*}[Y_{r}]^{\circ}) .T( m_{0, r} U_{t_{0},p}^{N!-r} e_{K})
=
\vol^{\circ}(K_{H, 0 , p})
\cdot (\jmath_{*}[Y_{r}]).T(  m_{0, r} U_{t_{0},p}^{N!-r} e_{K})$$
(see \eqref{vol circ} for $\vol^{\circ}(K_{H, 0 , p})$).

 Let $\det\colon \RH\to \RT$ be the determinant map.  For a compact open subgroup $K\subset \RH(\A^{\infty})$, 
 by  Shimura's reciprocity law we have an isomorphism of $G_{F}$-sets
$$\pi_{0}(Y_{K, \ol{F}})\cong \RT(F_{0})\bs \RT(\A^{\infty})/ \det K.$$
Now we have $\det K_{H,0,v}^{(r)}= 1+v^{r}\sO_{F_{0,v}}=1+w^{r}\sO_{F_{w}}\cong K_{T,v}^{(r)}$, where the last identification comes from the natural map $F_{w}^{\times}\subset F_{v}^{\times}\to T_{v}$. Thus
we deduce a natural surjection ${\rm p}\colon\pi_{0}(Y_{r, \ol{F}})\to \Gamma_{r}$. For each $\gamma\in \Gamma_{r}$, let $Y_{r, \gamma, \ol{F}}\subset Y_{r, \ol{F}}$ be the union of connected components in ${\rm p}^{-1}(\gamma)$; it arises as $Y_{r, \gamma}\times_{F_{r}}\ol{F}$ for an $F_{r}$-subvariety $$Y_{r, \gamma}\subset Y_{r, F_{r}}.$$ Then for any $\gamma_{0}\in \Gamma_{r}$, we have
$$Z_{K}^{\dag, N} = \vol^{\circ}(K_{H, 0 , p})
\cdot \RN_{F_{r}/F} (\jmath_{*}[Y_{r, \gamma_{0}}]).T(  m_{0, r} U_{t_{0},p}^{N!-r} e_{K}),$$
which belongs to $\RN_{F_{r}/F}( \RZ^{n}(X_{ K, F_{r}})_{\Z_{p}})$ since $\vol^{\circ}(K_{H, 0 , p})$ is a $p$-unit. 
\end{proof}

\subsection{The distribution and its spectral expansion}
\lb{sec: J def}
From now until the end of the paper, we suppose that every place $v\vert p$ of $F_{0}$ splits in $F$, and that $K_{p}\subset \G(F_{0, p})$ is a relative selfdual hyperspecial subgroup (as defined in \S~\ref{sss rel sd}).

If $L$ is an algebraic extension of $\Q_p$, we say that $\pi\in \widetilde \sC(\G)(L)$ is \emph{$K_p$-ordinary} if $\pi$ is ordinary and has $K_p$-fixed vectors. (By Lemma \ref{nv ord unr}, this is equivalent to the condition that ${\rm BC}(\pi)$ is $K_p\ts K_p$-ordinary.) This condition defines an ind-subscheme
$$ \sC(\RH\bs \G)_{K_{p}}^{\ord}\subset  \sC(\RH\bs \G).$$

\subsubsection{Height pairings} 
Considering the setup of \S~\ref{sec: ht gl cyc}, we denote by 
\beq \lb{ht coh}
h\colon  H^{1}_{f}(F, M_{\tord}^{K_{p}})\ts H^{1}_{f}(F, M^{\tord}_{K_{p}})
 \to \Gamma_{F_{0}, L}\eeq
the pairing induced by the family 
$h_{X_{K}, \lm}\colon H^{1}_{f}(F, M^{K}_{\tord})^{\ot 2}\to \Gamma_{F_{0}, L}$ for $K=K^{p}K_{p}$,
where
$$\lm\colon \Gamma_{F, L}\to \Gamma_{F_{0}, L}$$
is the natural surjection.
 It is well-defined by the projection formula (Lemma \ref{proj formula}).
 Note that the conditions of \S~\ref{pair def} for the definition of $h$ (as well as for the definition of the pairing $h_{\pi}$ from \S~\ref{ht intro}) are satisfied by Proposition \ref{wt-mon}.  
We also have a pairing (abusively denoted by the same name)
\beq\lb{ht coh2} h\colon  H^{1}_{f}(F, M_{\tord}^{K_{p}})\ts H^{1}_{f}(F, M_{\tord}^{K_{p}})  \to \Gamma_{F_{0}, L}
\eeq
obtained from \eqref{ht coh} by composing with the map induced by \eqref{inv to dir} on the second factor.

 For a non-archimedean place $w$ of $F$, we  denote by $h_{w}$ the corresponding local pairings \eqref{dec ht} on  pairs of (limits of) cycles with disjoint supports in $\RZ^{n}_{\tord}(X_{K, F_{w}})^{0}_{L}$ (if $w\vert p$) or $\RZ^{n}_{\rm temp}(X_{K, F_{w}})^{0}_{L}$ (if $w\nmid p$). For $w\nmid p$, this requires a projection formula for $w$-local heights, which is equivalent  to Lemma \ref{bc ht} (2).

\subsubsection{Definition of the distribution}
For $S$ a finite set of non-archimedean places of $F_{0}$ and $?\in \{{\rm temp}, \tord\}$, denote  by
 $$\sH(\G(\A^{S}),L)^{\circ}_{K_{S}\text{-?}}\subset \sH(\G(\A^{S}),L)_{K_{S}}$$
the subalgebra of measures $f^{S}=f^{S\infty}f_{\infty}$ such that $f_{\infty}\in L f_{\infty}^{\circ}$ (where $f_{\infty}^{\circ}= \eqref{f inf})$ and $M^{\oplus}.\hT(f^{S}e_{K_{S}})\subset M^{K_{S}}_{ {\rm ?}}$.

Define first
\beq\lb{j12}
\sJ_{K_{p}}\colon (\sH(\G(\A^{p}),L)^{\circ}_{K_{p}\text{-temp}})^{2}
  &\to \Gamma_{F_{0},L}\\
(f_{1}^{p}, f_{2}^{p}) &\mapsto  h(Z_{K_{p}}^{\ord}. \hT(f_{1}^{p}), Z^{\ord}_{K_{p}}.\hT(f_{2}^{p})),
\eeq
where the right-hand side uses the pairing \eqref{ht coh2}.
\begin{definition}\lb{artf}
For any $f^{p}\in \sH(\G(\A^{p}),L)^{\circ}_{K_{p}}$ that can be written as 
\beq \lb{dec fp}
f^p=f_{1}^{p}\star f_{2}^{p, \vee}\eeq
with $f_{i}^{p} \in \sH(\G(\A^{p}),L)^{\circ}_{K_{p}\text{-temp}}$,  we  define the \emph{arithmetic relative-trace distribution} by\footnote{The abuse of notation with respect to \eqref{j12}
should cause no confusion.}
$$\sJ_{K_{p}}(f^{p})\coloneqq  \sJ_{K_{p}}(f_{1}^{p}, f_{2}^{p}).$$
\end{definition}
\begin{remark}
The definition  is independent of the decomposition  \eqref{dec fp}. Indeed, let $K^{p}\subset \G(\A^{p\infty})$ be such that $f^{p}\in \sH(\G(\A^{p\infty}))_{K^{p}}$, and let $S$ be a finite set of finite places of $F_{0}$, not above $p$, such that  $K^{S}\coloneqq K^p\cap \G(\A^{Sp\infty})$ is a maximal  hyperspecial subgroup. Let $e^{\rm temp}_{K}\in \sH(\G(\A^{S\infty}))_{K^{S}}$ be an element acting as the  idempotent projection $M^{K}\to M^{K}_{\rm temp}$. Then by the projection formula (Lemma \ref{proj formula}), for each decomposition \eqref{dec fp} we have
$$\sJ_{K_{p}}(f_{1}^{p}, f_{2}^{p})= h(Z_{K_{p}}^{\ord}. \hT(f_{1}^{p}), Z^{\ord}_{K_{p}}.\hT(f_{2}^{p}e^{\rm temp}_{K}))= h(Z_{K_{p}}^{\ord}. \hT(f^{p}), Z^{\ord}_{K_{p}}.T(e^{\rm temp}_{K})),$$
which shows that  $\sJ_{K_{p}}(f^{p})$ is well-defined.
\end{remark}

Let now
\beq\lb{fp G}
f_{p, K_{p}, N}\coloneqq  \prod_{v\vert p}q_{v}^{rd(n)}\cdot m_{0,r, p} U_{t_{0}, p}^{N!-r}e_{K_{p}},
\eeq
where $1\leq r\leq N!$.
By the definitions, we have  
$$\sJ^{}_{K_{p}}(f^{p}_{1}, f_{2}^{p})
= \lim_{N\sto \infty}  h(Z.\hT(f_{1}^{p}  f_{p, K_{p}, N})
 , Z.\hT^{}(f_{2}^{p} f_{p, K_{p}, N})).
$$
It is independent of the choice of $r\leq N!$.

\subsubsection{Spectral expansion} \lb{J sp exp}
Let $\pi\in \wt \sC(\RH\bs\G)_{K_{p}}^{\ord}(L)$.
Denote by 
$$h_{M_{\pi}} 
\colon 
H^{1}_{f}(F, {M}_{\pi}^{K_{p}})\times  H^{1}_{f}(F, {M}^{\pi^{\vee}}_{K_{p}})\to \Gamma_{F_{0}, L}$$
 the restriction of $h$.
For any $f^{p}
\in  \sH(\G(\A^{p}),L)^{\circ}$, we define
$$\sJ_{\pi, K_{p}}(f^{p}) \coloneqq  h_{M_{\pi}}(Z^{\ord}_{\pi, K_{p}}. \hT(f^{p}) , Z^{\ord}_{\pi^{\vee}, K_{p}})
= 
\Tr_{(,)_{\pi^{K_{p}}}}^{h_{\pi}\circ\left(Z_{\pi, K_{p}}^{\ord}\boxtimes Z_{\pi^{\vee},K_{p}}^{\ord}\right)} (\hT(f^{p})),$$
where the pairing ${(,)_{\pi^{K_{p}}}}$ is the restriction of $(, )_{\pi}=\eqref{pair pi}$ to $\pi^{K_{p}}\times \pi^{\vee, K_{p}}$.
Then it is clear that if $f^{p}$ is as in Definition \ref{artf}, we have
$$\sJ_{K_{p}}(f^{p})=\sum_{\pi\in \sC(\RH\bs \G)_{K_{p}}^{\ord}} \sJ_{\pi ,K_{p}}^{}(f^{p}),$$
where for a Galois orbit $\pi=\{\pi^{\sg}\} \in \sC(\RH\bs\G)_{K_{p}}^{\ord}$ of isomorphism classes of representations, we put $\sJ_{\pi, K_{p}}\coloneqq \sum \sJ_{\pi^{\sg}, K_{p}}$.

\subsection{Vanishing of local heights at split places}
\lb{sec: J dec}
We will complete the arithmetic relative-trace formula by finding a geometric expansion for the distribution $\sJ_{K_{p}}$. Each term in the expansion will  be a sum over all  nonsplit finite places of $F_{0}$. The goal of this subsection is to show the preliminary result that $\sJ_{K_{p}}$ has a decomposition as a sum over nonsplit places, by proving some vanishing results for local height pairings at split ($p$-adic and non-$p$-adic) places.

\subsubsection{Decomposition over all places}
Let $v$ be a non-archimedean place of $F_{0}$. We define
$$  \sJ_{K_{p}}^{(v), N}(f_{1}^{p}, f_{2}^{p})\coloneqq \sum_{w\vert v} h_{w}(Z_{K_{p}}^{\dag, N}.\hT(f^{p}_{1}), Z_{K_{p}}^{\dag, N,}.T^{}(f_{2}^{p}))$$
for any $f_{1}^{p}$, $f_{2}^{p}\in \sH(\G(\A^{p}),L)^{\circ}_{K_{p}\text{-temp}}$ (respectively $\sH(\G(\A^{p}),L)^{\circ}_{K_{p}\text{-t-ord}}$ if $v\vert p$)  such that the two cycles involved have disjoint supports.
 Here, the sum ranges over the (one or two) places of $F$ above $v$.

It is then clear from the definitions that for $f_{1}^{p}$, $f_{2}^{p}\in \sH(\G(\A^{p}),L)^{\circ}_{K_{p}\text{-t-ord}}$, we have a decomposition
\beq \lb{dec sJ}
\sJ^{}_{K_{p}}(f^{p}_{1}, f_{2}^{p}) = \lim_{N\sto \infty}\sum_{v\nmid \infty}  \sJ_{K_{p}}^{(v), N}(f_{1}^{p}, f_{2}^{p}).\eeq

In the rest of this subsection, we show the vanishing of the contribution at split ($p$-adic and non-$p$-adic) places. 

\begin{remark}
If $v\nmid p\infty$, we can more generally define
$$\sJ^{(v)}(f_{1}, f_{2})\coloneqq  \sum_{w\vert v}  h_{w}(Z.\hT(f_{1}), Z.T^{}(f_{2}))$$
for $f_{1}$, $f_{2}\in \sH(\G(\A), L)_{\rm temp}^{\circ}$ such that the two cycles involved  have disjoint supports; then 
\beq \lb{sJv rel}
\sJ_{K_{p}}^{(v), N}(f_{1}^{p}, f_{2}^{p})= \sJ^{(v)}(f_{1}^{p}f_{p, K_{p}, N}, f_{2}^{p}f_{p, K_{p},  N}).\eeq
\end{remark}

\subsubsection{A regularity condition}  

\begin{definition}\lb{reg pair}
 Let $v$ be a finite place of $F_0$, and let $K_v\subset G_v$ be a compact open subgroup. We say that a pair $(f_{1,v},f_{2,v})\in \sH(G_{v})_{K_{v}}^{2}$ 
 is $K_{v}$-\emph{regular} if $f_{1,v}$ has regular support and $f_{2,v}=e_{K_{v}}$.

 If $S$ is a finite set of finite places of $F_{0}$ and $v\notin S$ is another finite place of $F_{0}$, we  say that a pair $(f^{S}_1,f^{S}_2)\in \sH(\G(\A^{S}), L)_{K^{S}}^{\circ}$ is \emph{$K^{S}$-regular} at $v$ if we can write $K^{S}=K^{Sv}K_{v}$ and  $f_{i}^S=f_{i,v}\ot f_{i}^{Sv}$ with $(f_{1, v}, f_{2,v})$ $K_{v}$-regular.
\end{definition}

\subsubsection{Auxiliary Shimura varieties} \lb{sec: aux sh}
We will use some auxiliary Shimura varieties associated to the following type of data. 
 \begin{definition}\lb{def CM type} Let $v$ be a finite place of $F_0$ and $w$ a place of $F$ above $v$. An \emph{auxiliary datum} relative to $w$ is a pair $(\Phi, u)$ where $\Phi$ is a CM type of $F$ that contains $\vphi_0$, and $u$ is a place of  $E=E_\Phi F$ above $w$ satisfying: 
 \begin{itemize}
     \item $u$  is unramified over $w$;
     \item if $v$ splits in $F$, the `matching condition' \eqref{cond CM} holds.
 \end{itemize}
\end{definition}

Denote by  $$d_{F/F_0} \coloneqq N_{F_0/\Q}(D_{F/F_0})$$ the absolute norm of the relative discriminant of $F/F_0$.

\begin{remark}\lb{rem exs CM}
For every finite place $v\nmid d_{F/F_0}$ of $F_0$ and every place $w$ of $F$ above $v$, as observed in \cite[p.851]{LL1}, there exists an auxiliary datum relative to~$w$. 
(For a weaker condition, see \cite[Definition 1.1, Remark 1.2]{LL2}.)
\end{remark}

For an auxiliary datum $(\Phi, u)$ (relative to some finite place $w$ of $F$), we consider the auxiliary Shimura variety 
\beq\lb{aux X'}
X'_{ u}\coloneqq X'_{\Phi, K/E_{u}} \coloneqq  {\rm Sh}_{K_{\wt \G}} (\wt\G)_{E_{u}},
\eeq
 and its integral model from \S \ref{ss:int mod}
 \beq\lb{aux int X'}
 \sX'_{u}\coloneqq  \CM_{K_{\wt\G}, \sO_{E,u}}.
 \eeq
Note that the notation hides the dependence on $\Phi$, for the sake of lightness. Recall also from \S\ref{ss:RSZ int} that $X'_{ u}$ further depends on the choice of a compact open subgroup $K_{\RZ^\Q}\subset \prod'_\ell {\RZ^\Q}(\Q_\ell)$ that we take to be maximal at the rational prime underlying~$v$. 
We observe  that $X'_{u}$ is of the form $X_{F_{w}}\times_{F_{w}} A$ for some finite \'etale $F_w$-algebra $A$.

We denote by
$$\sZ'_{u}\coloneqq  \wt{\jmath}_{*}\left(\vol(K_{\RH}) [ \wt\CM_{K_{\wt \RH}}]\right)$$ the $\sO_{E,u}$-integral model of the arithmetic diagonal cycle, where $\wt{\jmath}$ is as in  \S~\ref{ss:ADC int mod}. 

\subsubsection{Local heights at split places} 
The following lemma will be useful for considerations both at places above $p$ and away from $p$.  
 \begin{lemma}\lb{disjt sp fib} Let $v$
 be a split  place of $F_{0}$.
 Let $K=\prod_{v}K_v$ be an open subgroup of $\G(\A^{\infty})$, and let  $f_1,f_2\in \sH(\G(\A), L)_{K}^{\circ}$. Suppose that: 
\begin{itemize}
\item $(f_{1}, f_{2})$ is $K_{v'}$-regular at some finite place $v'\neq v$;
\item the subgroup $K_{v}=K_{n, v}\ts K_{n+1,v}$ satisfies either of the following conditions:
\begin{enumerate}
\item[(a)]    for some labelling $\{\nu, \nu'\} =\{n, n+1\}$, the subgroup  $K_{\nu,v}$ is maximal hyperspecial and $K_{\nu', v}$ is the principal congruence subgroup of level $m\in \Z_{\geq 0}$ (cf. \S~\ref{ss:int model spl});
\item[(b)] for both $\nu=n, n+1$, the subgroup $K_{\nu, v}$ is Iwahori (that is, $G_{\nu, v}$-conjugate to the standard Iwahori $\Iw_{\nu, v, 0}$).
\end{enumerate}
\end{itemize}
Then 
the following statements hold.

\begin{itemize}
\item[(i)] The cycles $Z.\hT(f_1)$ and $Z. \hT(f_2)$ have disjoint support (on the generic fiber).
\item [(ii)]  Let $(\Phi,u)$ be an auxiliary datum relative to a place $w\vert v$ of $F$, and consider the auxiliary objects associated in \S~\ref{sec: aux sh}. Abusing notation,  we still let $\hT(f_{i})$ denote the (flat) correspondence on the integral model $\sX_{u}'$.
Then the cycles 
$\sZ'_{u}.\hT(f_{1})$ and $ \sZ_{u}'. \hT(f_{2})$ have disjoint supports in $\sX'_{u}$. 
\end{itemize}
\end{lemma}
\begin{proof}
Part $(i)$ follows from \cite[Theorem 8.5~(i)]{RSZ3}. (The result in {\it loc. cit.} only treats the auxiliary Shimura variety attached to $\wt \G$; but it implies the desired result  for $\G$.)

For $(ii)$ case (a),
 the integral model with Drinfeld $m$-level structure  at one factor and with hyperspecial level ($m=0$) at the other factor is regular. The proof of \cite[Theorem 8.5~(ii)]{RSZ3} (only the case $f_2=e_{K}$ was considered there) still applies to show that  the cycles $ \sZ_{u}'. \hT(f_1)$ and  $\sZ_{u}'. \hT(f_2)$ have disjoint supports in $\sX'_{u}$.

 In case (b), the integral model is the resolution given in  \S \ref{ss:int mod} of the moduli scheme and the Hecke correspondences are  obtained  by base change and hence remain finite flat. The cycles are obtained  by strict transforms. Hence it suffices to show the disjointness before the resolution, which again follows from \cite[Theorem 8.5~(ii)]{RSZ3}. (Strictly speaking, the result of {\it loc. cit.} concerns  the case of Drinfeld $m$-levels rather than Iwahori level. However, we may pull back the cycles to the moduli scheme with Drinfeld level for $m=1$  and then apply that result.) 
\end{proof}

\subsubsection{Vanishing at split places away from $p$}

\begin{proposition} \lb{vanish ht spl}
Let $v\nmid p d_{F/F_0}$ 
be a split place of $F_{0}$.  Let $f_{1}$, $f_{2}\in \sH(\G(\A), L)_{\rm temp}^{\circ}$ be  as in Lemma \ref{disjt sp fib}. Assume furthermore that either $K_{v}$ is hyperspecial, or that for each $w\vert v$ there exists an auxiliary datum $(\Phi, u)$ relative to~$w$ such  that $\hT(f^{v}_{1})$, $\hT(f^{v}_{2})$ annihilate $H^{2n}(\sX_{u}', L(n))$. 
Then
$$\sJ^{(v)}(f_{1}, f_{2})=0.$$
\end{proposition} 
\begin{proof} 
We show that  $ h_{w}(Z.\hT(f_{1}), Z.T^{}(f_{2}))=0$ for each of the two places $w|v$.
By  Lemma \ref{bc ht} (2),
 it  suffices to show the vanishing of the local height after pull-back to the auxiliary Shimura variety $X'_u$ attached to some auxiliary datum $(\Phi, u)$. (If $K_v$ is hyperspecial, we take $(\Phi, u)$ to be arbitrary and in this  case the integral model $\sX'_u$ is smooth. Otherwise, we take  $(\Phi, u)$ to satisfy the given cohomological vanishing condition.)
Under our assumption, Proposition \ref{B13}  further reduces the question to the vanishing of the arithmetic intersection pairing on the model  $\sX_{u}'$ over $\sO_{E,u}$. This last vanishing follows from Lemma \ref{disjt sp fib} $(ii)$.
\end{proof}

\subsubsection{Vanishing at $p$-adic places}
\begin{proposition} \lb{vanish ht p} 
Let $v$ be a place of $F_{0}$ above $p$ (hence split and not dividing $d_{F/F_0}$). If $n>1$, assume $p>2n$.  
Let $f^{p}=f_{1}^{p}\star f_{2}^{p, \vee} \in \sH(\G(\A^{p}), L)^{\circ}_{K_{p}\text{-}{\tord}}$, and assume that the pair $(f_{1}^{p}, f_{2}^{p})$   has regular support. Then 
$$\lim_{N\sto \infty} \sJ^{(v), N}_{K_{p}}(f^{p})=0.$$
\end{proposition}
\begin{proof}
Write $Z_{1}\coloneqq Z_{K_{p}}^{\dag, N} .\hT(f_{1}^{p})$, $Z_{2}\coloneqq Z_{K_{p}}^{\dag, N} .T^{}(f_{2}^{p})$, and let $K^{p}$ be such that $f_{1}^p$, $f_{2}^p$ are right-$K^{p}$-invariant.
We  will show  that 
$$h_{X_{K,F_{w}}}(Z_{1}, Z_{2})\quad\in p^{N!-C} \sO_{F_{w}}^{\times}\hat{\ot}\sO_{L}$$ for some constant $C$; after taking limits, this implies the desired vanishing.  Up to multiplying by a nonzero scalar, we may assume that $f_{i}^{p}\in \sH(\G(\A^{p}), \sO_{L})^{\circ}$. 

By Lemma \ref{N rel}, for some constant $C'$ cancelling the denominators of $f_{i}^p$, and for any sufficiently large $r\leq N!$, we have $  Z_{i}=\RN_{F_{r}/F}(Z_{i,r})$ for some $Z_{i,r}\in p^{-C'} \RZ^{n}(X_{ K, F_{r}})_{\sO_{L}}$. Denote by $F_{w,r}$ the localization of $F_{r}$ at its unique place above $w$, and by $h_{X_{K, F_{w,r}}} \coloneqq h_{X_{K,F_{w,r}}, \lm_{F_{w,r}}'}$ the corresponding height pairing, where $\lm_{F_{w,r}}'\colon F_{w,r}^{\ts}\hat\ot L\to F_{w,r}^{\ts}\hat{\ot}L$ is the identity map.
First, we show that 
\beq\lb{vanish p1} h_{X_{K,F_{w,r}}}(Z_{1}, Z_{2,r}) \quad \in \sO_{F_{w,r}}^{\ts}\hat{\ot} L.\eeq
By Lemma \ref{bc ht} (1) (which applies thanks to the observation made after \eqref{aux X'}), it is enough to show the same result for the corresponding height pairing of  arithmetic diagonal cycles on the auxiliary Shimura variety \eqref{aux X'} (for some arbitrary auxiliary datum). 
This follows from Lemma \ref{disjt sp fib} (ii),  Proposition \ref{crys prop}, and Remark \ref{rmk crys}. 

By the integrality results of \cite[Proposition II.1.11]{N95}, we have in fact 
\beq\lb{Cr} h_{X_{K,F_{w,r}}}(Z_{1}, Z_{2,r}) \quad \in p^{-C''_{r}-C'} \sO_{F_{w,r}}^{\ts}\hat{\ot} \sO_{L},\eeq
for a constant $C''_{r}$ that, similarly to \cite[Proof of Proposition 4.35]{DL},
 can be bounded as follows. Let $\sM\coloneqq  M^{\tord}_{K, L}\cap H^{2n-1}(X_{K, \ol{F}}, \sO_{L}(n))/({\rm tors})$, and denote by 
 $${\rm N}_{\infty}H^{1}_{f}(F_{w,r},\sM)
 \coloneqq \ 
 \bigcap_{s\geq r} 
 {\rm Im}\left[\Tr_{F_{w,s}/F_{w,r}} \colon H^{1}_{f}(F_{w,s}, \sM) \to  H^{1}_{f}(F_{w,r}, \sM)\right].$$
 Then $p^{C''_{r}}\leq c_{r}''\coloneqq  |  H^{1}_{f}(F_{w,r}, \sM)/ {\rm N}_{\infty}  H^{1}_{f}(F_{w,r}, \sM)|$. However $c_{r}''$ is bounded independently of $r$: this follows  by the same argument as for \cite[Lemma 4.37]{DL} from the fact that $ M^{\tord}_{ K, L}$, as a representation of $G_{F_{w}}$,  is crystalline, Panchishkin-ordinary, and pure of weight~$-1$ (Proposition \ref{wt-mon}). Thus in \eqref{Cr} we may replace $C''_{r}+C'$ by a constant $C''$.

Finally, by  Lemma \ref{bc ht} (1), we have
$$p^{C''}\cdot h_{X_{K,F_{w}}} (Z_{1}, Z_{2}) =
p^{C''} {\rm N}_{F_{w,r}/F_{w}}
  (h_{X_{K}, F_{w,r}}(Z_{1, F_{w,r}}, Z_{2,r})) \quad \in  {\rm N}_{F_{w,r}/F_{w}}(\sO_{F_{w,r}}^{\ts}\hat{\ot} \sO_{L}). $$
By the definition of $F_{w,r}$ and local class field theory, $  {\rm N}_{F_{w,r}/F_{w}}(\sO_{F_{w,r}}^{\ts}\hat{\ot} \sO_{L})\subset   p^{r-C'''}(\sO_{F_{w}}^{\ts}\hat{\ot} \sO_{L})$ for some constant $C'''$. This completes the proof. 
\end{proof}

\subsection{The arithmetic relative-trace formula}
\lb{sec: J geom}
The previous subsection shows that, for suitable $f_{1}^{p}$, $f_{2}^{p}$,  we have a decomposition
$$\sJ^{}_{K_{p}}(f^{p}_{1}, f_{2}^{p}) = \lim_{N\sto \infty}\sum_{v\nmid \infty}  \sJ_{K_{p}}^{(v), N}(f_{1}^{p}, f_{2}^{p}),$$
where the sum runs over places that are nonsplit or divide $d_{F/F_0}$. 
We state a geometric expansion of   $ \sJ_{K_{p}}^{(v), N}$ (in fact,  $\sJ^{(v)}$)
for inert places $v$ with a special parahoric level. When $F/F_{0}$ is unramified, we then  deduce  a geometric expansion of $\sJ_{K_{p}}$, thus completing the corresponding RTF.

\subsubsection{Local arithmetic intersection numbers and geometric expansions at inert places} 
Let $v\nmid 2p$ be an inert finite place of $F_{0}$ and let $w$ be the unique place of $F$ above $v$. Let $\e=\e_v\in \{0, 1\}$ have the same parity as $v(\nrm)$. 
Consider the vertex parahoric subgroup $K_v=K_{n,v}\times K_{n+1,v}$ of type $(t,t+\e)$ for $0\leq t\leq n$, defined in \S~\ref{sss:vert}.

We define for $\delta\in  \G^{V(v)}_{\rm rs}(F_{0,v})$,
\begin{equation}\lb{eq: def int RZ}
\sJ_{\delta, v}(e_{K_v})\coloneqq - (\delta\cdot\CN_{n,v},\CN_{n,v})_{\CN_{(n,n+1),v}}\, \lm(\vpi_{w}),
\end{equation}
where, in the right hand side, $(-,-)$ denotes the arithmetic intersection number on the (relative) unitary Rapoport--Zink space $\CN_{(n,n+1),v}\coloneqq \CN_{n,v}\times_{\Spf \sO_{ \breve{F}_v}}\CN_{n+1,v}  $ (resp. the small resolution in \cite{zhiyu}) if $K_v$ is hyperspecial (resp. $K_v$ is a vertex parahoric of type $(t,t+\epsilon)\neq (0,0),(n,n+1)$), relative to the quadratic field extension $F_{w}/F_{0,v}$. Since  $\sJ_{\delta, v}(e_{K_v})$ only plays an intermediate role, we refer to \cite{afl-pf2} (resp. \cite{zhiyu}) for the unexplained notation in the hyperspecial  (resp. vertex parahoric) case.

Recall the matching of global orbits $\ul{\delta}$ of \eqref{dec Brs gl}, and the characteristic function $\one_{V'}$ of those orbits matching one from a given $V'\in \sV^{\circ}$ from  \S~\ref{63 notn}. 
{Recall also the local  unitary-group orbital integrals $J_{\delta,v}$  defined in \S~\ref{loc J def}.}

\begin{proposition} \lb{factor J inert}  Let $v\nmid 2p d_{F/F_0}$ be an inert finite place of $F_{0}$.  Let $K=K_{v}K^{v}\subset \G(\A^{\infty})$ be a compact open subgroup,   let 
$$f_{1}=f_{1,v}f_{1}^{v}, \quad f_{2} =f_{2,v}f_{2}^{v}\quad \in \sH(\G(\A), L)_{K{\textup{-temp}}}^{\circ},$$ and let $f= f_{1}\star f_{2}^{\vee}$.
 Suppose that:
 \begin{enumerate}
\item 
$(f_1,f_2)$ is $K$-regular at a place different from  $v$;
 \item 
  $f_{1,v}=f_{2,v}=e_{K_v}$ where $K_v$ is a vertex parahoric subgroup of type $(t,t+\epsilon_v)$ (cf. \S~\ref{sss:vert});
  \item 
   $K_{v}$ is hyperspecial or, for some auxiliary datum $(\Phi, u)$ relative to the unique place of $F$ above~$v$, the correspondences $\hT(f_{1})$, $\hT(f_{2})$ annihilate $H^{2n}(\sX'_{u}, L(n))$. 
  \end{enumerate}
Then 
\begin{align*}
 \sJ^{(v)}(f_{1} , f_{2})  =& \sum_{\delta\in \RB^{V(v)}_{\rm rs}(F_{0})}  J^{v}_{\delta}(f^{v}) \sJ_{\delta, v}(f_{v})\\
 =&  \sum_{\gamma\in \RB'_{\rm rs}(F_{0})}  \one_{V(v)}(\gamma) J^{v}_{\ul{\delta}(\gamma)}(f^{v}) \sJ_{\ul{\delta}(\gamma),v}(f_{v}).
\end{align*}
\end{proposition} 
\begin{proof}It suffices to show the first equality. Similar to the proof of Proposition \ref{vanish ht spl}, by  the base change property of Lemma \ref{bc ht}(2)
we have 
$$
\sJ^{(v)}(f_{1} , f_{2})=\frac{1}{\deg(X'_{u}/X_w)} h_{u}(Z'.\hT(f_{1}), Z'.T^{}(f_{2})),
$$
where $h_{u}$ denotes the local height on $X'_{u}$ over $\sO_{E,u}$.  Under our assumption, by Proposition \ref{B13} we have
$$
h_{u}(Z'.\hT(f_{1}), Z'.T^{}(f_{2}))= 
-
(\sZ_{u}'.\hT(f_{1}),\sZ_{u}' .T^{}(f_{2}) )
 \, \lm(\Nm_{E_{u}/F_w}\vpi_{u}).
$$
Since $\lm_{|F_{w}^{\times}}$ is necessarily unramified and $E_{u}/F_w$ is an unramified extension, we have $$\lm(\Nm_{E_{u}/F_w}\vpi_{u})=\deg(E_{u}/F_w)\lm(\vpi_{w}).$$

In the hyperspecial case, by \cite[Theorem 8.15]{RSZ3}  (the statement there is for the sum over all places of $E$ above $w$, but the proof contains the formula for each place $u$),  we obtain
$$
 -(\sZ_{u}'.\hT(f_{1}),\sZ_{u}' .T^{}(f_{2-}) )
 =\deg(X'_{u}/X_w) \sum_{\delta\in \RB^{V(v)}_{\rm rs}(F_{0})}  J^{v}_{\delta}(f^{v}) \sJ_{\delta, v}(f_{v}).
$$
The vertex parahoric case is similar \cite{zhiyu,LMZ}, and we omit the details.
\end{proof}

\subsubsection{The arithmetic relative-trace formula}
We are ready to deduce the following relative-trace formula for $\sJ_{K_{p}}$.

\begin{theorem}[Arithmetic relative-trace formula]
 \lb{ARTF thm}
Suppose that:
\begin{itemize}
\item
$F/F_{0}$ is unramified,
\item $p>2n$ if $n>1$,
 \item all places $v\vert 2p$ of $F_{0}$ are split in $F$.
\end{itemize}
Suppose also that there is a finite set $S$ of places of $F_{0}$, not above $p$ or $\infty$,  and a compact open subgroup $K^{p}=\prod_{v\nmid p}K_v \subset \G(\A^{p\infty})$ satisfying:
\begin{itemize}
\item  $K_{v}$ is relative selfdual hyperspecial (\S~\ref{sss rel sd}) for $v \notin S$,
\item for every split place $v\in S$,  $K_{v}=K_{n, v}\ts K_{n+1 ,v}$ where either at least one of the factors is maximal hyperspecial, or both are Iwahori,
\item for every inert $v\in S$, $K_v$ is a vertex parahoric subgroup of type $(t,t+\epsilon_v)$ (cf. \S\ref{sss:vert}).
\end{itemize}
For $i=1, 2$, let $ f_{i}^{p}=f_{i}^{S p} \ot \ot_{v\mid S } f_{i,v} \in \sH(\G(\A^p), L)_{K^{p}\textup{-t-ord}}^{\circ}$   
satisfy the following properties:
\begin{itemize}
\item for every inert $v$,  $f_{1,v}=f_{2,v}=e_{K_v}$,
\item  for two (necessarily split) places $v\in S$, the pair $(f_{1,v},f_{2,v})$ is $K_{v}$-regular (in the sense of Definition \ref{reg pair});
\item  for every finite place $v\in S$ and every place $w$ of $F$ above $v$, there exists an auxiliary datum $(\Phi, u)$ (Definition \ref{def CM type})
such that $\hT(f_{i}^{S p})$ annihilates $H^{2n}(\sX'_{u},L(n))$,
where  $\sX'_{u}$ is the integral model defined by \eqref{aux int X'}.
\end{itemize}
 Let $f^{p}\coloneqq f_{1}^{p}\star f_{2}^{p, \vee}\in \sH(\G(\A^p), L)_{K^{p}}^{\circ}$.
Then we have a spectral and a geometric expansion\footnote{The geometric expansion implicitly includes the assertion that the integrand is integrable for the Radon measure $I^{\ord}_{\gamma,p, K'_{p}}$.}
\begin{align*} \sJ_{K_{p}} (f^{p})=&\sum_{\pi\in \sC(\RH\bs\G)_{K_{p}}^{\ord}} \sJ_{\pi, K_{p}}(f^{p})\\
= &\int_{\RB'_{\rs}(F_{0})^{\circ}} \sum_{v\nmid p\infty \atop \rm nonsplit}
\one_{V(v)}(\gamma) J^{vp}_{\ul{\delta}(\gamma)}(f^{vp}) \sJ_{\ul{\delta}(\gamma),v}(f_{v}) \, d I^{\ord}_{\gamma,p, K'_{p}},
\end{align*}
where   $d I^{\ord}_{\gamma,p, K'_{p}}= dI^{\ord}_{\gamma,p, K'_{p}}(\one_{p}) $ is as in \eqref{muK def}  for   $K_{p}'= \RG'(\sO_{F_{0, p}})$.
\end{theorem}

\begin{proof} 
The spectral expansion was noted in \S~\ref{J sp exp}. We establish the geometric expansion. 
By \eqref{dec sJ}, we have
$$\sJ_{K_p}(f^{p})=\lim_{N\sto \infty} \sum_{v\nmid \infty} \sJ^{(v), N}_{K_{p}} (f_{1}^{p}, f_{2}^{p}).$$
By Propositions \ref{vanish ht spl}, \ref{vanish ht p}, only the terms corresponding to nonsplit places $v\nmid p$  contribute. (We use the `second' place of regular support  to apply Proposition \ref{vanish ht spl} to the `first' one.) By  \eqref{sJv rel} and Proposition \ref{factor J inert}, we then have
$$\sJ_{K_p}(f^{p})=\lim_{N\sto \infty} \sum_{v\nmid p\infty  \atop \ {\rm nonsplit}}
 \sum_{\gamma \in \RB_{\rm rs}'(F_{0})} 
\one_{V(v)}(\gamma) J^{vp}_{\ul{\delta}(\gamma)}(f^{vp}) \sJ_{\ul{\delta}(\gamma),v}(f_{v}) \cdot J_{\ul{\delta}(\gamma)} (f_{p, K_{p}, N} \star f_{p, K_{p}, N}^{\vee}).$$
The asserted form of the geometric expansion then follows, via Lemma \ref{match split} and Lemma  \ref{I Pi dag},   from the definition  of $d I^{\ord}_{\gamma,p, K'_{p}}$.
\end{proof}

\part*{Epilogue}

\section{Comparison of RTFs and proof of the main theorem} \lb{sec:11}

In this concluding section, we compare the arithmetic distribution $\sJ_{K_{p}}$ with the derivative   $\partial\sI_{K_{p}'}$  of the analytic distribution, and  deduce our main theorem.

Throughout this section we assume:
\begin{itemize}
\item
$F/F_{0}$ is unramified,
\item $p>2n$ if $n>1$,
 \item all places $v\vert 2p$ of $F_{0}$  split in $F$.
\end{itemize}

In \S~\ref{sec:121}, we prove our  comparison theorem. In \S~\ref{sec:122}, we construct appropriate matching Hecke measures. In \S~\ref{sec:123}, we deduce  the main theorems from the comparison of the relative traces  of the chosen Hecke measures. {(The proof of Theorem \ref{main thm} is formally analogous to that of the Ichino--Ikeda conjecture in \S~\ref{sec: II}.)} 

\subsection{Comparison of relative-trace formulas}\lb{sec:121}
Recall the fixed $\nrm \in F_0^\ts$. For every inert finite place $v$ of $F_0$, we let $\e_v\in\{0,1\}$ have the same parity as $v(\nrm)$ as in \S~\ref{sss:ver level}.

The comparison will be based on the following local result.
Recall from \S\ref{eq: def int RZ} the definition of the local arithmetic intersection numbers $\sJ_{\delta, v}(e_{K_v})$ on RZ spaces in the vertex parahoric case.

\begin{theorem}[Arithmetic Fundamental Lemma and Arithmetic Transfer, {\cite{afl-pf1, afl-pf2, zhiyu,LMZ}}] \lb{afl thm}
Let $v$ be an odd inert place of $F_{0}$ and assume that either of the following conditions on $K_{v}\subset G_{v}$, $K_{v}'\subset G_{v}'$ holds:
\begin{enumerate}
\item  $K_{v}$ is hyperspecial, and $K_{v}'=\G'(\sO_{F_{0,v}})$;
\item  $K_{v}=K_{n,v}\times K_{n+1,v}$ is a vertex parahoric subgroup of type $(t,t+\epsilon_v)$ (cf. \S~\ref{sss:vert}), and $K_v'\subset G_v'$ is related to $K_v$ as in Proposition \ref{zhiyu tr}.
\end{enumerate}
Suppose that $\gamma\in B_{{\rm rs}, v}'$ matches an orbit $\delta=\ul{\delta}(\gamma)\in B_{{\rm rs}, v, V_{v}}$ for the hermitian pair $V_{v}$ with $\epsilon(V_{v})=-1$ (cf. \eqref{eq: HW}). Then
\beqq  \sJ_{{\delta},v}(e_{K_v})=\partial \sI_{\gamma,v}(e_{K'_v}).\eeqq
\end{theorem}
\begin{proof} 
By the definitions, the identity is equivalent to 
\beq\lb{afl 11}- (\delta\cdot\CN_{n,v},\CN_{n,v}) = 
{\partial \sI_{\gamma}(e_{K'_v}) /  \lm(\vpi_{w})}
 = \left.{d\over ds}\right|_{s=0} I_{\gamma,v}^{\C}(e_{K'_{v}}, |\cdot  |_{F_{0,v}}^{s}) / (-\log q_{0,v}^{2})
\eeq
(where $w$ is the place of $F$ above $v$, and the `division' in the second term has the obvious meaning). 

In the hyperspecial case where the cardinality of the residue field of $F_{0,v}$ is at least $n$, the identity \eqref{afl 11} is the Arithmetic Fundamental Lemma conjecture proved in \cite{afl-pf1, afl-pf2}. In the hyperspecial case where the cardinality of the residue field of $F_{0,v}$ is small, and in the vertex parahoric case (an instance of the Arithmetic Transfer conjecture), the identity \eqref{afl 11} is proved by Z. Zhang \cite{zhiyu} when $v$ is unramified over~$\Q$ and by Luo--Mihatsch--Z. Zhang \cite{LMZ} in the ramified case. 

 There are two points where the formulation in those works appears different. First, they consider a version with derivatives of `inhomogeneous' orbital integrals; this is verified to be equivalent to the above homogeneous version as in \cite[Proposition 14.1 (ii)]{RSZ-reg}. Second, their identity apparently differs from ours by a sign~$-1$: the reason is that their orbital integral contains a transfer factor defined as in \S~2.4 \emph{ibid.}; under our assumptions on $\gamma$ and $v$, that transfer factor (in its inhomogeneous version), evaluated at a preimage 
$\gamma'\in G'_{{\rm rs}, v}$ of $\gamma$, differs from our $\kappa_{v}(\gamma')$ by~$-1$. 
\end{proof}

We can now make the global comparison.

 \begin{theorem}[Comparison of RTFs]
 \lb{comp rtf}
Let $S$, $K^{p}=\prod_{v\nmid p} K_v$, and $$f^{p}\in\sH(\G(\A^{p}),L)_{K_{p}}^{\circ}$$ be as in Theorem \ref{ARTF thm}. Write $S=S^{\rm spl}\sqcup S^{\rm in}$ as a union of sets of split and inert places.

 Let  $K_{p}'\coloneqq \RG'(\sO_{F_{0,p}})$ and let $K'^p=\prod_{v\nmid p} K_v'\subset \G'(\A^{p\infty})$  be a compact open subgroup satisfying:
  \begin{itemize}
\item   for every $v \notin S $, $K'_{v}=\G'(\sO_{F_0,v})$ is hyperspecial;
\item for every inert $v\in S$, {$K_v'\subset G_v'$ is related to $K_v$ as in Proposition \ref{zhiyu tr}.}
\end{itemize}

Let 
$$f'^{p}=f'^{S p}\otimes  f'_{S^{\rm spl}\infty} \otimes f'_{S^{\rm in}}
\in \sH(\G'(\A^{p}), L)_{K'^{p}{\textup{-rs, qc}}}^{\circ} $$ 
be a quasicuspidal  {Gaussian with globally  regular semisimple support} whose factors satisfy the following properties:
\begin{itemize}
\item $f'^{S^{\rm spl}p}=\ot_{v}f'_{v}$ with   $f'_{v}=e_{K'_v}$;
\item $f'_{S^{\rm spl}\infty}$ matches $f_{S^{\rm spl} \infty}$.
\end{itemize}

  Then $\sI_{K'_{p}}(f'^{p}, \one)=0$ and 
$$\sJ_{K_{p}}(f^p)=\partial \sI_{K'_{p}}(f'^p).$$
\end{theorem}

\begin{proof}

We first show that $f^p$ and $f'^p$ match under the assumptions. Since $V$ is incoherent, this implies that $\sI_{K'_{p}}(f'^{p}, \one)=0$ by Proposition \ref{prop:de I}.

By our conditions and the Jacquet--Rallis Fundamental lemma (Proposition \ref{JRFL}), $f^{S^{\rm in}p}$ and $f'^{S^{\rm in}p}$ match.
Proposition \ref{zhiyu tr} on transfer at vertex parahoric levels shows that  $f_v$ and $f_v'$ match at the places in $S_{\rm in}$ as well.

Next we compare the geometric expansions of both sides of the desired equality, given by Theorem \ref{ARTF thm} and Proposition \ref{prop:de I} (3) respectively.
By the identity of Theorem \ref{afl thm}, 
these are equal term by term.
The proof is complete.
 \end{proof}

\subsection{Test Hecke measures} \lb{sec:122}
We find  some $f^{p}\in \sH(\G(\A^{p}), L)^{\circ}$, $f'^{p}\in  \sH(\G'(\A^{p}), L)^{\circ}$ to which the comparison may be applied, and that isolate a given pair of representations over $L$. {Recall the notion of local relative characters $J_{\pi_v}$ from \S~\ref{loc J def} and Remark \ref{rat matching}.}

\begin{lemma} \lb{lem key}
Let $\pi\in \sC(\RH\bs {\G})_{K_{p}}^{\ord, \rm st}(L)$ and let $\Pi={\rm BC}(\pi)$. Assume  that:
\begin{itemize}
\item  for every place $v$ of $F_0$ that is split in $F/F_0$, at least one of $\pi_{n,v}$ and $\pi_{n+1,v}$ is unramified;
\item  for every place $v$ of $F_0$ that is inert in $F/F_0$, $\pi_{n,v}$ and $\pi_{n+1,v}$ are either unramified or  almost unramified, and
if $\pi_{n,v}$ is not unramified then  $\pi_{n+1,v}$ is not   unramified.
 \end{itemize}
 
 Then there exist:
 \begin{itemize}
\item   a finite set $S$ of places of $F_{0}$, not above $p$ or $\infty$,  
 \item open compact subgroups $K^{p}=\prod_{v\nmid p} K_v\subset \G(\A^{p\infty})$ and  $K'^{p}=\prod_{v\nmid p} K'_v\subset \G'(\A^{p\infty})$, 
 \item Hecke measures  $f_{1}^{p}, f_{2}^{p}
 \in \sH(\G(\A^p), L)_{K_{p}}^{\circ}$ and $f'^{p}\in  \sH(\G'(\A^p), L)_{K'_{p}\textup{-rs, qc}}^{\circ}$,
 \end{itemize}
such that, upon setting  $f^{p}\coloneqq f_{1}^{p}\star f_{2}^{p,\vee}$:
\begin{itemize}
\item $(S, K^{p},   f_{1}^{p}, f_{2}^{p}, K'^{p},  f'^{p})$ satisfy the conditions of Theorem \ref{ARTF thm} and of Theorem \ref{comp rtf};
\item $M^{\oplus, *}.\hT(f_i^{p}e_{K_{p}})\subset \bigoplus_{\pi'\in{\rm BC}^{-1}(\Pi)} M_{\pi'}^{K_{p}}$;
\item $\Pi'(f'^{p}e_{K'_{p}})=0$ for every $\Pi' \in \sC(\G') -\{\Pi\}$;
\item $\ot_{v\nmid p}J_{\pi_{v}}(f^p ) = \ot_{v\nmid p} I_{\Pi, v}(f'^{p})\neq 0$.
\end{itemize}
\end{lemma}

\begin{proof}
We construct $f_{1}^p$, $f_{2}^{p}$, $f'^{p}$ as  products whose various factors take care of the required conditions.

\smallskip

\paragraph{\em Regularity of the supports} Let  $v_{+}$, $v_{-}$ be  two split places of $F_0$ at which $\Pi$  
 is an unramified regular principal series (cf. Lemma
\ref{cebo}). 
Let $$f_{v_{\pm}}' \in \sH(G_{v_{\pm}}', L)$$ be an element with $\pm$-regular support such that $I_{\Pi_{v_{\pm}}}(f_{v_{\pm}}', \one)\neq 0$ as provided by  Proposition \ref{test fin} \eqref{test fin 3}; we take any sufficiently small $K'_{v_\pm}$
such that $f_{v_{\pm}}' \in \sH(G_{v_{\pm}}', L)_{K_{v_{\pm}}'}$.
Upon a choice of a basis of $V_{v_{\pm}}$, we have the matching $f_{v_{\pm},1}\in \sH(G_{v_{\pm}}, L)$; by that proposition, we may arrange that $f_{v_{\pm},1}$ is bi-invariant under an Iwahori subgroup $K_{v_{\pm}}\subset G_{v_{\pm}}$.
We put $f_{v_{\pm},2}\coloneqq  e_{K_{v_{\pm}}}$. 
Thus $f_{v_{\pm}}\coloneqq f_{v_\pm,1}\star f_{v_{\pm},2}^{\vee}=f_{v_\pm,1}$ still matches $f_{v_{\pm}}'$ and has $K_{v_{\pm}}$-regular support.
Any global Hecke measure with component $f'_{v_{+}}\ot f'_{v_{-}}$ at $ \{v_{+}, v_{-}\}$  has  globally regular semisimple support (Definition \ref{def: weak support}) since $\G'_{\rs}=\G'_{\reg^{+}}\cap \G'_{\reg^{-}}$.

\smallskip

\paragraph{\em Choices at places of ramification} Let $S^{R}$ be the finite set  of places $v\nmid v_{+} v_{-} p\infty$ of $F_0$  
 where at least one of $\pi_{n,v},\pi_{n+1,v}$ is ramified. For every {\em split}  $v\in S^{R}$, we let $K_v=K_{n,v}\times K_{n+1,v}$ such that $\pi^{K_v}\neq 0$ and $K_{\nu,v}$ is  hyperspecial if $\pi_{\nu,v}$ is unramified.  
 Then we pick any 
$f_{v, 1}$, $ f_{v, 2}$, $f_{v}\coloneqq  f_{v, 1}\star f_{v,2}^{\vee}\in \sH(G_{v}, L)_{K_{v}}$ such that $J_{\pi_{v}}(f_{v})\neq 0$. For every {\em inert}  $v\in S^{R}$,  we let $K_v$ be a vertex parahoric subgroup of the type (depending on $\pi_v$) specified in Lemma \ref{lem dang},
 and we let $f_{v}=e_{K_v}$. 
We put $f_{S^{R},i}=\ot_{v\in S^{R}} f_{v,i}$ and we let $f'_{S^{R}}\in \sH(G'_{S^{R}})$  purely match $f_{S^{R}} \coloneqq  f_{S^{R},1} \star f_{S^{R},2}^{\vee}$.

\smallskip

 \paragraph{\em Isolation of $\pi$ and $\Pi$} 
 Now  take $S=S^{R}\cup \{v_{+}, v_{-}\}$. For $v\notin Sp$, we let $K_v$, $K_{v}'$ be hyperspecial, and form $K=\prod_{v}K_v$, $K'=\prod_{v}K_{v}'$.  
Consider the split Hecke algebras 
\beqq\BT=\BT^{\rm spl}&\coloneqq \bigotimes_{v\nmid Sp\atop \text{split}} \sH(G_{v}, L)_{K_v}& \subset \sH(\G(\A^{Sp}), L)^{\circ}_{K^{Sp}}\\
\BT'=\BT'^{\rm spl}&\coloneqq \bigotimes_{v\nmid Sp\atop \text{split}} \sH(G'_{v}, L)_{K'_v}\ot_{L} \sH(G'_{\infty}, L)^{\circ}&\subset  \sH(\G'(\A^{Sp}), L)^{\circ}_{K'^{Sp}}. \eeqq
Let $f_{\pi, 1}=f_{\pi, 2}\in \BT$ be  an element acting as the idempotent projection from   $M^{\oplus, *}_{K}$ onto $\bigoplus_{\pi'\in {\rm BC}^{-1}(\Pi)}M_{\pi'}^{K}$, which exists by Lemma \ref{isolate pi} (for $\Sigma$ the finite set of representations  occurring in $M^{\oplus, *}_{K}$).   Let $f_{\pi}'\in { \BT}' $ be an element supported at the finite places and matching $f_{\pi}\coloneqq f_{\pi, 1} \star f_{\pi, 2}^{\vee}$.

Let $f_{\Pi}'\in \BT'$ be an element such that {$\Pi(f'_{\Pi})=\Pi(\ot_{v\nmid Sp \textup{ split}} e_{K'_v})$} and that for each $\iota\colon L\into \C$, $R(f_{\Pi}'^{\iota})$ sends $\sA(\G')^{K'}$ into $\Pi^{\iota,K'}$, which exists by Proposition \ref{test gauss}; let $f_{\Pi,1}\in \BT$ be a matching element and let $f_{\Pi, 2}$ be the unit of $\BT$. 

\smallskip

\paragraph{\em Annihilation of absolute cohomology}
For every place $v\in S$, and every place $w\vert v$ of $F$, we fix an arbitrary auxiliary datum $(\Phi, u)$ relative to~$w$, which exists  by Remark \ref{rem exs CM} since  $d_{F/F_0}=1$, and we consider the scheme $\sX'_{u}$ from \eqref{aux int X'}. By the vanishing theorem of Proposition \ref{thm:van}  (1) (applied to  the maximal ideal $\fkm$ of $\BT$ corresponding to the eigensystem attached to $\pi$), there exists $f_{\{v\},1}=f_{\{v\}, 2}\in \BT$ which annihilates $H^{2n}(\sX'_{u},L(n))$,
and acts by a non-zero scalar on the line $\bigotimes_{v\nmid Sp\atop \text{split}}\pi_v^{K_v}$.  Let $f_{\{ v\}}'\in \BT'$ be an element supported at the finite places and matching $f_{\{v\}}\coloneqq f_{\{v \}, 1} \star f_{\{v\}, 2}^{\vee}$. 
 
\smallskip

\paragraph{\em Assembly} For $i=1,2, \emptyset$, we define 
$$f_{i}^{Sp} =f_{\pi, i} f_{\Pi, i}\ot \ot_{ v\in S} f_{\{v\}, i} \in \BT, \qquad
 f'^{Sp}= f_{\pi}' f_{\Pi}' \ot\ot_{v\in S} f'_{\{ v\}} \in \BT',$$ 
 viewed naturally as  elements in $ \sH(\G(\A^{Sp}), L)^{\circ}_{K^{Sp}}$, $\sH(\G'(\A^{Sp}), L)^{\circ}_{K'^{Sp}}$. Then we define 
 $$
 f_i^p=   f_{S, i} f_{i}^{Sp }, \qquad f'^{p}=   f'_{S} f'^{Sp }.
 $$

Then it is easy to see that, by construction, $f_i^p$ satisfies the  required conditions. To check the  condition on relative characters, we use 
$$
\ot_{v\nmid p}J_{\pi_{v}}
(f^p )= \ot_{v\notin Sp} J_{\pi_{v}}(f^{Sp})\prod_{v\in S}J_{\pi_{v}}(f_{v} ).
$$
The product over $v\in S$ does not vanish by construction; the first factor is  the product of  $ \ot_{v\notin Sp}J_{\pi_{v}}(e_{K^{Sp}})\neq 0$ and of 
 the eigenvalue of $ f^{Sp}$ acting on the line $\pi^{K^{Sp}}$, which is a non-zero scalar.
 \end{proof}

\subsection{Proofs of the main theorems} \lb{sec:123}
We first reduce the identity
\beq\lb{eq: ref aggp}
{h_{\pi} (Z_{\pi}(\phi), Z_{\pi^{\vee}}(\phi'))}
 =  e_{p}(\RM_{\Pi})^{-1} \cdot {1\over 4} \partial \sL_{p}  (\RM_{\Pi}) 
\cdot \alpha(\phi, \phi')\eeq
 of Theorem \ref{main thm} to the   factorization 
\beq\lb{RTF fact} \sJ_{\pi, K_{p}}(f^{p})= {1\over 4} \partial \sL_{p}(\RM_{\Pi})\cdot  \ot_{v\nmid p}J_{\pi_{v}}(f^p ).\eeq

\begin{lemma} \lb{reduce to J}
Let $\pi \in \sC(\G)_{K_{p}}^{\ord}$, and let $\Pi={\rm BC}(\pi)$, $L=\Q_{p}(\pi)$. The following are equivalent:
\begin{enumerate}
\item\lb{D for all} For every $\phi \in \pi$, $\phi'\in \pi^{\vee}$, the identity \eqref{eq: ref aggp} holds.
\item\lb{D for one} For some $\phi \in \pi$, $\phi'\in \pi^{\vee}$ such that $\alpha(\phi, \phi')\neq 0$, the identity \eqref{eq: ref aggp} holds.
\item \lb{RTF for all}
For every $f^{p}\in \sH(\G(\A^p), L)^{\circ}$, the factorization \eqref{RTF fact} holds.
\item \lb{RTF for one} For some  $f^{p}\in \sH(\G(\A^p), L)^{\circ}$ such that $\ot_{v\nmid p}J_{\pi_{v}}(f^p )\neq 0$, the factorization \eqref{RTF fact} holds.
\end{enumerate}
\end{lemma}

\begin{proof}
It is trivial that \eqref{D for all} implies \eqref{D for one}, and \eqref{RTF for all} implies \eqref{RTF for one}. The two converse implications follow from multiplicity one and the non-vanishing of $\alpha$. 

We prove that \eqref{RTF for all} is equivalent to  \eqref{D for all}. It is clear that \eqref{D for all} is equivalent to 
\beq \lb{rtr id} \Tr_{(,)_{\pi}}^{h\circ Z_{\pi}\boxtimes Z_{\pi^{\vee}}}( \tau) =  
 e_{p}(\RM_{\Pi})^{-1} \cdot {1\over 4} \partial \sL_{p}  (\RM_{\Pi})\cdot \Tr_{(,)_{\pi}}^{\alpha}(\tau) \eeq
 for all $\tau\in \End(\pi)$, and equivalently for some $\tau$ such that $ \Tr_{(,)_{\pi}}^{\alpha}(\tau)\neq 0$. Thus it is enough to show that \eqref{RTF fact} is equivalent to  \eqref{rtr id} for some such $\tau$. 

 Choose a factorization $(,)_{\pi}=(,)_{\pi^{p}}(,)_{\pi_{p}}$. For any $N\geq 1$, let $f_{p, K_{p},N}\coloneqq  \eqref{fp G}\in \sH(G_{p},L)$,  let $f_{p,K_{p}, N}^{\star}\coloneqq f_{p,K_{p}, N}\star f_{p,K_{p}, N}^{\vee}$,   and for $?\in \{\emptyset, \vee\} $, let 
$$\pi_{p}^{?}(f_{p,K_{p}})\coloneqq \lim_{N\sto \infty}\pi_{p}^{?}(f_{p, K_{p},N})\in \End(\pi_{p}).$$
(This does not depend on the integer $1\leq r\leq N!$ implicit in \eqref{fp G}.)
 Let 
 $$\pi_{p}(f^{\star}_{p, K_{p}}) \coloneqq \pi_{p}(f_{p, K_{p}}) \circ(\pi_{p}^{\vee}( f_{p, K_{p}}))^{\vee},$$
  where $(-)^{\vee}$ denotes the transpose with respect to  $(,)_{\pi_{p}}$. Then by the definition in \S~\ref{J sp exp}, we have 
\beq\lb{sJ tr}
\Tr_{(,)_{\pi}}^{h\circ Z_{\pi}\boxtimes Z_{\pi^{\vee}}}(\pi^{p}(f^{p}) \pi_{p}(f^{\star}_{p, K_{p}})) =  \sJ_{\pi, K_{p}}(f^{p}).
\eeq

On the other hand, it is clear from the definitions that 
\beq\lb{J tr}
\Tr_{(,)_{\pi}}^{\alpha}(\pi^{p}(f^{p})\pi_{p} (f^{\star}_{p, K_{p}})) = \ot_{v\nmid p} J_{\pi_{v}} (f^{p}) \cdot \lim_{N\sto \infty}J_{\pi_{p}}(f_{p, K_{p},N}^{\star}).
\eeq
Now by Lemma \ref{match split}, $f^{\star}_{p, K_{p},N}$ matches the function $f'_{p,K_{p}',N}$ attached to $U_{t_{p}}^{N!}$ as in Lemma \ref{I Pi dag}. By the definitions and  Corollary \ref{Iord ep}, we then have 
\beq \lb{J p}
\lim_{N\sto \infty} J_{\pi_{p}}(f_{p, K_{p},N}^{\star})= \lim_{N\sto \infty} I_{\Pi_{p}}(f'_{p,K_{p}', N}) =e_{p}(\RM_{\Pi}).\eeq
(Recall that $e_{p}(\RM_{\Pi})$ is the product of the factors $e(\Pi_{v}, \one_{v})$ of \eqref{e v ord}.) Thus by \eqref{sJ tr}, \eqref{J tr}, \eqref{J p}, the identity \eqref{RTF fact} for $f^{p}$ is equivalent to  \eqref{rtr id} for    $\tau=\pi^{p}(f^{p}) \pi_{p}(f_{p, K_{p}}^{\star}).$
This completes the proof.
\end{proof}

We may now prove Theorem \ref{main thm} based on the comparison of relative-trace formulas in Theorem \ref{comp rtf}.

\begin{proof}[Proof of Theorem \ref{main thm}] By Lemma \ref{reduce to J}, it suffices to prove 
\beq\lb{des fact}
 \sJ_{\pi, K_{p}}(f^{p})= {1\over 4} \partial \sL_{p}(\RM_{\Pi})\cdot  \ot_{v\nmid p}J_{\pi_{v}}(f^p )
\eeq
for any   $f^p$ such that  $ \ot_{v\nmid p}J_{\pi_{v}} (f^p )\neq 0$.

Let $S$, $K^{p}$, $f^{p}$, $f'^{p}$ be as in Lemma \ref{lem key}.  By construction, $\ot_{v\nmid p}J_{\pi_{v}}(f^p)\neq 0$, the elements  $f^{p}$ and $f'^{p}$ match (geometrically), and  Theorem  \ref{comp rtf} is applicable and it gives
$$\sJ_{K_{p}}(f^p)=\partial \sI_{K'_{p}}(f'^p).$$
By Theorem \ref{ARTF thm} and  Proposition \ref{prop:de I} (2), we have an equality of spectral expansions
$$
\sum_{\pi\in \sC(\RH\bs\G)_{K_{p}}^{\ord}} \sJ_{\pi, K_{p}}(f^{p})
=\sum_{\Pi\in \sC(\G')^{{\rm her, ord}, V}_{K_{p}}} \partial\sI^{}_{\Pi, K_{p}} (f'^{p}),$$
but by construction only the terms corresponding to $\pi$ and $\Pi$ may be nonzero. We deduce that 
 $$\sJ_{\pi, K_{p}}(f^{p})= {1\over 4} \partial \sL_{p}(\RM_{\Pi})\cdot  \ot_{v\nmid p}I_{\Pi_{v}}(f'^p ),$$
 which is equivalent to the desired factorization \eqref{des fact} by the  (spectral) matching of $f^{p}$ and~$f'^{p}$. 
 \end{proof}

\begin{proof}[Proof of Theorem \ref{th BBK}]
Let $(V, \pi)$ be the distinguished element in the Vogan $L$-packet of  $\Pi$, in the sense of Proposition \ref{cor C her}. Since $\eps(\Pi)=-1$, the pair $V$ is incoherent.

We need to verify that, for each finite place $v$ of $F_0$ that is inert in $F$, the pair $(V_v,\pi_v=\pi_{n,v}\boxtimes\pi_{n+1,v})$ satisfies the conditions of Theorem \ref{main thm}.

There are three cases (in all of which both $\Pi_{n,v}$, $\Pi_{n+1, v}$ have trivial central character).
\begin{itemize} 
\item Both $\Pi_{n,v}$ and $\Pi_{n+1,v}$ are  unramified. Then $\e(V_v)=+1$  and $\pi_{v}$
is unramified.
\item $\Pi_{n,v}$ is unramified and $\Pi_{n+1,v}$ has conductor~$1$. 
{Then by Lemma \ref{lem: GGP au} and Remark \ref{rmk: GGP au}, we have that $\e(V_v)=(-1)^n$ and $\pi_{n+1,v}$ is almost unramified and not unramified.}
\item both $\Pi_{n,v}$ and   $\Pi_{n+1, v}$ have conductor~$1$.
{Then by Lemma \ref{lem: GGP au} and Remark \ref{rmk: GGP au}, we have that $\e(V_v)=-1$ and $\pi_{v}$ is almost unramified.}
\end{itemize}

Then the main implication of Theorem \ref{th BBK}   follows immediately from Theorem \ref{main thm} applied to~$\pi$. 
The strengthened implication of Theorem \ref{th BBK} then follows from \cite{LTXZZ}  (or \cite{Lai-Ski} under a different condition), as observed in Remark \ref{C and D}.
\end{proof}

  \bigskip

\
  
   \bigskip

\addtocontents{toc}{\medskip}

\printindex

\end{document}